\documentclass[authoryear]{elsarticle}
\usepackage{listings}
\usepackage[utf8]{inputenc}
\usepackage{amsmath,amssymb,amsfonts,amsthm,bbm,mathrsfs,verbatim,bm} 
\usepackage[misc]{ifsym}
\usepackage{graphics}                 
\usepackage{graphicx}
\usepackage{caption}
\usepackage{subcaption}
\usepackage{float}
\usepackage{array}
\usepackage{hyperref}                 
\usepackage{makeidx}
\usepackage{algorithm}
\usepackage{graphicx}
\usepackage[shortlabels]{enumitem}
\usepackage{amsmath}
\usepackage[nottoc]{tocbibind}
\usepackage{tcolorbox}
\usepackage{mdframed}
\usepackage{pdfpages}
\usepackage[utf8]{inputenc}
\usepackage[english]{babel}
\graphicspath{ {images/} } 

 \usepackage{blindtext}
\usepackage[utf8]{inputenc}
\usepackage{natbib}
\usepackage{lipsum}
\usepackage{xr}
\usepackage{blindtext}
\usepackage[T1]{fontenc}
\usepackage[utf8]{inputenc}
\usepackage{blindtext}
\usepackage[T1]{fontenc}
\usepackage[utf8]{inputenc}

\usepackage{tikz}
\usetikzlibrary{positioning}
\newdimen\nodeDist
\usepackage[utf8]{inputenc}
\usepackage[english]{babel}
 \usepackage{amsmath}
\usepackage{caption}
\usepackage{algorithm}
\usepackage{algorithmic}
\newtheorem{remark}{Remark}

\usepackage{xcolor}
\usepackage{float}
\usepackage{bm}

\newcommand{\vs}{\boldsymbol{s}}

\DeclareMathOperator*{\TSR}{TSR}
\DeclareMathOperator*{\LR}{LR}
\DeclareMathOperator*{\TR}{TR}

\title{BART-based inference for Poisson processes}
\author[1]{Stamatina Lamprinakou\corref{cor1}}
\ead{s.lamprinakou18@imperial.ac.uk}
\cortext[cor1]{Corresponding author}
\author[1]{Mauricio Barahona}
\author[1]{Seth Flaxman}
\author[1]{Sarah Filippi}
\author[1]{Axel Gandy}
\author[1]{Emma McCoy}
\affiliation[1]{Department of Mathematics, Imperial College London, London, United Kingdom}

\begin{document}
\begin{abstract}
The effectiveness of Bayesian Additive Regression Trees (BART) has been demonstrated in a variety of contexts including non-parametric regression and classification. A BART scheme for estimating the intensity of inhomogeneous Poisson processes is introduced. Poisson intensity estimation is a vital task in various applications including medical imaging, astrophysics and network traffic analysis. The new approach enables full posterior inference of the intensity in a non-parametric regression setting. The performance of the novel scheme is demonstrated through simulation studies on synthetic and real datasets up to five dimensions, and the new scheme is compared with alternative approaches.
\end{abstract}

\maketitle

\section{Introduction}

The Bayesian Additive Regression Trees (BART) model is a Bayesian framework, which uses a sum of trees to predict the posterior distribution of a response $y$ given a $p$-dimensional covariate $X$ and priors on the
function relating the covariates to the response.  \citet{chipman2010bart} proposed an inference procedure using Metropolis Hastings within a Gibbs Sampler, whereas \citet{lakshminarayanan2015particle} used a Particle Gibbs Sampler to increase mixing when the true posterior consists of deep trees or when the dimensionality of the data is high. Several theoretical studies of BART models~\citep{rockova2017posterior, rockova2018theory, linero2018bayesian} have recently established optimal posterior convergence rates. The BART model has been applied in various contexts including non-parametric mean regression \citep{chipman2010bart}, classification \citep{chipman2010bart,zhang2010bayesian,kindo2016multinomial}, variable selection\citep{chipman2010bart, bleich2014variable,linero2018bayesian2}, estimation of monotone functions \citep{10.1214/21-BA1259}, causal inference \citep{hill2011bayesian}, survival analysis \citep{sparapani2016nonparametric}, and heteroskedasticity \citep{bleich2014bayesian, doi:10.1080/10618600.2019.1677243}.  \citet{linero2018bayesian} illustrated how the BART model suffers from a lack of smoothness and the curse of dimensionality, and overcome both potential shortcomings by considering a sparsity assumption similar to \citep{linero2018bayesian2} and treating decisions at branches probabilistically.

The original BART model \citep{chipman2010bart} assume that the response has a Gaussian distribution and the majority of applications have used this framework. \citet{murray2017log} adapted the BART model to count data and categorical data  via a log-linear transformation, and provided an efficient MCMC sampler. Our focus is on extending this methodology to estimate the intensity function of inhomogeneous Poisson processes.   
 
 The question of estimating the intensity of Poisson processes has a long history, including both frequentist and Bayesian methods. Frequentist methods include fixed-bandwidth and adaptive bandwidth kernel estimators with edge correction \citep{diggle1983statistical}, and  wavelet-based methods \citep[e.g.\ ][]{fryzlewicz2004haar,patil2004counting}. Bayesian methods include using a sigmoidal Gaussian Cox process model for intensity inference \citep{adams2009tractable}, 
 a Markov random field (MRF) with Laplace prior \citep{doi:10.1198/016214504000000188},
 variational Bayesian intensity inference \citep{lloyd2015variational}, and non-parametric Bayesian estimations of the intensity via piecewise functions with either random or fixed partitions of constant intensity \citep{arjas1994nonparametric, heikkinen1998non, gugushvili2018fast}.  

 In this paper, we introduce an extension of the BART model~\citep{chipman2010bart} for Poisson Processes whose intensity at each point 
 is estimated via a tiny ensemble of trees. Specifically, the logarithm of the intensity at each point is modelled via a sum of trees (and hence the intensity is a product of trees). 
 This approach enables full posterior inference of the intensity in a non-parametric regression setting. Our main contribution is a novel BART scheme for estimating the intensity of an inhomogeneous Poisson process. The simulation studies demonstrate that our algorithm is competitive with the Haar-Fisz algorithm in one dimension, kernel smoothing in two dimensions, and outperforms the kernel approach for multidimensional intensities. The simulation analysis also demonstrates that our proposed algorithm is competitive with the inference via spatial log-Gaussian Cox processes. We also demonstrate its ability to track varying intensity in synthetic and real data.  
 
 The outline of the article is as follows. Section 2 introduces our approach for estimating the intensity of a Poisson process through the BART model, and Section 3 presents the proposed inference algorithm. Sections 4 and 5 present the application of the algorithm to synthetic data and real data sets, respectively. Section 6 provides our conclusions and plans for future work.

\section{The BART Model for Poisson Processes}
%
 Consider an inhomogeneous Poisson process defined on a $d$-dimensional domain $S \subset \mathbb{R}^d$, $d \geq 1$, with intensity $\lambda: S \rightarrow \mathbb{R}^+$. 
 For such a process, the number of points within a subregion $B \subset S$ has a Poisson distribution with mean $\lambda_B=\int_{B}\lambda(s) \, \mathrm{d}s$,
 and the number of points in disjoint subregions are independent~\citep{daley2003elementary}. The homogeneous Poisson process is a special case with constant intensity $\lambda(s)=\lambda_0, \, \forall s\in S$.

To estimate the intensity of the inhomogeneous Poisson process, we use $m$ partitions of the domain $S$, each associated with a tree $T_h, \, h=1,\ldots, m$. The partitions are denoted $T_h=\{\Omega_{ht}\}_{t=1}^{b_h}$, where $b_h$ is the number of terminal nodes in the corresponding tree $T_h$, and each leaf node $t$ corresponds to one of the subregions $\Omega_{ht}$ of the partition $T_h$. Being a partition, every tree covers the full domain, i.e.\ $\mathop{S=\cup_{t=1}^{b_h}\Omega_{ht}}$ for every $h$. 
  Each subregion $\Omega_{ht}$ has an associated parameter $\lambda_{ht}$, and hence each tree $T_h$ has an associated vector of leaf intensities $\Lambda_h=(\lambda_{h1},\lambda_{h2},..,\lambda_{hb_h})$.

We model the intensity  of $s \in S$ as: 
\begin{align} 
\log(\lambda(s))&=\sum_{h=1}^m \sum_{t=1}^{b_h} \log\left(\lambda_{ht}\right) \, I(s \in \Omega_{ht})  \label{eq:model}  \\
T_h &\sim \text{heterogeneous Galton-Watson process for a partition of $S$ } \label{eq:GW} \\ 
\lambda_{ht}|T_h &\overset{\text{iid}}\sim  \text{Gamma}(\alpha,\beta)  \label{eq:Gamma}
\end{align}
where $I(\cdot)$ denotes the indicator function.  Equivalently, \eqref{eq:model}  can be expressed as
\begin{align} \label{eq:lambda_product}
\lambda(s) &=\prod_{h=1}^m\prod_{t=1}^{b_h} \lambda_{ht}^{I(s\in \Omega_{ht})}.
\end{align}

Given a fixed number of trees, $m$, the parameters of the model are thus the regression trees $T=\{T_h\}_{h=1}^m$ and their corresponding intensities $\Lambda=\{\Lambda_h\}_{h=1}^m$. 
Following~\cite{chipman2010bart}, we assume that the tree components $(T_h,\Lambda_h)$ are independent of each other, and that the terminal node parameters of every tree are independent, so that the prior can be factorized as:
\begin{equation}
P(\Lambda,T)=\prod_{h=1}^{m}P(\Lambda_h,T_h)=\prod_{h=1}^{m}P(\Lambda_h|T_h)P(T_h)=\prod_{h=1}^m \left[\prod_{t=1}^{b_h}P(\lambda_{ht}|T_h) \right] P(T_h).
\label{eq:factorized_P}
\end{equation}

\paragraph*{Prior on the trees} 
The trees $T_h$ of the BART model are stochastic regression trees
generated through a heterogeneous Galton-Watson (GW) process~\citep{harris1963theory,rockova2018theory}. The GW process is the simplest branching process concerning the evolution of a population in discrete time. Individuals (tree nodes) of a generation (tree depth) give birth to a random number of individuals (tree nodes), called offspring, mutually independent and all with the same offspring distribution that may vary from generation (depth) to generation (depth).
In our case, we use the prior introduced by \citet{chipman1998bayesian}, that is a GW process in which each node has either zero or two offspring 
and the probability of a node splitting depends on its depth in the tree. Specifically, a node $\eta \in T_h$ splits into two offsprings with probability
\begin{align}
\label{eq:splitting}
p_\text{split}(\eta)=\frac{\gamma}{\left(1+d(\eta)\right)^\delta},
\end{align}
where $d(\eta)$ is the depth of node $\eta$ in the tree, and $\gamma \in(0,1)$ and $\delta \geq 0$ are parameters of the model. Classic results from the theory of branching processes show that $\gamma \leq 0.5$ guarantees that the expected depth of the tree is finite. %
In our construction, each tree $T_h$ is associated with a partition of $S$. Namely, if node $\eta$ splits, we select uniformly at random one of the $d$ dimensions of the space of the Poisson process, followed by uniform selection from the available split values associated with that dimension respecting the splitting rules higher in the tree.

\paragraph*{Prior on the leaf intensities} 
Our choice of a Gamma prior for the leaf parameters $\lambda_{ht}$ builds upon previous work by~\citet{murray2017log}, who used a mixture of Generalized Inverse Gaussian (GIG) distributions as the prior on leaf parameters in a BART model for count regression. Here we impose a Gamma prior (a special case of GIG) on the leaf parameters, which simplifies the model and leads to a closed form of the conditional integrated likelihood below (see Section \ref{sIA}) as the Gamma distribution is the conjugate prior for the Poisson likelihood. We discuss the selection of it hyperparameters $\alpha$ and $\beta$ in Section \ref{sec:hyperparameters}.

\section{The Inference Algorithm}\label{sIA} 

Given a finite realization of an inhomogeneous Poisson process with $n$ sample points $\vs = s_1,\ldots,s_n \in S \subset \mathbb{R}^d$, we seek to infer the parameters of the model $(\Lambda,T)$ by sampling from the posterior $P(\Lambda,T|\vs)$.   

Before presenting the sampling algorithm we summarize a preliminary result.
To simplify our notation, let us define 
\[g(s_i;T_h,\Lambda_h) 
=\prod_{t=1}^{b_h}\lambda_{ht}^{I(s_i\in \Omega_{ht})},\] 
so that Eq.~\eqref{eq:lambda_product} becomes 
$\lambda(s_i) =\prod_{h=1}^m g(s_i;T_h,\Lambda_h)$.

Let us choose any arbitrary tree $T_h$ in our ensemble $T$, and let us denote the set with the rest of the trees as $T_{(h)}=\{T_j\}_{j=1,j\neq h}^{m}$ and their leaf parameters as $\Lambda_{(h)}=\{\Lambda_j\}_{j=1,j\neq h}^{m}$. 
The intersection of all the partitions associated with the trees in $T_{(h)}$ gives us a global partition $\{\overline{\Omega}_{k}^{(h)}\}_{k=1}^{K(T_{(h)})}$ with $K(T_{(h)})$ subregions~\citep{rockova2017posterior}.   

Then we have the following result.

\begin{remark}

\begin{enumerate}[label=(\roman*)]
\item \label{rem1}
The conditional likelihood of the realization  is given by
\begin{align} \label{equ2}  
P(\vs |\Lambda,T)&=c_h \, \prod_{t=1}^{b_h} \lambda_{ht}^{n_{ht}} \, e^{-\lambda_{ht}c_{ht}}, \\
\text{with}\quad  c_h &=\prod_{i=1}^n \prod_{j=1, j \neq h}^m g (s_i;T_j,\Lambda_j), 
\nonumber\\
c_{ht}& =\sum_{k=1}^{K(T_{(h)})}\overline{\lambda}_k^{(h)}|\overline{\Omega}_{k}^{(h)}\cap \Omega_{ht}|, \nonumber
\end{align}
where 
$\overline{\lambda}_{k}^{(h)}=\prod_{t=1,t \neq h}^{m} \prod_{l=1}^{b_t}\lambda_{tl}^{I(\Omega_{tl} \cap \overline{\Omega}_{k}^{(h)}\neq 0)}$,
$n_{ht}$ is the cardinality of the set $\{i:s_i\in \Omega_{ht}\}$, and
$ |\overline{\Omega}_{k}^{(h)} \cap \Omega_{ht}|$ is the volume of the region $ \overline{\Omega}_{k}^{(h)}\cap \Omega_{ht}$.
\item \label{rem2} For a tree $h$, the conditional integrated likelihood obtained by integrating out $\Lambda_h$ is
\begin{align}
\label{equ3}
P(\vs |T_h,T_{(h)},\Lambda_{(h)})= c_h  \left( \frac{\beta^{\alpha}}{\Gamma(\alpha)}  \right)^{b_h}\prod_{t=1}^{b_h}\frac{\Gamma(n_{ht}+\alpha)}{(c_{ht}+\beta)^{n_{ht}+\alpha}}.
\end{align} 
\end{enumerate}
\end{remark}
A proof can be found in ~\ref{app:remark1}~and~\ref{app:remark2}.

We now summarize our sampling algorithm. To sample from $P(\Lambda, T|\vs)$, we implement a
Metropolis-Hastings within block Gibbs sampler (Algorithm~\ref{alg:MH}), which requires $m$ successive draws from  $(T_h,\Lambda_h)|T_{(h)},\Lambda_{(h)},\vs$. 
Note that
\begin{align} 
P(T_h,\Lambda_h&|T_{(h)},\Lambda_{(h)},\vs)  \nonumber 
=P(T_h|T_{(h)},\Lambda_{(h)},\vs) \,
P(\Lambda_h|T_h,T_{(h)},\Lambda_{(h)},\vs) \nonumber \\
&\propto P(T_h|T_{(h)},\Lambda_{(h)},\vs) \, P(\vs|\Lambda,T)P(\Lambda_h|T_h)  \nonumber \\
&= P(T_h|T_{(h)},\Lambda_{(h)},\vs)  \, P(\vs|\Lambda,T)\, \, \prod_{t=1}^{b_h}P(\lambda_{ht}|T_h) \nonumber \\
&= P(T_h|T_{(h)},\Lambda_{(h)},\vs) \,
c_h \, \prod_{t=1}^{b_h} \lambda_{ht}^{n_{ht}} \, e^{-\lambda_{ht}c_{ht}}
\prod_{t=1}^{b_h} \frac{\beta^\alpha}{\Gamma(\alpha)} \lambda_{ht}^{\alpha-1} \, e^{-\beta \lambda_{ht}} \nonumber \\
&\propto P(T_h|T_{(h)},\Lambda_{(h)},\vs) \,
 \prod_{t=1}^{b_h}\lambda_{ht}^{n_{ht}+\alpha-1}e^{-(c_{ht}+\beta)\lambda_{ht}} \\
 \label{eq:conditional_draw}
\end{align}
which follows directly from Bayes' rule and Eqs.~\eqref{eq:factorized_P}~and~\eqref{eq:Gamma}. 

From~\eqref{eq:conditional_draw}, it is clear that a draw from $(T_h,\Lambda_h)|T_{(h)},\Lambda_{(h)},\vs$ can be achieved in ($b_h$+1) successive steps consisting of:
\begin{itemize}
\item{sampling $T_h|T_{(h)},\Lambda_{(h)},\vs$ using Metropolis-Hastings (Algorithm~\ref{alg2}})
\item{sampling $\lambda_{ht}|T_h, T_{(h)},\Lambda_{(h)},\vs$  from a Gamma distribution with shape $n_{ht}+\alpha$ and rate $c_{ht}+\beta$  for $t=1,..,b_h$.}
\end{itemize} 
These steps are implemented through Metropolis-Hastings in Algorithm~\ref{alg:MH}. 
Note also that  \[P(T_h|T_{(h)},\Lambda_{(h)},\vs) \propto P(\vs|T_h,T_{(h)},\Lambda_{(h)}) \, P(T_h),\] so that the conditional integrated likelihood~\eqref{equ3} is required to compute the Hastings ratio.

\begin{algorithm}[H] 
\caption{Metropolis-Hastings within Gibbs sampler} 
\begin{algorithmic}
\label{alg:MH}
\FOR{$v=1,2,3,..$ }
\FOR{$h=1$ to $m$ }
\STATE{Sample $T_h^{(v+1)}|\vs, \{T_j^{(v+1)}\}_{j=1}^{h-1}, \{T_j^{(v)}\}_{j=h+1}^{m},\{\Lambda_j^{(v+1)}\}_{j=1}^{h-1}, \{\Lambda_j^{(v)}\}_{j=h+1}^m$ \\ \hspace*{0cm} using Algorithm~\ref{alg2}}
\FOR{$t=1$ to $b_h$}
\STATE {Sample $\lambda_{ht}^{(v+1)}|\vs,\{T_j^{(v+1)}\}_{j=1}^{h}, \{T_j^{(v)}\}_{j=h+1}^m, \{\Lambda_j^{(v+1)}\}_{j=1}^{h-1}, \{\Lambda_j^{(v)}\}_{j=h+1}^m$} \,  \, \text{from \mbox{Gamma}($n_{ht} + \alpha$, $c_{ht}+\beta$) }
\ENDFOR
\ENDFOR
\ENDFOR
\end{algorithmic}
\end{algorithm}

\begin{algorithm}[H] 
\caption{Metropolis-Hastings Algorithm for sampling from the posterior  $P(T_j|\vs,T_{(j)},\Lambda_{(j)})$ } 
\label{alg2}
\begin{algorithmic}
\STATE{Generate a candidate value $T_j^{*}$ with probability $q(T_j^{*}|T_j^{(v)})$.}
\STATE {Set $T_j^{(v+1)}=T_j^{*}$ with probability 
\begin{equation*}
\alpha(T_j^{(v)},T_j^{*})=\min \left\{1,\frac{q(T_j^{(v)}|T_j^{*})}{q(T_j^{*}|T_j^{(v)})} \frac{P(\vs|T_j^{*},T_{(j)},\Lambda_{(j)})}{P(\vs|T_j^{(v)},T_{(j)},\Lambda_{(j)})} \frac{P(T_j^{*})}{P(T_j^{(v)})} \right \}
\end{equation*}Otherwise, set $T_j^{(v+1)}=T_j^{(v)}$.
}
\end{algorithmic}
\end{algorithm}

 The transition kernel $q$ in Algorithm~\ref{alg2} is chosen from the three proposals: GROW, PRUNE, CHANGE \citep{chipman2010bart,kapelner2013bartmachine}. The GROW proposal randomly picks a terminal node, splits the chosen terminal into two new nodes and assigns a decision rule to it. The PRUNE proposal randomly picks a parent of two terminal nodes and turns it into a terminal  node by collapsing the nodes below it. The CHANGE proposal randomly picks an internal node and randomly reassigns to it a splitting rule. We describe the implementation of the proposals in \ref{ap:A}.

 For completeness, in the supplementary material, we present the full development of the algorithm for inference of the intensity of inhomogeneous Poisson processes via only one tree.

\subsection{Fixing the hyperparameters of the model}
\label{sec:hyperparameters}

\paragraph*{Hyperparameters of the Gamma distribution for the leaf intensities}
We use a simple data-informed approach to fix the hyperparameters $\alpha$ and $\beta$ of the Gamma distribution~\eqref{eq:Gamma}. We discretize the domain into $N_G$ subregions of equal volume ($N_G=(\lceil100^{1/d}\rceil)^d$ works well in practice up to 5 dimensions) and count the number of samples $s_i$ per subregion. We thus obtain the empirical densities in each of the subregions: $\xi_i, \, i=1,\ldots,N_G$.
Given the form of the intensity~\eqref{eq:lambda_product} as a product of $m$ trees, we consider the $m$-th roots $\Xi=\{\xi_i^{1/m}\}_{i=1}^{N_G}$ as candidates for the intensity of each tree. Taking the sample mean $\widehat{\mu}_{\Xi}$ and sample variance $\widehat{\sigma}^2_{\Xi}$, we choose the model hyperparameters $\alpha$ and $\beta$ to correspond to those of a Gamma distribution with the same mean and variance, i.e.,
 $\alpha=\widehat{\mu}^2_\Xi/\widehat{\sigma}^2_\Xi$ and $\beta=\widehat{\mu}_\Xi/\widehat{\sigma}^2_\Xi$, although fixing $\beta=1$ can also give good estimates of the intensity.
Although setting $N_G=(\lceil100^{1/d}\rceil)^d$  leads to convergence and good estimates of the intensity in our simulation studies below, there are other possibilities.  Alternatively, we can bin the data based on a criterion that takes into account the number of samples, $n$, and the number of dimensions, $d$. For example, the number of bins per dimension, $n_b$, can be computed as~\citep{inbook, wand1997data}: 
(i) $n_b=\lceil{n^{1/(d+1)}}\rceil$, 
(ii) $n_b=\lceil{n^{1/(d+2)}\rceil}$, or 
(iii) $n_b=\max_{k\in\{1,2,..,d\}} [\lceil{DR_k \cdot n^{1/(d+2)}/(2 \cdot \text{IQR}(\{s_{i,k}\})\rceil}]$,
where IQR denotes the interquartile range of the sample, $DR_k$ is the range of the domain in dimension $k$ (here we scale the initial domain to a unit hypercube so that $DR_k=1, \, \forall k$), and by extension $N_G={n_b}^{d}$.
In our simulation scenarios below, all these approaches lead to comparable convergence times and estimates of the intensity.

\paragraph*{Hyperparameters of the stochastic ensemble of regression trees}
The GW stochastic process that generates our tree ensemble has several hyperparameters. The parameters $(\gamma,\delta)$ control the shape of trees. The parameter $\gamma>0$ controls the probability that the root of a tree will split into two offspring, while the parameter $\delta>0$ penalizes against deep trees.  As noted in \citep{chipman2010bart}, for a sum-of-trees model, we want to keep the depth of the tree small whilst ensuring non-trivial trees, hence, in our simulation study we fix $\gamma=0.98$ and $\delta=2$. 
Second, each of the $d$ dimensions has to be assigned a grid of split values, from which the subregions of the partition are randomly chosen, yet always respecting the consistency of the ancestors in the tree (that is respecting the splitting rules higher in the tree). Here, we use a simple uniform grid for each of the $d$-dimensions~\citep{pratola2016efficient}: we normalize each dimension of the space from (0,1) and discretize each dimension into $N_d$ segments. ($N_d=100$ works well in practice and is used throughout our examples.) More sophisticated, data-informed grids are also possible, although using, e.g., the sample points as split values does not improve noticeably the performance in our examples. 
Finally, the number of trees $m$ also needs to be fixed as in \citet{chipman2010bart}. In our examples below, we have checked the performance of our algorithm with varying number of trees $m$ between 2 and 50. We find that good performance can be achieved with a moderate number of trees, $m$, between 3 and 10 depending on the particular example.  

\section{Simulation Study on Synthetic Data}
We carried out  a simulation study on synthetic data to illustrate the performance of Algorithm~\ref{alg:MH} to estimate first the intensity of one dimensional and two dimensional inhomogeneous Poisson processes and finally the intensity of multidimensional Poisson processes.    

 We simulate realizations of Poisson processes on the domain $[0,1)^d$ for $d\in \{1,2,3,4,5\}$ via thinning \citep{lewis1979simulation}.
The hyperparameters of the model (for the trees and the leaf intensities) are fixed as described in Section~\ref{sec:hyperparameters}.
 We initially randomly generate $m$ trees of zero depth. The probabilities of the proposals in Algorithm~\ref{alg2} are set to: $P(\text{GROW})=P(\text{PRUNE})=0.4$ and $P(\text{CHANGE})=0.2$. 
A set $\{z_i\}$ is defined by uniformly sampling points in the domain  $[0,1)^d$. 

We run 3 parallel chains of the same length. We discard their first halves treating the second halves as a sample from the target distribution. We assess chain convergence using the Gelman-Rubin convergence diagnostic~\citep{gelman1992inference} applied to the estimated intensity for each point of the set  $\{z_i\}$ , as well as trace plots and autocorrelation plots for some points of the testing set. 

At each state $t$ of a  simulated chain we estimate the intensity for each point $z_i$ by a product of trees denoted as 
\[\widehat{\lambda}^{(t)}(z_i)=\prod_{j=1}^{m}g(z_i;T^{(t)}_j,\Lambda^{(t)}_{j}).\]
The induced sequence $\{\widehat{\lambda}^{(t)}(\cdot)\}_{t=1}^\infty$ for the sequence of draws $\{(T_1^{(t)},\Lambda_1^{(t)}),..,(T_m^{(t)},\Lambda_m^{(t)})\}_{t=1}^\infty$ converges to $P(\widehat{\lambda}|\vs)$. 
We estimate the posterior mean $E[\widehat{\lambda}(\cdot) |s_1,..s_n] $, the posterior median of $\widehat{\lambda}(\cdot)$, and the highest density interval (hdi) using the function $hdi$ provided by the \textbf{R} package \textbf{bayestestR} \citep{makowski2019bayestestr}. To assess the performance of our algorithm, we compute the Average Absolute Error (AAE) of the computed estimate:
  \begin{align}
      \text{AAE}(\widehat{\lambda})= \frac{1}{N_z} \sum_{i=1}^{N_z}| \widehat{\lambda}(z_i)-\lambda(z_i) |
  \end{align}
  and the Root Integrated Square Error (RISE):
  \begin{align}
     \text{RISE}(\widehat{\lambda})=  \left( \frac{1}{N_z} \sum_{i=1}^{N_z}( \widehat{\lambda}(z_i)-\lambda(z_i) )^2 \right)^{1/2}
  \end{align}
where $N_z$ is the number of test points. 

In the spirit of Akaike information criterion (AIC) \citep{loader1999local}, we also introduce two diagnostics targetting the likelihood function to evaluate if increasing the number of trees leads to better intensity estimation: 
\begin{equation*}
D_g= 2\left( \log P(s_1,..,s_n) - k_g\right),
\end{equation*}and 
\begin{equation*}
D_l= 2\left( \log P(s_1,..,s_n) - k_l\right),
\end{equation*}
where $k_g$ is the number of global cells, and $k_l$ is the overall number of leaves in the ensemble. We estimate both diagnostics using the sequence of the draws $\left(T^{(w)},\Lambda^{(w)}\right)=\left\{\left(T_1^{(w)},\Lambda_1^{(w)}\right),...,\left(T_m^{(w)},\Lambda_m^{(w)}\right) \right \}$ after the burn-in period as 
\begin{equation*}
D_g \approx 2\frac{1}{N_w}\sum_{w=1}^{N_w}\left(\log P\left(s_1,..,s_n|\left(T^{(w)},\Lambda^{(w)}\right)\right)-k_g^{(w)}\right),
\end{equation*} and
\begin{equation*}
D_l \approx 2\frac{1}{N_w}\sum_{w=1}^{N_w}\left(\log P\left(s_1,..,s_n|\left(T^{(w)},\Lambda^{(w)}\right)\right)-k_l^{(w)}\right),
\end{equation*}
where $k_g^{(w)}$ and $k_l^{(w)}$ are the number of global cells and the overall number of leaves in the ensemble associated to  the $w_{th}$ draw, respectively.

AIC has been shown to be asymptomatically equal to leave-one-out cross validation (LOO-CV) \citep{stone1977asymptotic, gelman2014understanding}. According to \citet{leininger2017bayesian},  the computational burden required for leave-one-out cross validation considering a point pattern data is impractical. We introduce a leave-partition-out (LPO) method, assuming that the initial process $N(t)$ is obtained by combining independent processes $\{N_i(t)\}_{i=1}^{N_p}$, as follows 
\begin{equation} \label{LPO}
D_{LPO}=\sum_{i=1}^{N_P}\log P\left(N_i(t)|N(t)-\{N_i(t)\}\right)
\end{equation} where $P\left(N_i(t)|N(t)-\{N_i(t)\}\right)$ is the leave-partition-out predictive intensity given the process $N(t)$
 without the $i_{th}$ partition, $N_i(t)$. We can evaluate \ref{LPO} as follows,
 \begin{equation*}
    D_{LPO}=\sum_{i=1}^{N_P}\log \left( \frac{1}{N_w}\sum_{w=1}^{N_w}  P\left(N_i(t)| \left(T^{(w,i)},\Lambda^{(w,i)}\right)\right)\right)
\end{equation*} where $(T^{(w,i)},\Lambda^{(w,i)})$ is the sequence of draws $\left\{\left(T_1^{(w,i)},\Lambda_1^{(w,i)}\right),...,\left(T_m^{(w,i)},\Lambda_m^{(w,i)}\right) \right \}$ after the burn-in period leaving out the partition $N_i(t)$. 
We assume that each event of $N(t)$ is coming from $N_i(t)$ with probability $p_i$. The bias of the method is introduced by randomly splitting the process into individual processes. We can get the LOO-CV by LPO, defining appropriately the parameter $N_p$. As higher the number $N_p$ is, as less biased the method is. In the simulation scenarios, we consider that $p_i=0.1$, $i=1,...,N_p$ and $N_p=10$ for computational reasons. The diagnostics show that tiny ensembles of trees provide good estimates in our simulation scenarios.

To confirm the proposed diagnostics, we use $p$-thinning \cite[Chapter 6]{illian2008statistical} with $p=0.8$ to create training and test datasets in two of the simulation scenarios. We employ Root Standardized Mean Square Error (RSMSE) and Rank Probability Score (RPS) with the test data set comparing observed counts in disjoint equal volume subregions $\{S_i\}_{i=1}^{N_s}$ as follows:

  \begin{align}
     \text{RSMSE}(\widehat{N})=  \left( \frac{1}{N_s} \sum_{i=1}^{N_s}\frac{(\widehat{N}(S_i)-N(S_i) )^2}{\widehat{N}(S_i)} \right)^{1/2}
  \end{align} and
  
  \begin{align}
     \text{RPS}(N(S_j))= \sum_{u=0}^{N(S_j)-1}F(u)^2 + \sum_{u=N(S_j)}^{\infty} \left( F(u)-1\right)^2 ,  
  \end{align}where $F$ is the Poisson distribution with parameter $m=\int_{S_j}\hat{\lambda}(s)ds$, $N(S_i)$ the actual number of testing points in $S_i$ and $\widehat{N}(S_i)$ the estimated number of testing points in $S_i$ given by  \begin{align}
     \widehat{N}(S_i)=  \int\limits_{S_i} \frac{1-p}{p} \widehat{\lambda}(s)ds \simeq \frac{1-p}{p}\frac{1}{N_{z}^i}\sum\limits_{z_j\in S_i} \widehat{\lambda}(z_j)|S_i|
  \end{align} with $N_{z}^i$ being the number of points $\{z_j\}$ falling in $S_i$ and estimating the intensity at each points $s$ , $\widehat{\lambda}(s)$, via the posterior mean $E[\widehat{\lambda}(\cdot) |s_1,..s_n] $.

For one dimensional processes, we compare the results of Algorithm~\ref{alg:MH} to the Haar-Fisz algorithm~\citep{fryzlewicz2004haar}, a wavelet based method for estimating the intensity of one dimensional Poisson Processes that outperforms well known competitors. We apply the Haar-Fisz algorithm to the counts of points falling into 256 consecutive intervals using the \textbf{R} package \textbf{haarfisz} \citep{haarfisz}. Our algorithm is competitive with the Haar algorithm for smooth intensity functions and is not strongly out-performed by the Haar-Fisz algorithm when the underlying intensity is a stepwise function. 

For two-dimensional processes, we compare the results of our algorithm with fixed-bandwidth estimators and log-Gaussian Cox processes (LGCP) with intensity $\lambda(\vs)=\exp{(a + u(\vs)) }$ where $u$ is a Gaussian process with exponential covariance function. We used a discretization version of the LGCP model defined on a regular grid over space which we implemented using Stan-code \citep{gelman2015stan}. As noted in \citet{davies2018fast}, the choice of the kernel is not of primary importance, we choose a Gaussian kernel for its wide applicability.
In our tables of results, the smoothing bandwidth, sigma, selected using likelihood cross-validation \citep{loader1999local} denoted by (LCV), and we have also included other values of sigma to demonstrate the sensitivity to bandwidth choice. The kernel estimators, and the bandwidth value given by likelihood cross-validation, were computed using the \textbf{R} package \textbf{spatstat} \citep{spatstat}.  Our algorithm outperforms the maximum likelihood approach using linear conditional intensity, as expected. Our algorithm outperforms kernel smoothing and LGCP for stepwise functions and is competitive with them for a smooth intensity.

 Finally, we examine the performance of our algorithm for multidimensional intensities by generating realizations of Poisson Processes on the domain $[0,1)^d$ for $d\in\{3,5\}$ via thinning. Future work includes the study of intensities in higher dimensions ($d>5$). We compare our intensity estimates with kernel smoothing estimators having isotropic standard deviation matrices with diagonal elements equal to $h$ and the methodology for applying maximum likelihood to point process models with linear conditional intensity \citep{peng2003multi}. We select the bandwidth $h$ using likelihood cross-validation \citep{loader1999local} denoted by (LCV).

\subsection{One dimensional Poisson Process with stepwise intensity}
\begin{figure}[H] 
 \begin{subfigure}{8cm}
    \centering\includegraphics[width=8cm]{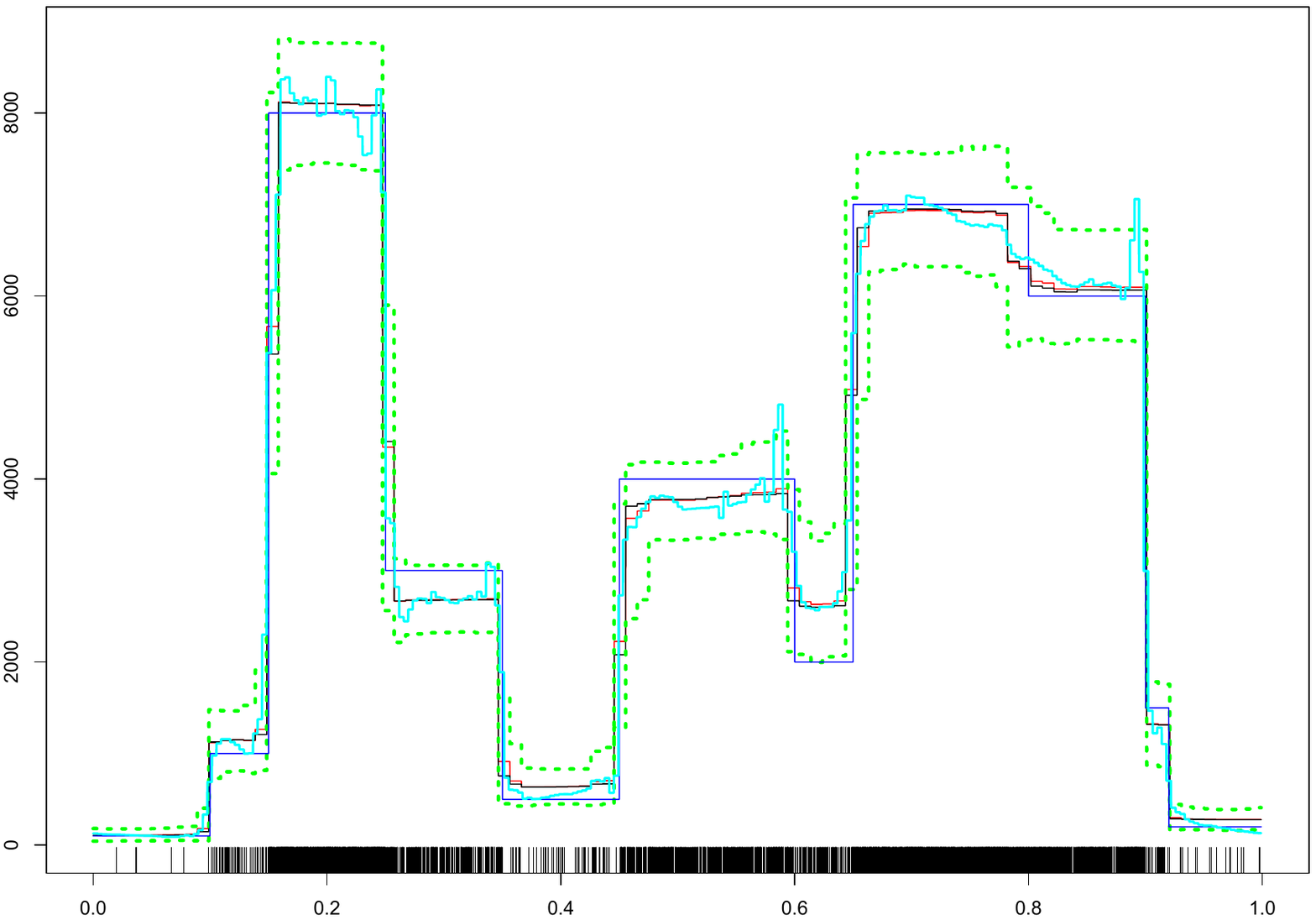}
    \caption{5 Trees }
  \end{subfigure}
\begin{subfigure}{8cm}
    \centering\includegraphics[width=8cm]{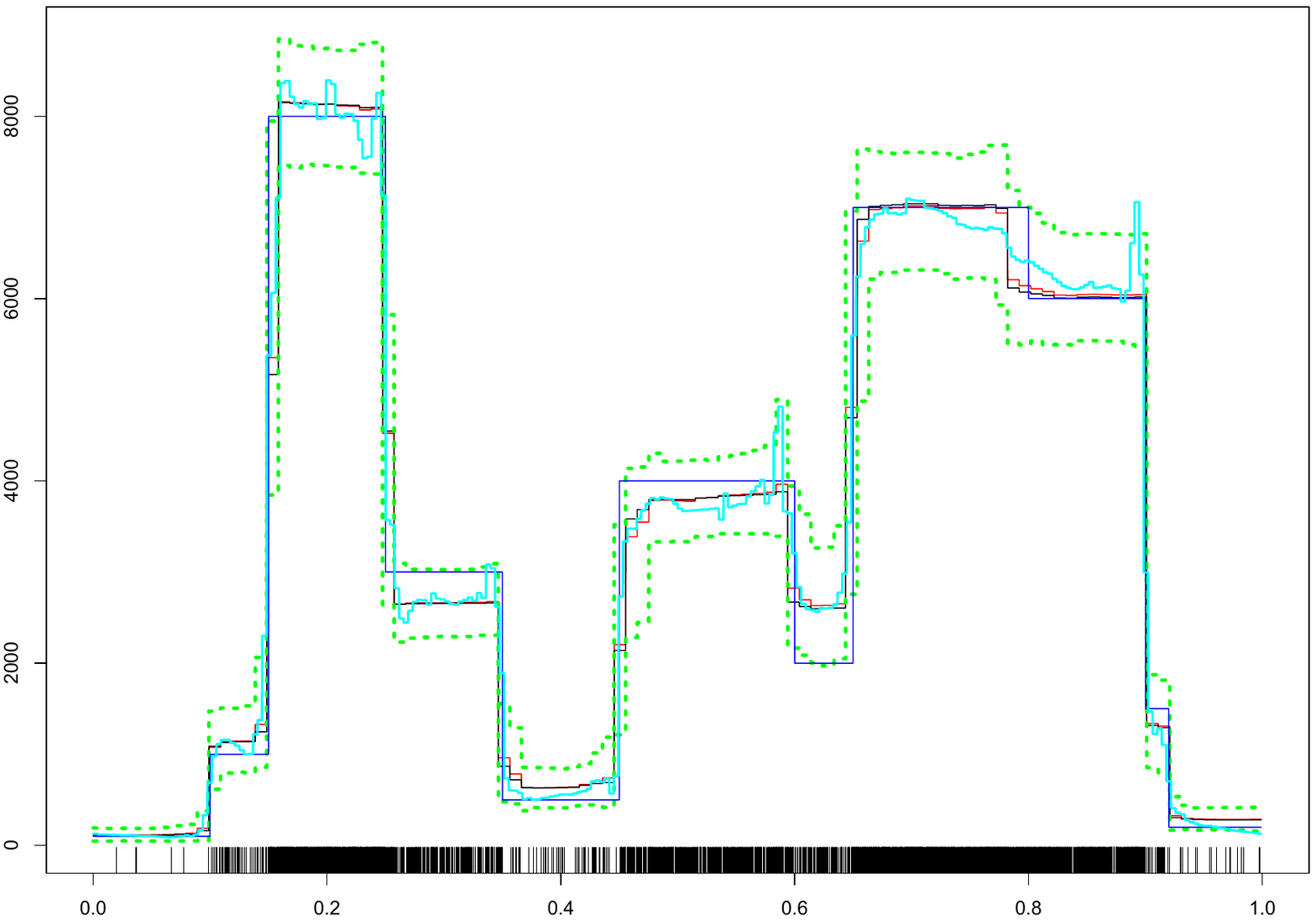}
    \caption{7 Trees }
  \end{subfigure}
\caption{The original intensity (blue curve), the posterior mean (red curve), the posterior median (black curve), the 95\% hdi interval of the estimated intensity illustrated by the dotted green lines and the Haar-Fisz estimator (cyan curve). The rug plot on the bottom displays the 3590 event times. }
\label{fig1}
\end{figure}

 Our first example is a one dimensional Poisson Process with piecewise constant intensity with several steps (Fig.~\ref{fig1}). We run 3 parallel chains of the same length for 200000 iterations for 2-10 trees, 100000 for 12 trees, 50000 iterations for 15 trees and 30000 iterations for 20 trees.  
 
 Our algorithm detects the change points and provides good estimates of the intensity and is competitive in terms of AAE with the Haar-Fisz algorithm, but does not perform as well in terms of RISE (see Fig.~\ref{fig1} and Tables \ref{table1}-\ref{tableH1}). We have found the metrics and convergence diagnostics in a set of uniformly chosen points without excluding the points close to jumps. Due to inferring the intensity via a product of stepwise functions, it is expected that the proposed algorithm will provide estimates with higher variability close to jumps. The proposed algorithm outperforms the Haar-Fisz algorithm without considering the points close to jumps. Tables \ref{table1b}-\ref{tableH1b} show the metrics for various number of trees without considering the points in a distance $=\pm 0.02$ from the jumps.

The diagnostics $D_g$, $D_l$ and $D_{LPO}$ obtain their highest values for 7, 4 and 8 trees, respectively. 
 The analysis demonstrates only small differences between log-likelihood values as the number of trees increases, supporting results found in previous BART studies that the method is robust to the choice of $m$.  The average RSMSE and RPS on testing points over 7 different splits of the original data set (Tables \ref{table41}-\ref{table41b}) provide evidence that ensembles with more than seven trees do not improve the fit of the proposed algorithm.

\begin{table}[H]
\begin{tabular}{ p{1.5 cm} p{1.5 cm} p{1.5 cm}  p{1.5 cm}  p{1.5 cm}  p{1.5 cm} p{1.5 cm} }
 \hline
 \multicolumn{7}{c}{Proposed BART Algorithm}  \\
 \hline
 Number of trees & $N_s=1$ & $N_s=10$ & $N_s=25$ &  $N_s=50$ & $N_s=75$ &  $N_s=75$ \\
2          & 17.05 &  5.24 &  3.10 &  2.08 &  1.66 &  1.44 \\ 
3          & 17.12 &  5.25 &  3.11 &  2.08 &  1.65 &  1.43 \\ 
4          & 16.98 &  5.28 &  3.09 &  2.07 &  1.64 &  1.42 \\ 
5          & 16.94 &  5.26 &  3.04 &  2.06 &  1.63 &  1.41 \\ 
7          & 17.01 &  5.22 &  3.00 &  2.04 &  1.62 &  1.40 \\ 
8          & 17.09 &  5.20 &  2.99 &  2.03 &  1.62 &  1.40 \\ 
9          & 17.07 &  5.20 &  2.98 &  2.02 &  1.61 &  1.39 \\ 
10          & 20.10 &  5.78 &  3.12 &  2.12 &  1.63 &  1.43 \\ 
12          & 17.74 &  5.37 &  2.99 &  2.06 &  1.63 &  1.39 \\ 
15          & 17.03 & 5.15 &  2.94 &  2.01 & 1.60 & 1.38 \\ 
20          & 17.08 &  5.16 & 2.93 & 2.00 &  1.60 &  1.38 \\ 
 \hline
\end{tabular}
\caption{The average RPS on testing points over 7 different splits of the original data set in Fig.~\ref{fig1}.}
\label{table41}
\end{table}

\begin{table}[H]
\begin{tabular}{ p{1.5 cm} p{1.5 cm} p{1.5 cm}  p{1.5 cm}  p{1.5 cm}  p{1.5 cm} p{1.5 cm} }
 \hline
 \multicolumn{7}{c}{Proposed BART Algorithm}  \\
 \hline
 Number of trees & $N_s=1$ & $N_s=10$ & $N_s=25$ &  $N_s=50$ & $N_s=75$ &  $N_s=75$ \\
 2          & 0.95 & 1.13 & 1.06 & 1.02 & 0.99 & 1.00 \\ 
3          & 0.95 & 1.13 & 1.06 & 1.02 & 0.98 & 0.99 \\ 
4          & 0.95 & 1.14 & 1.06 & 1.02 & 0.98 & 0.98 \\ 
5          & 0.94 & 1.13 & 1.04 & 1.02 & 0.98 & 0.98 \\ 
7          & 0.95 & 1.13 & 1.04 & 1.01 &0.97 & 0.97 \\ 
8          & 0.95 & 1.12 & 1.03 & 1.01 & 0.97 & 0.96 \\ 
9          & 0.95 & 1.13 & 1.03 & 1.01 & 0.97 & 0.97 \\ 
10          & 1.10 & 1.20 & 1.06 & 1.04 & 0.98 & 0.98 \\ 
12          & 0.98 & 1.18 & 1.04 & 1.02 & 0.98 & 0.96 \\ 
15          & 0.95 & 1.12 & 1.02 & 1.00 & 0.97 & 0.96 \\ 
20          & 0.95 & 1.12 & 1.02 & 1.00 & 0.97 & 0.96 \\ 
 \hline
\end{tabular}
\caption{The average RSMSE on testing points over 7 different splits of the original data set in Fig.~\ref{fig1}.}
\label{table41b}
\end{table}

\begin{table}[H]
\begin{tabular}{ p{1.5cm} p{2cm} p{2cm} p{2cm} p{2cm} p{1.5cm} p{1.5cm} p{1.5cm}}
\hline
 \multicolumn{8}{c}{Proposed BART Algorithm} \\ 
 \hline
 Number of  trees & AAE for Posterior Mean & AAE for Posterior Median & RISE for Posterior Mean & RISE for Posterior Median & $D_g$& $D_l$ & $D_{LPO}$\\
 \hline
 3 & 308.87 & 320.84 & 603.54 & 633.48 & 54095.1 & 54090 & -339.7 \\
 4 & 287.89 & 283.03 & 580.69 & 587.08 & 54096.5 & 54090 & -368\\ 
 5 & 289.27 &281.13 &580.55 &586.24 & 54098 & 54088.4 & -352.5\\
7 &281.59 & 274.88 & 588.7 &592.11 & 54098 & 54082.7 & -263.5\\
8 &280.62 & 274.07 & 588.73 &591.29 & 54097.9 & 54079.5 & -261.5\\
9 & 282.78 & 276.99 & 593.93 & 595.23 & 54096.9 & 54075.2 & -327.9 \\
10 & 283.79 & 279.07 & 593.95 & 595.41 & 54095.6 & 54071.6 & -322.6\\
20 & 297.21 & 287.86 & 599.77 & 595.04 & 54082.9 & 54029.7 & -436 \\
 \hline
\end{tabular}
\caption{Average Absolute Error and Root Integrated Square Error for  various number of trees for the data in Fig.~\ref{fig1}.}
\label{table1}
\end{table}

\begin{table}[H]
\begin{tabular}{ p{2cm} p{2cm} p{2cm} p{2cm} p{2cm}  }
\hline
 \multicolumn{5}{c}{Proposed BART Algorithm} \\ 
 \hline
 Number of  trees & AAE for Posterior Mean & AAE for Posterior Median & RISE for Posterior Mean & RISE for Posterior Median\\
 \hline
 4 & 144.48 & 139.58 & 181.21 & 174.82\\
5 & 144.55 & 139.02 & 180.74 & 176.19\\
7 & 124.53 & 123.2 & 175.74 & 172.4 \\
 \hline
\end{tabular}
\caption{Average Absolute Error and Root Integrated Square Error for the data in Fig.~\ref{fig1} without considering points close to steps.}
\label{table1b}
\end{table}

\begin{table}[H]
\begin{tabular}{ p{3cm}p{3cm}  }
 \hline
 \multicolumn{2}{c}{Haar-Fisz Algorithm}  \\
 \hline
 AAE  &  RISE   \\
 \hline
141.95 & 192.6   \\
 \hline
\end{tabular}
\caption{Average Absolute Error and Root Integrated Square Error for Haar-Fisz estimator for the data in Fig.~\ref{fig1} without considering points close to steps.}
\label{tableH1b}
\end{table}

\begin{table}[H]
\begin{tabular}{ p{3cm}p{3cm}  }
 \hline
 \multicolumn{2}{c}{Haar-Fisz Algorithm}  \\
 \hline
 AAE  &  RISE   \\
 \hline
272.3 & 476.9   \\
 \hline
\end{tabular}
\caption{Average Absolute Error and Root Integrated Square Error for Haar-Fisz estimator for the data in Fig.~\ref{fig1}.}
\label{tableH1}
\end{table}

\subsection{Two-dimensional Poisson process with stepwise intensity function}
\begin{figure}[H] 
  \begin{subfigure}{8 cm}
    \centering\includegraphics[width=8cm]{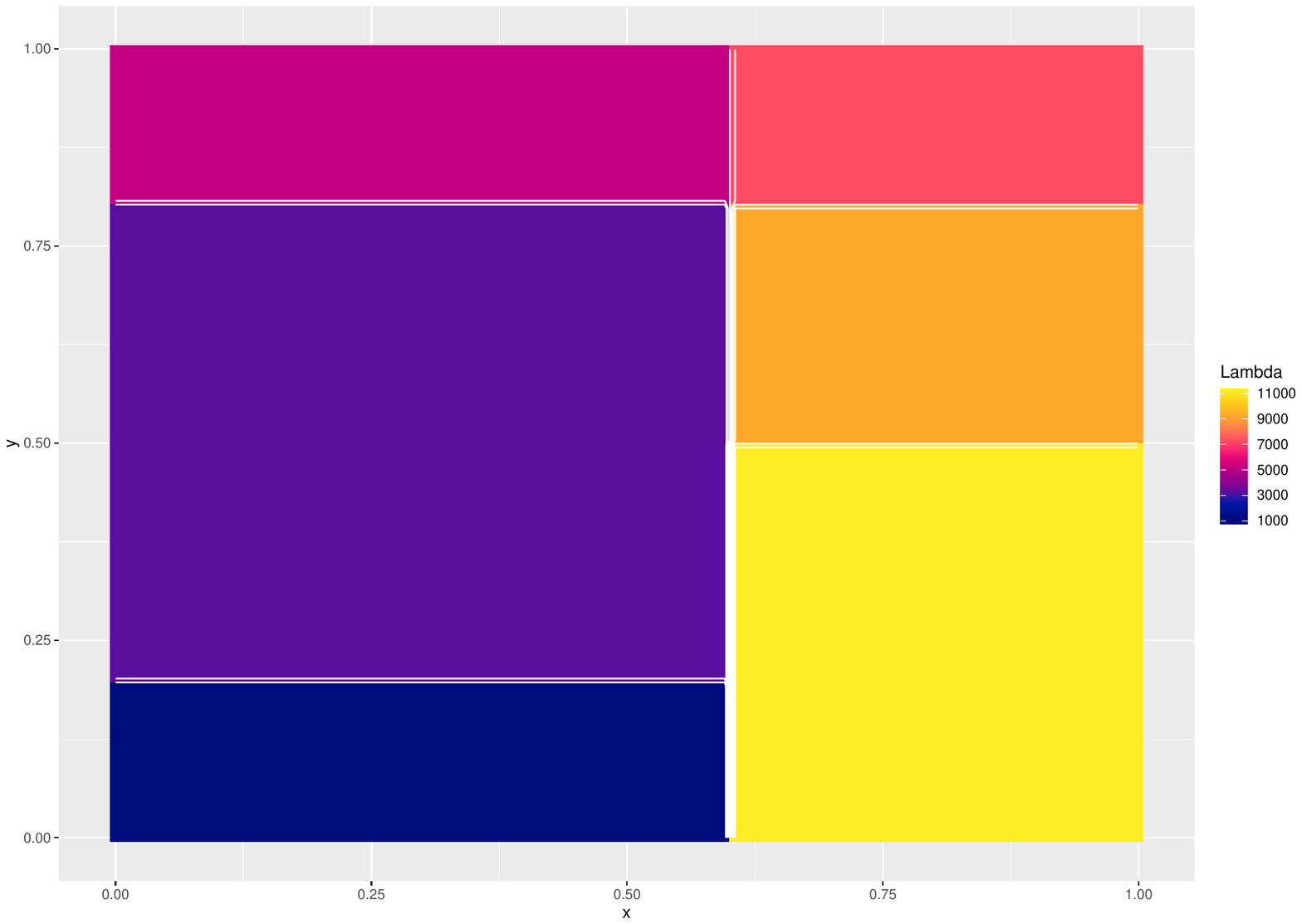}
    \caption{ Original Intensity}
  \end{subfigure}
 \begin{subfigure}{8cm}
    \centering\includegraphics[width=8cm]{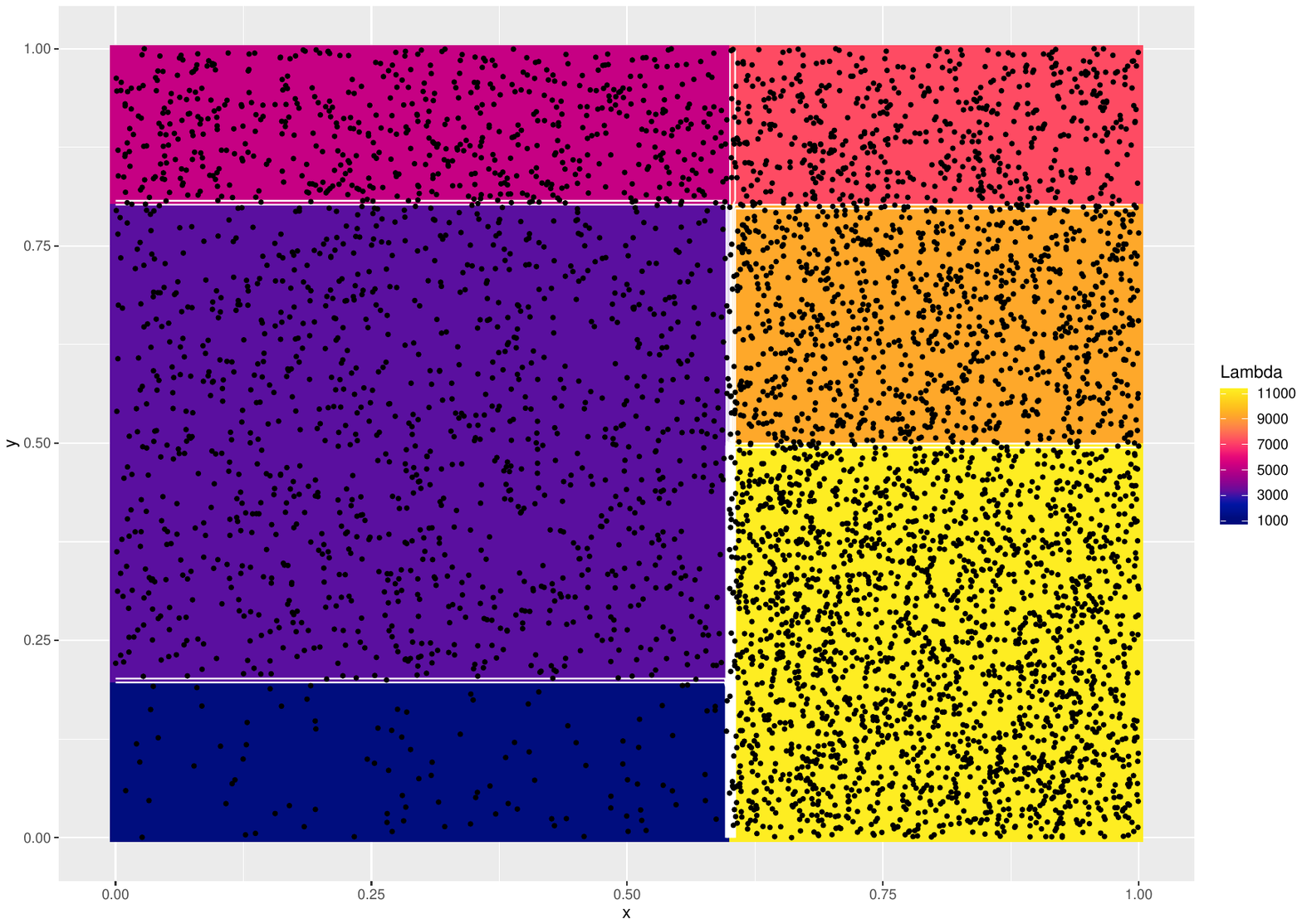}
    \caption{Realization of Process consisting of 5579 points}
  \end{subfigure}
 \begin{subfigure}{8 cm}
    \centering\includegraphics[width=8cm]{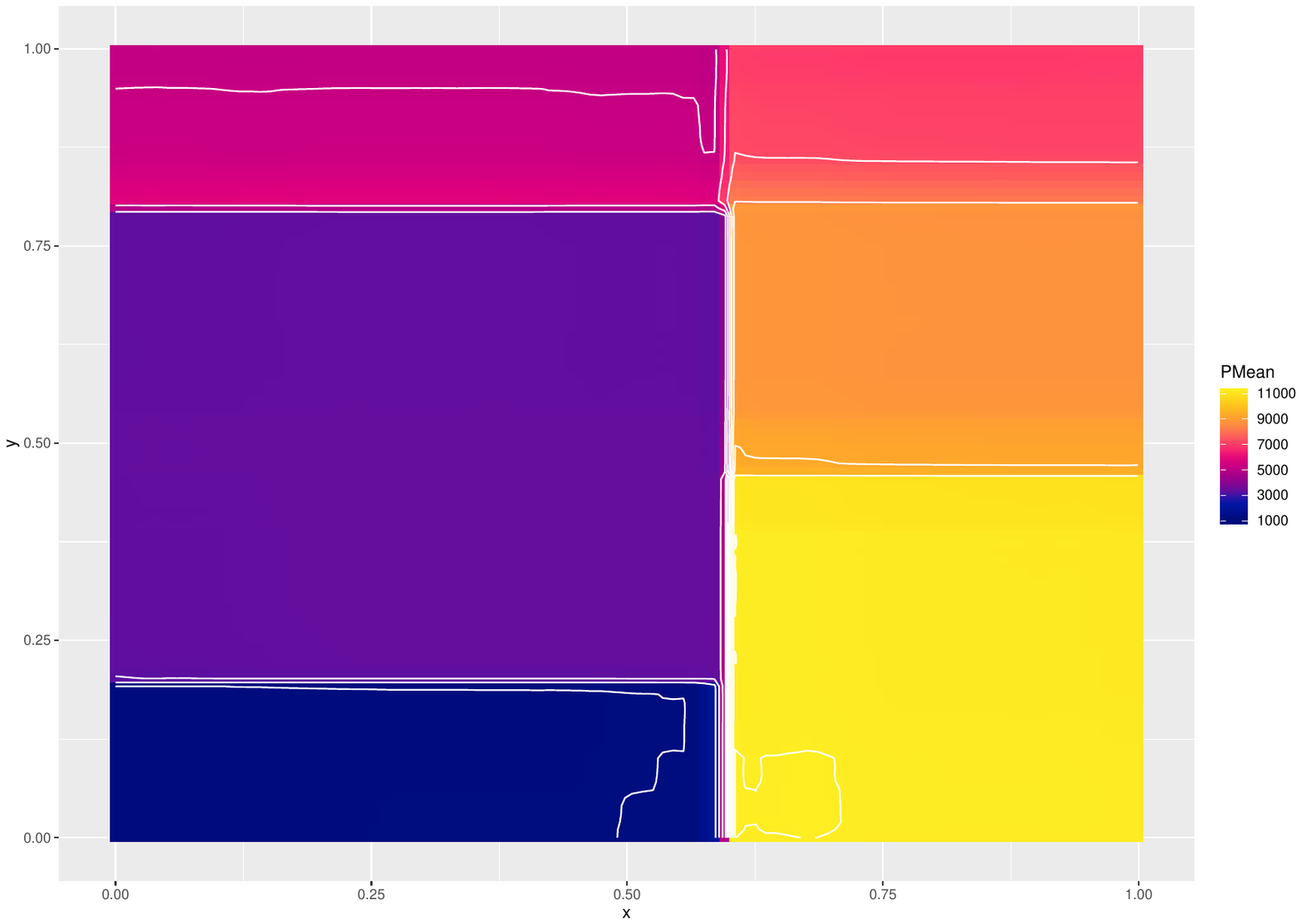}
    \caption{ Posterior Mean}
  \end{subfigure}
 \begin{subfigure}{8 cm}
    \centering\includegraphics[width=8cm]{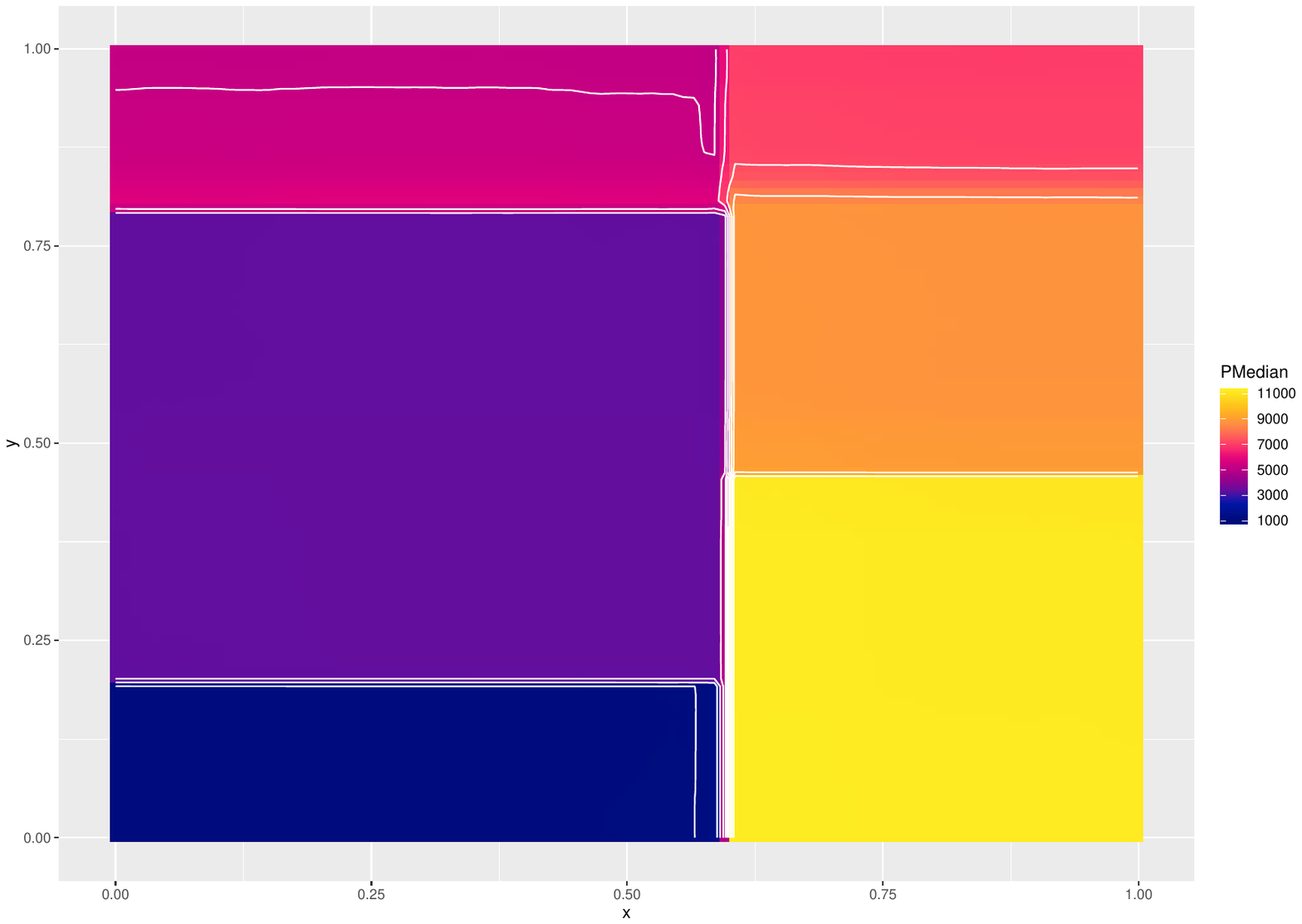}
    \caption{ Posterior Median}
  \end{subfigure}
\caption{Original Intensity, posterior mean and posterior median for 4 trees.}
\label{fig3}
\end{figure}

\begin{figure}[H] 
  \begin{subfigure}{8 cm}
    \centering\includegraphics[width=8cm]{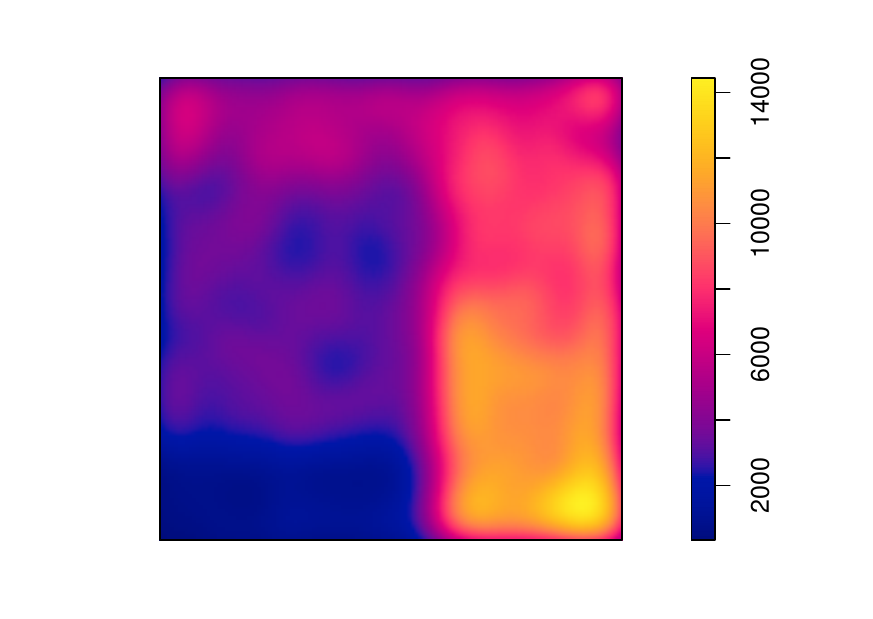}
    \caption{Kernel Smoothing with $h=0.039$.}
  \end{subfigure}
 \begin{subfigure}{8cm}
    \centering\includegraphics[width=8cm]{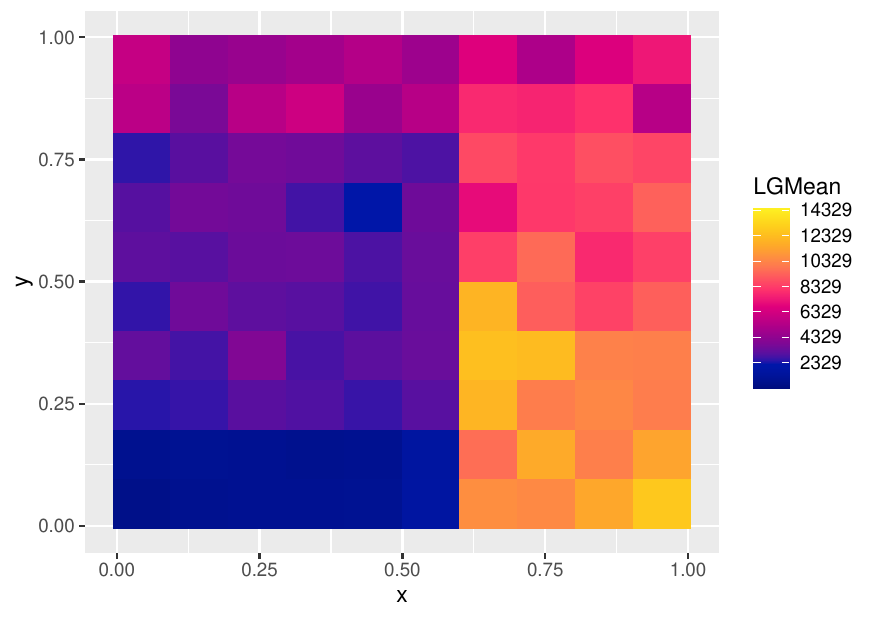}
    \caption{LGCP}
  \end{subfigure}
\caption{Kernel estimator and inference with spatial log-Gaussian Cox processes.}
\label{fig3l}
\end{figure}

To demonstrate the applicability of our algorithm in a two-dimensional setting, Figures \ref{fig3}-\ref{fig3l} and Tables \ref{table3}-\ref{table4l} reveal that our algorithm outperforms kernel smoothing and inference with spatial log-Gaussian Cox processes for stepwise intensity functions. We run 3 parallel chains of the same length for 100000 iterations for 3-6 trees. The convergence criteria indicate convergence of the simulated chains for the majority of points. As may be expected, the simulation study shows that points close to jumps are estimated with less reliability.  The algorithm converges less well at these points, as demonstrated by the Gelman-Rubin diagnostic (see supplementary material). The diagnostics $D_g$, $D_l$ and $D_{LPO}$ obtain their highest values for three trees, respectively. The diagnostics indicate that small ensembles of trees can provide a good estimate of the intensity.

\begin{table}[H]
\begin{tabular}{ p{2 cm} p{2 cm} p{2 cm} p{2 cm} p{2 cm} p{1.5 cm} p{1.5 cm} p{1.5 cm}}
 \hline
 \multicolumn{8}{c}{Proposed BART Algorithm}  \\
 \hline
 Number of  trees & AAE for Posterior Mean & AAE for Posterior Median &RISE for Posterior Mean & RISE for Posterior Median & $D_g$ & $D_l$ &$D_{LPO}$  \\
 \hline
 3 & 224.1 & 230.3 & 419.2 & 453.7 & 87227.2 & 87232.3 & 505.1\\ 
 4 & 208.7  &213  &410.2 &447.9 & 87223.7 & 87230.5  & 491.2\\
 5& 216.8  &212.9  &389.5 &410.9 & 87211.6 & 87220.6  & 406\\
 6 & 228.9  &221.9  &395.8 &412.8 & 87197.5 & 87214.7  &463.9\\
 \hline
\end{tabular}
\caption{Average Absolute Error, Root Integrated Square Error and diagnostics for various trees for the data in Figure~\ref{fig3}.}
\label{table3}
\end{table}

\begin{table}[H]
\begin{tabular}{ p{4cm} p{3cm} p{3cm}  }
 \hline
 \multicolumn{3}{c}{Kernel Smoothing}\\ 
 \hline
 Bandwidth (sigma) & AAE  &  RISE  \\
 \hline
0.027 &   763.8& 1041.3 \\ 
0.038&  662.7&956.8  \\ 
0.047 (LCV) & 636.7 & 960.6\\ 
0.067 & 672.8& 1042.5\\ 
 \hline
\end{tabular}
\caption{Average Absolute Error and Root Integrated Square Error for fixed bandwidth estimators for the data in Figure~\ref{fig3}.}
\label{table4}
\end{table}

\begin{table}[H]
\begin{tabular}{ p{4cm} p{3cm} p{3cm}  }
 \hline
 \multicolumn{3}{c}{Inference with spatial log-Gaussian Cox processes}\\ 
 \hline
 grid & AAE  &  RISE  \\
 \hline
$10\times 10$ &568& 751\\
 $20\times 20$ & 678& 953 \\
 \hline
\end{tabular}
\caption{Average Absolute Error and Root Integrated Square Error with LGCP for the data in Figure~\ref{fig3}.}
\label{table4l}
\end{table} 

\subsection{Inhomogeneous three-dimensional Poisson Process with Gaussian intensity} 

\begin{figure}[H] 
  \begin{subfigure}{6cm}
    \centering\includegraphics[width=6cm]{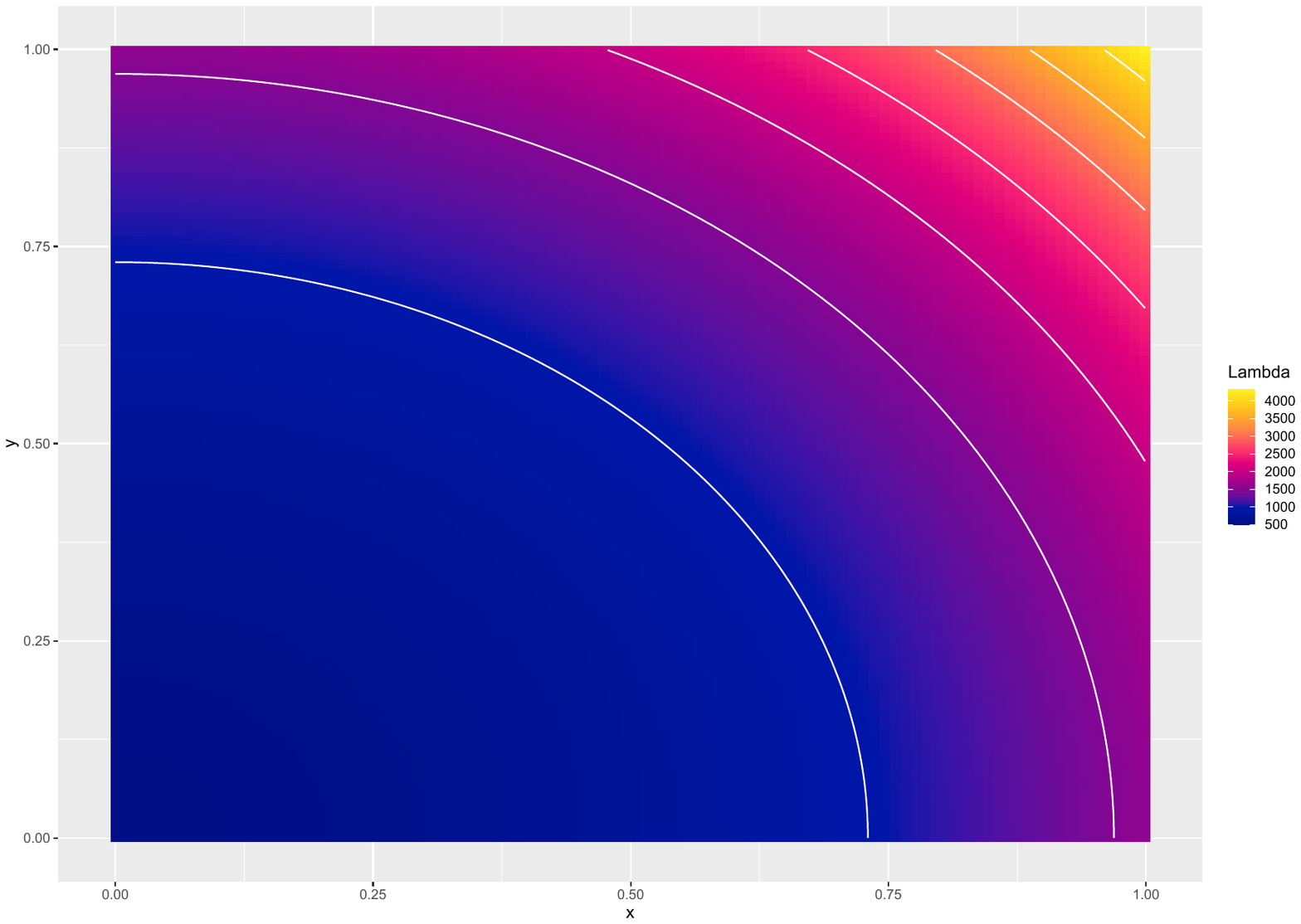}
    \caption{Original Intensity}
  \end{subfigure}
  \begin{subfigure}{6cm}
    \centering\includegraphics[width=6cm]{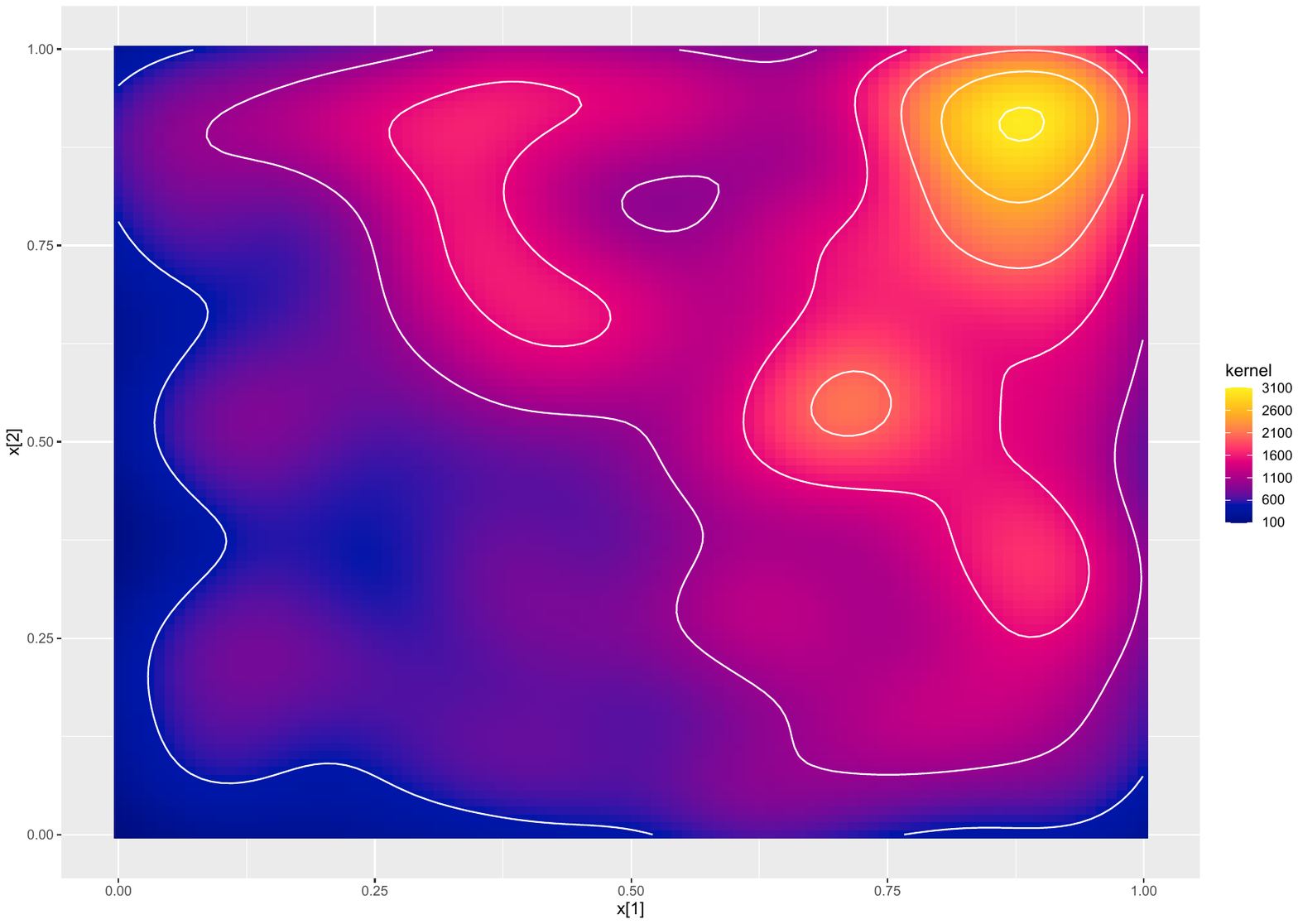}
    \caption{Kernel Estimator with $h$=0.073}
  \end{subfigure}
  
  \begin{subfigure}{6cm}
    \centering\includegraphics[width=6cm]{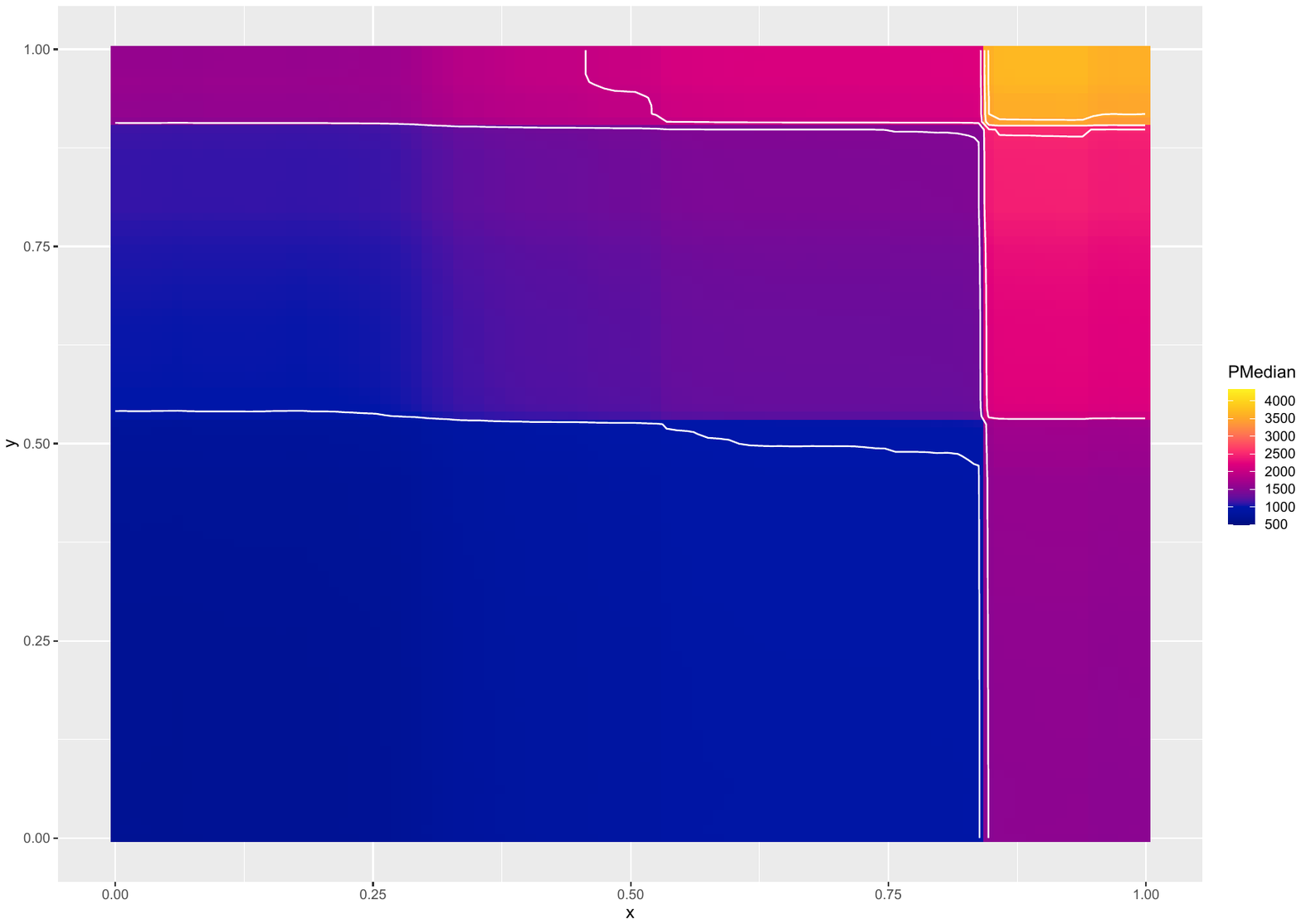}
    \caption{Posterior Median for 8 Trees}
  \end{subfigure}
  \begin{subfigure}{6cm}
    \centering\includegraphics[width=6cm]{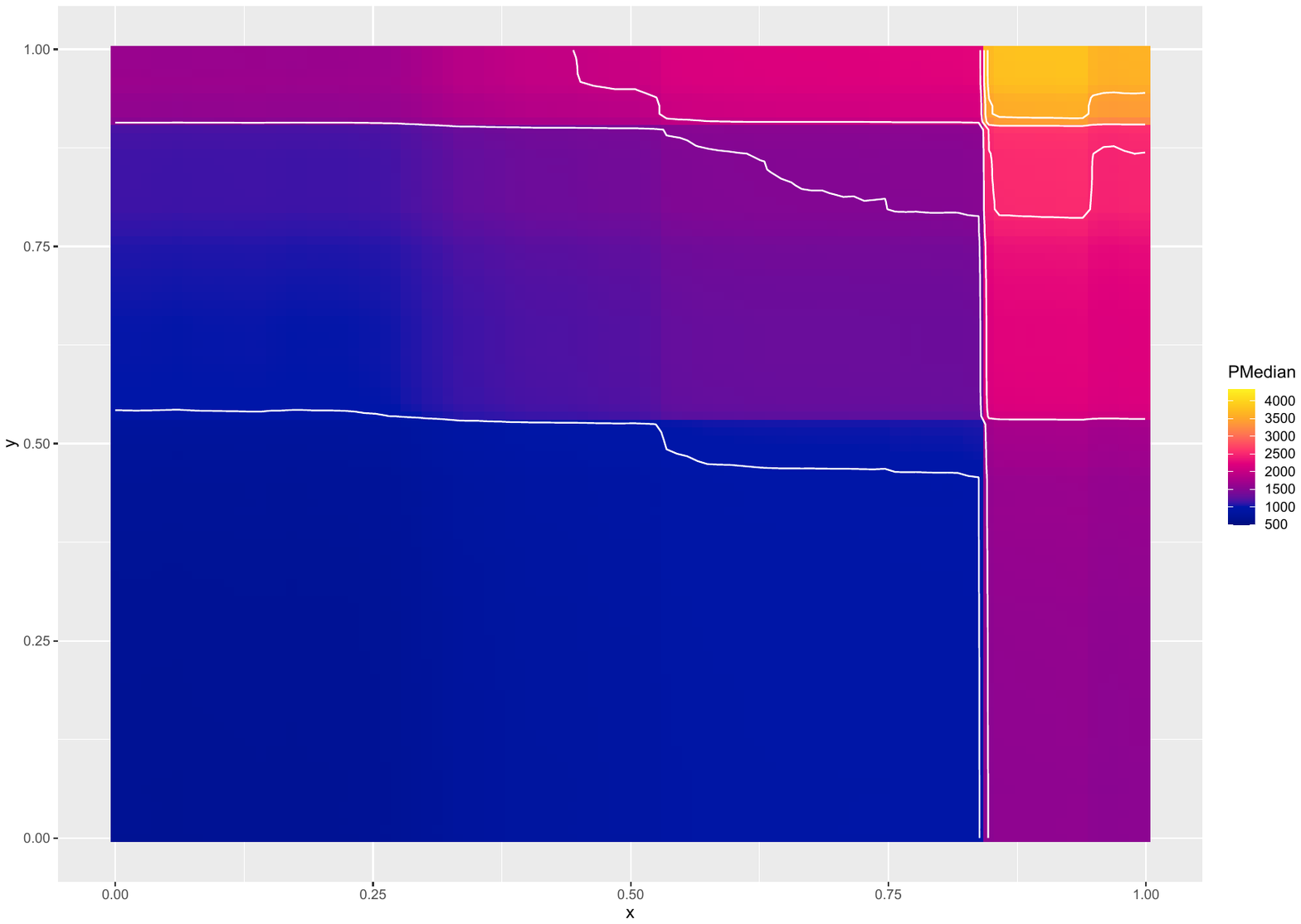}
    \caption{Posterior Median for 10 Trees}
  \end{subfigure}

\caption{Kernel estimator and Posterior Median for 8 and 10 Trees with $\text{x}[3]=0.4$.}
\label{figMD1}
\end{figure}

\begin{figure}[H] 
  \begin{subfigure}{6cm}
    \centering\includegraphics[width=6cm]{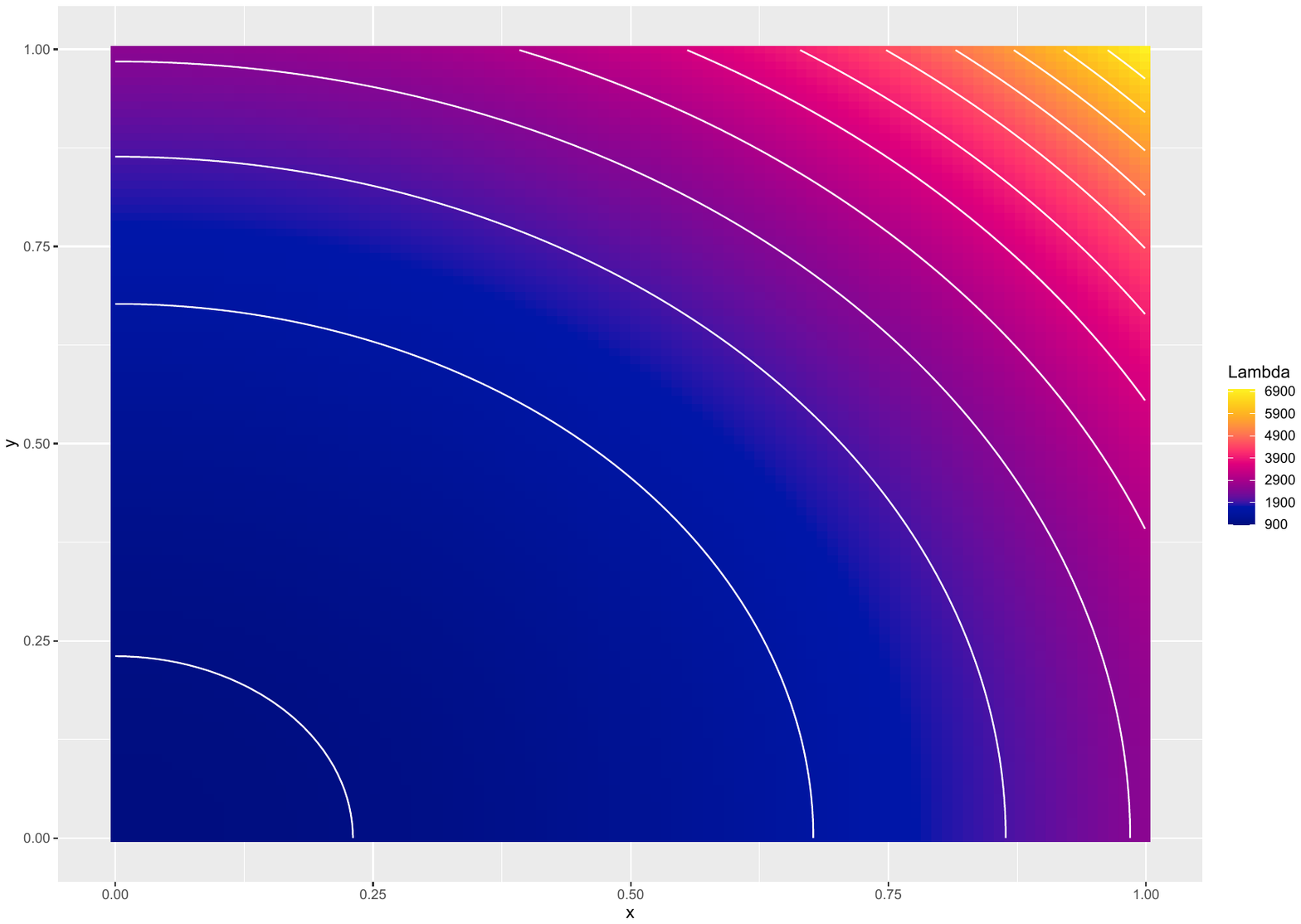}
    \caption{Original Intensity}
  \end{subfigure}
 \begin{subfigure}{6cm}
    \centering\includegraphics[width=6cm]{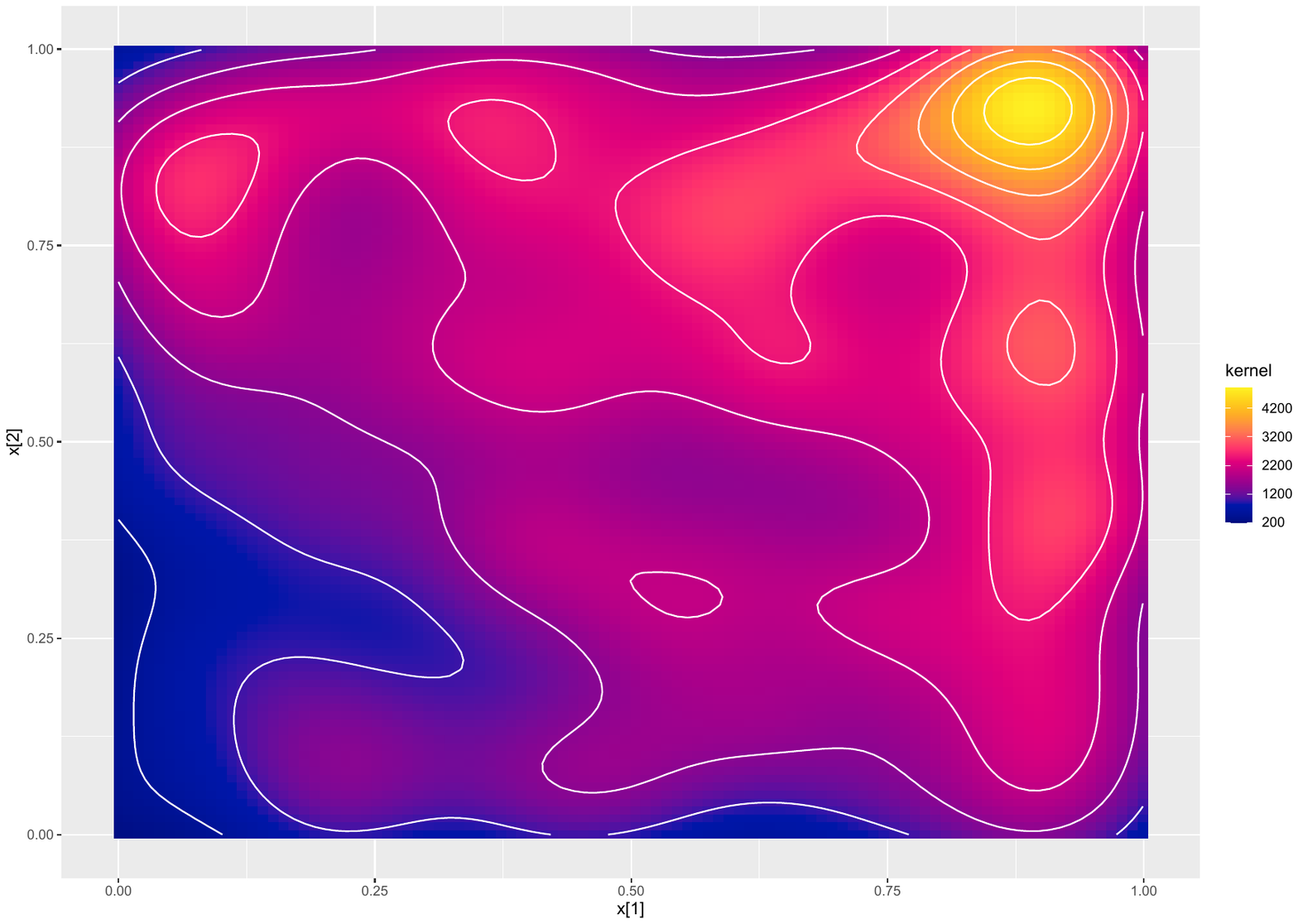}
    \caption{Kernel Estimator with $h$=0.073}
  \end{subfigure}

  \begin{subfigure}{6cm}
    \centering\includegraphics[width=6cm]{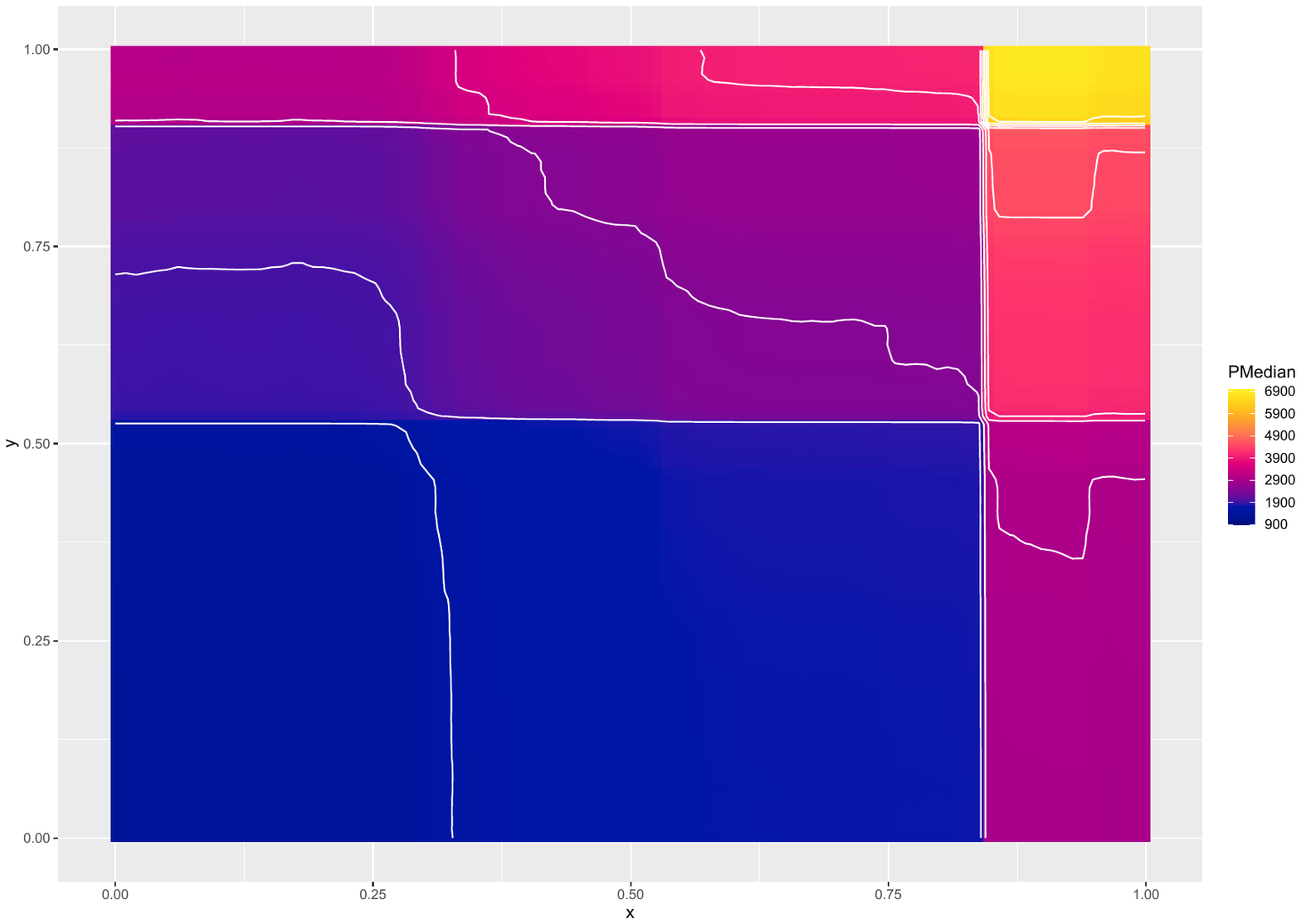}
    \caption{Posterior Median for 8 Trees}
  \end{subfigure}
  \begin{subfigure}{6cm}
    \centering\includegraphics[width=6cm]{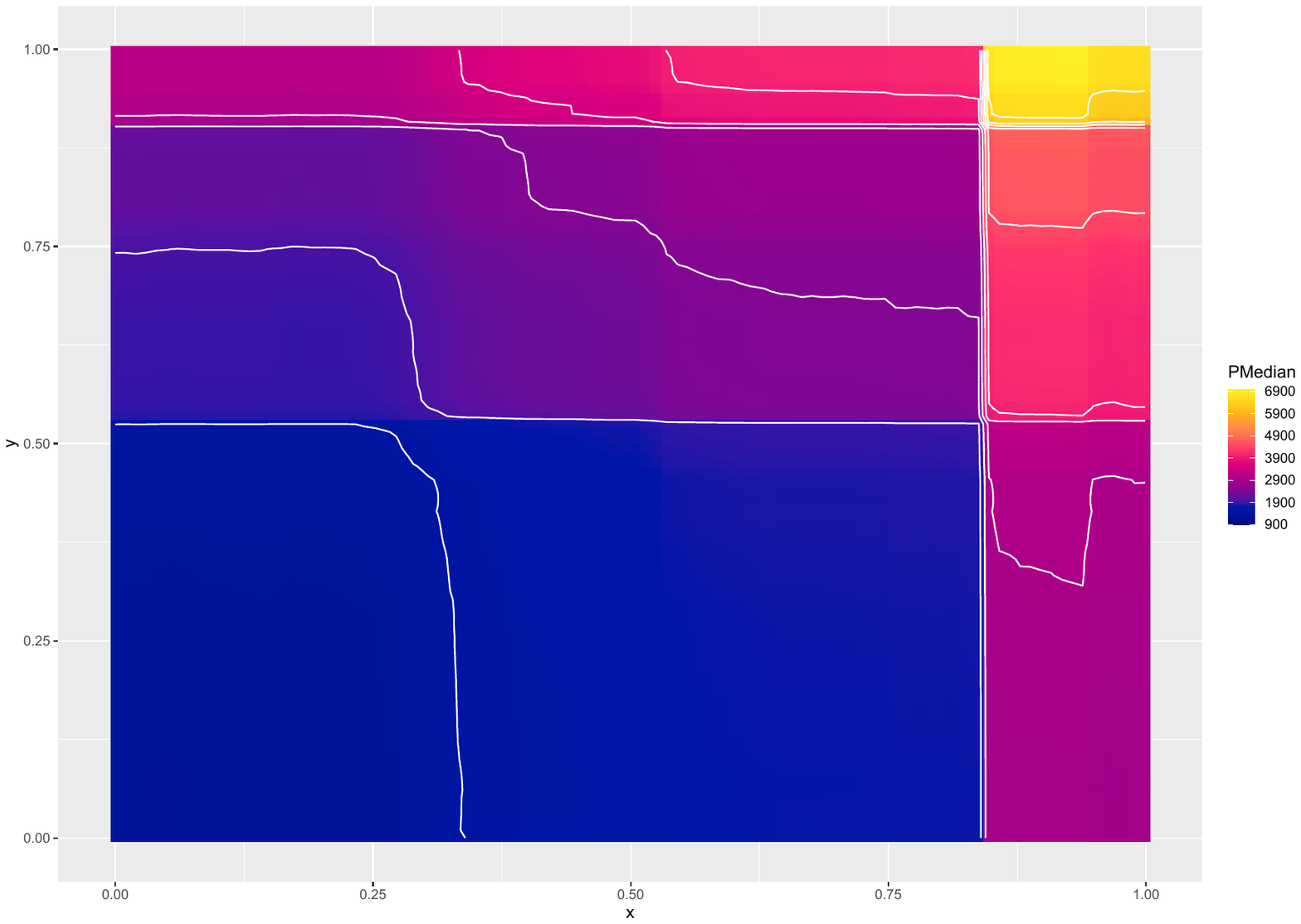}
    \caption{Posterior Median for 10 Trees}
  \end{subfigure}

\caption{Kernel Estimator and Posterior Median for 8 and 10 Trees with $\text{x}[3]=0.8$.}
\label{figMD2}
\end{figure}

Our first example for multidimensional intensities is a three-dimensional Poisson process with intensity $\lambda(x)=500e^{\text{x}^T\text{x}}$ for $\text{x}\in[0,1)^3$. We generated a realization of 1616 points via thinning. We run 3 parallel chains of the same length for 100000 iterations for 3-10 trees and 30000 iterations for 12 trees. Tables~\ref{Table_MD1}~and~\ref{Table_MD2} illustrate the statistics of our algorithm and kernel smoothing. Figures~\ref{figMD1}~and~\ref{figMD2} show our estimators and the kernel estimator with $h$=0.073 for 8 Trees and 10 Trees with fixed third dimension ($\text{x}[3]$) at 0.4 and 0.8, respectively. 

The diagnostics $D_g$, $D_l$ and $D_{LPO}$ get their highest values with
4 trees, respectively. We observe that the diagnostic
$D_l$ slightly differs between 4 and 8 trees. The diagnostic $D_g$ is similar between 4 and 5 trees. The estimate of the average logarithm of Poisson process likelihood does not change significantly from 4 trees to 12 trees. Specifically, we observe its maximum equal to 10536.3 at 12 trees, while its minimum to 10531.9 at 4 trees. In addition, the estimated average number of leaves in a tree of an ensemble is about 3 for $4-12$ trees. That explains why we observe higher values of diagnostics for a small number of trees. The metrics $AAE$ and $RISE$ are optimised with 12 trees. However, it should be noted that only small variations in the metrics are seen between 4 and 12 trees. The diagnostics provide evidence that increasing the number of trees does not improve the fit of the proposed model.

\begin{table}[H]
\begin{tabular}{ p{1.2cm}p{1.6cm}p{1.6cm}p{1.8cm}p{1.8cm}p{1.6cm}p{1.6cm} p{1.6cm}}
 \hline
 \multicolumn{8} {c} {Proposed BART Algorithm}  \\
 \hline
 Number of  trees & AAE for Mean & AAE for Median & RISE for Mean & RISE for Median & $D_g$ & $D_l$ & $D_{LPO}$ \\
 \hline
 4 & 247.6 & 254.9 & 360.7 & 376.3 & 20993.7 & 21040.6 &-1409 \\
 5& 250.2 & 258.3 & 364.3 & 380.4 & 20992.1 & 21039.8 &  -1492 \\
 6 & 247.6 & 254.7 & 360.8 & 375.4 & 20979.5 & 21038.6 & -1529\\ 
 8 & 234.8 & 239.4 & 341 & 352.4 & 20938.6 & 21032.8 & -1515\\
 10 & 226.8 & 229.4 &330.5 & 338.5 &20883 & 21026.9 & -1539\\
 12 & 221.6 & 222.3 & 320.6 & 326.4 & 20810.4 & 21020.7 & -1609\\
 \hline
\end{tabular}
\caption{Average Absolute Error, Root Integrated Square Error and diagnostics for  various number of trees.}
\label{Table_MD1}
\end{table}

\begin{table}[H]
\begin{tabular}{ p{3.2cm}p{2.2cm}p{2.2cm} }
 \hline
 \multicolumn{3} {c} {Kernel Smoothing}  \\
 \hline
 $h$ & AAE &  RISE \\
 \hline
 0.053 & 480.8 & 667.5 \\
 \hline
  0.073 (LCV) & 415.86 & 645.16       \\
\hline
0.08 &417.7 & 661.4   \\
\hline
0.085 &423.2 & 676.2  \\
\hline
0.1 &450.3 & 727.6  \\
\hline
0.3 &890.4 & 1236  \\
\hline
\end{tabular}
\caption{Average Absolute Error and Root Integrated Square Error for  various isotropic variance matrices.}
\label{Table_MD2}
\end{table}

\subsection{Inhomogeneous five dimensional Poisson Process with sparsity assumption} 

Here, we demonstrate the performance of our algorithm to detect the dimensions that contribute most in the intensity of $s\in \mathcal{S}$ in a noisy environment. Consider a five dimensional inhomogeneous Poisson process with intensity function of $\mathbf{x}=\left(x_1,x_2,x_3,x_4,x_5\right) \in [0,1)^5$ depending on 3 of 5 dimensions:  
\begin{align*}
    \lambda(\mathbf{x}) &= \left(2\mathbbm{1}\left(x_1<0.2\right) + 10 \mathbbm{1}\left(x_1\geq0.2\right)\right) * \left(3 \mathbbm{1}\left(x_2< 0.5\right) + 15\mathbbm{1}\left(x_2\geq0.5\right)\right) \\&* \left(3 \mathbbm{1}\left(x_3< 0.8\right) + 30\mathbbm{1}\left(x_3\geq0.8\right)\right)
\end{align*}We generate a realization of 669 points via thinning. We run 3 parallel chains of the same length for 100000 iterations for 4-8 trees, 50000 iterations for 10 trees, 30000 iterations for 12 trees and 10000 iterations for 15 trees. The convergence criterion is smaller than 1.1 for the majority of testing points. 

Table \ref{tableC5D} shows the metrics and diagnostics $D_g$ and $D_l$ of the estimated intensity over various numbers of trees. The diagnostics $D_g$ and $D_l$ obtain their highest values with 4 trees, and the diagnostic $D_l$ shows only small differences between 4 and 5 trees. We note that (i) the average number of leaves in a tree of the ensemble is about 2.2 for 4-5 trees, and (ii) the estimated logarithm of Poisson process likelihood for 4 and 5 trees are 4271.5 and 4271.8, respectively. The diagnostic $D_{LPO}$ gets its highest value with 5 trees. The $p$-thinning approach confirms the diagnostics, and indicates that increasing the number of trees does not improve the fit of the proposed model to the data.

\begin{figure}[tbp] 
  \begin{subfigure}{8cm}
    \centering\includegraphics[width=6cm]{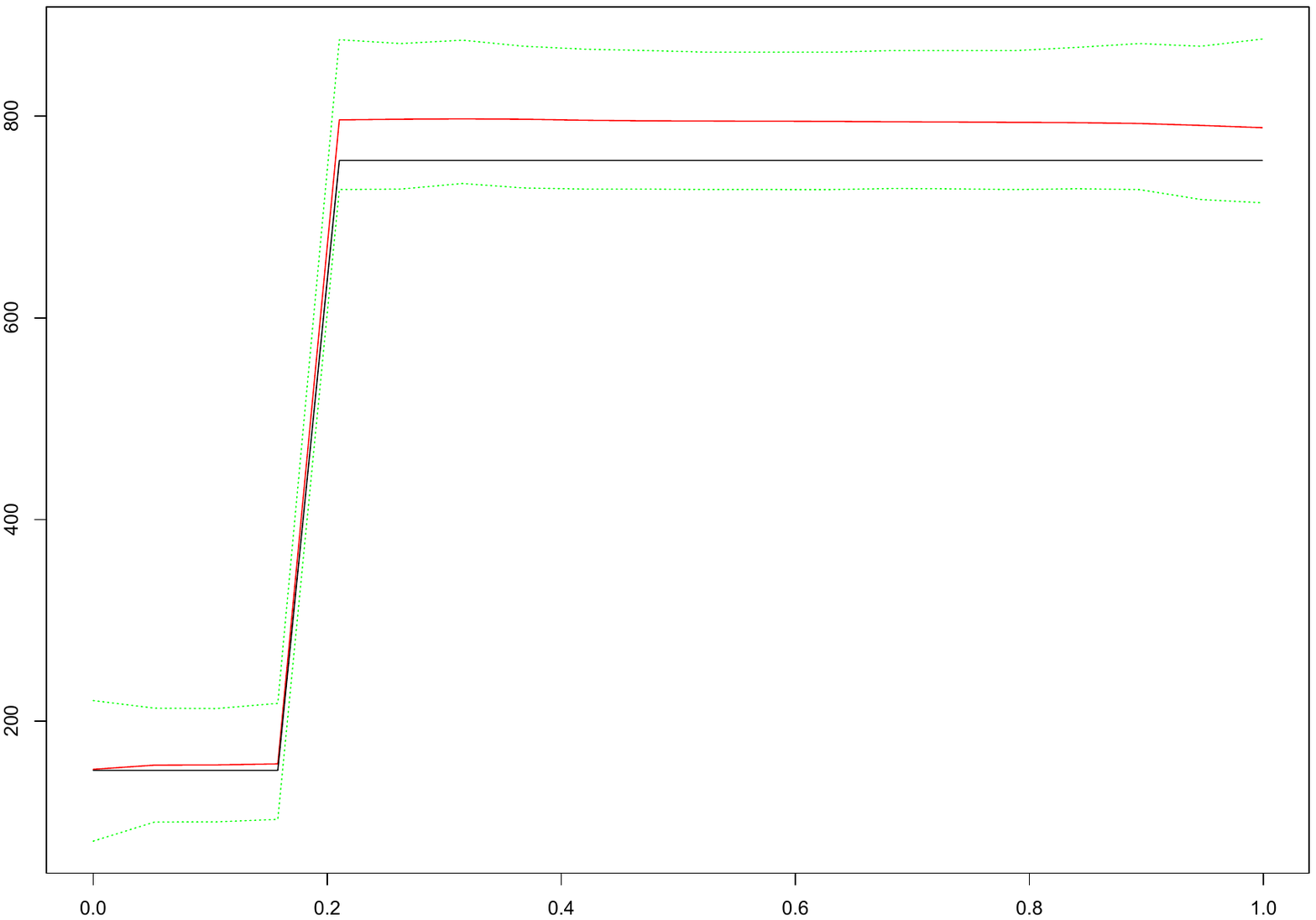}
    \caption{ The posterior mean of $\lambda(x_1)$ (red line; 95\% CI (green line)) and the true expected value of $\lambda(x_1)$ (black line).} 
  \end{subfigure}
  \begin{subfigure}{8cm}
    \centering\includegraphics[width=6cm]{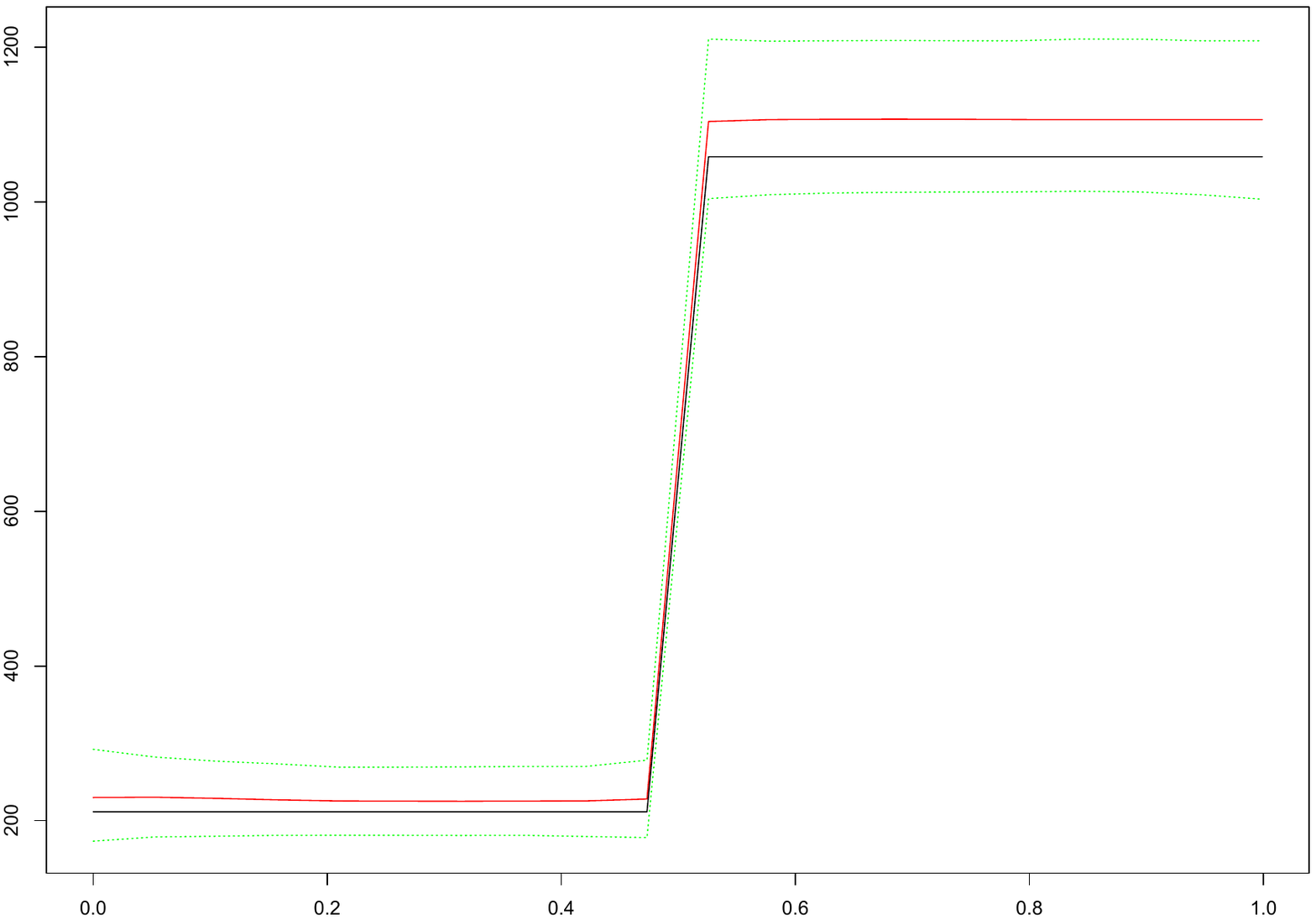}
   \caption{The posterior mean of $\lambda(x_2)$ (red line; 95\% CI (green line)) and the true expected value of $\lambda(x_2)$ (black line).}
  \end{subfigure}
  \begin{subfigure}{8cm}
    \centering\includegraphics[width=6cm]{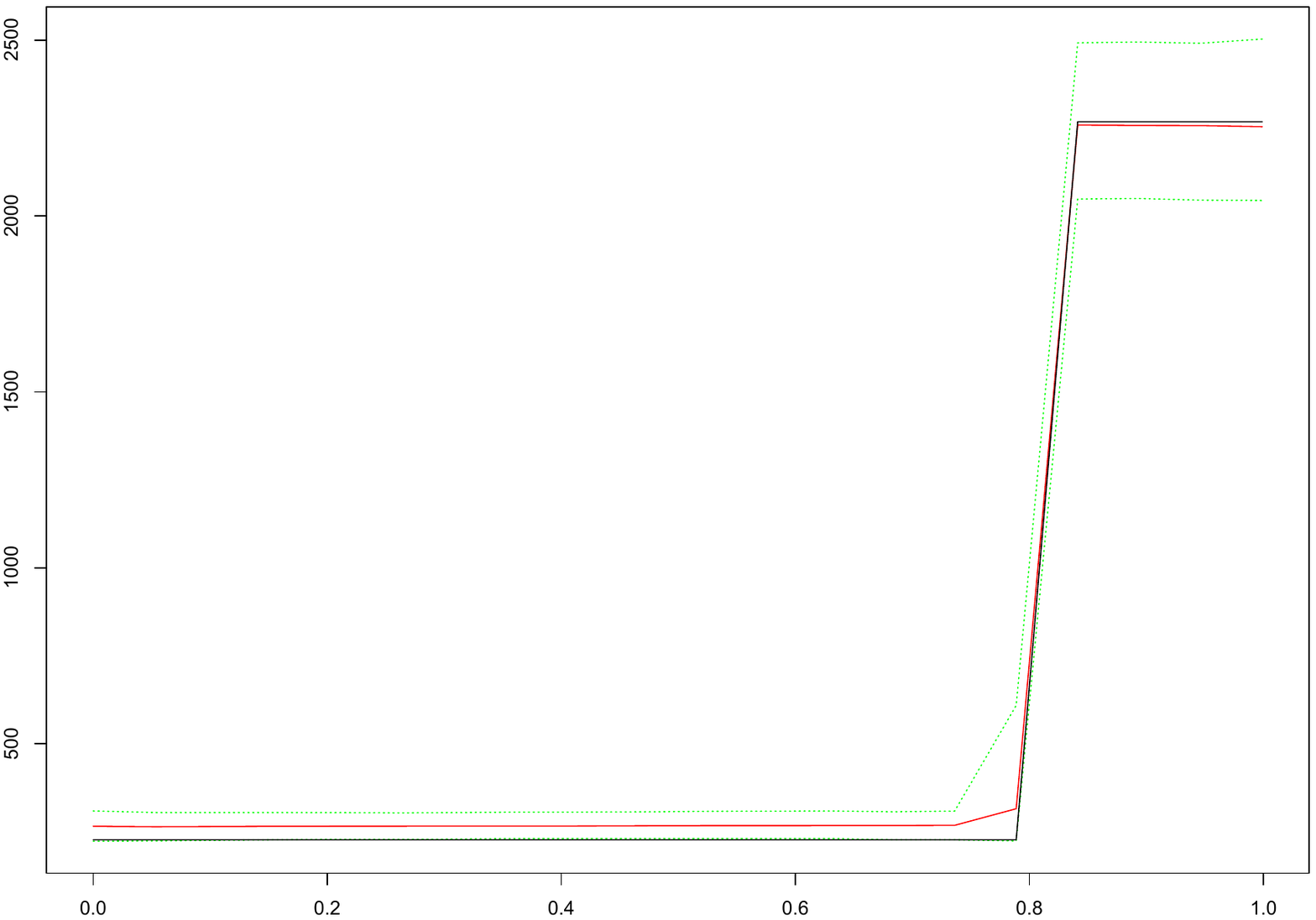}
   \caption{The posterior mean of $\lambda(x_3)$ (red line; 95\% CI (green line)) and the true expected value of $\lambda(x_3)$ (black line).}
  \end{subfigure}
  \begin{subfigure}{8cm}
    \centering\includegraphics[width=6cm]{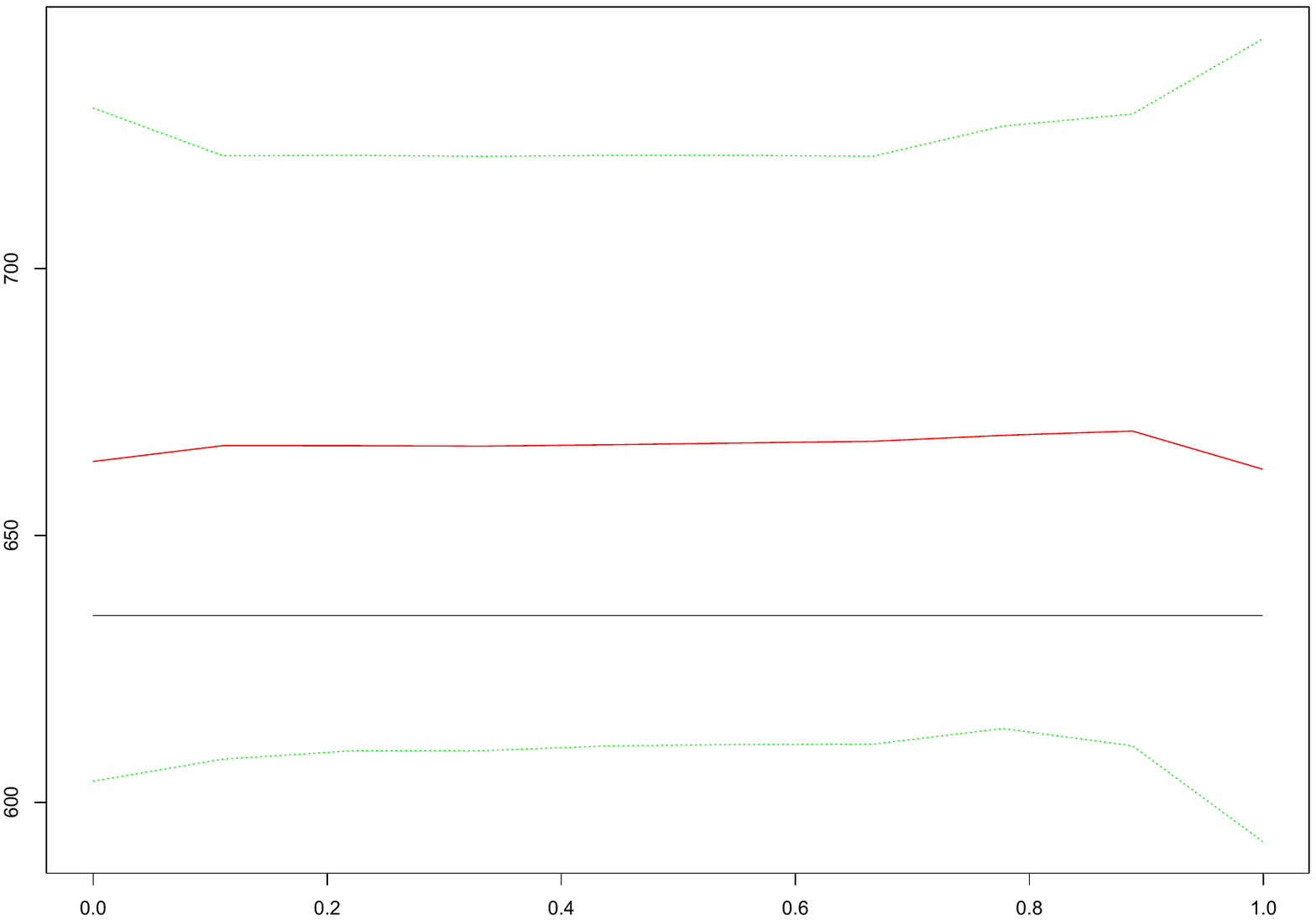}
   \caption{The posterior mean of $\lambda(x_4)$ (red line; 95\% CI (green line)) and the true expected value of $\lambda(x_4)$ (black line).}
  \end{subfigure}
  \begin{subfigure}{8cm}
    \centering\includegraphics[width=6cm]{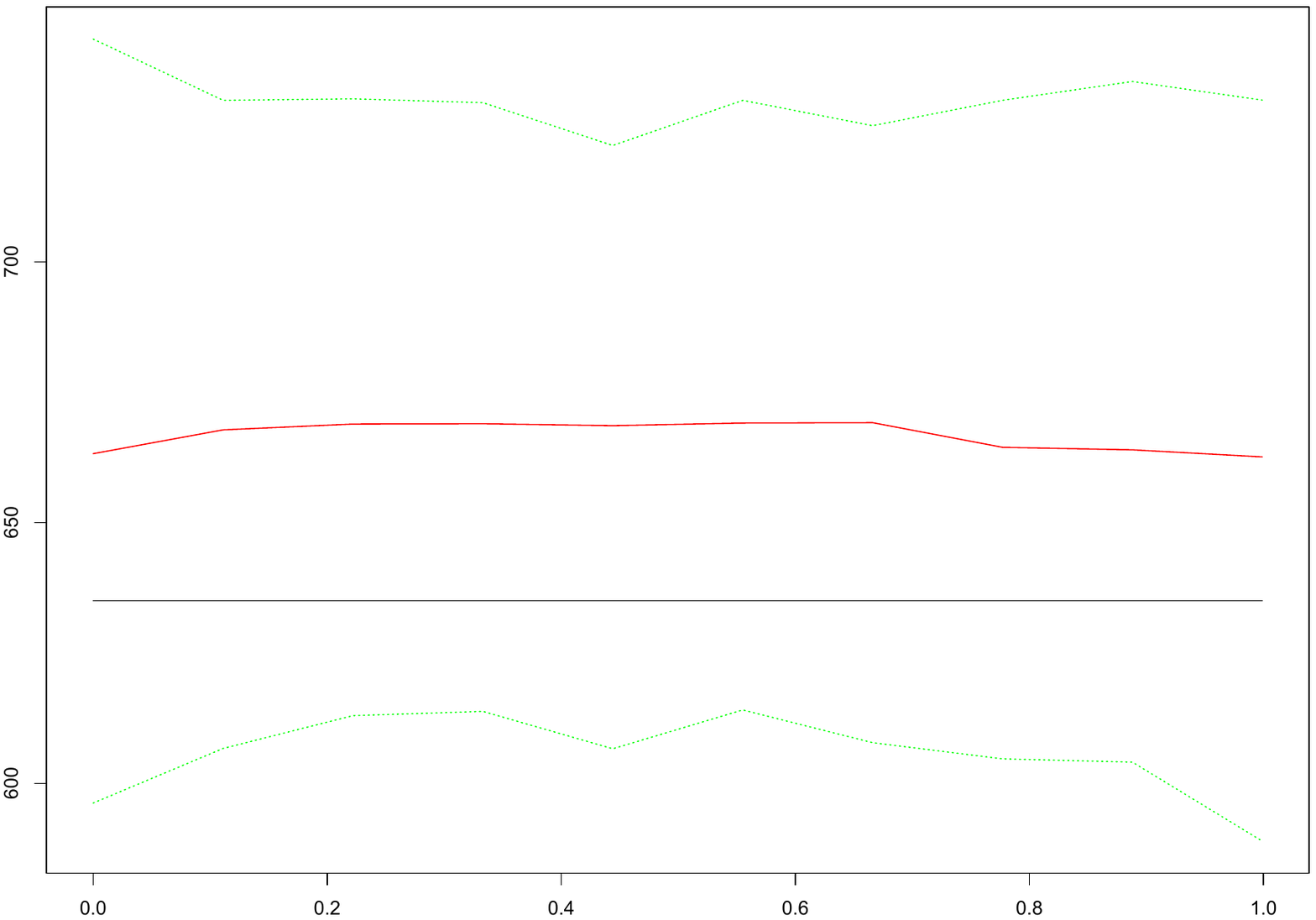}
   \caption{The posterior mean of $\lambda(x_5)$ (red line; 95\% CI (green line)) and the true expected value of $\lambda(x_5)$ (black line).}
  \end{subfigure}
  \caption{Posterior marginal intensities considering 4 trees. }
  \label{FigC5M}
\end{figure}

  Table \ref{table5F} demonstrates the frequency of times we meet each dimension in the decision rules of a tree. Table \ref{table5R} shows how likely each dimension is to be involved in the root's decision rule. The results illustrate that the important covariates $x_1$, $x_2$ and $x_3$ are more likely to be involved in the decision rules of a tree than the noisy dimensions $x_4$ and $x_5$. That indicates the algorithm prioritizes the dimensions that contribute most to the intensity. Figure \ref{FigC5M} shows that the mean of the posterior marginal intensities are similar to the expected marginal intensities given that $\{x_i\}_{i=1}^5$ are uniform independent covariates. 
 
 Tables \ref{tableC5D}, \ref{tableC5K} and \ref{tableC5L} show that our algorithm outperforms kernel smoothing and the maximum likelihood approach considering linear conditional intensity as expected.  The ability of our method to identify important features demonstrates an important advantage over other procedures.

 \begin{table}[H]
\begin{tabular}{ p{1.5 cm}p{1.6 cm}p{1.6cm}p{1.8 cm}p{1.8 cm}p{1.5 cm}p{1.5 cm}p{1.5 cm}}
 \hline
 \multicolumn{8}{c} {Proposed BART Algorithm}  \\
 \hline
 Number of  trees & AAE for Mean & AAE for Median & RISE for Mean & RISE for Median & $D_g$ & $D_l$ & $D_{LPO}$\\
 \hline
 4 & 48.36 & 45.47 & 159.95 & 170.35 & 8510.1 & 8525.4 & -485.9\\
 5& 49.18 & 44.54 & 158.82 & 169.07& 8486.1 & 8520.9 & -467.1\\
 6 & 50.59 & 45.05 & 161.36 & 170.61 & 8462.1 & 8519 & -477.4\\ 
 8 & 56.06 & 47.94 & 162.56 &164.46 &8349 & 8511.8 & -490.8   \\
 10 & 61.55 & 52.23 & 169.72 & 166.62 & 8141.8 & 8505.5 & -503.1\\
 12 & 67.01 & 57.06 & 180.53 & 175.23 & 7774.7 & 8499.8 &-522.2\\
 15 & 75 & 65.06 & 192.88 & 181.03 & 6813.1 & 8490 & -500.2\\
\hline
\end{tabular}
\caption{Average Absolute Error, Root Integrated Square Error and diagnostics for  various number of trees in the case of Inhomogeneous five dimensional Poisson Process with sparsity assumption.}
\label{tableC5D}
\end{table}

 \begin{table}[H]
\begin{tabular}{ p{4cm} p{3cm} p{3cm}  } \hline
 \multicolumn{3}{c}{Kernel Smoothing } \\
 \hline
 Bandwidth (sigma) & AAE  &  RISE  \\
 \hline
0.121 (LCV) &  407.1 & 888.1  \\
 \hline
\end{tabular}
\caption{Average Absolute Error and Root Integrated Square Error for fixed bandwidth estimators in the case of Inhomogeneous five dimensional Poisson Process with sparsity assumption.}
\label{tableC5K}
\end{table}

 \begin{table}[H]
\begin{tabular}{ p{3cm} p{3cm}  } \hline
 \multicolumn{2}{c}{Linear conditional intensity} \\
 \hline
 AAE  &  RISE  \\
 \hline
654.2 & 1076.5  \\
 \hline
\end{tabular}
\caption{Average Absolute Error and Root Integrated Square Error for linear conditional intensity in the case of Inhomogeneous five dimensional Poisson Process with sparsity assumption.}
\label{tableC5L}
\end{table}

\begin{table}[H]
\begin{tabular}{ p{2 cm} p{2 cm} p{2 cm} p{2 cm} p{2 cm} p{2 cm} }
 \hline
 \multicolumn{6}{c}{Proposed BART Algorithm}  \\
 \hline
 Number of  trees & $x_1$ & $x_2$ & $x_3$ & $x_4$ & $x_5$  \\
 \hline
 4 & 0.31  & 0.29  & 0.34 & 0.03 & 0.03   \\
 5& 0.35  & 0.29  & 0.26  & 0.05 & 0.06   \\
 \hline
\end{tabular}
\caption{How likely each dimension is to be involved in the root's decision rule.}
\label{table5R}
\end{table}

\begin{table}[H]
\begin{tabular}{ p{2 cm} p{2 cm} p{2 cm} p{2 cm} p{2 cm} p{2 cm} }
 \hline
 \multicolumn{6}{c}{Proposed BART Algorithm}  \\
 \hline
 Number of  trees & $x_1$ & $x_2$ & $x_3$ & $x_4$ & $x_5$  \\
 \hline
 4 & 0.35  & 0.36  & 0.37 & 0.06 & 0.07   \\
 5& 0.39  & 0.34  & 0.37  & 0.09 & 0.10   \\
 \hline
\end{tabular}
\caption{The frequency of times we meet each dimension in the decision rules of a tree.}
\label{table5F}
\end{table}

\begin{table}[H]
\begin{tabular}{ p{1.5 cm} p{1.5 cm} p{1.5 cm}  p{1.5 cm}  }
 \hline
 \multicolumn{4}{c}{Proposed BART Algorithm}  \\
 \hline
 Number of trees & $N_s=1$ & $N_s=32$ & $N_s=243$  \\
4          & 5.40 &  0.99 &  0.71  \\ 
5          & 5.40 & 1 &  0.71  \\ 
6          & 5.42 &  1 &  0.71  \\ 
8          & 5.42 &  1 &  0.71  \\ 
10          & 5.42 &  1 &  0.71  \\ 
15          & 5.43 &  1 &  0.71 \\ 
 \hline
\end{tabular}
\caption{The average RPS on testing points over 7 different splits of the original data set in the case of Inhomogeneous five dimensional Poisson Process with sparsity assumption.}
\label{table44}
\end{table}

\begin{table}[H]
\begin{tabular}{ p{1.5 cm} p{1.5 cm} p{1.5 cm}  p{1.5 cm}  }
 \hline
 \multicolumn{4}{c}{Proposed BART Algorithm}  \\
 \hline
 Number of trees & $N_s=1$ & $N_s=32$ & $N_s=243$  \\
4          & 0.64 &  0.95 &  1  \\ 
5          & 0.64 &  0.95 &  1  \\ 
6          & 0.65 &  0.95 &  1  \\ 
8          & 0.65 &  0.96 &  1.01  \\ 
10          & 0.65 &  0.96 &  1.01  \\ 
15          & 0.65 &  0.96 &  1.01  \\ 
 \hline
\end{tabular}
\caption{The average RSMSE on testing points over 7 different splits of the original data set in the case of Inhomogeneous five dimensional Poisson Process with sparsity assumption.}
\label{table44b}
\end{table}

\section{Intensity estimation for  Real Data}
 In this section, we first apply our algorithm to real data sets when modelled as realizations of inhomogeneous Poisson processes in one and two dimensions. To assess the performance of our algorithm, we break the domain $[0,1)^d$ into equal volume subareas $\{S_i\}_{i=1}^{N_S}$ and consider a set $\{z_i\}$ by uniformly sampling points in the domain  $[0,1)^d$. We compute the AAE of the estimated expected number of points falling into each of the subareas : 
 \begin{align}
      \text{AAE}(\widehat{N})= \frac{1}{N_S} \sum_{i=1}^{N_S}| \widehat{N}(S_i)-N(S_i) |
  \end{align}
  and Root Integrated Square Error (RISE):
  \begin{align}
     \text{RISE}(\widehat{N})=  \left( \frac{1}{N_S} \sum_{i=1}^{N_S}( \widehat{N}(S_i)-N(S_i) )^2 \right)^{1/2},
  \end{align}
where $N(S_i)$ is the actual number of points in $S_i$ and \begin{align}
     \widehat{N}(S_i)=  \int\limits_{S_i} \widehat{\lambda}(s)ds \simeq \frac{1}{N_{S_i}}\sum\limits_{z_j\in S_i} |S_i| \widehat{\lambda}(z_j)
  \end{align} with $N_{S_i}$ being the number of testing points $\{z_j\}$ falling in $S_i$. We apply the metrics AAE and RISE to compare our intensity estimates of one dimensional processes with those obtained by applying the Haar-Fisz algorithm for one dimensional data; and with kernel estimators for two-dimensional data. We observe that our algorithm, the Haar-Fisz algorithm and the kernel smoothing lead to similar results. As expected, the reconstructions of the intensity function are less smooth than those derived with kernel smoothing. The kernel estimator, as well as the bandwidth value given by likelihood cross-validation were computed using the \textbf{R} package \textbf{spatstat} \citep{spatstat}. We provide more simulation results in the supplementary material.

\subsection{Earthquakes Data}  
 This data set is available online from the Earthquake Hazards Program and consists of the times of 1088 earthquakes from 2-3-2020 to 1-4-2020. We consider the period from 27-2-2020 to 5-4-2020 to avoid edges. We run 3 parallel chains of the same length for 100000 iterations for 3-10 trees. The convergence criteria included in the supplementary material indicate that the considered chains have converged.

\begin{figure}[H] 
  \begin{subfigure}{8cm}
    \centering\includegraphics[width=8cm]{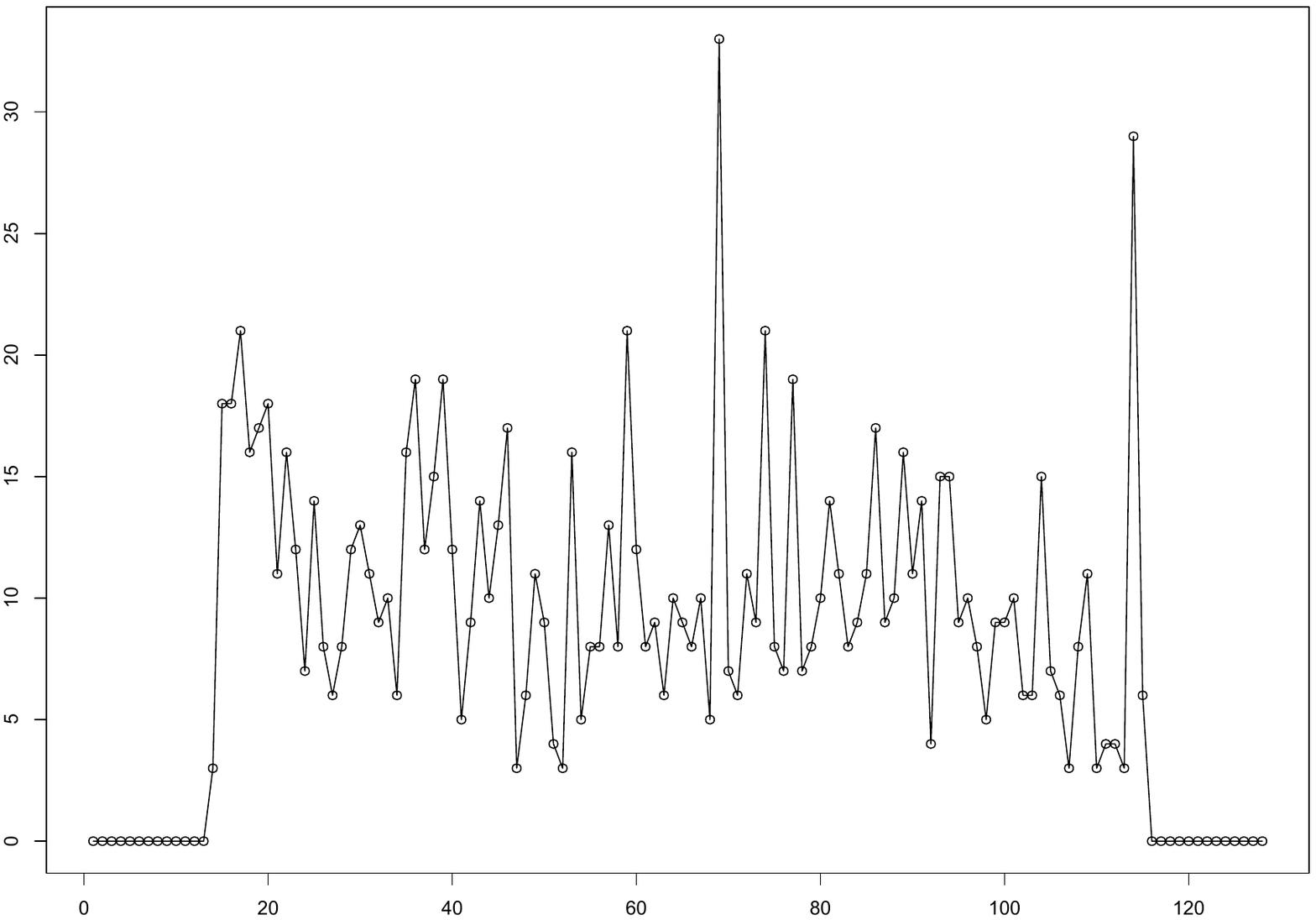}
    \caption{Aggregated Earthquakes in 128 consecutive intervals}
  \end{subfigure}
  \begin{subfigure}{8cm}
    \centering\includegraphics[width=8cm]{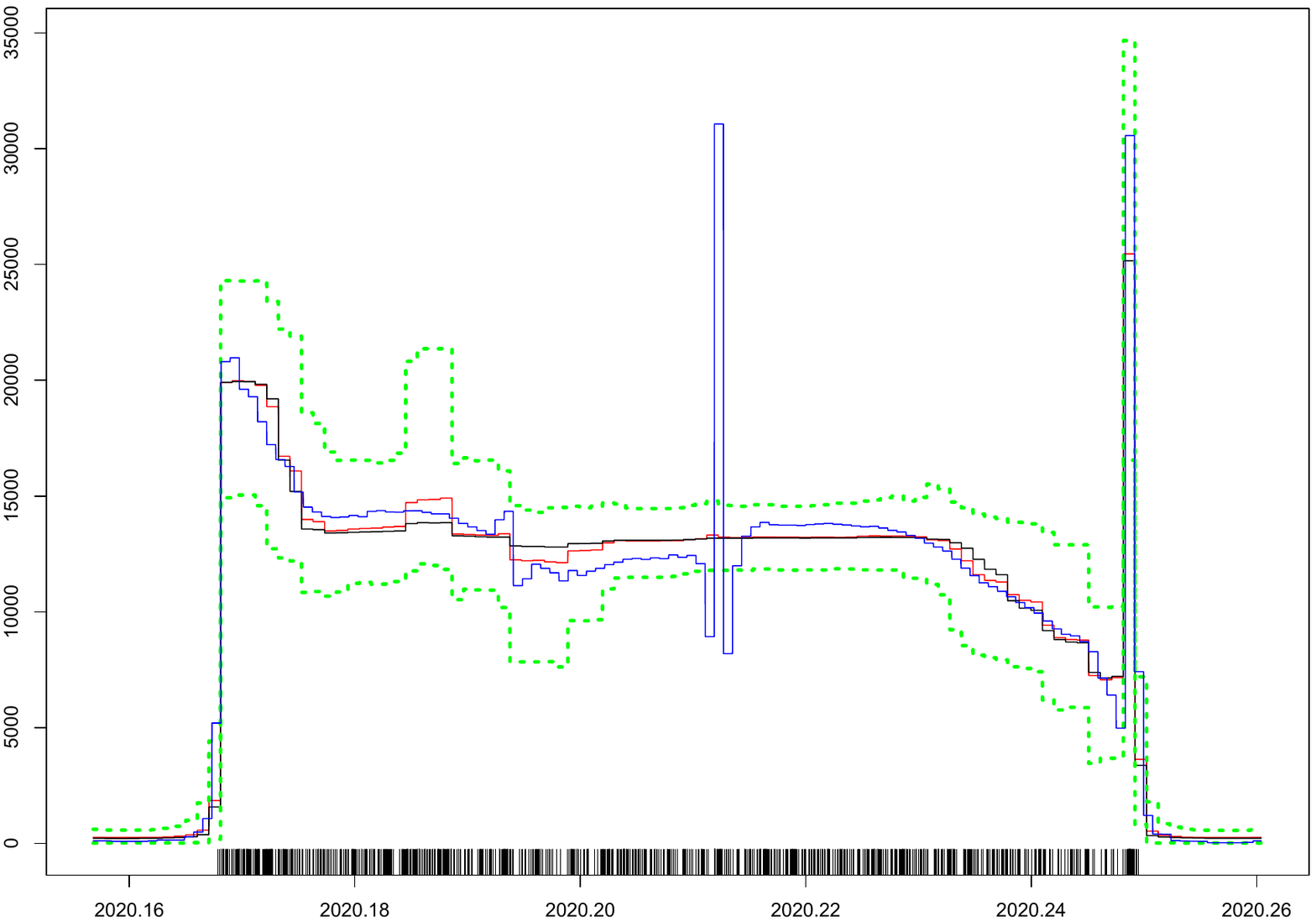}
    \caption{5 Trees and Haar-Fisz Algorithm}
  \end{subfigure}
\caption{Earthquakes Data: The posterior mean (red curve), the posterior median (black curve), the 95\% hdi interval of the estimated intensity illustrated by the dotted green lines and the intensity estimator of the Haar-Fisz Algorithm illustrated by the blue line. The rug plot on the bottom displays the event times.}
\label{earths}
\end{figure}

\begin{figure}[H] 
  \begin{subfigure}{8cm}
    \centering\includegraphics[width=8cm]{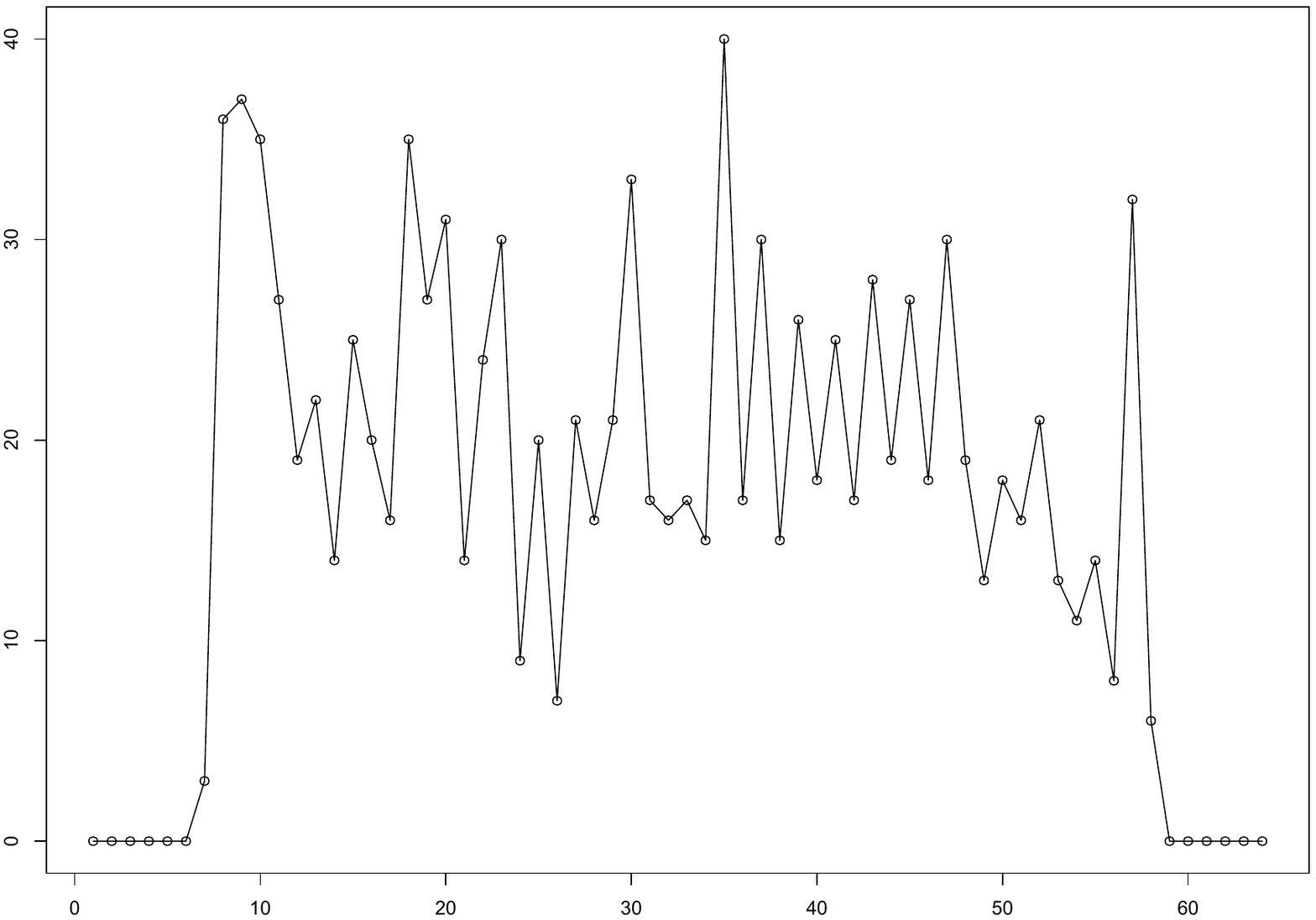}
    \caption{Aggregated Earthquakes in 64 consecutive intervals}
  \end{subfigure}
  \begin{subfigure}{8cm}
    \centering\includegraphics[width=8cm]{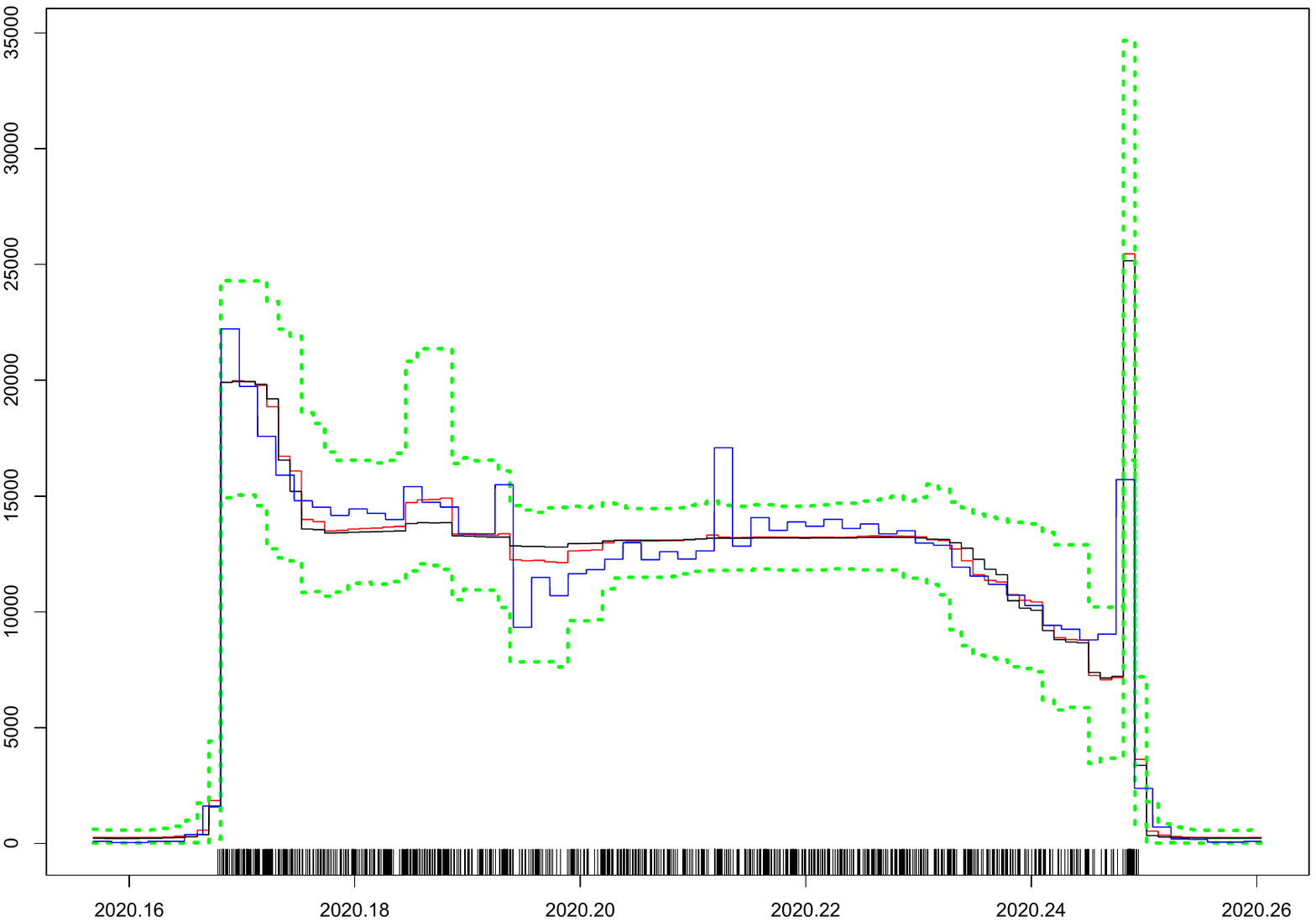}
    \caption{5 Trees and Haar-Fisz Algorithm}
  \end{subfigure}
\caption{Earthquakes Data: The posterior mean (red curve), the posterior median (black curve), the 95\% hdi interval of the estimated intensity illustrated by the dotted green lines and the intensity estimator of the Haar-Fisz Algorithm illustrated by the blue line. The rug plot on the bottom displays the event times.}
\label{earths64}
\end{figure}

 Figure \ref{earths} presents the Posterior Mean and the Posterior Median for 5 Trees, as well as the intensity estimate of the Haar-Fisz algorithm applied to the counts in 128 consecutive intervals of equal length. The deterministic discretized intensity of the \textbf{R} package \textbf{haarfisz} is divided by the duration of an interval. The differences between both algorithms are due to different assumptions; the Haar-Fisz algorithm considers the aggregated counts into disjoint subintervals of the domain, while the proposed algorithm the times of individual events. The most noticeable difference is observed between 2020.212 and 2020.213 (69th interval) where we see a jump in earthquakes from 5 to 33 and again to 7. The Haar-Fisz algorithm detects that peak as we feed it with that information, while the proposed algorithm does not indicate a sharp rise in the intensity in that period, treating it as an outlier. The intensity estimate of the Haar-Fitz algorithm applied in 64 consecutive intervals is closer to the proposed algorithm (see Figure \ref{earths64}), as expected. Similar to coarser binning, the proposed algorithm is less prone to overfitting to spikes in the data, which get filtered out. 
 The estimated AAE and RISE demonstrate good performance compared to the Haar-Fisz method. The simulation results illustrate that our algorithm can track the varying intensity of earthquakes. 
 
  The diagnostics $D_g$, $D_l$ and $D_{LPO}$ obtain their highest values at 9, 3 and 8 trees, respectively. The AIC diagnostics values between 3 and 9 trees show only small variations,  we choose 5 trees for the analysis, noting that the results will not vary significantly for other choices of $m$ in this region. 

\begin{table}[H]
\begin{tabular}{ p{1.8 cm} p{2.4 cm} p{2.4 cm} p{1.5 cm} p{1.5 cm} p{1.5 cm} p{1.5 cm}   p{1.5 cm}}
 \hline
 \multicolumn{8}{c}{Proposed BART Algorithm}  \\
 \hline
 Number \\of trees & AAE for Posterior Mean & AAE for Posterior Median & RISE for Posterior Mean & RISE for Posterior Median & $D_g$ & $D_l$ & $D_{LPO}$\\
 \hline
3 & 93.8 & 94.1 & 106.9 & 107.1 & 13570.1  & 13565.4 & -1194.7 \\
4 & 94 & 94.1 & 106.8 & 107 & 13570.6 & 13563.6  & -1163.8 \\
5 & 93.6 & 94 & 106.7 & 107 & 13570.1 & 13560.7   & -1150.5\\
6 & 93.8 & 94 & 106.9 & 107 & 13571.6 & 13559.6  & -1169.7 \\
8 & 93.5 & 94 & 106.6 & 107.1 & 13571.8 & 13554.9 & -1140.2  \\
9 & 93.4 & 94 & 106.7 & 107.3 & 13572.1 & 13552.6  & -1192.5 \\
10 & 93.4 & 94 & 106.8 & 107.3 & 13571.9 & 13549.8  & -1184.4 \\
 \hline
\end{tabular}
\caption{Average Absolute Error, Root Integrated Square Error and diagnostics for the data in Fig.~\ref{earths}.}
\label{tableE2}

\end{table}
\begin{table}[H]
\begin{tabular}{ p{3cm} p{3cm} p{3cm} }
 \hline
 \multicolumn{3}{c}{Haar-Fisz Algorithm}  \\
 \hline
 Subintervals & AAE  &  RMSE   \\
 \hline
128 & 94.1 & 107.8   \\
64 & 94 & 107\\
\hline
\end{tabular}
\caption{Average Absolute Error and Root  Mean Square Error for Haar-Fisz estimator for the data in Fig.~\ref{earths}}
\label{tableEH2}
\end{table}

\subsection{Lansing Data}  

\begin{figure}[H] 
  \begin{subfigure}{8cm}
    \centering\includegraphics[width=6cm]{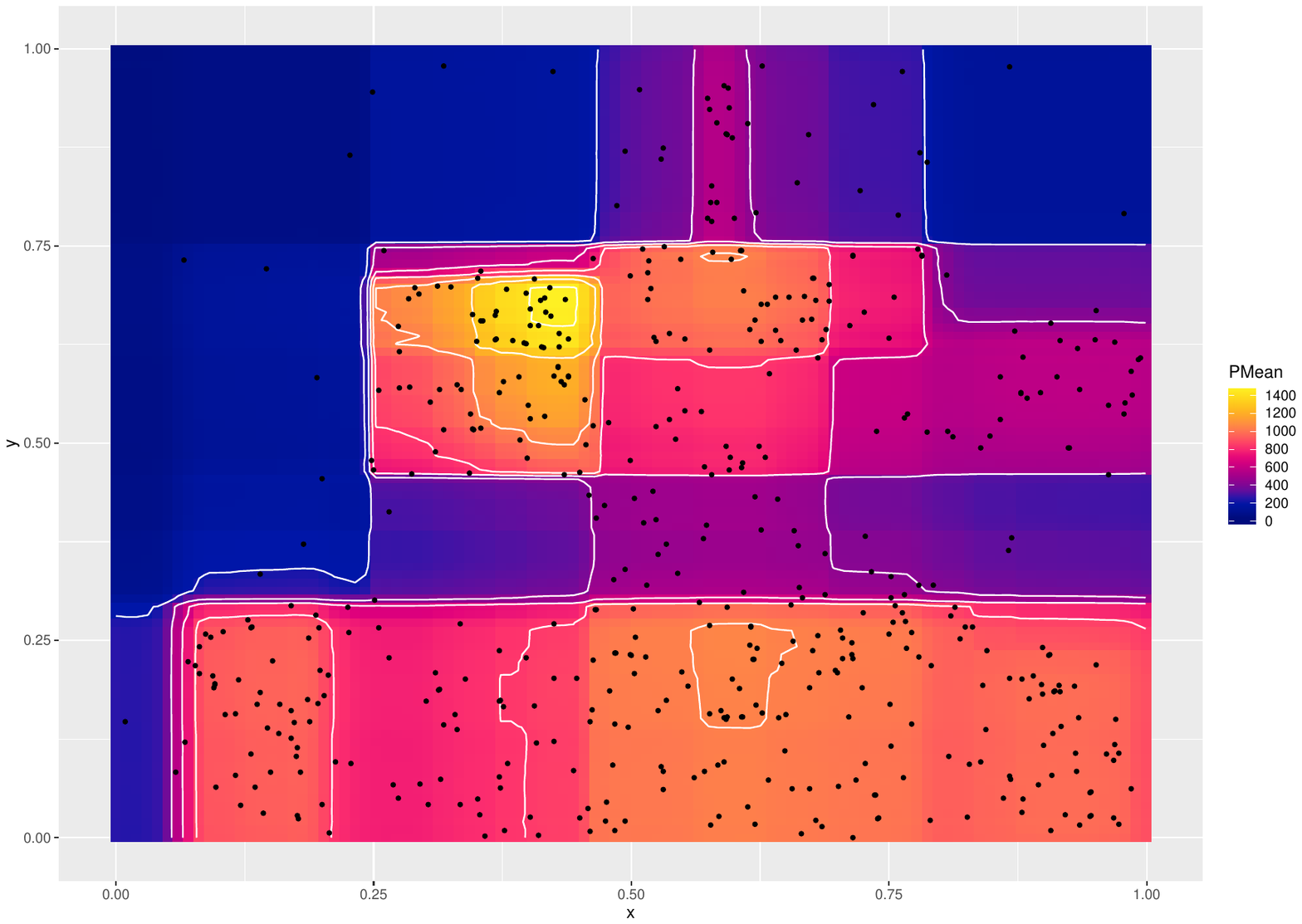}
    \caption{Posterior Mean for 5 Trees}
  \end{subfigure}
  \begin{subfigure}{8cm}
    \centering\includegraphics[width=6cm]{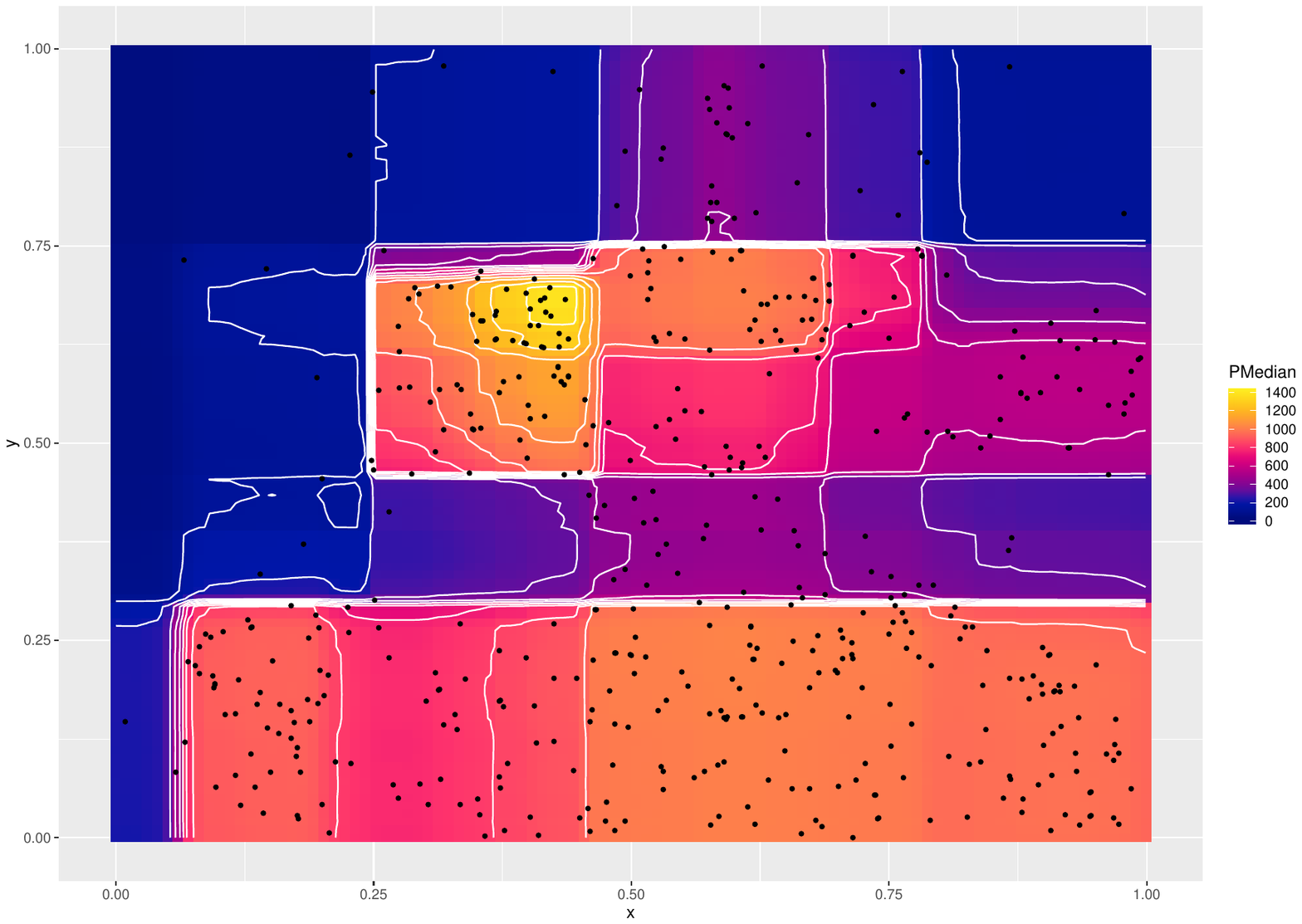}
    \caption{Posterior Median for 5 Trees}
  \end{subfigure}

 \begin{subfigure}{8cm}
    \centering\includegraphics[width=6cm]{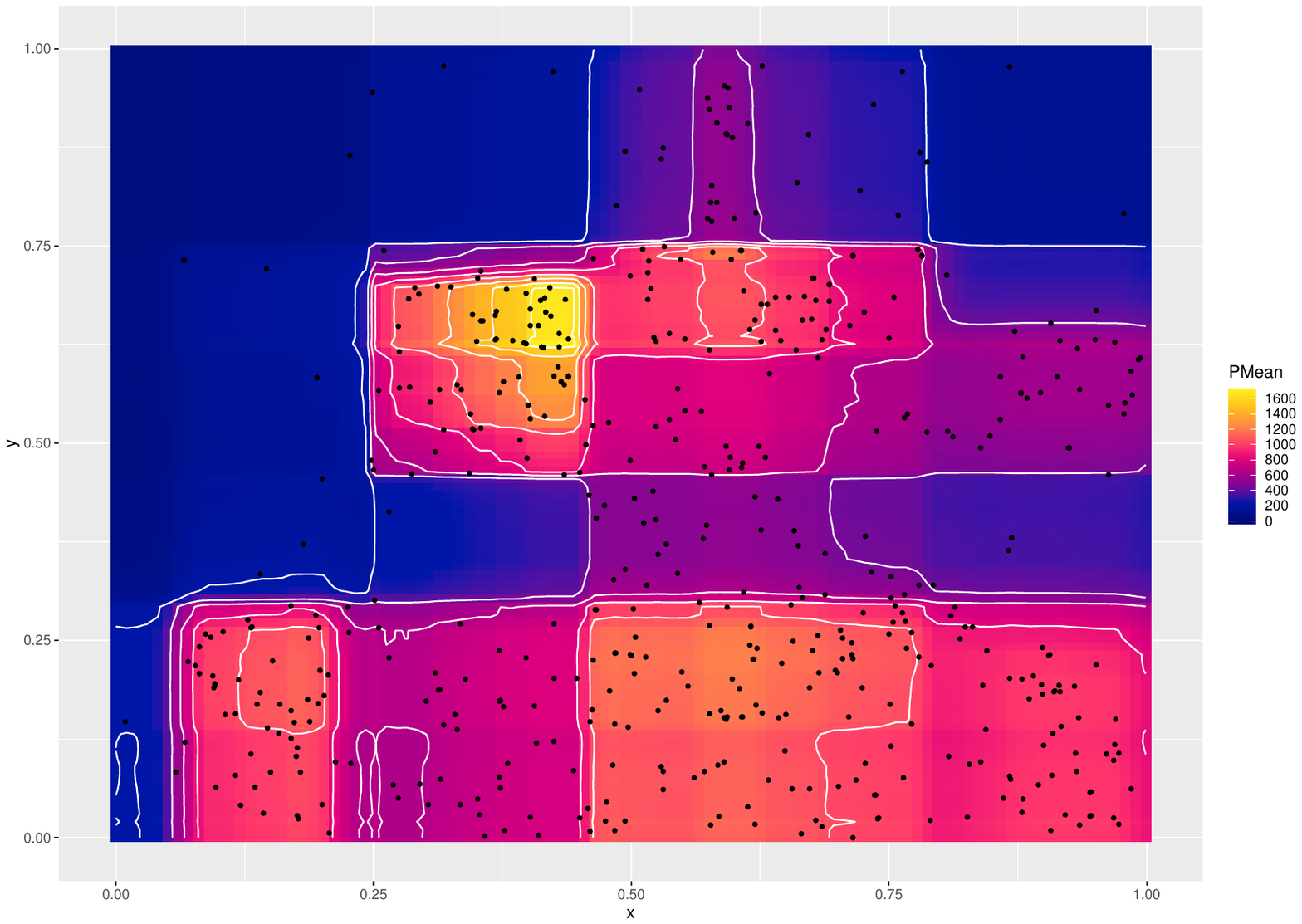}
    \caption{Posterior Mean for 10 Trees}
  \end{subfigure}
  \begin{subfigure}{8cm}
    \centering\includegraphics[width=6cm]{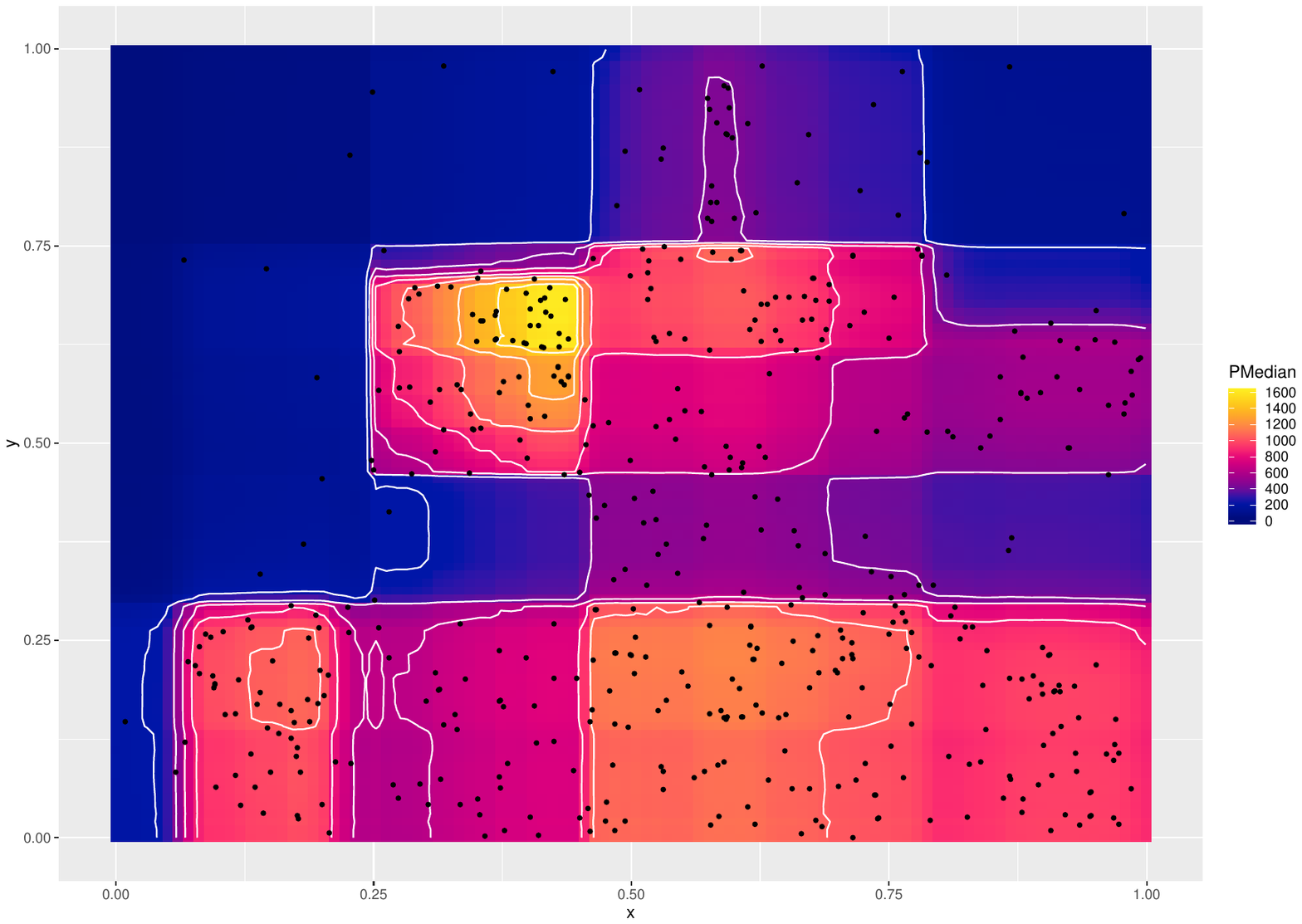}
    \caption{Posterior Median for 10 Trees}
  \end{subfigure}
\caption{Posterior Mean and Posterior Median for 5 and 10 Trees}
\label{fig6}
\end{figure}

\begin{figure}[H]
 \centering\includegraphics[width=6cm]{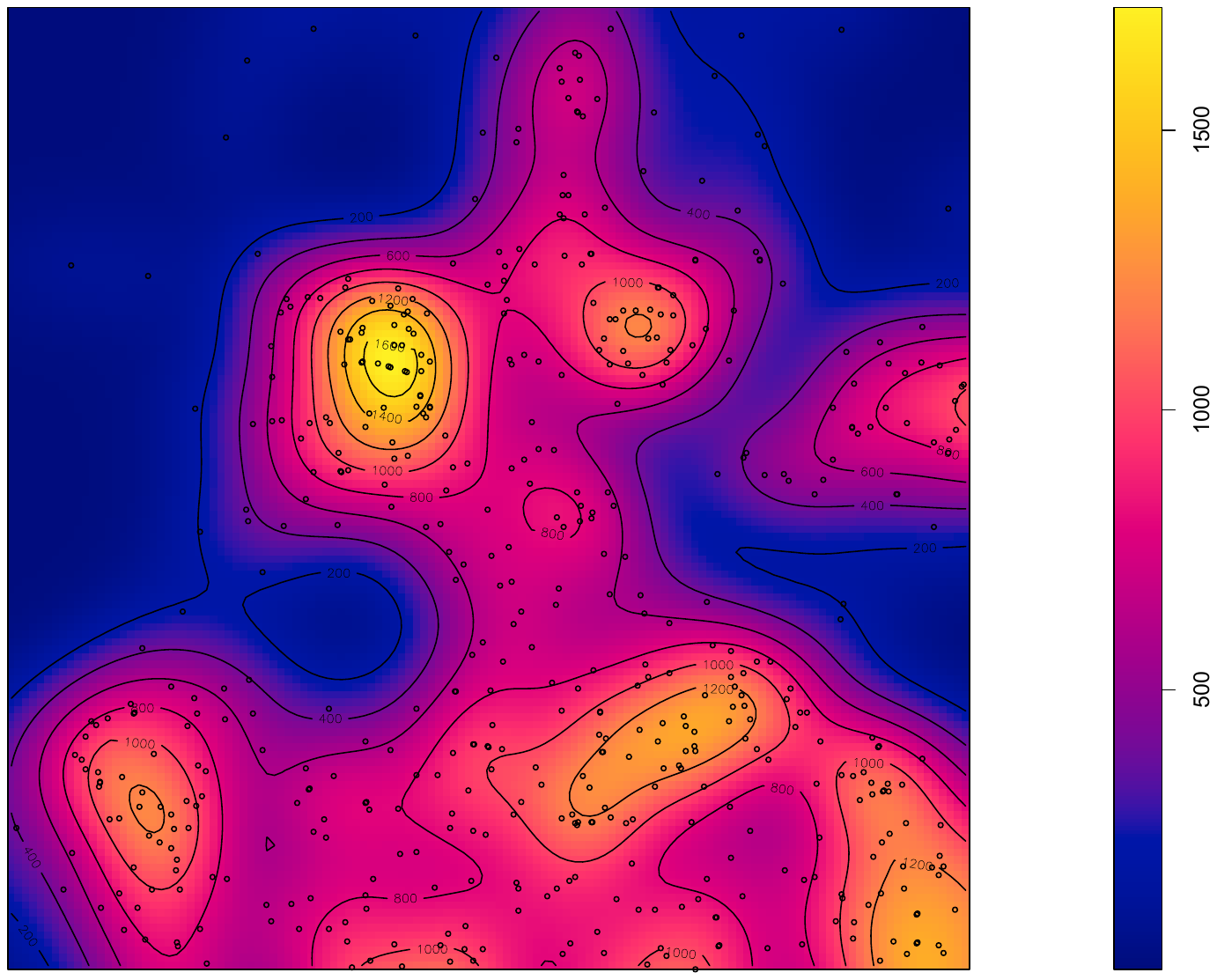}
   \caption{Fixed-bandwidth chosen using likelihood cross-validation.}
\label{fig7} 
\end{figure}

The \textbf{lansing} data set included in the \textbf{R} package \textbf{spatstat} describes the locations of different types of trees in the Lansing woods forest. Our attention is restricted to the locations of 514 maples that are presented with dots in Figures~\ref{fig6}-\ref{fig7}. We run 3 parallel chains of the same length for 200000 iterations for 3-10 trees and 100000 iterations for 12 trees. The diagnostic criteria included in the supplementary material indicate that the considered chains have converged for the majority of testing points. 

We compare our algorithm to a fixed bandwidth estimator using a Gaussian kernel. Our algorithm and the kernel estimator are consistent in the overall structure. The differences are due to the different nature of the methods. Given the tree locations, our algorithm recovers the spatial pattern of trees as rectangular regions of different intensities (Fig.~\ref{fig6}), whereas the kernel method produces a continuum with more localized peaks in space. As expected, the kernel estimator presented in~Figure~\ref{fig7} consists of smoother subregions with various intensities.  Tables \ref{tableM1}-\ref{tableM2}
show that our algorithm is competitive to kernel smoothing with fixed bandwidth chosen with likelihood cross-validation. In contrast to our method, kernel methods are highly sensitive to parameter (bandwidth) choice. 

The diagnostics $D_g$ and $D_l$ obtain their highest values at 4 and 10 trees, respectively.  

\begin{table}[H]
\begin{tabular}{ p{2cm} p{2cm} p{2cm} p{2cm} p{2cm} p{1.5 cm} p{1.5 cm} }
 \hline
 \multicolumn{5}{c}{Proposed BART Algorithm}  \\
 \hline
 Number of trees & AAE for Posterior Mean & AAE for Posterior Median & RMSE for Posterior Mean & RMSE for Posterior Median & $D_g$ & $D_l$\\
 \hline
 3 & 1.3 & 1.2 & 1.7 &  1.8 & 5686.5 & 5705.3 \\
 4 & 1.2 & 1.2 & 1.7 &  1.8 & 5683.8 & 5709.5 \\
 5 & 1.2 & 1.2 & 1.7 &  1.7 & 5672.4 & 5705.4 \\
 7 & 1.2 & 1.2 & 1.7 &  1.71 & 5643.5 & 5702 \\
 8 & 1.2 & 1.2 & 1.7 &  1.7 & 5634 & 5707.2 \\
 9 & 1.2 & 1.2 &  1.7 &  1.7 & 5614.3 & 5698.1 \\
 10& 1.2 & 1.2 &  1.6 &  1.7 & 5596.6 & 5699.8 \\
12 & 1.2 & 1.2 & 1.7 & 1.7 & 5558.2 & 5692.5  \\
 \hline
\end{tabular}
\caption{Average Absolute Error, Root Integrated Square Error with $N_S=225$ and diagnostics for the data in Fig.~\ref{fig6}.}
\label{tableM1}
\end{table}

\begin{table}[H]
\begin{tabular}{ p{2cm} p{2.4cm} p{2.4cm} p{2.5cm} p{2.5cm} p{2 cm}  }
 \hline
 \multicolumn{5}{c}{Proposed BART Algorithm}  \\
 \hline
 Number of trees & AAE for Posterior Mean & AAE for Posterior Median & RMSE for Posterior Mean & RMSE for Posterior Median\\
 \hline
 3 & 0.9 & 0.9 & 1.3 &  1.3 \\
 4 & 0.9 & 0.9 & 1.2 &  1.3  \\
 5 & 0.9 & 0.9 & 1.2 &  1.3  \\
 7 & 0.9 & 0.9 & 1.2 &  1.2  \\
 8 & 0.9 & 0.9 & 1.2 &  1.2  \\
 9 & 0.9 & 0.9 &  1.2 &  1.2  \\
 10& 0.9 & 0.9 &  1.2 &  1.2  \\
12 & 0.9 & 0.9 & 1.2 & 1.2   \\
 \hline
\end{tabular}
\caption{Average Absolute Error and Root Integrated Square Error with $N_S=400$ for the data in Fig.~\ref{fig6}.}
\label{tableM11}
\end{table}

\begin{table}[H]
\begin{tabular}{ p{4cm} p{3cm} p{3cm}  } \hline
 \multicolumn{3}{c}{Kernel Smoothing} \\
 \hline
 Bandwidth (sigma) & AAE  &  RISE  \\
 \hline
0.05 (LCV) for $N_S=225$ & 1.03 & 1.42  \\
0.05 (LCV) for $N_S=400$ & 0.82 & 1.13  \\
  \hline
\end{tabular}
\caption{Average Absolute Error and Root Integrated Square Error for fixed bandwidth estimators for data in Fig.~\ref{fig7}.}
\label{tableM2}
\end{table}

 \section{Discussion and Future Work}
 In this article, we have studied how the Bayesian Additive Regression Trees (BART) model can be applied to estimating the intensity of Poisson processes. The BART framework provides a flexible non-parametric approach to capturing non-linear and additive effects in the underlying functional form of the intensity. Our numerical experiments show that our algorithm provides good approximations of the intensity with ensembles of less than 10 trees. This enables our algorithm to detect the dimensions contributing most to the intensity. The ability of our method
to identify important features demonstrates an important advantage over other procedures. 

Our approach enables full posterior inference of the intensity in a non-parametric regression setting. In addition, the method extends easily to higher dimensional settings. The simulation study on 
synthetic data sets shows that our algorithm can detect change points and provides good estimates of the intensity via either the posterior mean or the posterior median. Our algorithm is competitive with the Haar-Fisz algorithm and kernel methods in one and two dimensions and inference using spatial log-Gaussian Cox processes. The strength of our method is its performance in higher dimensions, and we demonstrate that it outperforms the kernel approach for multidimensional intensities.  We also demonstrate that our inference for the intensity is consistent with the variability of the rate of events in real and synthetic data. The convergence criteria included in the  supplementary material indicate good convergence of the considered chains. We ran each chain for at least 100000 iterations to increase our confidence in the results. However, our algorithm works well with considerably fewer iterations (around 10000). 
 The BART model assumes independence of the underlying tree structure. The alternative method of \citep{doi:10.1198/016214504000000188} makes use of a locally dependent Markov Random Field, and one way of extending our model in this direction is to consider neighbouring intensities following \citet{10.1214/21-BA1259}. 
 
 Our method has only considered the standard priors commonly used in BART procedures, an interesting avenue of future research would be to implement different prior assumptions.  In addition, we have fixed the parameters for the Galton-Watson prior on the trees, and further work on sensitivities to hyperparameter selection and alternative methods for inference of the hyperparameters is of interest.
 Currently, our model is limited to non-homogeneous Poisson Process and we believe the flexibility of the BART approach could be extended to more general point processes.

\bibliography{references}

\appendix


\section{Metropolis Hastings Proposals}
\label{ap:A}
We describe the proposals of Algorithm~\ref{alg2}.
The Hastings ratio can be expressed as the product of three terms \citep{JSSv070i04}:
\begin{itemize}
\item{Transition Ratio: 
\[\TR=\frac{q(T_j^{(t)}|T_j^{*})}{q(T_j^{*}|T_j^{(t)})}\] 
}
\item{Likelihood Ratio: 
\[\LR=\frac{P(\vs|T_j^{*},T_{(j)},\Lambda_{(j)})}
{P(\vs|T_j^{(t)},T_{(j)},\Lambda_{(j)})}\] 
}
\item{Tree Structure Ratio: 
\[ \TSR=\frac{P(T_j^{*})}{P(T_j^{(t)})} \] 
}
\end{itemize}

\subsection{GROW Proposal}
This proposal randomly picks a terminal node, splits the chosen terminal into two new nodes and assigns a decision rule to it.  

 Let $\eta$ be the randomly picked terminal node in tree $T_j^{(t)}$. We denote the new nodes as $\eta_L$ and $\eta_R$. We now derive the expressions for the transition ratio ($\TR$), tree structure ratio ($\TSR$) and likelihood ratio ($\LR$).  

\paragraph*{Transition Ratio}
 It holds that: 
\begin{enumerate}[label=(\roman*)]
\item{
$q(T_j^{*}|T_j^{(t)}) =$ P(GROW)\\ 
\hspace*{20mm}$\times$ P(selecting a leaf $\eta$ to grow from)\\ 
\hspace*{20mm}$\times$ P(selecting an available dimension $j$ to split on)\\
\hspace*{20mm}$\times$ P(selecting the slitting value given the chosen dimension to split on) 
\\ \hspace*{18mm}= P(GROW) $\frac{1}{b_j}\frac{1}{\text{card}(k_{\eta})}\frac{1}{\text{card}(\tau_{\eta})}$\\
 where $b_j$ is the number of terminal nodes in the tree $T_j^{(t)}$, $k_h$ the set of all available dimensions to split the node $\eta$, $\tau_{\eta}$ the set of all available splitting values given the chosen dimension for splitting the node $\eta$ and card($S$) the cardinality of a set $S$.
}
\item{
$q(T_j^{(t)}|T_j^{*})=$ P(PRUNE)\\ 
\hspace*{20mm}$\times$ P(selecting a node $\eta$ having two terminal nodes to prune from)\\ 
\hspace*{18mm}= P(PRUNE) $\frac{1}{w^{*}}$ \\  
where $w^*$ is the number of internal 
nodes with two terminal nodes as children in the tree $T_j^{*}$. 
}
\end{enumerate}
 Hence the transition ratio is given by 
\begin{equation*}
\TR=\frac{P(\text{PRUNE})\frac{1}{w^{*}}}{P(\text{GROW})\frac{1}{b_j}\frac{1}{\text{card}(k_{\eta})}\frac{1}{\text{card}(\tau_{\eta})}}.
\end{equation*}

\paragraph*{Tree Structure Ratio:}
 The difference between the structures of the proposed tree $T_{j}^{(t)}$ and the tree $T_j^{*}$ is the two offsprings $\eta_{L}$ and $\eta_{R}$. Thus the tree structure ratio is: 
\begin{align*}
\TSR&=\frac{P(T_j^{*})}{P(T_j^{(t)})}=
\frac{(1-p_{\text{SPLIT}}(\eta_{L})) \, (1-p_{\text{SPLIT}}(\eta_{R})) \, p_{\text{SPLIT}}(\eta) \, p_{RULE}(\eta)}{(1-p_{\text{SPLIT}}(\eta))} \\
&=\frac{ \left(1-\frac{\gamma}{(1+d(\eta_L))^{\delta}}  \right)  \left(1-\frac{\gamma}{(1+d(\eta_R))^{\delta}}  \right) \frac{\gamma}{(1+d(\eta))^{\delta}}\frac{1}{\text{card}(k_{\eta})}\frac{1} {\text{card}(\tau_{\eta})} }{1-\frac{\gamma}{(1+d(\eta))^{\delta}}},
\end{align*}
 where $p_{\text{SPLIT}}(\eta)$ is the splitting probability for a node $\eta$ and $p_{RULE}(\eta)$ the distribution of decision rule associated to node $\eta$.

\paragraph*{Likelihood Ratio}
 The likelihood ratio is an application of equation \ref{equ3} twice, that is once considering the proposed tree, $T_j^{*}$  (numerator)  and the other considering the tree of the current iteration $t$, $T_{j}^{(t)}$ (denominator), which can be simplified as follows
\begin{align*} 
\LR&=\frac{\beta^{\alpha}}{\Gamma(\alpha)}\frac{\frac{\Gamma(n_{j\eta_L}+\alpha)}{(c_{j\eta_L}+\beta)^{n_{j\eta_L}+\alpha}}\frac{\Gamma(n_{j\eta_R}+\alpha)}{(c_{j\eta_R}+\beta)^{n_{j\eta_R}+\alpha}}}{\frac{\Gamma(n_{j\eta}+\alpha)}{(c_{j\eta}+\beta)^{n_{j\eta}+\alpha}}}\\
 &=\frac{\beta^{\alpha}}{\Gamma(\alpha)} \frac{\Gamma(n_{j\eta_L}+\alpha)\Gamma(n_{j\eta_R}+\alpha)}{\Gamma(n_{j\eta}+\alpha)} \frac{(c_{j\eta}+\beta)^{n_{j\eta}+\alpha}}{(c_{j\eta_L}+\beta)^{n_{j\eta_L}+\alpha}(c_{j\eta_R}+\beta)^{n_{j\eta_R}+\alpha}}\end{align*}

\subsection{PRUNE Proposal}
This proposal randomly picks a parent of two terminal nodes and turns it into a terminal  node by collapsing the nodes below it.  

 Let $\eta$ be the picked parent of two terminal nodes, $y$ and $c$ the dimension and splitting value of the rule linked to the node $\eta$.  

\paragraph*{Transition Ratio}
 It holds that: 
\begin{enumerate}[label=(\roman*)]
\item{
$q(T_j^{*}|T_j^{(t)}) =$ P(PRUNE)\\ 
\hspace*{20mm} $\times$ P(selecting a parent of two terminal nodes to prune from)\\  \hspace*{18mm}= P(PRUNE) $\frac{1}{w}$ \\   where $w$ is the number of nodes with two 
terminal nodes as children in the tree $T_j^{(t)}$. 
}
\item{
$q(T_j^{(t)}|T_j^{*})=$ P(GROW)\\ 
\hspace*{20mm}$\times$ P(selecting the node $\eta$ to grow from)\\ \hspace*{20mm}$\times$ P(selecting the dimension $y$)\\ \hspace*{20mm}$\times$ P(selecting the slitting value c given the chosen dimension $y$) 
\\ \hspace*{18mm}= P(GROW) $\frac{1}{w^*}\frac{1}{\text{card}(k_{\eta})}\frac{1}{\text{card}(\tau_{\eta})}$\\
 where $w^*$ is the number of terminal nodes in the tree $T_j^{*}$, $k_h$ the set of all available dimensions to split the node $\eta$ and $\tau_{\eta}$ the set of all available splitting values given the chosen dimension y for splitting the node $\eta$.
}
\end{enumerate}
 Hence the transition ratio is given by 
\begin{equation*}
\TR=\frac{P(\text{GROW})\frac{1}{w^*}\frac{1}{\text{card}(k_{\eta})}\frac{1}{\text{card}(\tau_{\eta})}}{P(\text{PRUNE})\frac{1}{w}}.
\end{equation*}

\paragraph*{Tree Structure Ratio}
 The proposed tree differs by not having the two children nodes $\eta_{L}$ and $\eta_{R}$. Thus the tree structure ratio is: 
\begin{align*}
\TSR=\frac{P(T_j^{*})}{P(T_j^{(t)})}&=
\frac{(1-p_{\text{SPLIT}}(\eta))} {(1-p_{\text{SPLIT}}(\eta_{L})) \, (1-p_{\text{SPLIT}}(\eta_{R})) \, p_{\text{SPLIT}}(\eta) \, p_{RULE}(\eta)}\\
&=\frac{1-\frac{\gamma}{(1+d(\eta))^{\delta}}}{ \left(1-\frac{\gamma}{(1+d(\eta_L))^{\delta}}  \right)  \left(1-\frac{\gamma}{(1+d(\eta_R))^{\delta}}  \right) \frac{\gamma}{(1+d(\eta))^{\delta}}\frac{1}{\text{card}(k_{\eta})}\frac{1}{\text{card}(\tau_{\eta}) }}
\end{align*}

\paragraph*{Likelihood Ratio}
 Similar to the GROW proposal, the likelihood ratio can be written as follows
\begin{align*} 
\LR&= \left(\frac{\beta^{\alpha}}{\Gamma(\alpha)}  \right)^{-1}\frac{\frac{\Gamma(n_{j\eta}+\alpha)}{(c_{j\eta}+\beta)^{n_{j\eta}+\alpha}}}{\frac{\Gamma(n_{j\eta_L}+\alpha)}{(c_{j\eta_L}+\beta)^{n_{j\eta_L}+\alpha}}\frac{\Gamma(n_{j\eta_R}+\alpha)}{(c_{j\eta_R}+\beta)^{n_{j\eta_R}+\alpha}}} 
\\&= \left(\frac{\beta^{\alpha}}{\Gamma(\alpha)}  \right)^{-1}\frac{\Gamma(n_{j\eta}+\alpha)}{\Gamma(n_{j\eta_L}+\alpha) \Gamma(n_{j\eta_R}+\alpha)}\frac{(c_{j\eta_L}+\beta)^{n_{j\eta_L}+\alpha}(c_{j\eta_R}+\beta)^{n_{j\eta_R}+\alpha}}{(c_{j\eta}+\beta)^{n_{j\eta}+\alpha}}\end{align*}

\subsection{CHANGE Proposal}
 This proposal randomly picks an internal node and randomly reassigns to it a splitting rule.
 
 Let $\eta$ be the picked internal node having rule $y<c$ and children denoted as $\eta_R$ and $\eta_L$. We assume that $\widetilde{y}<\widetilde{c}$ is its new assigned rule in the proposed tree, $T_j^{*}$. Following \citet{JSSv070i04}, for simplicity we are restricted to picking an internal node having two terminal nodes as children. 

\paragraph*{Transition Ratio}
 It holds that: 
\begin{enumerate}[label=(\roman*)]
\item{
$q(T_j^{*}|T_j^{(t)})=$ P(CHANGE) \\   
\hspace*{20mm}$\times$ P(selecting an internal node $\eta$ to change) \\ 
\hspace*{20mm}$\times$ P(selecting the new available dimension $\widetilde{y}$ to split on) \\ 
\hspace*{20mm}$\times$ P(selecting the new splitting value $\widetilde{c}$ given the chosen dimension $\widetilde{y}$)
}
\item{
$q(T_j^{(t)}|T_j^{*})=$ P(CHANGE) \\   
\hspace*{20mm}$\times$ P(selecting the node $\eta$ to change)\\ 
\hspace*{20mm}$\times$ P(selecting the dimension $y$ to split on) \\  \hspace*{20mm}$\times$ P(selecting the splitting value c given the chosen dimension $y$)
}
\end{enumerate}
 Thus the Transition Ratio is 
\begin{equation*}
\TR=\frac{P(\text{selecting $c$ to split on given the chosen dimension $y$})}{
P(\text{selecting $\widetilde{c}$ to split on given the chosen dimension $\widetilde{y}$})}
\end{equation*}

\paragraph*{Tree Structure Ratio} 
 The two trees differ in the splitting rule at node $\eta$. Thus we have that 
\begin{align*}
\TSR &=\frac{P(T_j^{*})}{P(T_j^{(t)})}=\frac{p_{\text{SPLIT}}(\eta) \, p_{\text{RULE}}(\eta|T_j^{*})}{p_{\text{SPLIT}}(\eta) 
\, p_{\text{RULE}}(\eta|T_j^{(t)})}\\
&= \frac{P(\text{selecting $\widetilde{y}$}) \, P(\text{selecting $\widetilde{c}$ given $\widetilde{y}$)}}{P(\text{selecting $y$}) \, P(\text{selecting  $c$ given $y$})}\\
&=\frac{P(\text{selecting $\widetilde{c}$ given $\widetilde{y}$)}}{P(\text{selecting  $c$ given $y$})}.
\end{align*}
It then follows that $\TR \times \TSR=1$, and hence only the likelihood ratio needs to be found to obtain the Hastings ratio.  

\paragraph*{Likelihood Ratio}
 Let $n_L^{*}=n_{j\eta_L}^{(T_j^{*})}$,  $n_R^{*}=n_{j\eta_R}^{(T_j^{*})}$,  $c_L^{*}=c_{j\eta_L}^{(T_j^{*})}$,  $c_R^{*}=c_{j\eta_R}^{(T_j^{*})}$,  $n_L^{(t)}=n_{j\eta_L}^{(T_j^{(t)})}$,  $n_R^{(t)}=n_{j\eta_R}^{(T_j^{(t)})}$,  $c_L^{(t)}=c_{j\eta_L}^{(T_j^{(t)})}$ and $c_R^{(t)}=c_{j\eta_R}^{(T_j^{(t)})}$, where $(T_j^{*})$ and $(T_j^{(t)})$ indicate that the corresponding quantities are related to the tree $T_j^{*}$ and $T_J^{(t)}$ respectively. Following the previous proposals, the likelihood ratio is
\begin{align*} 
\LR&=\frac{\frac{\Gamma(n_L^{*}+\alpha)}{(c_L^{*}+\beta)^{n_L^{*}+\alpha}}\frac{\Gamma(n_R^{*}+\alpha)}{(c_R^{*}+\beta)^{n_R^{*}+\alpha}}}{\frac{\Gamma(n_L^{(t)}+\alpha)}{(c_L^{(t)}+\beta)^{n_L^{(t)}+\alpha}}\frac{\Gamma(n_R^{(t)}+\alpha)}{(c_R^{(t)}+\beta)^{n_R^{(t)}+\alpha}}}  
&=\frac{(c_L^{(t)}+\beta)^{n_L^{(t)}+\alpha} \, (c_R^{(t)}+\beta)^{n_R^{(t)}+\alpha}}
{(c_L^{*}+\beta)^{n_L^{*}+\alpha} \, (c_R^{*}+\beta)^{n_R^{*}+\alpha}} \,
\frac{\Gamma(n_L^{*}+\alpha) \, \Gamma(n_R^{*}+\alpha)}{\Gamma(n_L^{(t)}+\alpha) \, \Gamma(n_R^{(t)}+\alpha)}.\end{align*}


\section{The Poisson Process conditional likelihood}
\label{app:remark1}
 Let us consider a finite realization of an inhomogeneous Poisson process with $n$ points $\vs$. Given the tree components $(T,\Lambda)$, and approximating the intensity of a point 
$s_i\in S$ by a product of $m$ trees $\lambda(s_i)=\prod_{j=1}^m g(s_i;T_j,\Lambda_j)$, the likelihood is: 
\begin{align}
P(\vs|\Lambda,T) &=
 \prod_{i=1}^n\lambda(s_i)\exp\left(-\int_{S}\lambda(s)ds\right) \nonumber\\
&= \prod_{i=1}^n\prod_{j=1}^mg(s_i;T_j,\Lambda_j)
\exp \left(-\int_{S}\prod_{j=1}^{m}g(s;T_j,\Lambda_j)ds \right).
 \label{e10}
\end{align}
 The first term of the above equation can be written as follows
\begin{align*}\prod_{i=1}^n\prod_{j=1}^mg(s_i;T_j,\Lambda_j)&=\prod_{i=1}^n\prod_{j=1,j\neq h}^mg(s_i;T_j,\Lambda_j)g(s_i;T_h,\Lambda_h)\\
&= \prod_{i=1}^n\prod_{j=1,j \neq h}^mg(s_i;T_j,\Lambda_j) \left(\prod_{i=1}^ng(s_i;T_h,\Lambda_h)  \right)=c_h\prod_{t=1}^{b_h}\lambda_{ht}^{n_{ht}}
\end{align*}
 where $c_h=\prod_{i=1}^n \prod_{j=1,j \neq h}^mg (s_i;T_j,\Lambda_j)$ and $n_{ht}$ is the cardinality of the set $\{i:s_i\in \Omega_{ht}\}$.  

 The exponential term of~\eqref{e10} can be expressed as:
\begin{align*}
\exp \left(-\int_{S}\prod_{j=1}^mg(s;T_j,\Lambda_j)ds  \right) &=\exp \left(-\int_{S}\prod_{j=1,j \neq h}^mg(s;T_j,\Lambda_j)g(s;T_h,\Lambda_h)  \right)\\
&=\exp \left(-\int_{S}\prod_{j=1,j \neq h}^m g(s;T_j,\Lambda_j) \left(\sum_{t=1}^{b_h}\lambda_{ht}I(s \in \Omega_{ht})  \right)ds  \right)\\
&=\exp \left(-\int_{S}\sum_{t=1}^{b_h}\lambda_{ht}\prod_{j=1,j\neq h}^mg(s;T_j,\Lambda_j)I(s\in\Omega_{ht})ds  \right)\end{align*}
 Tonelli’s theorem allows the change of order between summation and integral.
\begin{align*}
\exp \left(-\int_{S}\prod_{j=1}^mg(s;T_j,\Lambda_j)ds \right) 
&= \exp \left(-\sum_{t=1}^{b_h}\lambda_{ht}\int_{S}\prod_{j=1,j \neq h}^mg(s;T_j,\Lambda_j)I(s \in \Omega_{ht})ds   \right) \\  &= \exp \left( -\sum_{t=1}^{b_h}\lambda_{ht}c_{ht}  \right) 
\end{align*}
 where 
\begin{equation*}
c_{ht}= \int_{S} \left( \prod_{j=1,j \neq h}^{m}g(s;T_j,\Lambda_j)  \right)I(s\in \Omega_{ht})ds .
\end{equation*}

 Let $T_{(h)}=\{T_j\}_{j=1,j\neq h}^{m}$ be an ensemble of trees not including the tree $T_h$ that defines  the global partition 
$\{\overline{\Omega}_{k}^{(h)}\}_{k=1}^{K(T_{(h)})}$ by merging all cuts in $\{T_j\}_{j=1,j\neq h}^{m}$. Giving, 
\begin {equation*}{\prod_{j=1,j\neq h}^{m}g(s;T_j,\Lambda_j)=\sum_{k=1}^{K(T_h)}\overline{\lambda}_k^{(h)}I(s \in \overline{\Omega}_k^{(h)}}) \end{equation*} where
\begin{equation*} \overline{\lambda}_{k}^{(h)}=\prod_{t=1,t\neq h}^{m}\prod_{l=1}^{b_t}\lambda_{tl}^{I(\Omega_{tl}\cap\overline{\Omega}_{k}^{(h)}\neq 0)},\end{equation*}
leading to the following expression for $c_{ht}$,
\begin{align*}
c_{ht}&=\int_{S} \left( \prod_{j=1,j \neq h}^{m}g(s,T_j,\Lambda_j)  \right)I(s\in \Omega_{ht})ds=
\int_{S} \left(\sum_{k=1}^{K(T_{(h)})}\overline{\lambda}_k^{(h)}I(s \in \overline{\Omega}_k^{(h)})   \right)I(s \in \Omega_{ht})ds \\ 
&= \sum_{k=1}^{K(T_{(h)})}\overline{\lambda}_k^{(h)}\int_{S}I(s \in \overline{\Omega}_{k}^{(h)}\cap \Omega_{ht})ds 
=\sum_{k=1}^{K(T_{(h)})}\overline{\lambda}_k^{(h)}|\overline{\Omega}_{k}^{(h)}\cap \Omega_{ht}|,
\end{align*}
 where $ |\overline{\Omega}_{k}^{(h)} \cap \Omega_{ht}|$ is the volume of the region $ \Omega_{k}^{(h)}\cap \Omega_{ht}$. 
 Hence the conditional likelihood can be written as follows 
\begin{equation*} 
P(\vs|\Lambda,T)=c_h\prod_{t=1}^{b_h}\lambda_{ht}^{n_{ht}}e^{-\lambda_{ht}c_{ht}}.
\end{equation*}

\section{The conditional integrated likelihood}
\label{app:remark2}
The conditional integrated likelihood is given by
\begin{align*}
P(\vs|T_h,T_{(h)},\Lambda_{(h)})&=\int_0^\infty P(\vs,\Lambda_h|T_h,T_{(h)},\Lambda_{(h)})d\Lambda_{h}\\
&=\int_0^\infty P(\vs|\Lambda,T) \, P(\Lambda_{h}|T_{h},T_{(h)},\Lambda_{(h)})d\Lambda_{h}\\
&= c_h \int_0^\infty \hdots \int_0^\infty \prod_{t=1}^{b_h}\lambda_{ht}^{n_{ht}}e^{-\lambda_{ht}c_{ht}}\prod_{t=1}^{b_h}\frac{\beta^{\alpha}}{\Gamma(\alpha)}e^{-\beta\lambda_{ht}}\lambda_{ht}^{\alpha-1}d\lambda_{h1}\hdots d\lambda_{hb_h}  \\
&=c_h  \left( \frac{\beta^{\alpha}}{\Gamma(\alpha)}  \right)^{b_h}\prod_{t=1}^{b_h}\int_0^\infty \lambda_{ht}^{n_{ht}+\alpha-1}e^{-(c_{ht}+\beta)\lambda_{ht}}d\lambda_{ht}\\&=c_h \left(\frac{\beta^{\alpha}}{\Gamma(\alpha)}  \right)^{b_h}\prod_{t=1}^{b_h}\frac{\Gamma(n_{ht}+\alpha)}{(c_{ht}+\beta)^{n_{ht}+\alpha}}
\end{align*}

\newpage



\maketitle
\begin{center}
 {\Large  Supplementary Material}
\end{center}

\section{The model for the case of one tree}
 The proposed model for considering only one tree can be written as follows
\begin {align*}
\lambda(s_i) &=g(s_i;T,\Lambda) =\sum_{k=1}^{b}\lambda_{k} \, I(s_i \in \Omega_{k}) \\
T&\sim \text{heterogeneous Galton-Watson process for a partition of $S$} \\
\lambda_{k}|T& \sim \mbox{Gamma}(\alpha,\beta)
\end{align*}
underpinned by a tree-shaped partition $T=\{\Omega_{k}\}_{k=1}^{b}$ where $b$ is the number of terminal nodes in the tree $T$. Each leaf node k associated to region $\Omega_{k}$ is linked with a 
parameter $\lambda_{k}$. All parameters $\lambda_{k}$ are collected in the vector $\Lambda=(\lambda_{1},\lambda_{2},..,\lambda_{b})$.  
 The parameters of the model are 
\begin{enumerate}
\item{the regression tree $T$}
\item{the parameters $\Lambda=(\lambda_{1},\lambda_{2},..,\lambda_{b})$.}
\end{enumerate}
 We assume that the leaf parameters are independent, i.e., $P(\Lambda|T)=\prod_{k=1}^{b}P(\lambda_k|T).$

\subsection{Poisson Process conditional likelihood}
 The conditional likelihood of a finite realization of an inhomogeneous Poisson process with $n$ points $s_1, \ldots , s_n$ is derived by describing $\lambda(s)$ using one tree $(\Lambda,T)$ as: $\lambda(s)= g(s;T,\Lambda)$.
\begin{align}
 \label{e2} 
 P(s_1, \ldots , s_n|\Lambda,T)=\prod_{i=1}^n\lambda(s_i)\exp \left(-\int_{S}\lambda(s)ds  \right)=\prod_{i=1}^ng(s_i;T,\Lambda)\exp  \left(-\int_{S}g(s;T,\Lambda)ds  \right) .
\end{align}

 The first term of the above equation can be written as follows
\begin{align*}\prod_{i=1}^ng(s_i;T,\Lambda)&=\prod_{k=1}^b\lambda_{k}^{n_k}
\end{align*}
 where $n_{k}$ is the cardinality of the set $\{i:s_i\in \Omega_{k}\}$.  

 The exponential term of~\eqref{e2} can be expressed as follows
\begin{align*}
\exp \left(-\int_{S}g(s;T,\Lambda)ds  \right) &=\exp \left(-\int_{S}\sum_{k=1}^b \lambda_{k}I(s \in \Omega_{k})ds  \right)\\
&=\exp \left(-\sum_{k=1}^b\lambda_k\int_{S}I(s \in \Omega_{k})ds  \right)=\exp \left(-\sum_{k=1}^b\lambda_k|\Omega_k|  \right)\end{align*}

 Hence the conditional likelihood can be written as 
\begin{equation} \label{eq10}
P(s_1, \ldots , s_n|\Lambda,T)=\prod_{k=1}^{b}\lambda_{k}^{n_{k}} \, e^{-\lambda_{k}|\Omega_{k}|},
\end{equation}
 where $|\Omega_k|$ is the volume of the region $\Omega_k$.

\subsection{Inference Algorithm}
 Inference on the model parameters $(\Lambda,T)$ induces sampling from the posterior $P(\Lambda,T|s_1,...,s_n)$. A Metropolis Hastings within Gibbs sampler (Algorithm~\ref{alg3}) is  proposed for sampling from the posterior $P(\Lambda,T|s_1, \ldots , s_n)$. Noting that, 
\begin{equation*} 
P(\Lambda,T|s_1, \ldots , s_n)=P(\Lambda|T,s_1, \ldots , s_n) \, P(T|s_1,\ldots, s_n) 
\end{equation*} and 
\begin{equation*} P(\Lambda|T,s_1, \ldots , s_n)\propto P(s_1,...,s_n|\Lambda,T) \, P(\Lambda|T) \propto \prod \limits_{k=1}^{b}\lambda_{k}^{n_{k}+\alpha-1}e^{-(|\Omega_{k}|+\beta)\lambda_{k}}, \end{equation*}
a draw from $(T,\Lambda)|s_1, \ldots , s_n$ can be achieved in ($b$+1) successive steps:
\begin{itemize}
\item{sample $T|n,s_1, \ldots , s_n$ through Metropolis-Hastings Algorithm summarized in Algorithm \ref{alg4}}
\item{sample $\lambda_{k}|T,n,s_1, \ldots , s_n$ from a Gamma distribution with shape equal to $n_{k}+\alpha$ and rate equal to $|\Omega_k|+\beta$  for $k=1,..,b$.}
\end{itemize} 
Noting that 
\[P(T|s_1, \ldots , s_n)\propto P(s_1, \ldots , s_n|T) \, P(T)\] 
the integrated likelihood (integrating out the parameters $\Lambda$) is: 
\begin{align}
 P(s_1, \ldots , s_n|T)&=\int P(s_1, \ldots , s_n,\Lambda|T)d\Lambda=\int P(s_1, \ldots , s_n|\Lambda,T)P(\Lambda|T)d\Lambda  \nonumber \\
&=  \left( \frac{\beta^{\alpha}}{\Gamma(a)}  \right)^{b}\prod_{k=1}^b \int \lambda_k^{n_k+\alpha-1}e^{-(|\Omega_k|+\beta)\lambda_k}d\lambda_k \nonumber\\
&= \left( \frac{\beta^{\alpha}}{\Gamma(a)}  \right)^{b} \prod_{k=1}^b \frac{\Gamma(n_k+\alpha)}{(\beta+|\Omega_k|)^{n_k+\alpha}}.
\label{eq:integrated_onetree}
\end{align}
In the tree sampling Algorithm~\ref{alg4}, the transition kernel $q$ is chosen from the three proposals: GROW, PRUNE, CHANGE~\citep{chipman2010bart,JSSv070i04}, and Eq.~\eqref{eq:integrated_onetree} allows us to compute the Metropolis Hastings ratio to accept or reject the proposal.

\begin{algorithm}[H] 
\caption{Proposed Algorithm: Metropolis Hastings within Gibbs sampler} 
\label{alg3}
\begin{algorithmic}
\FOR{$t=1,2,3,..$ }
\STATE{Sample $T^{(t+1)}|s_1, \ldots , s_n$}
\FOR{$k=1$ to $b$}
\STATE {Sample $\lambda_{k}^{(t+1)}|s_1, \ldots , s_n,T^{(t+1)}$}
\ENDFOR
\ENDFOR
\end{algorithmic}
\end{algorithm}

\begin{algorithm}[H] 
\caption{Metropolis Hastings Algorithm for sampling from the posterior  $P(T|s_1, \ldots , s_n)$ } 
\label{alg4}
\begin{algorithmic}

\STATE{Generate a candidate value $T^{*}$ with probability $q(T^{*}|T^{(t)})$.}
\STATE {Set $T^{(t+1)}=T^{*}$ with probability 
\begin{equation*}
\alpha(T^{(t)},T^{*})= \min\left(1,\frac{q(T^{(t)}|T^{*})}{q(T^{*}|T^{(t)})} \frac{P(s_1, \ldots , s_n|T^{*})}{P(s_1, \ldots , s_n|T^{(t)})} \frac{P(T^{*})}{P(T^{(t)})} \right)
\end{equation*}
Otherwise, set $T^{(t+1)}=T^{(t)}$.
}
\end{algorithmic}
\end{algorithm}

\section{Simulation results on synthetic data with various number of sampling iterations}
 In this appendix we show that our algorithm works equally well for 10000 iterations by running three parallel chains, examining their convergence and assessing the performance of our algorithm via AAE and RMSE of computed estimates over various number of iterations. We also check the convergence of chains using the Gelman-Rubin criterion in all cases.
 
\subsection{One dimensional Poisson Process with stepwise intensity}
\begin{figure}[H] 
 \begin{subfigure}{8cm}
    \centering\includegraphics[width=6cm]{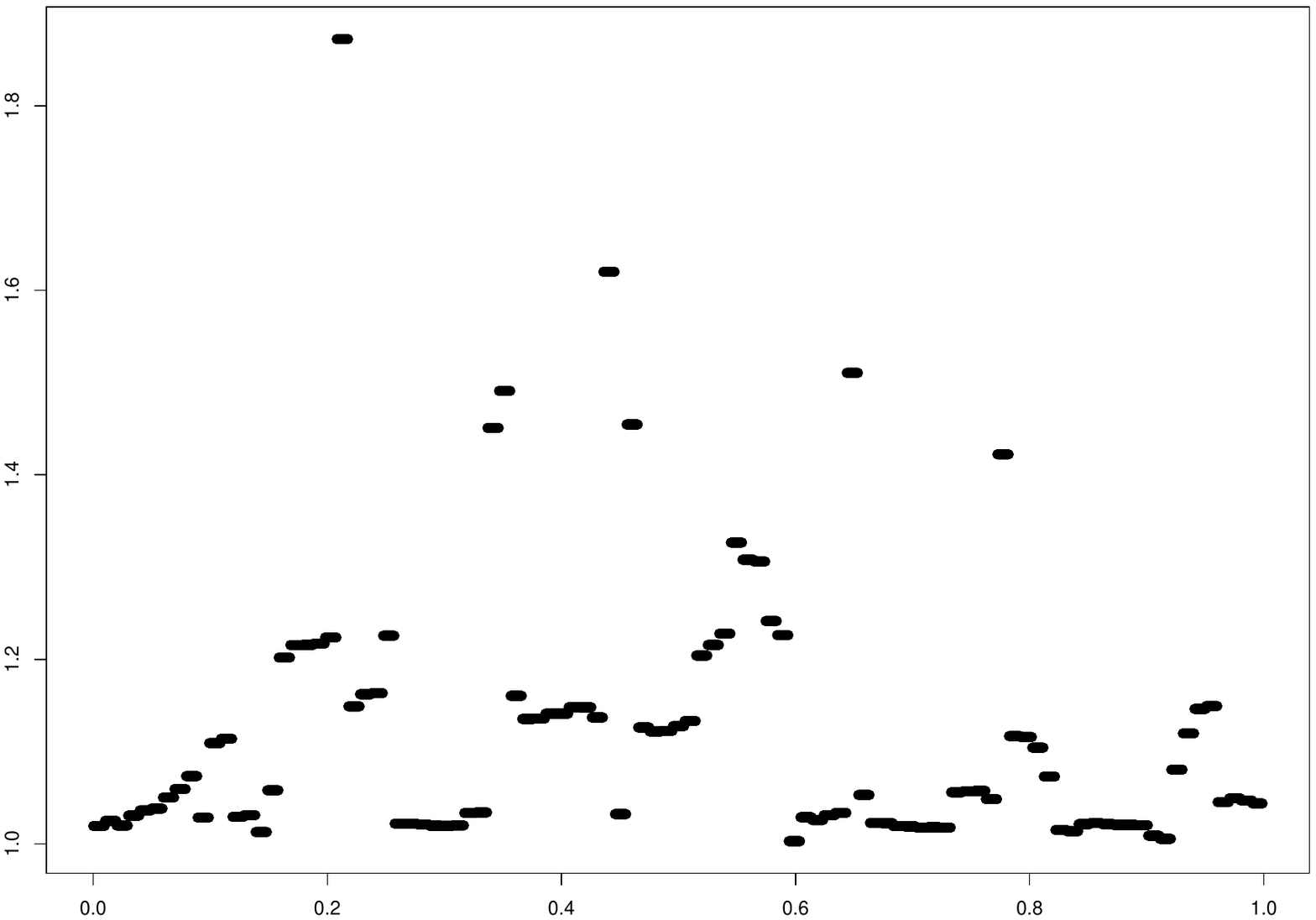}
    \caption{5 Trees and 10000 iterations }
  \end{subfigure}
\begin{subfigure}{8cm}
    \centering\includegraphics[width=6cm]{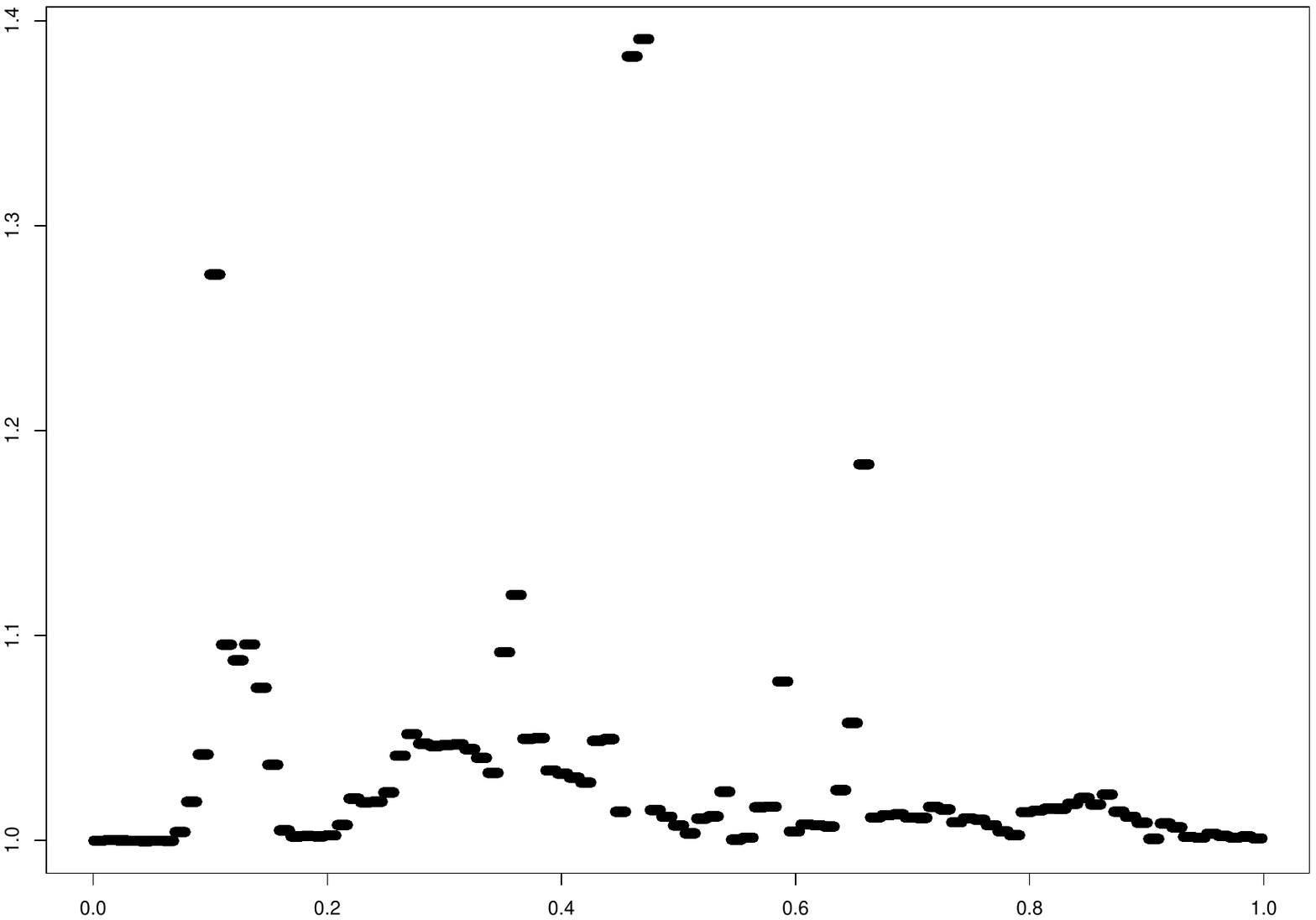}
    \caption{7 Trees and 10000 iterations }
  \end{subfigure}
  \begin{subfigure}{8cm}
    \centering\includegraphics[width=6cm]{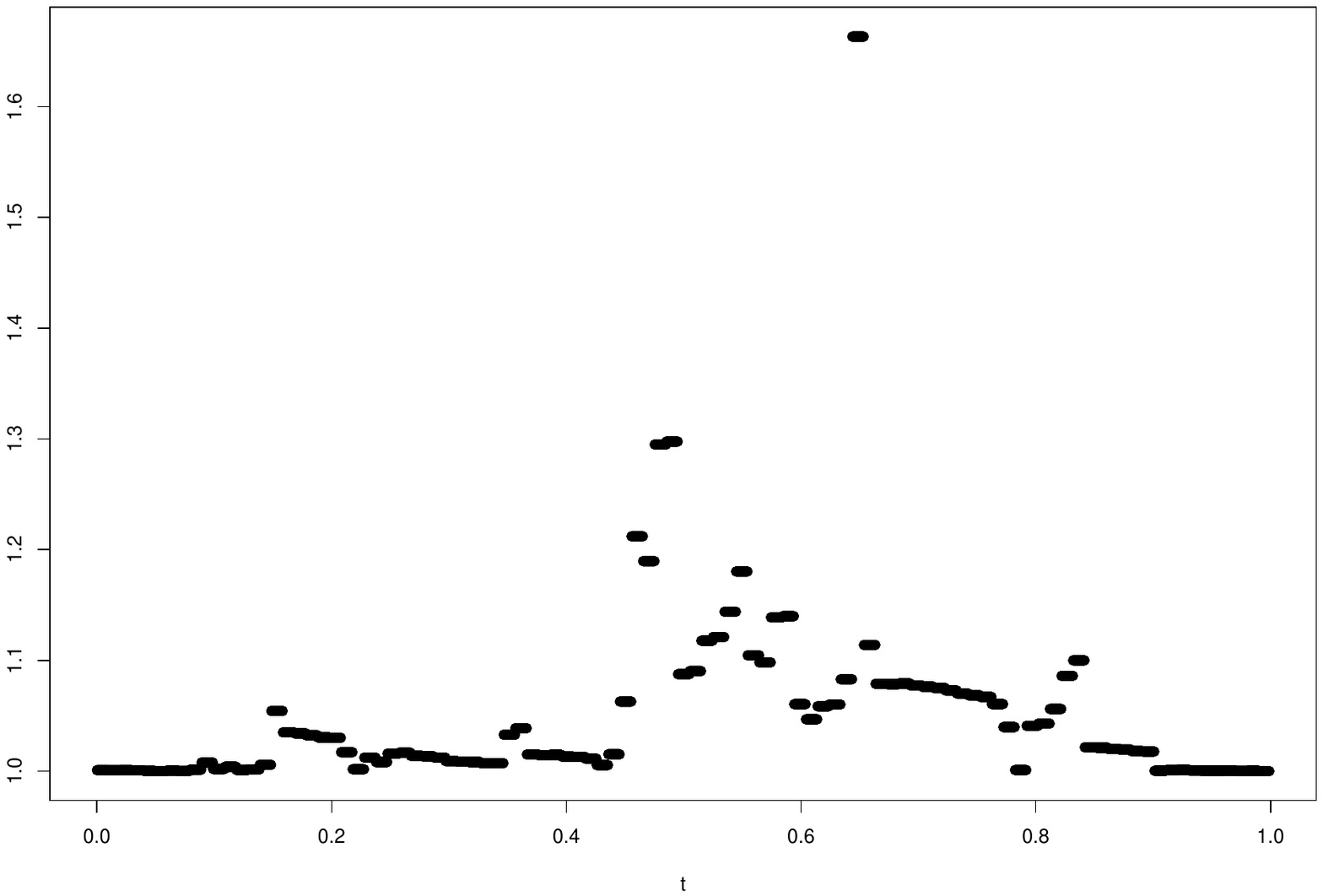}
    \caption{5 Trees and 50000 iterations }
  \end{subfigure}
\begin{subfigure}{8cm}
    \centering\includegraphics[width=6cm]{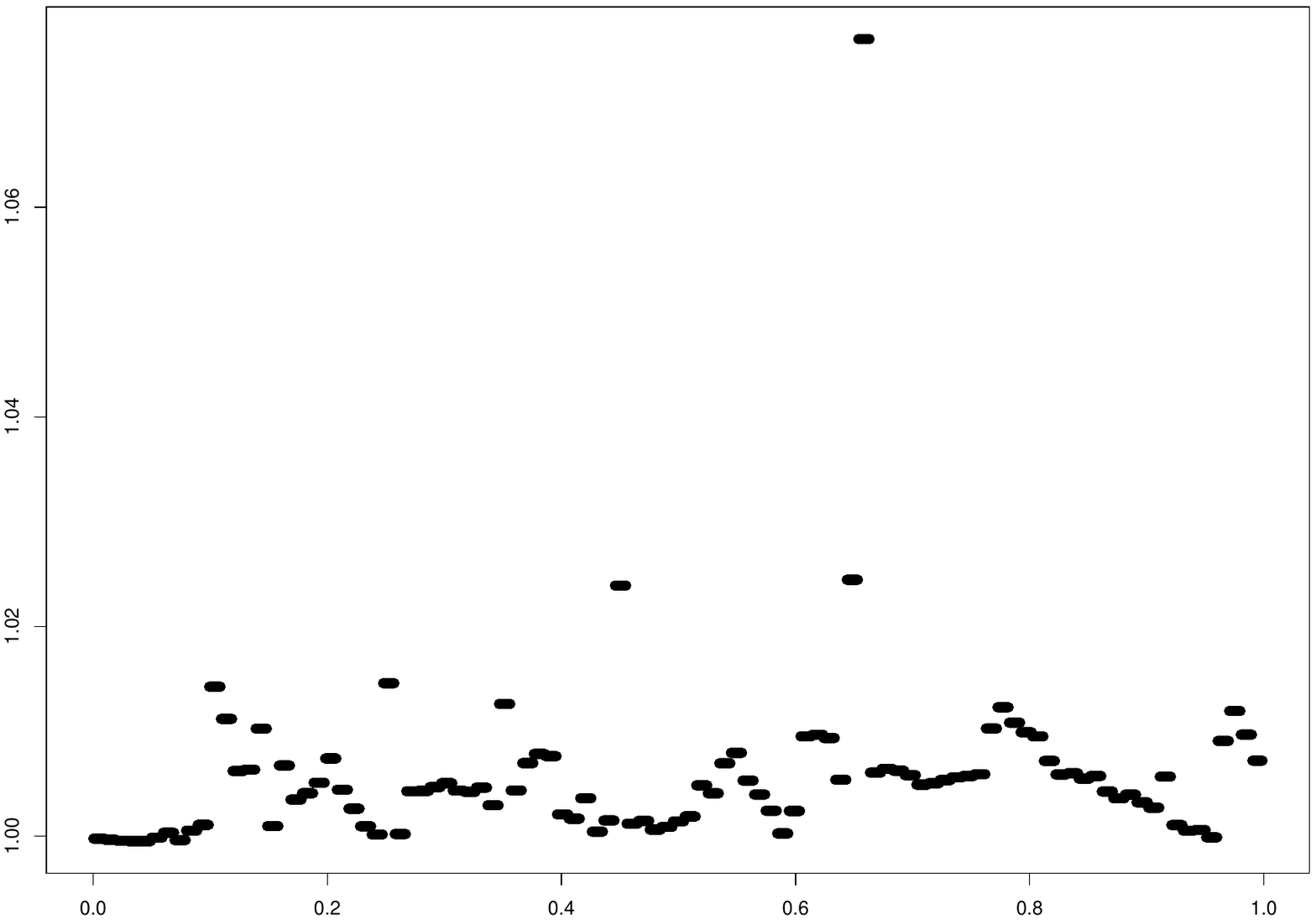}
    \caption{7 Trees and 50000 iterations }
  \end{subfigure}
  \begin{subfigure}{8cm}
    \centering\includegraphics[width=6cm]{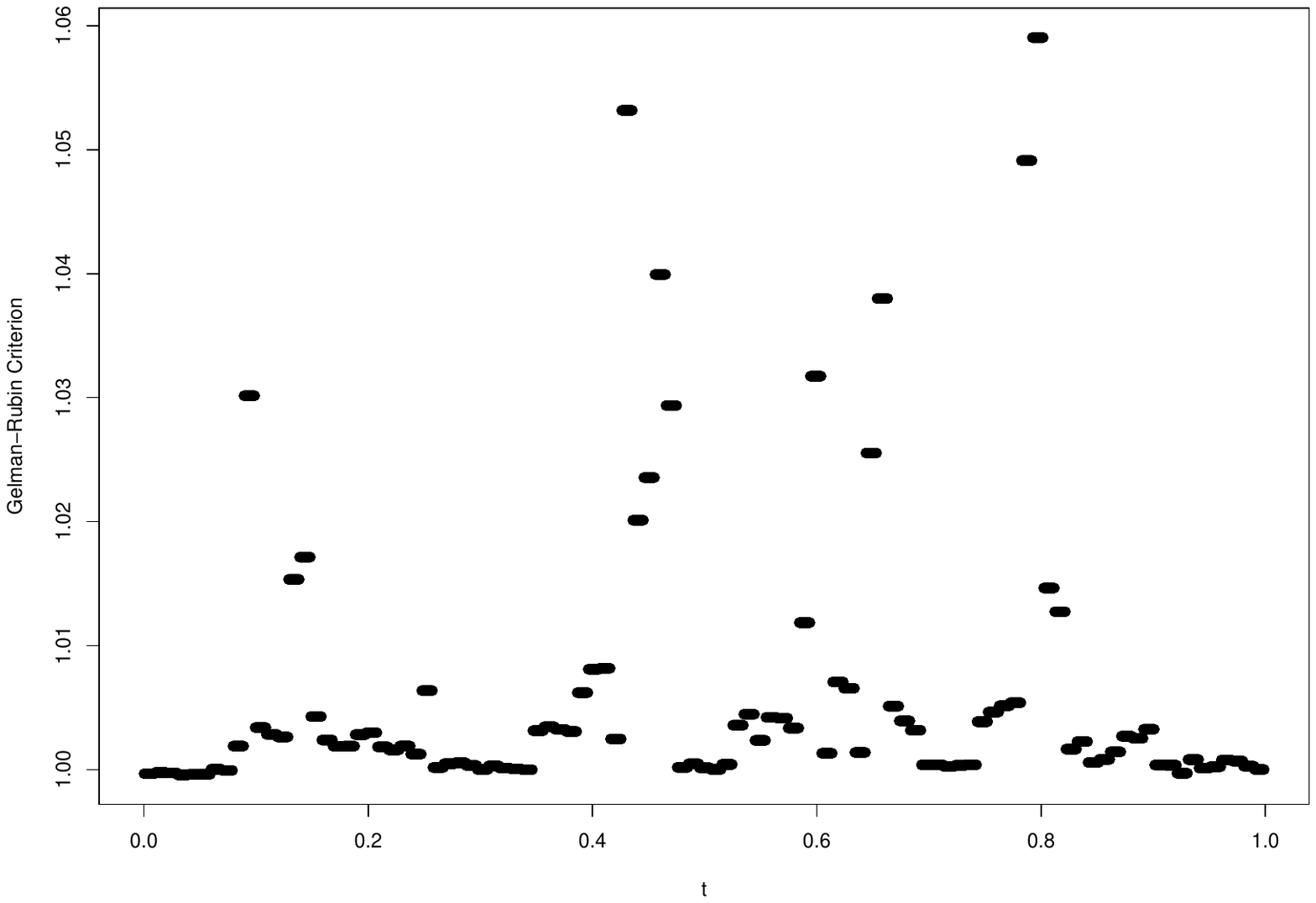}
    \caption{5 Trees and 200000 iterations }
  \end{subfigure}
\begin{subfigure}{8cm}
    \centering\includegraphics[width=6cm]{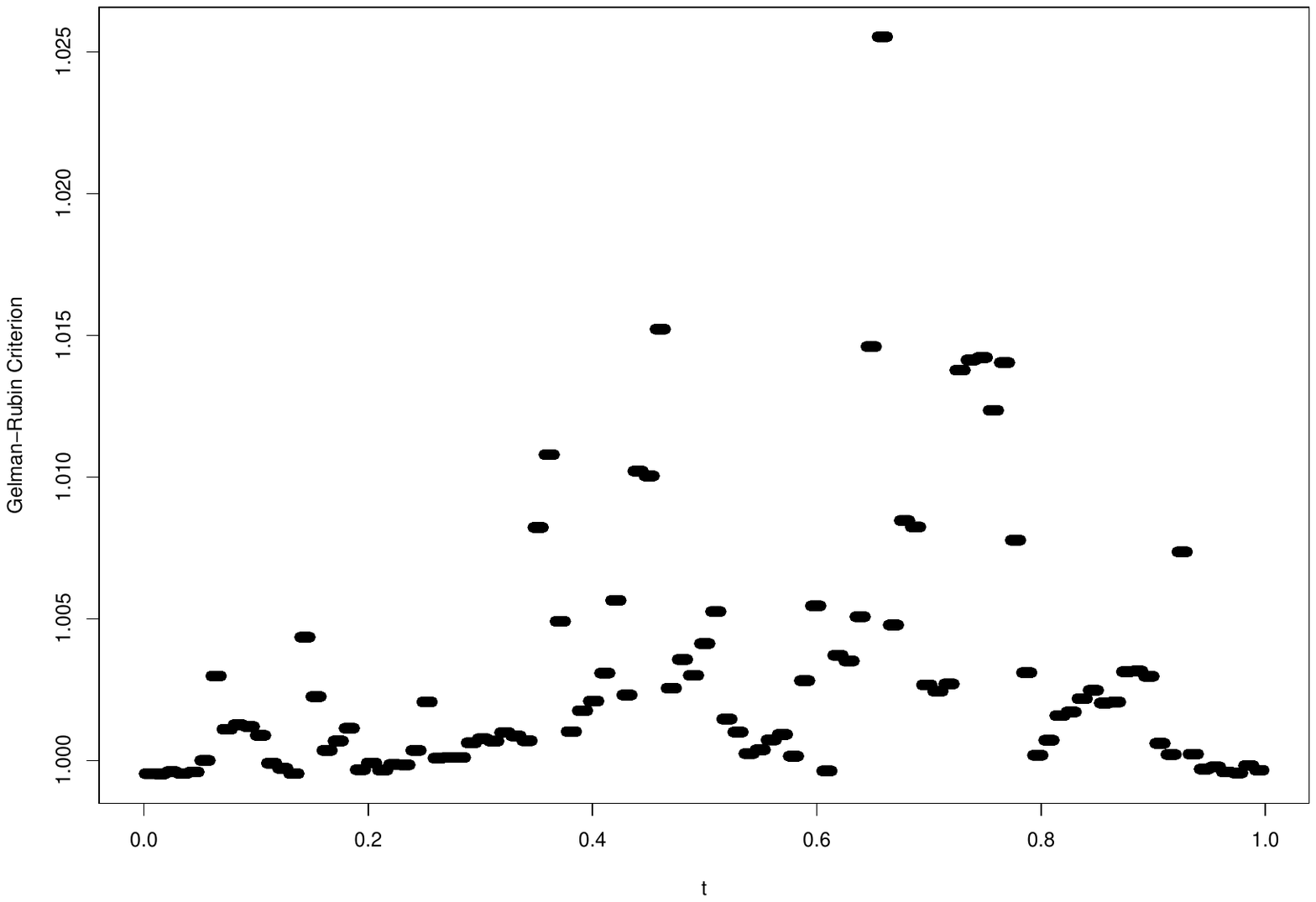}
    \caption{7 Trees and 200000 iterations }
  \end{subfigure}
\caption{The Gelman-Rubin Criterion for various number of iterations and trees.}
\label{figE1}
\end{figure}

Table~\ref{Table_E1} shows that there are no significant difference in errors increasing the number of iterations from 10000 to 200000.
Figure \ref{figE1} reveals that the chains work less well at points close to jumps for small number of iterations.

\begin{table}[H]
\begin{tabular}{ |p{1.2cm}||p{2cm}|p{2.2cm}|p{2.2cm}|p{2.2cm}|p{2.3cm}| }
 \hline
 \multicolumn{6}{|c|}{Proposed Algorithm}  \\
 \hline
 Number of  trees & Number of Iterations &AAE for Posterior Mean & AAE for Posterior Median & RMSE for Posterior Mean & RMSE for Posterior Median \\
 \hline
 5 & 10000 &284.61 & 274.3&588.88 &590.5   \\
 & 50000 &289.11 & 284.56&575.11 &579.17   \\
 & 200000 &279.88&269.81&572.94 &576.94   \\
 \hline
7 &10000&265.22&257.49 & 572.33&576.58\\
&50000&276.19&267.75 & 580.35&584.47\\
&200000&278.37&269.78 & 582.82&584.1\\
 \hline
\end{tabular}
\caption{Average Absolute Error and Root Mean Square Error for  various number of iterations and trees.}
\label{Table_E1}
\end{table}

\subsection{One dimensional Poisson Process with continuously varying intensity}
\begin{figure}[H] 
 \begin{subfigure}{8cm}
    \centering\includegraphics[width=6cm]{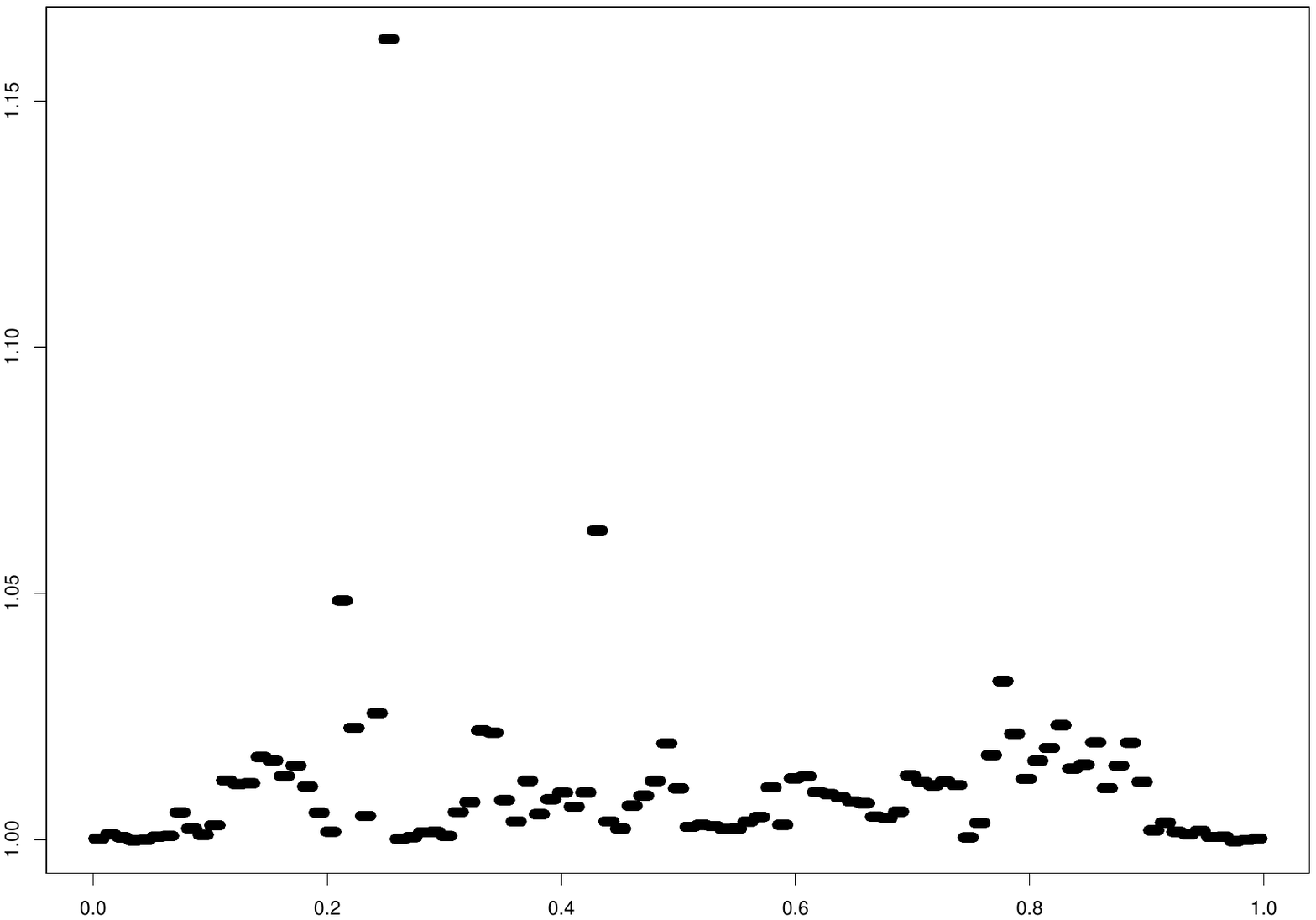}
    \caption{5 Trees and 10000 iterations }
  \end{subfigure}
\begin{subfigure}{8cm}
    \centering\includegraphics[width=6cm]{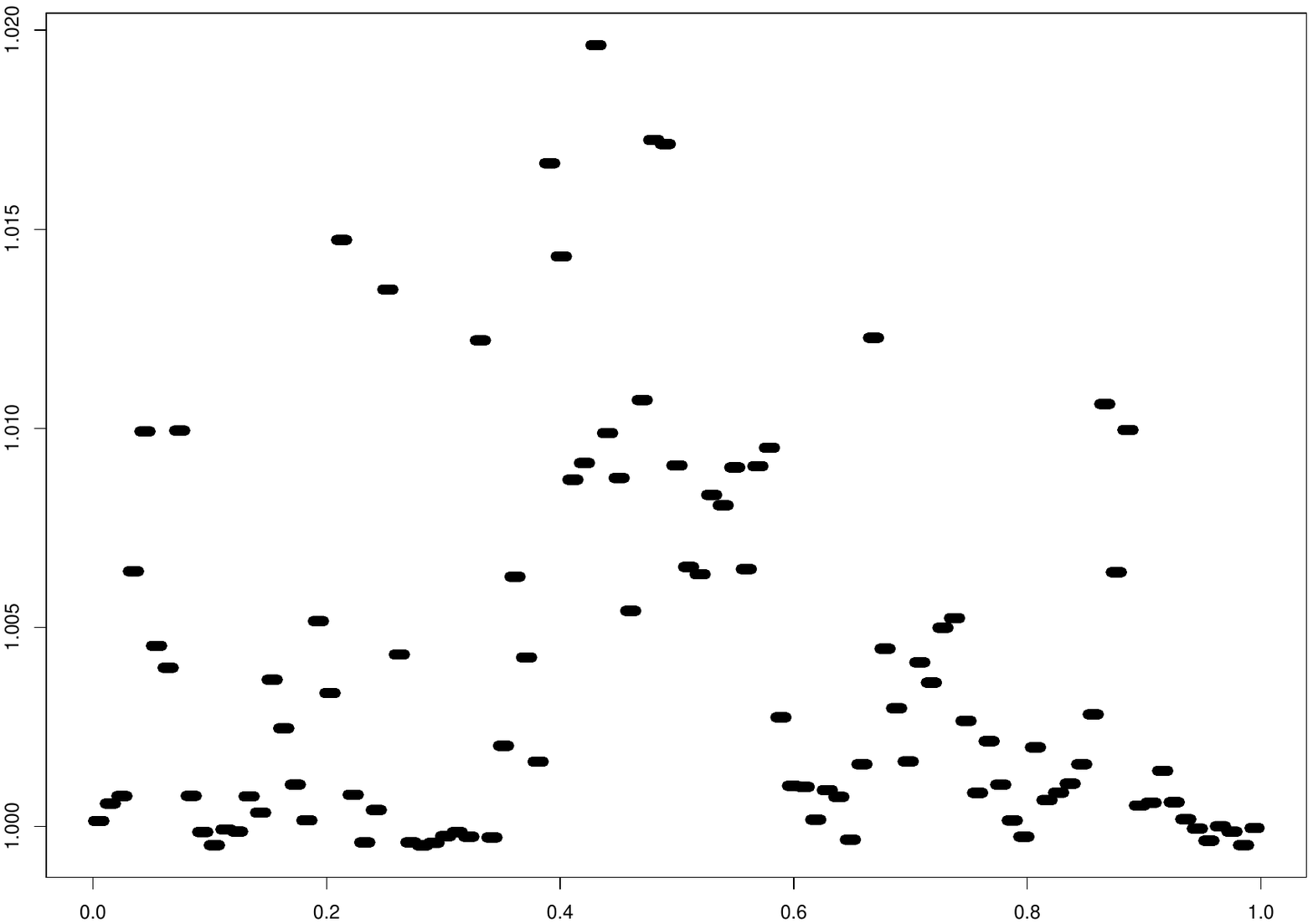}
    \caption{10 Trees and 10000 iterations }
  \end{subfigure}
  \begin{subfigure}{8cm}
    \centering\includegraphics[width=6cm]{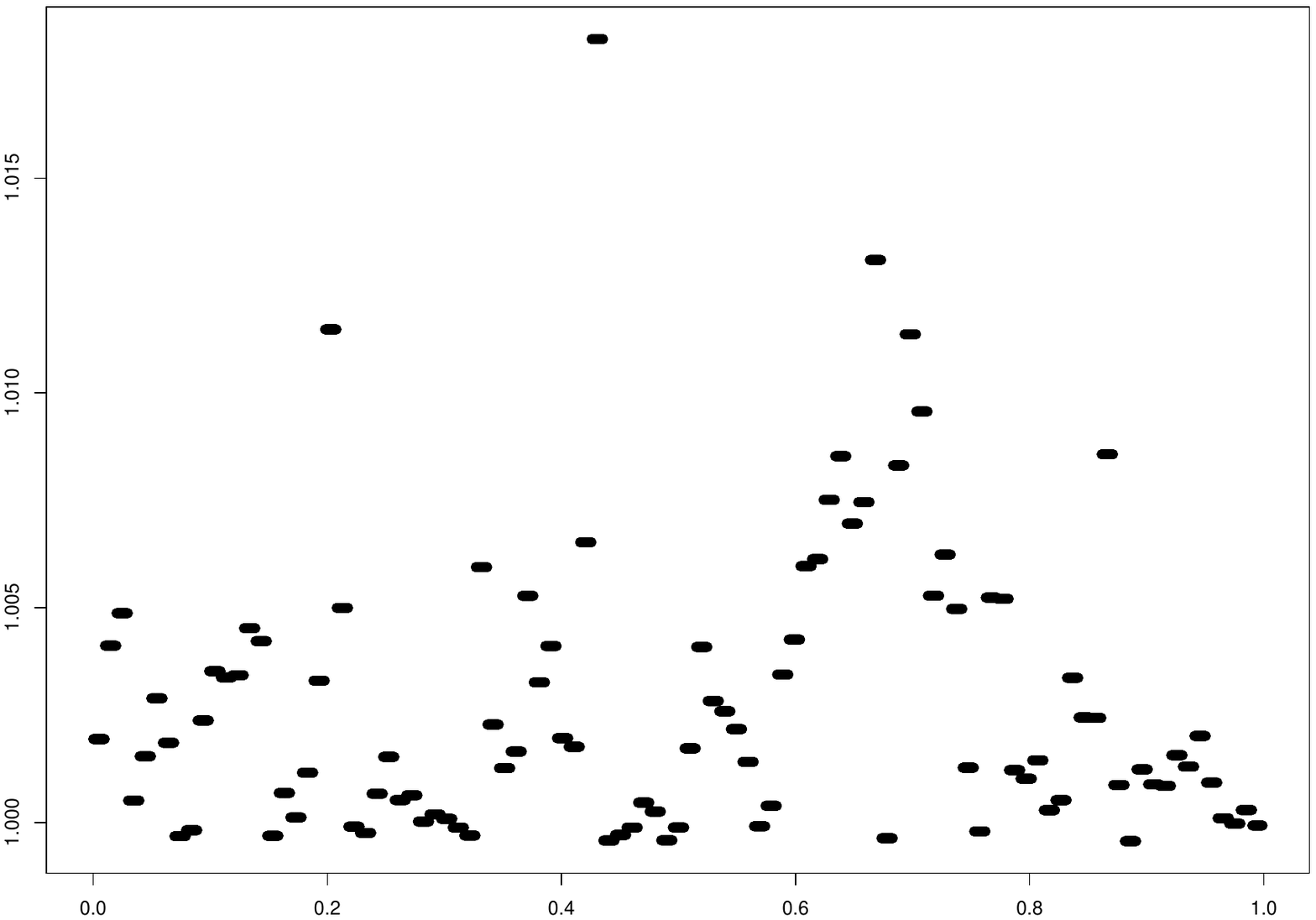}
    \caption{5 Trees and 50000 iterations }
  \end{subfigure}
\begin{subfigure}{8cm}
    \centering\includegraphics[width=6cm]{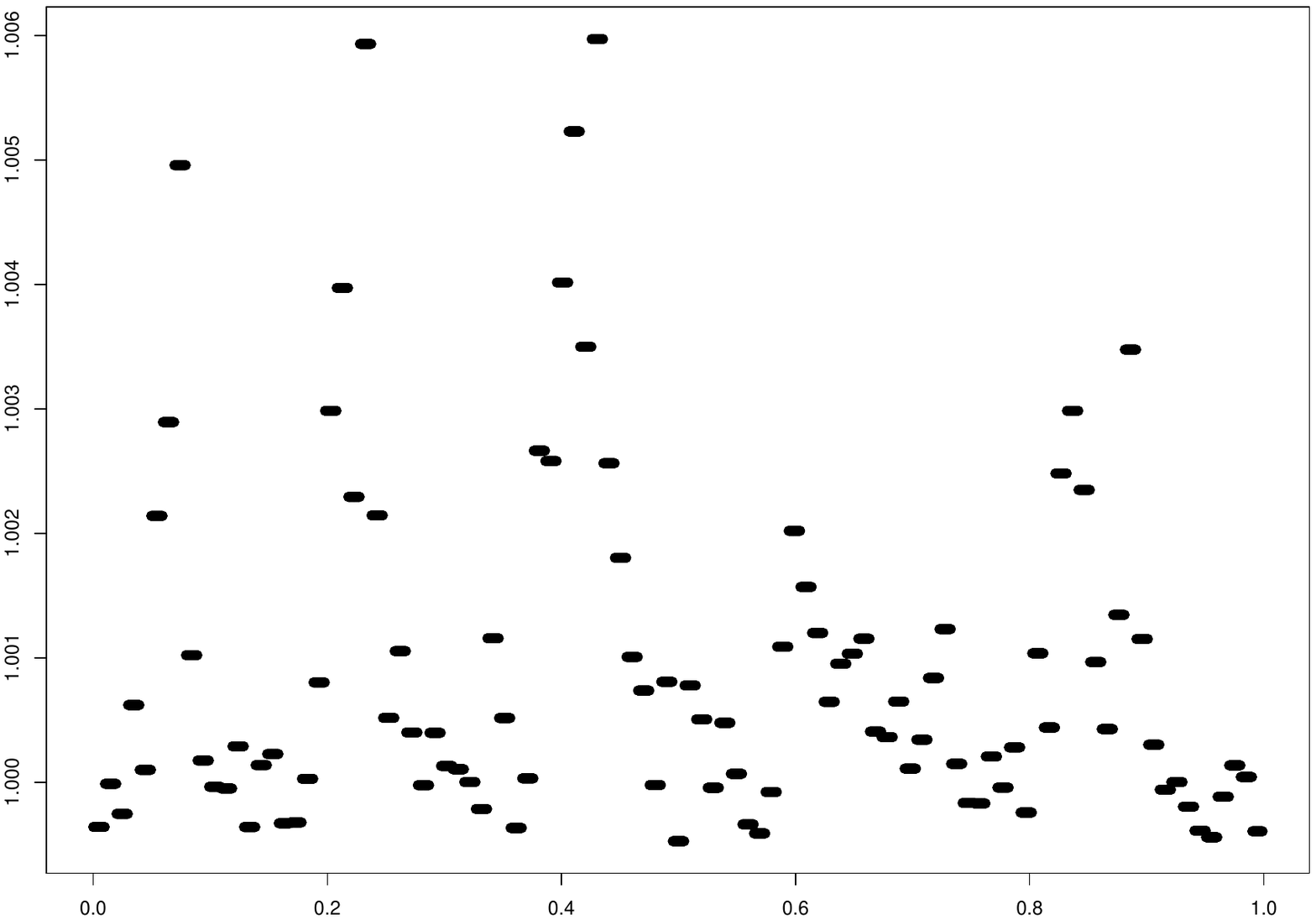}
    \caption{10 Trees and 50000 iterations }
  \end{subfigure}
  \begin{subfigure}{8cm}
    \centering\includegraphics[width=6cm]{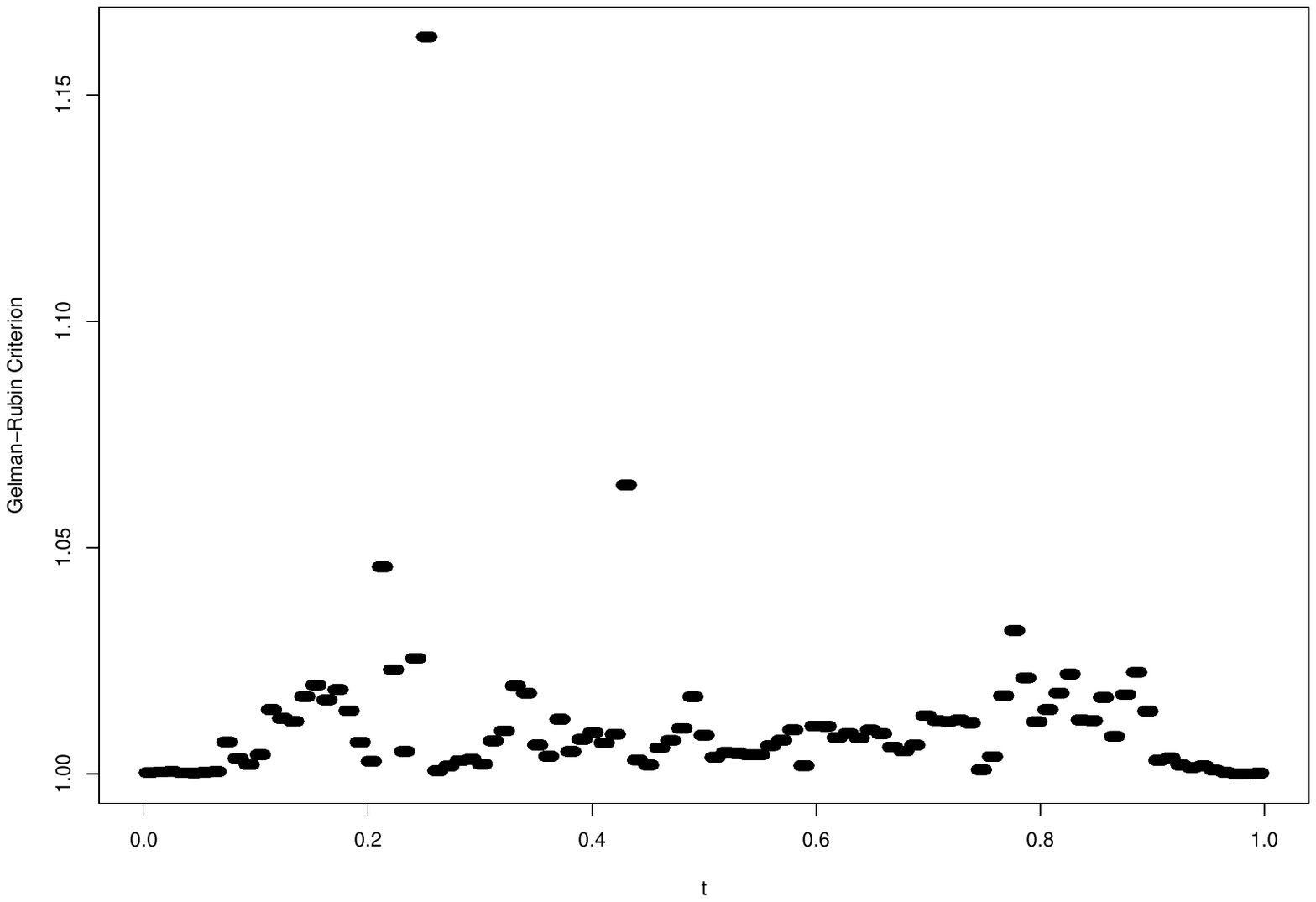}
    \caption{5 Trees and 100000 iterations }
  \end{subfigure}
\begin{subfigure}{8cm}
    \centering\includegraphics[width=6cm]{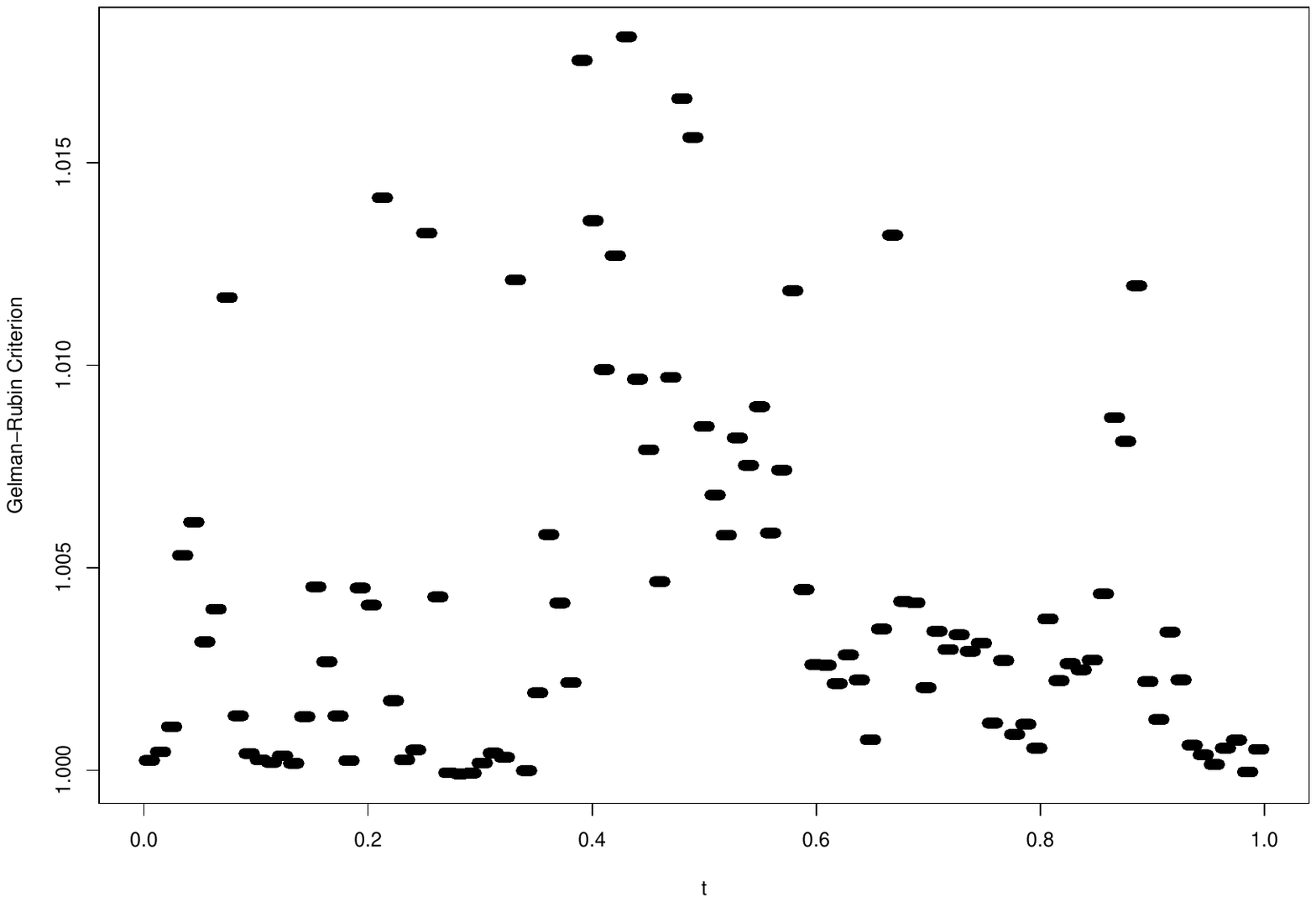}
    \caption{10 Trees and 100000 iterations }
  \end{subfigure}
\caption{The Gelman-Rubin Criterion for various number of iterations and trees.}
\label{figE2}
\end{figure}

Table \ref{Table_E2} shows that increasing the number of iterations does not change essentially the error for the synthetic data presented in Section~\ref{sec:cosine_intenstity}.
The convergence criterion indicates that even for small number of iterations, the chains converge for 10 trees. For 5 trees they converge for the majority of the range (Figure \ref{figE2}).

\begin{table}[H]
\begin{tabular}{ |p{1.2cm}||p{2cm}|p{2.2cm}|p{2.2cm}|p{2.2cm}|p{2.3cm}| }
 \hline
 \multicolumn{6}{|c|}{Proposed Algorithm}  \\
 \hline
 Number of  trees & Number of Iterations &AAE for Posterior Mean & AAE for Posterior Median & RMSE for Posterior Mean & RMSE for Posterior Median\\
 \hline
 5 & 10000 & 6.27 & 6.71 & 9.83 & 10.62   \\
 & 50000 &6.16 & 6.51&9.63 &10.42   \\
 & 100000 &6.14&6.38&9.52 &10.17   \\
 \hline
7 &10000&5.99& 6.03 &9.54&9.95\\
&50000&6.04&6.1 & 9.49&9.88\\
&100000&5.95&6.01 & 9.39&9.8\\
 \hline
\end{tabular}
\caption{Average Absolute Error and Root Mean Square Error for  various number of iterations and trees.}
\label{Table_E2}
\end{table}

\subsection{Two dimensional Poisson process with stepwise intensity function}

\begin{figure}[H] 
 \begin{subfigure}{8cm}
    \centering\includegraphics[width=6cm]{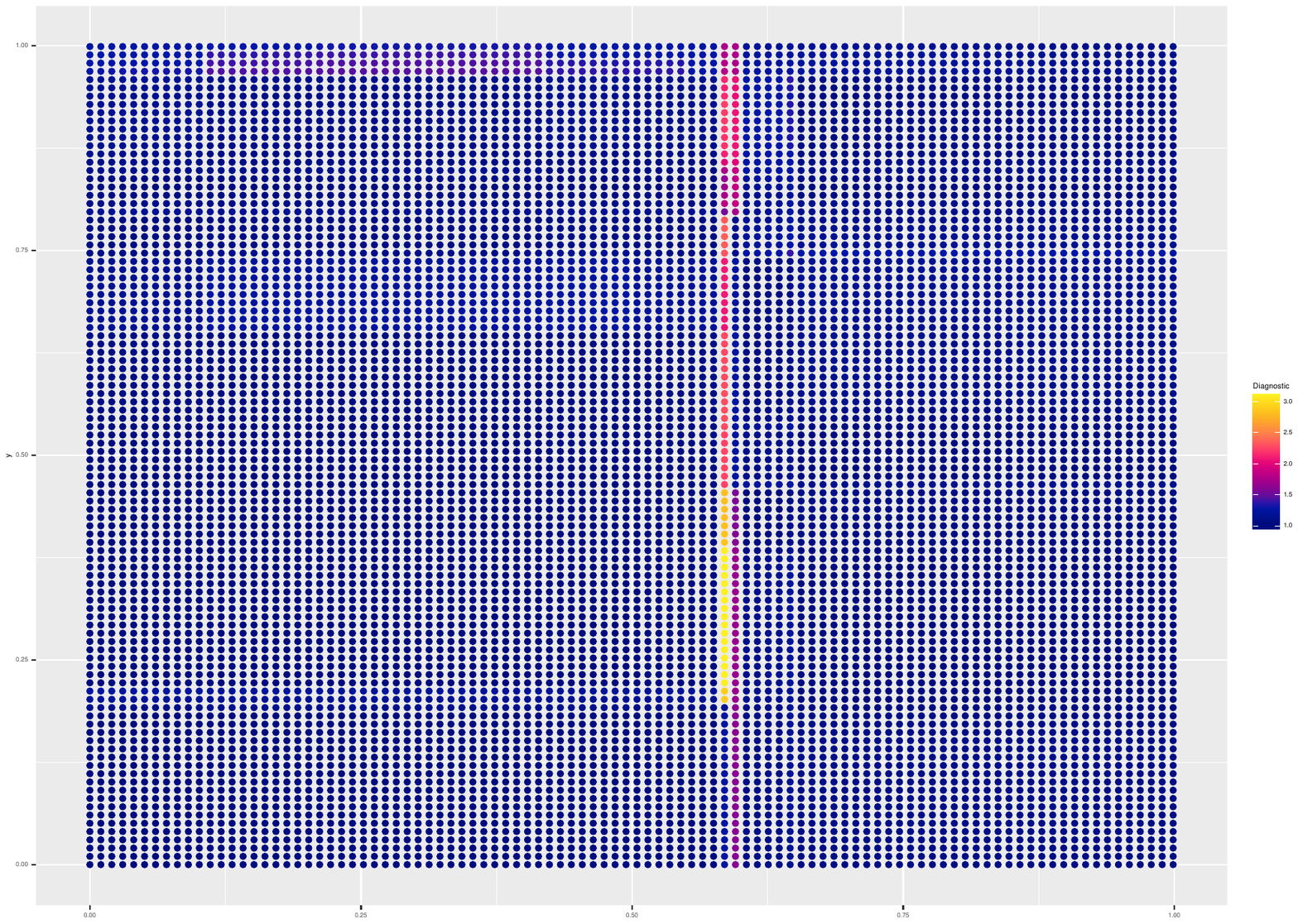}
    \caption{4 Trees and 10000 iterations }
  \end{subfigure}
  \begin{subfigure}{8cm}
    \centering\includegraphics[width=6cm]{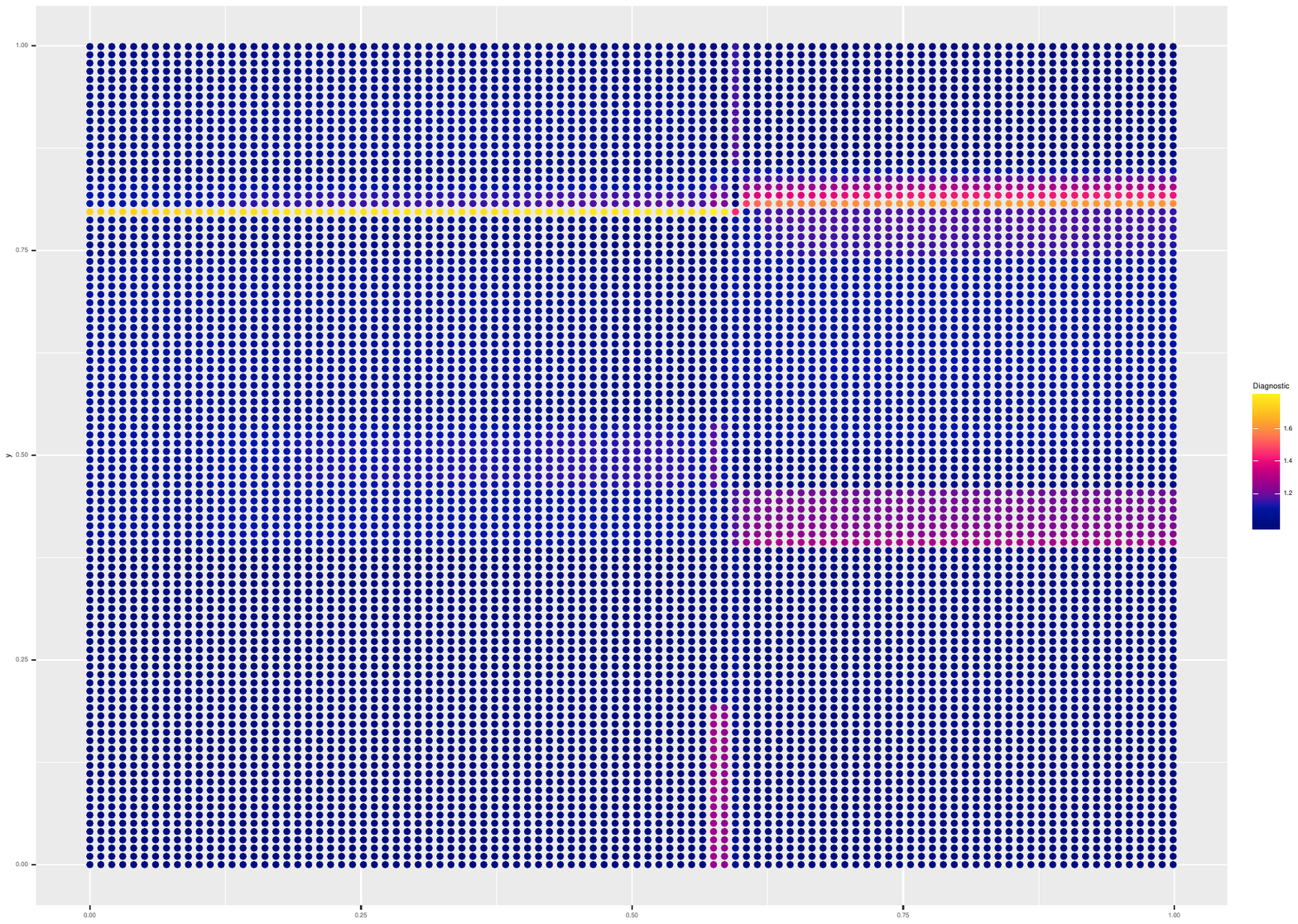}
    \caption{4 Trees and 50000 iterations }
  \end{subfigure}
\begin{subfigure}{8cm}
    \centering\includegraphics[width=6cm]{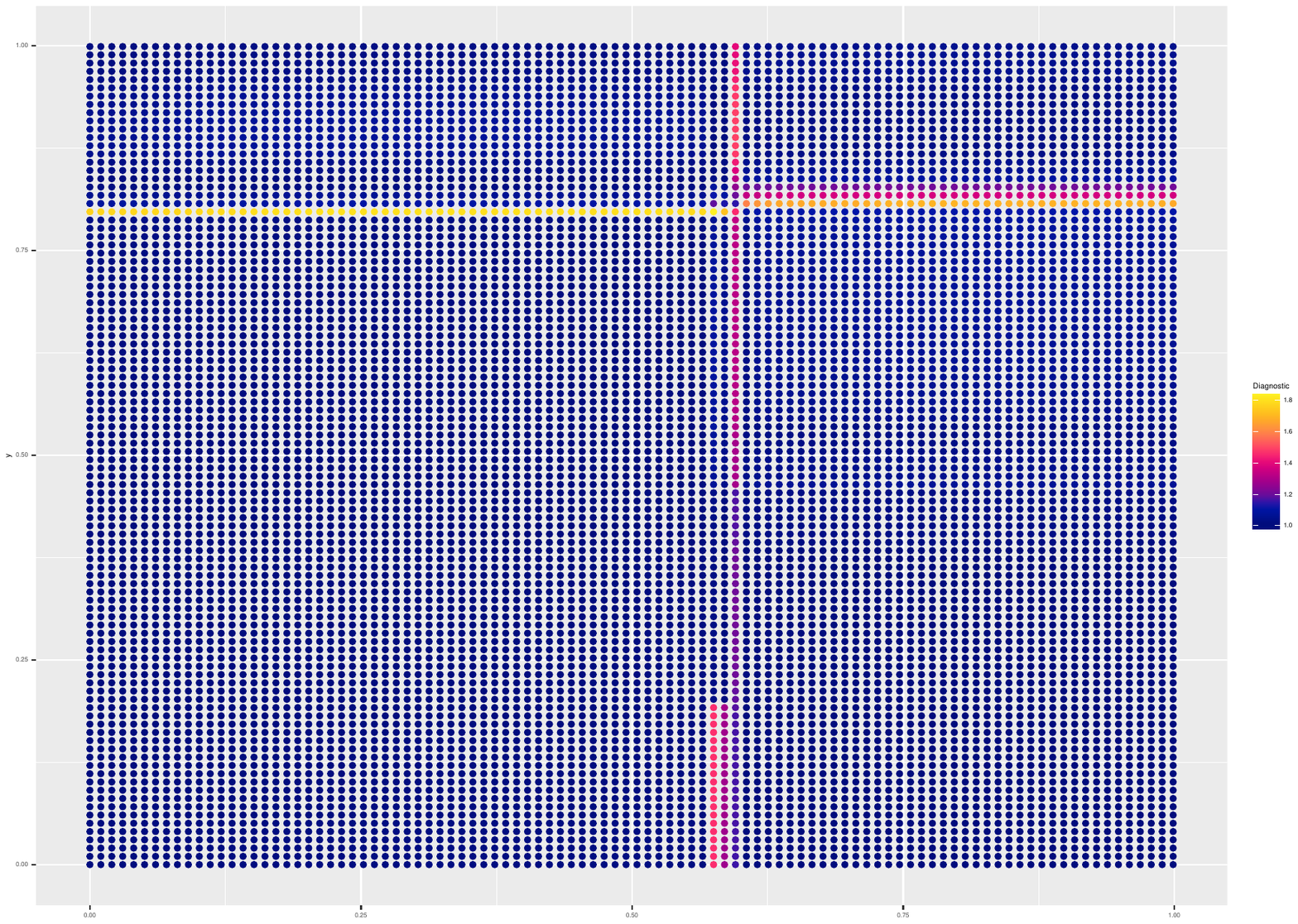}
    \caption{4 Trees and 100000 iterations }
  \end{subfigure}

\caption{The Gelman-Rubin Criterion for 4 trees and various number of iterations.}
\label{figE3}
\end{figure}

Likewise, we do not observe significant improvement in AAE and RMSE beyond 10000  iterations (see Table \ref{Table_E3}). Moreover, increasing the number of iterations does not fix the convergence issues at points close to jumps (see Figure \ref{figE3}).

\begin{table}[H]
\begin{tabular}{ |p{1.2cm}||p{2cm}|p{2.2cm}|p{2.2cm}|p{2.2cm}|p{2.3cm}| }
 \hline
 \multicolumn{6}{|c|}{Proposed Algorithm}  \\
 \hline
 Number of  trees & Number of Iterations &AAE for Posterior Mean & AAE for Posterior Median & RMSE for Posterior Mean & RMSE for Posterior Median\\
 \hline
 4 & 10000 & 241.82 & 240.1 & 464.99 & 489.93   \\
 & 50000 &209.95 & 209.58&392.43 &418.37   \\
 & 100000 &208.74&213.04&410.19 &447.86   \\
 \hline
\end{tabular}
\caption{Average Absolute Error and Root Mean Square Error for 4 Trees and various number of iterations.}
\label{Table_E3}
\end{table}

\subsection{Inhomogeneous two dimensional Poisson Process with Gaussian intensity} 
\begin{figure}[H] 
 \begin{subfigure}{8cm}
    \centering\includegraphics[width=6cm]{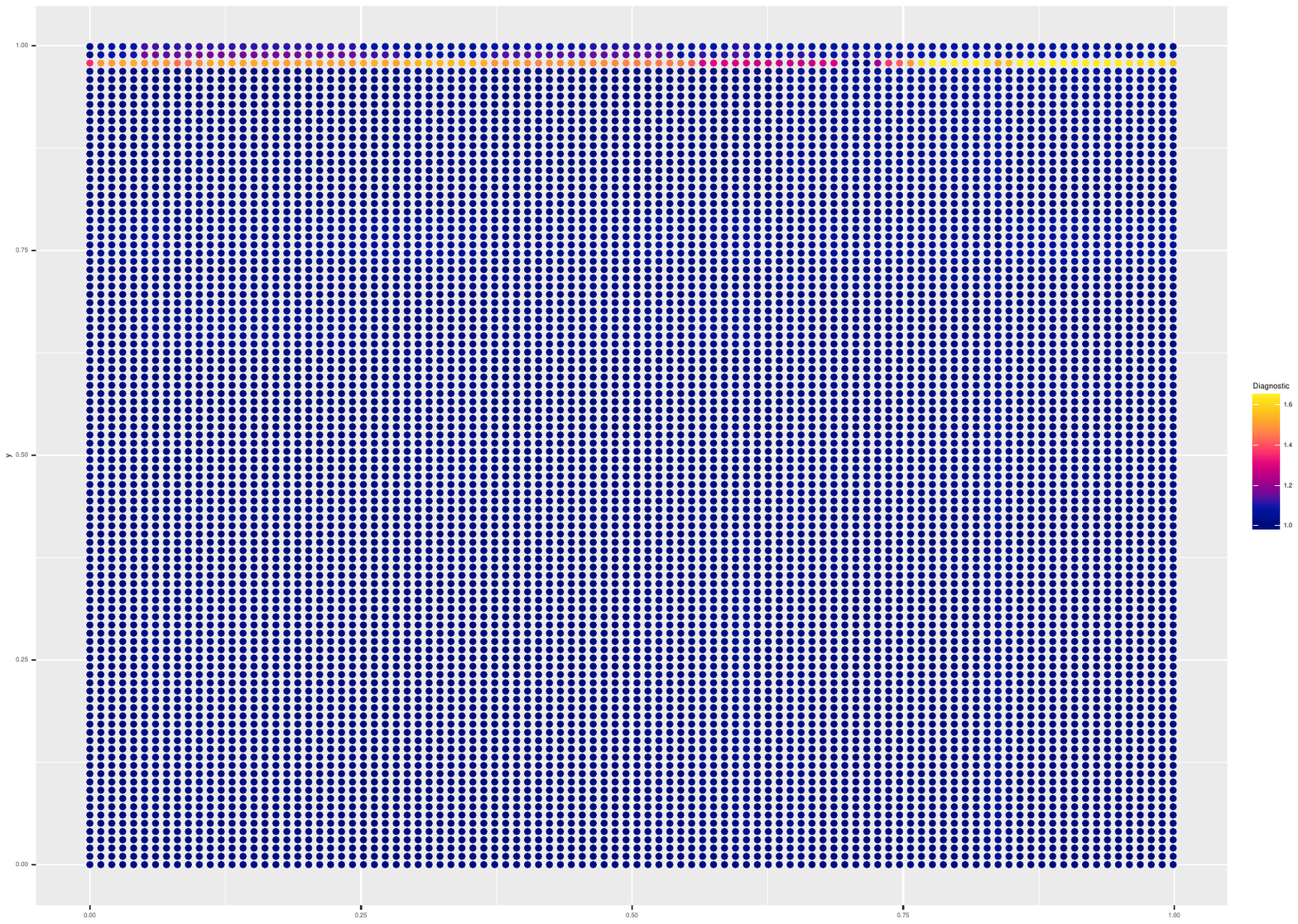}
    \caption{8 Trees and 10000 iterations }
  \end{subfigure}
\begin{subfigure}{8cm}
    \centering\includegraphics[width=6cm]{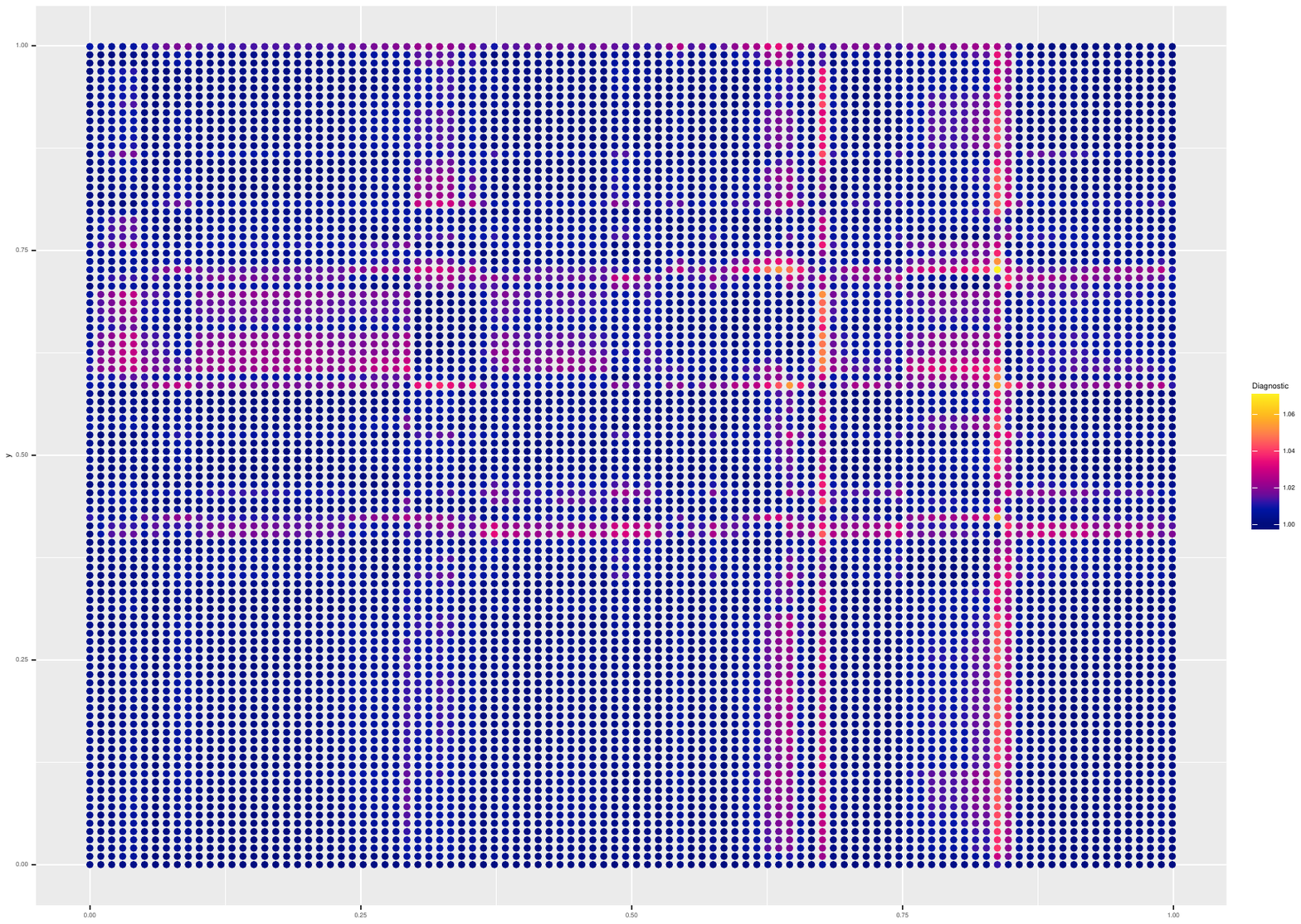}
    \caption{10 Trees and 10000 iterations }
  \end{subfigure}
  \begin{subfigure}{8cm}
    \centering\includegraphics[width=6cm]{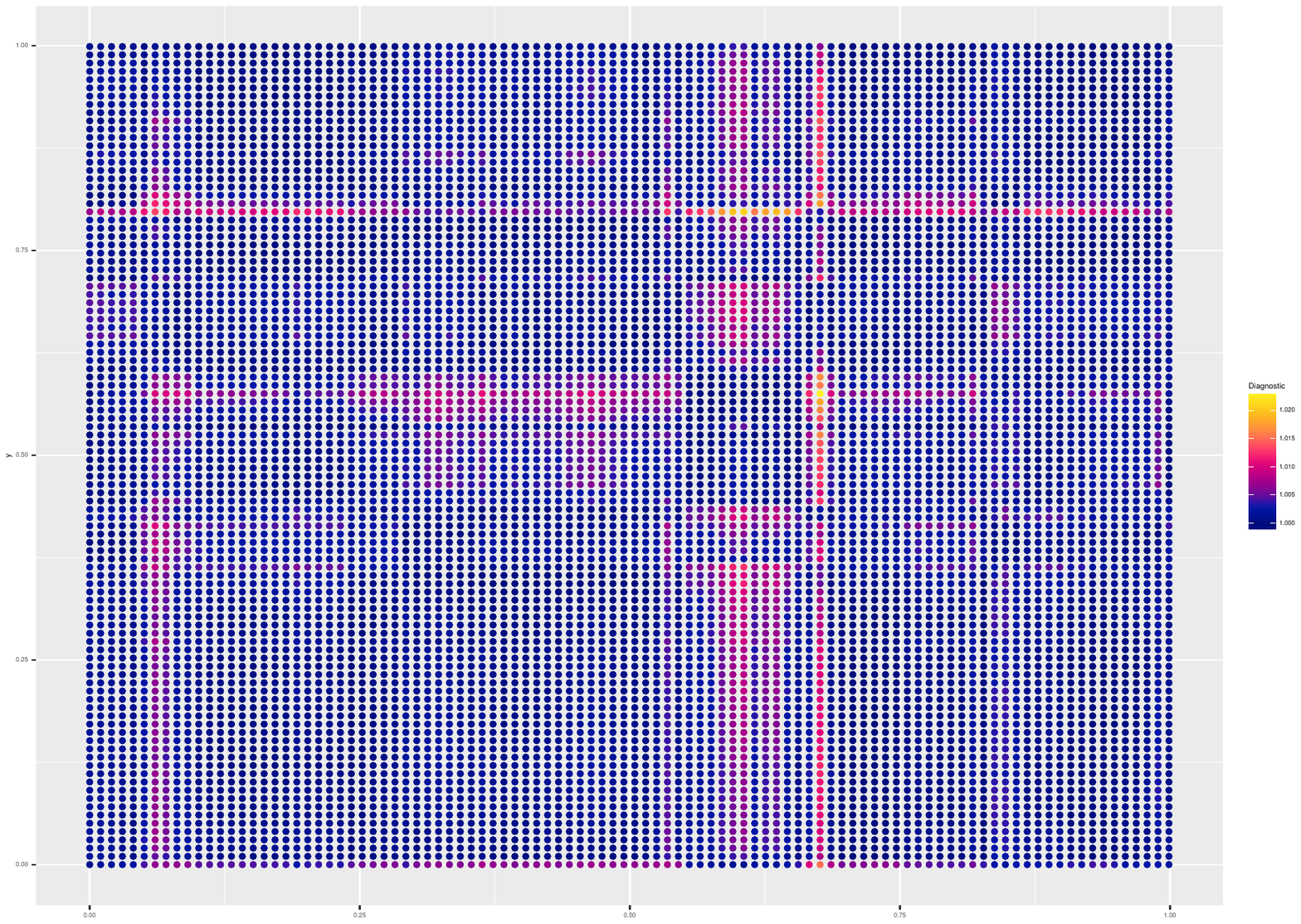}
    \caption{8 Trees and 50000 iterations }
  \end{subfigure}
\begin{subfigure}{8cm}
    \centering\includegraphics[width=6cm]{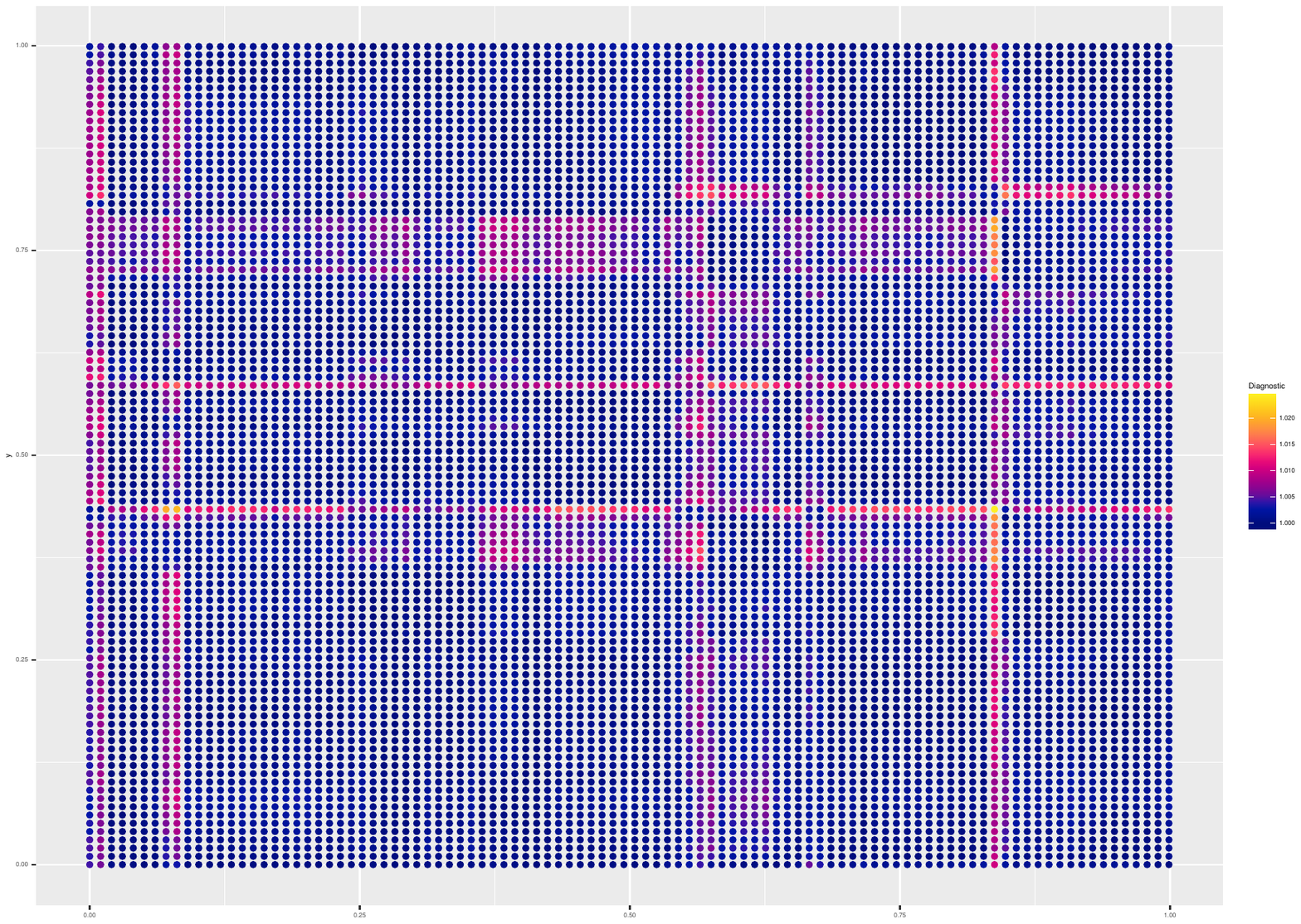}
    \caption{10 Trees and 50000 iterations }
  \end{subfigure}
  \begin{subfigure}{8cm}
    \centering\includegraphics[width=6cm]{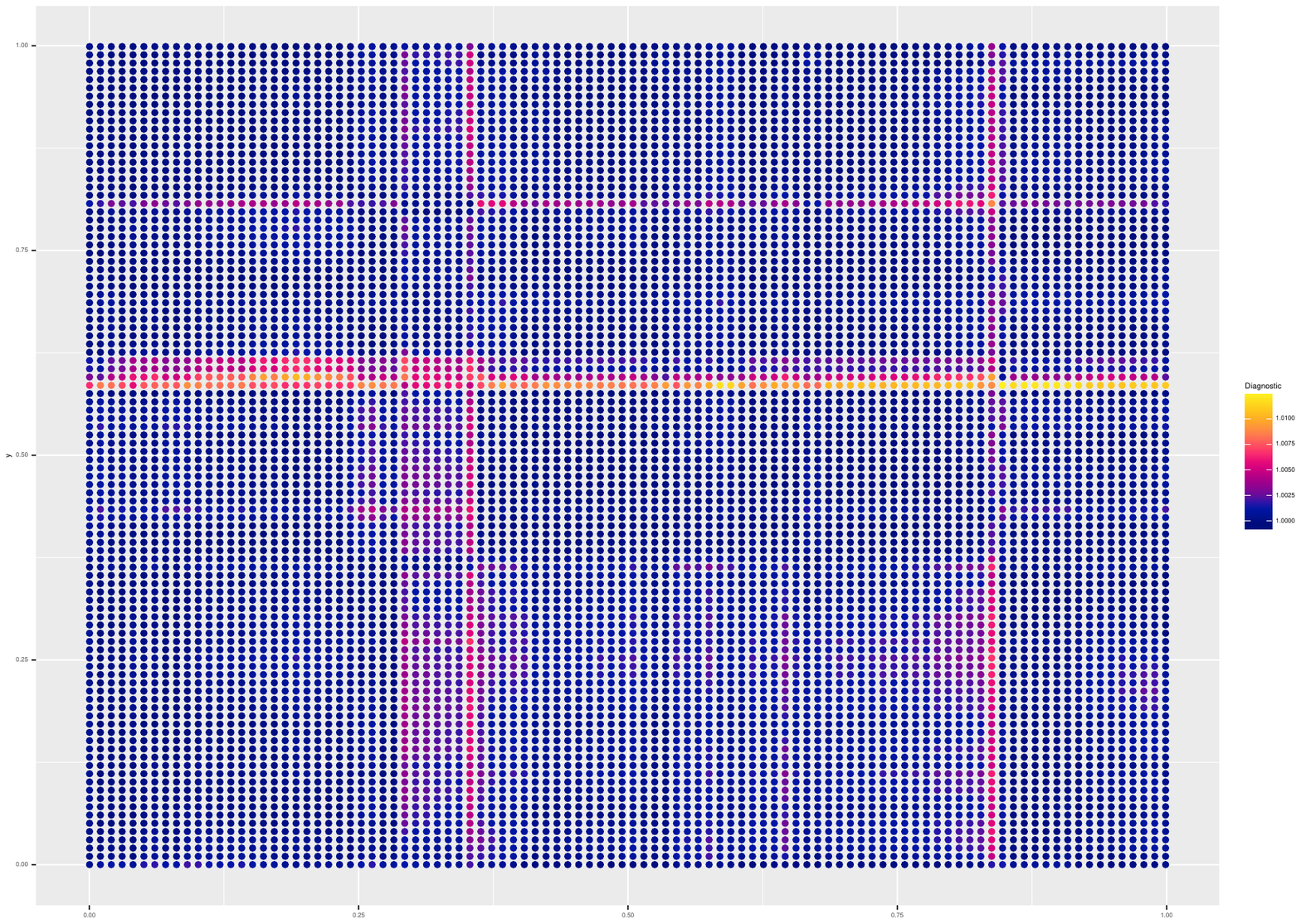}
    \caption{8 Trees and 200000 iterations }
  \end{subfigure}
\begin{subfigure}{8cm}
    \centering\includegraphics[width=6cm]{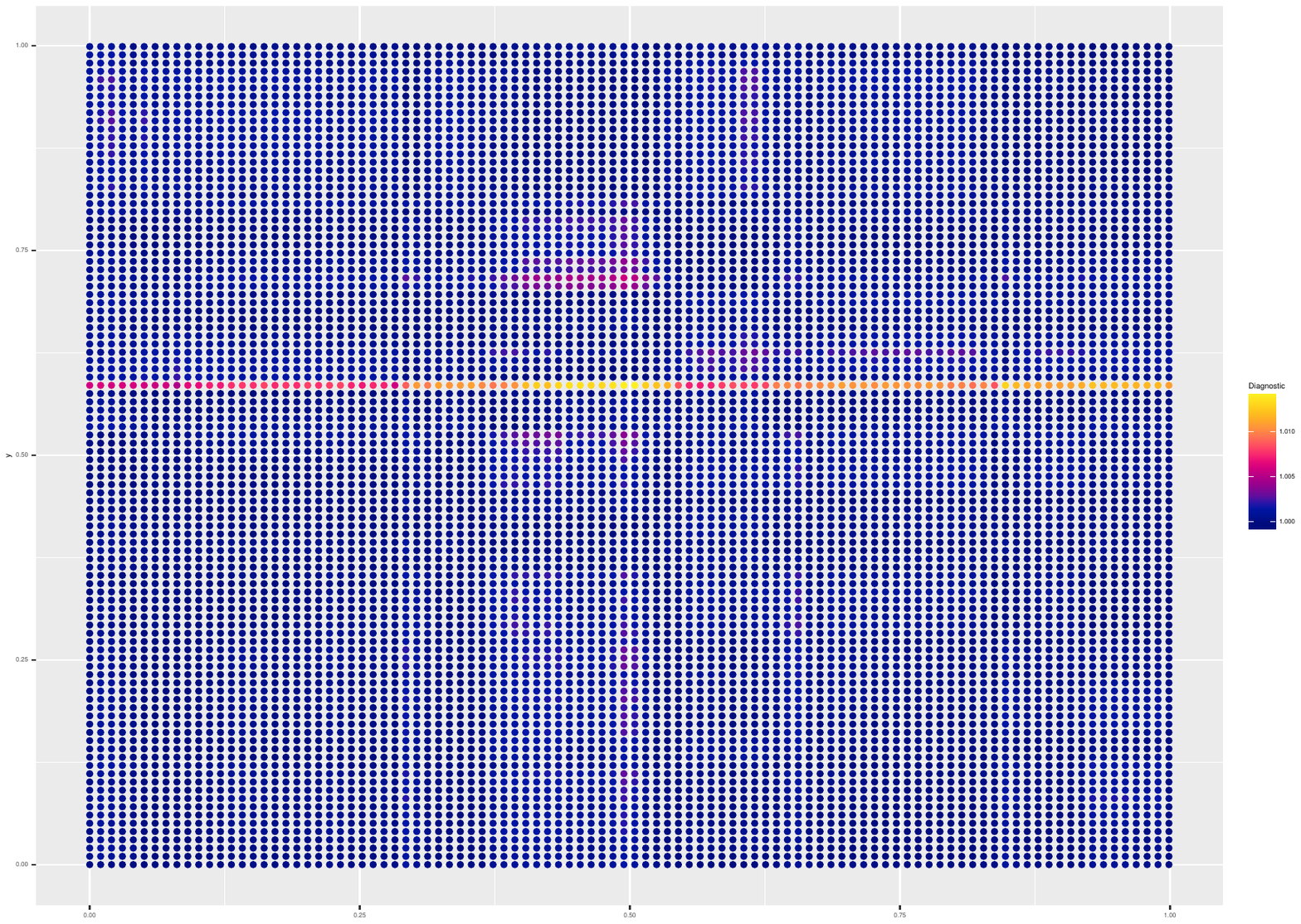}
    \caption{10 Trees and 200000 iterations }
  \end{subfigure}
\caption{The Gelman-Rubin Criterion for various number of iterations and trees.}
\label{figE4}
\end{figure}

Similarly to all the above scenarios, the error with 10000 iterations are already comparable to those obtained with a larger number of iterations (see Table \ref{Table_E4}).
Figure~\ref{figE4} shows that the chains converge for 10 Trees even if we consider a relatively  small number of iterations. The same holds for the majority of testing points for 8 Trees. The algorithm only provides less accurate estimations for the testing points close to the upper end of the domain for 8 Trees and relatively small number of iterations.

\begin{table}[H]
\begin{tabular}{ |p{1.2cm}||p{2cm}|p{2.2cm}|p{2.2cm}|p{2.2cm}|p{2.3cm}| }
 \hline
 \multicolumn{6}{|c|}{Proposed Algorithm}  \\
 \hline
 Number of  trees & Number of Iterations &AAE for Posterior Mean & AAE for Posterior Median & RMSE for Posterior Mean & RMSE for Posterior Median\\
 \hline
 8 & 10000 & 173.02 & 175.61 & 247.5 & 255.81   \\
 & 50000 &169.54 &170.5&242.03 &250.74   \\
 & 200000 &177.44&175.62&255.23 &258.88   \\
 \hline
10 &10000&168.91& 168.78 &242.62&249.38\\
&50000&177.72&173.93 &254.67&256.32\\
&200000&176.52&174.02 & 253.14&255.92\\
 \hline
\end{tabular}
\caption{Average Absolute Error and Root Mean Square Error for  various number of iterations and trees.}
\label{Table_E4}
\end{table}

\section{Intensity estimation for Real Data}
\subsection{Coal Data}   
\begin{figure}[H] 
 \begin{subfigure}{8 cm}
    \centering\includegraphics[width=8cm]{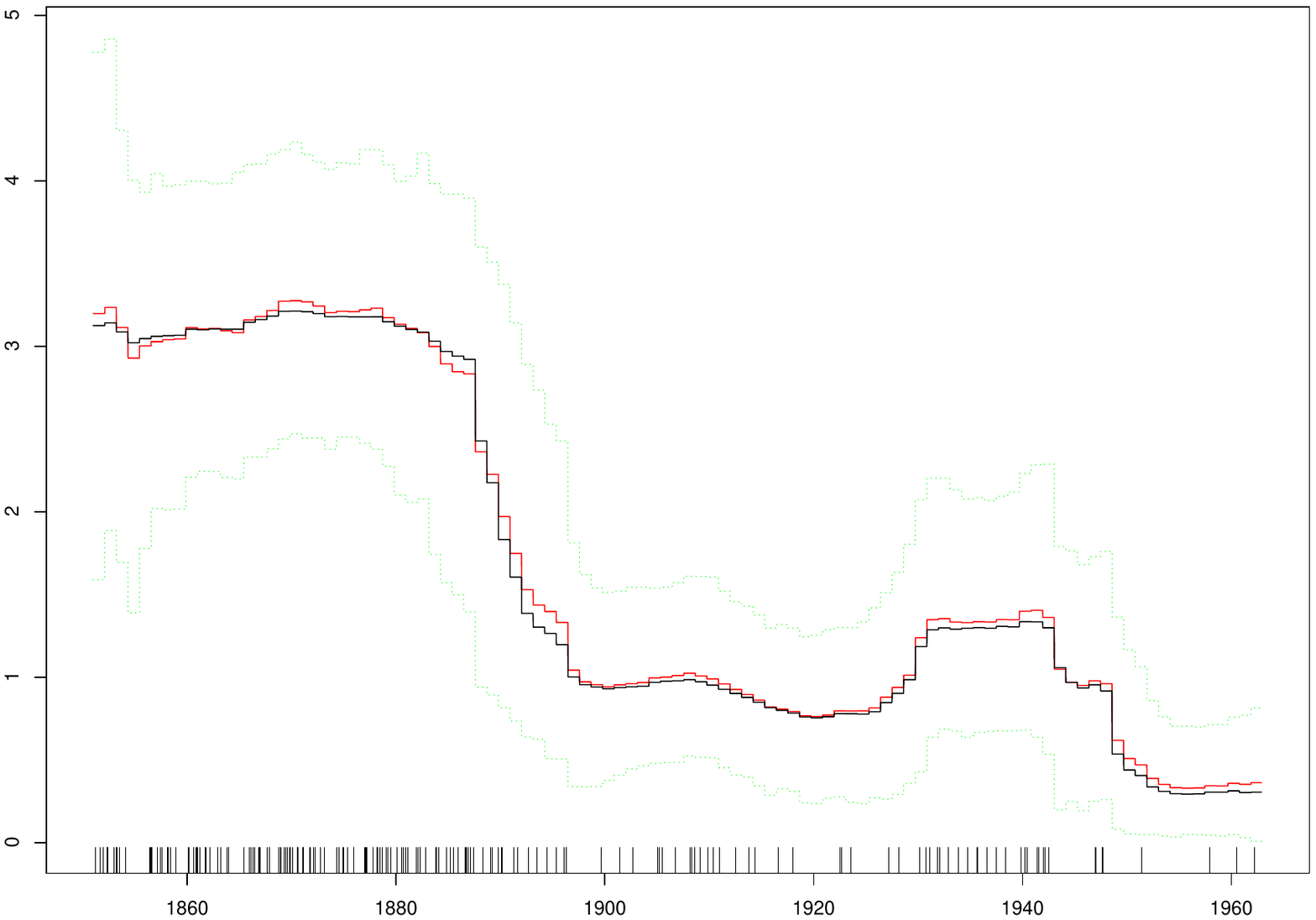}
    \caption{ 8 Trees}
  \end{subfigure}
  \begin{subfigure}{8cm}
    \centering\includegraphics[width=8cm]{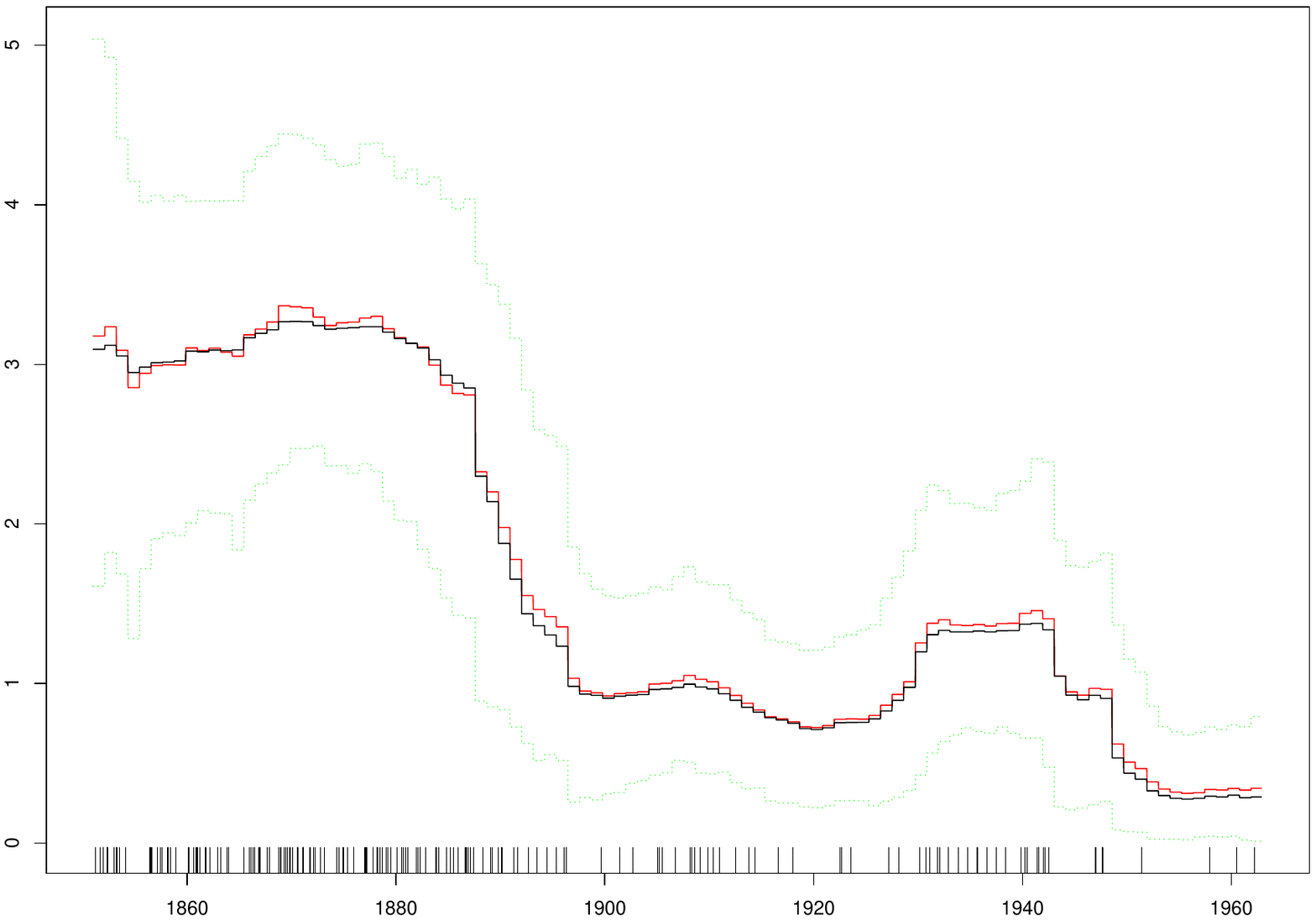}
    \caption{10 Trees}
  \end{subfigure}
\caption{Coal Data: The posterior mean (red curve), the posterior median (black curve), the 95\% hdi interval of the estimated intensity illustrated by the dotted green lines. The rug plot on the bottom displays the event times.}
\label{coal}
\end{figure}

 The first real data set under consideration is composed of the dates of 191 explosions which caused at least 10 occurrences of death from  March 22, 1962 until March 15, 1981. The data set is available in the \textbf{R} package
\textbf{boot} \citep{boot} as \textbf{coal}. 
Figure \ref{coal} illustrates the Posterior Mean and the Posterior Median for 8 and 10 Trees. We observe that our algorithm captures the fluctuations of the rate of accidents in the period under consideration. The diagnostic criteria included in the Supplementary Material indicate that the considered chains have converged. See  \citet{adams2009tractable}, \citet{gugushvili2018fast} and \citet{lloyd2015variational} for alternative analyses. 

\subsection{Redwoodfull Data} 
\begin{figure}[H] 
  \begin{subfigure}{8cm}
    \centering\includegraphics[width=6cm]{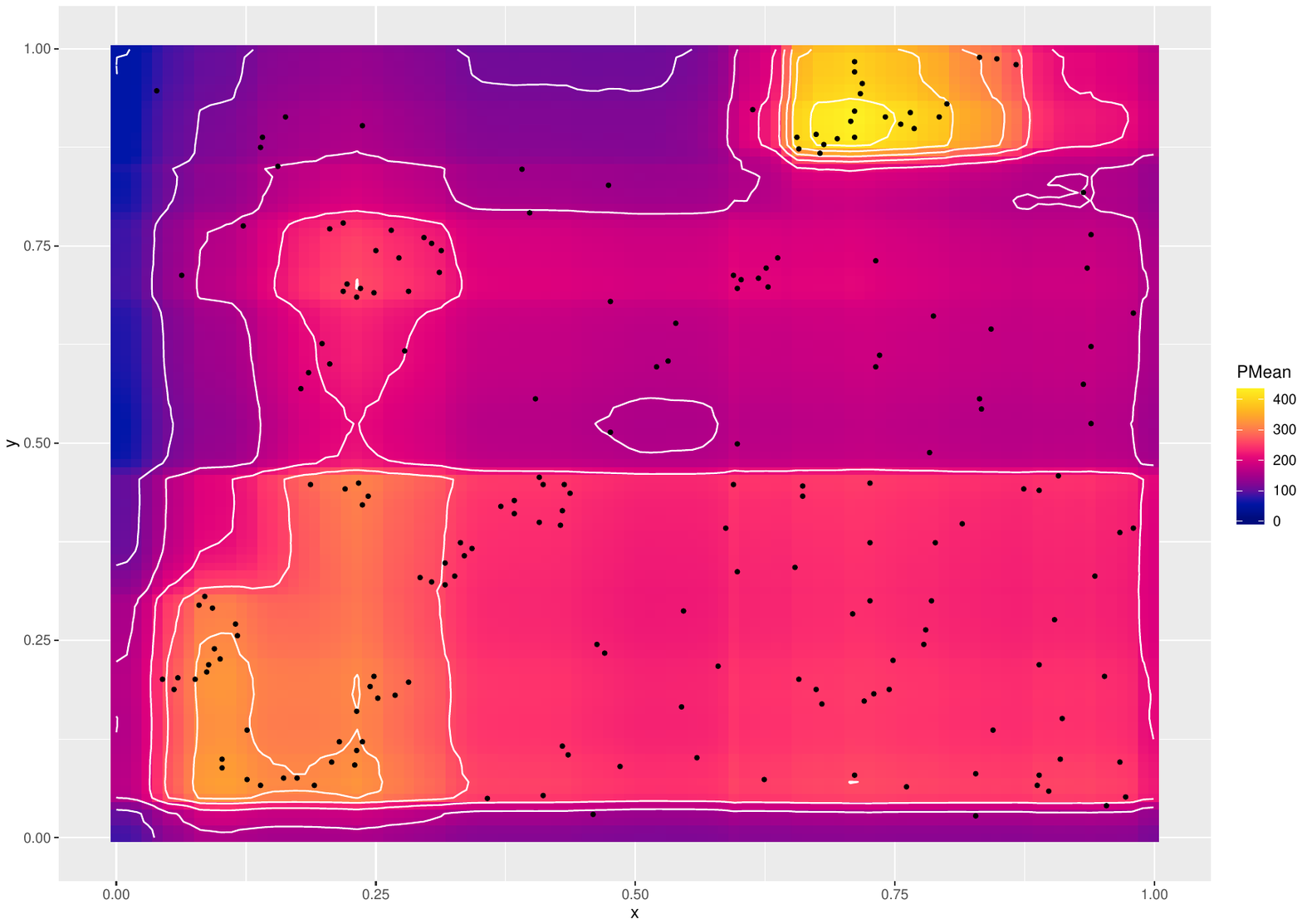}
    \caption{Posterior Mean for 5 Trees}
  \end{subfigure}
  \begin{subfigure}{8cm}
    \centering\includegraphics[width=6cm]{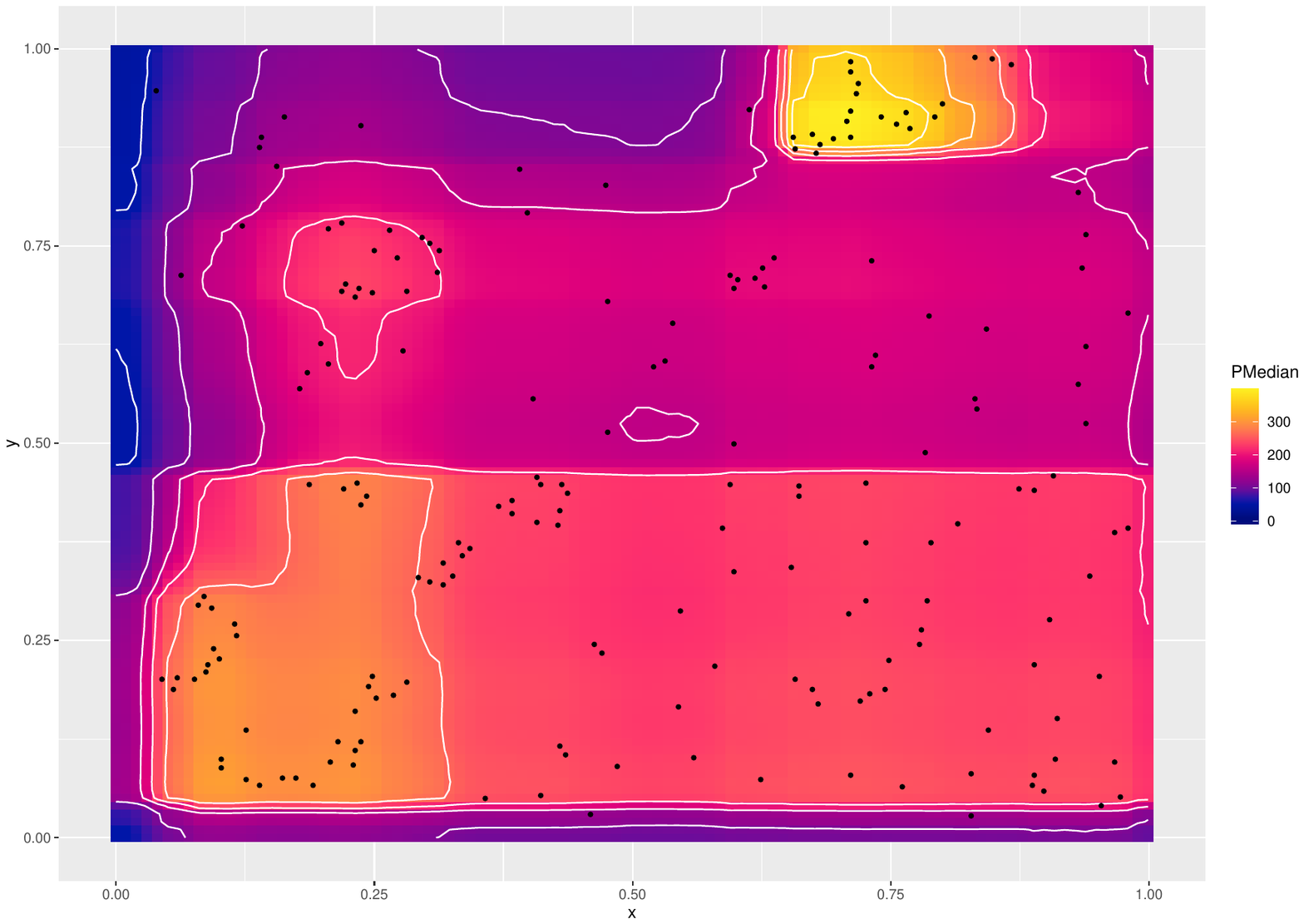}
    \caption{Posterior Median for 5 Trees}
  \end{subfigure}

 \begin{subfigure}{8cm}
    \centering\includegraphics[width=6cm]{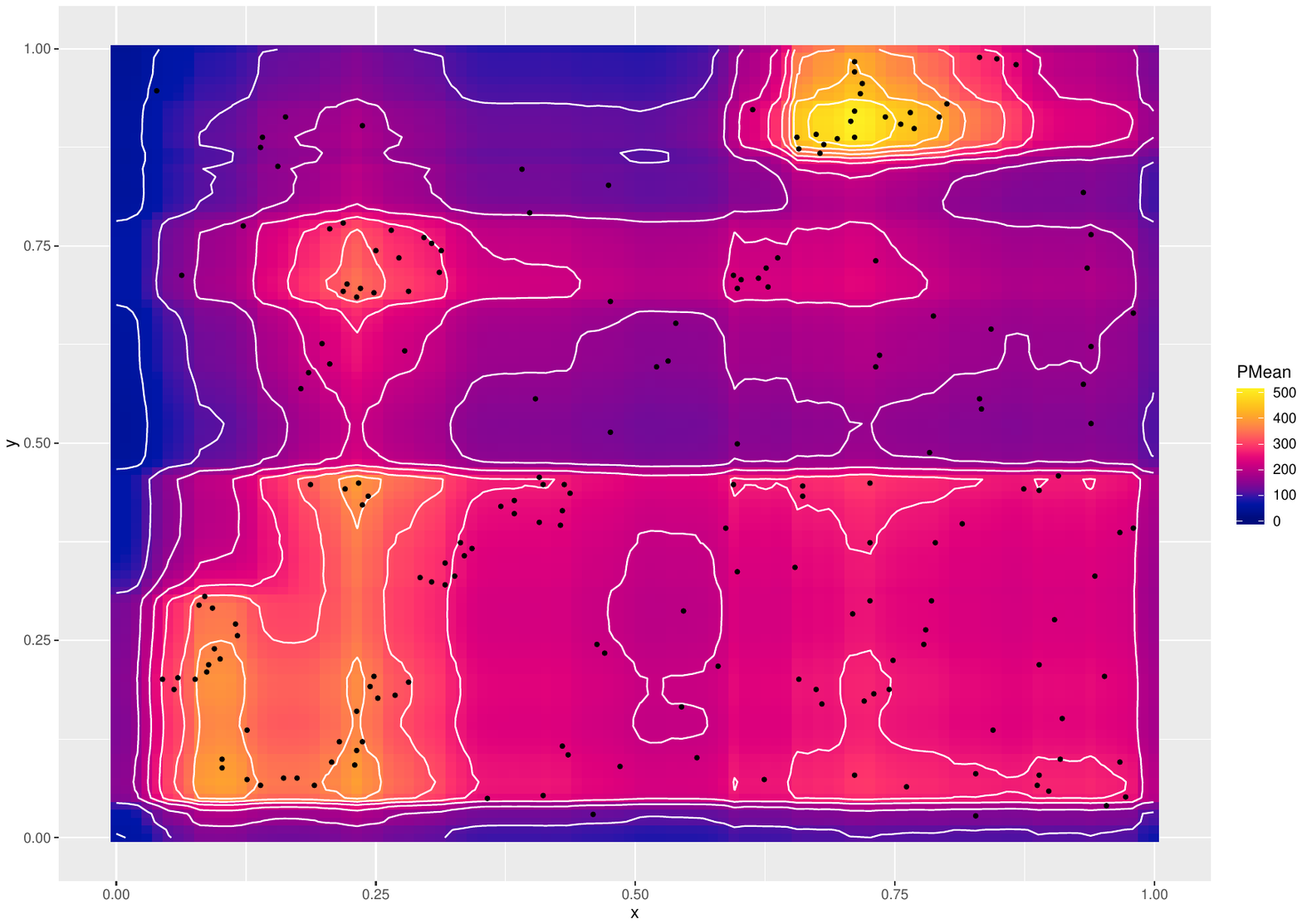}
    \caption{Posterior Mean for 10 Trees} 
\label{Red_mean}
  \end{subfigure}
  \begin{subfigure}{8cm}
    \centering\includegraphics[width=6cm]{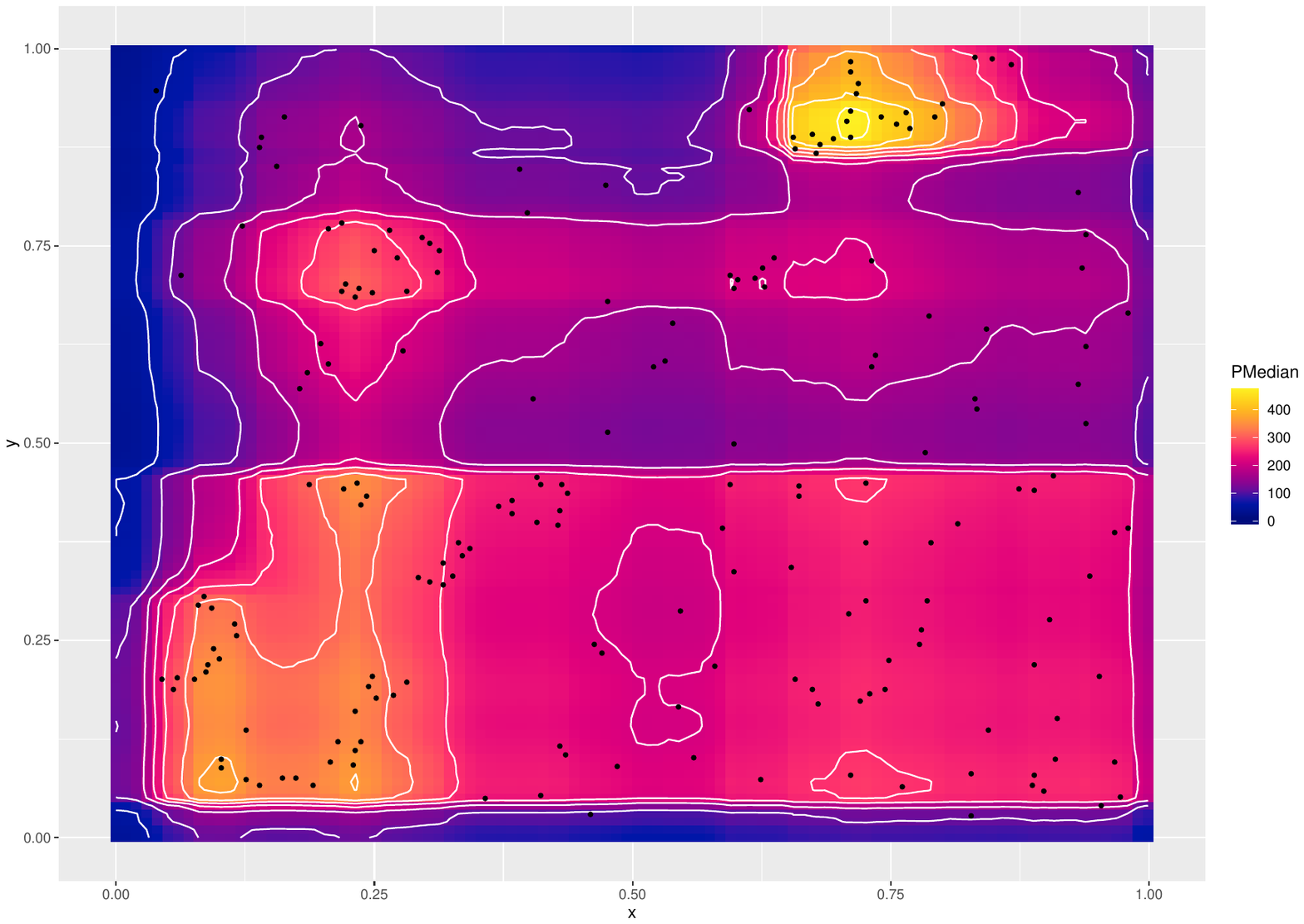}
    \caption{Posterior Median for 10 Trees}
\label{Red_median}
  \end{subfigure}
\caption{Posterior Mean and Posterior Median for 3, 5 and 10 Trees}

\end{figure}

\begin{figure}[H]
  \begin{subfigure}{8cm}
 \centering\includegraphics[width=6cm]{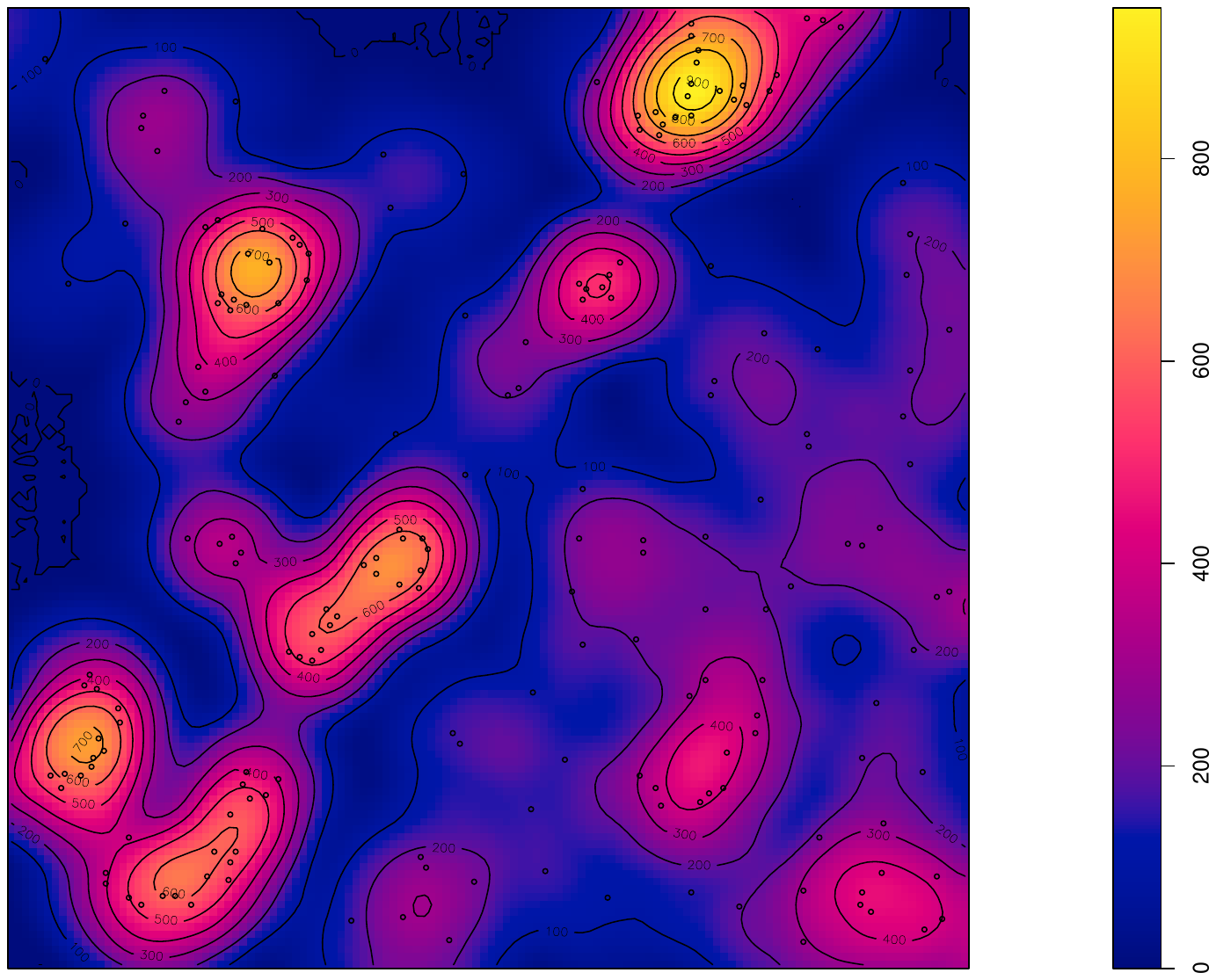}
    \caption{Fixed-bandwidth estimator}
\label{Red_kernel}
  \end{subfigure}
  \begin{subfigure}{8cm}
    \centering\includegraphics[width=6cm]{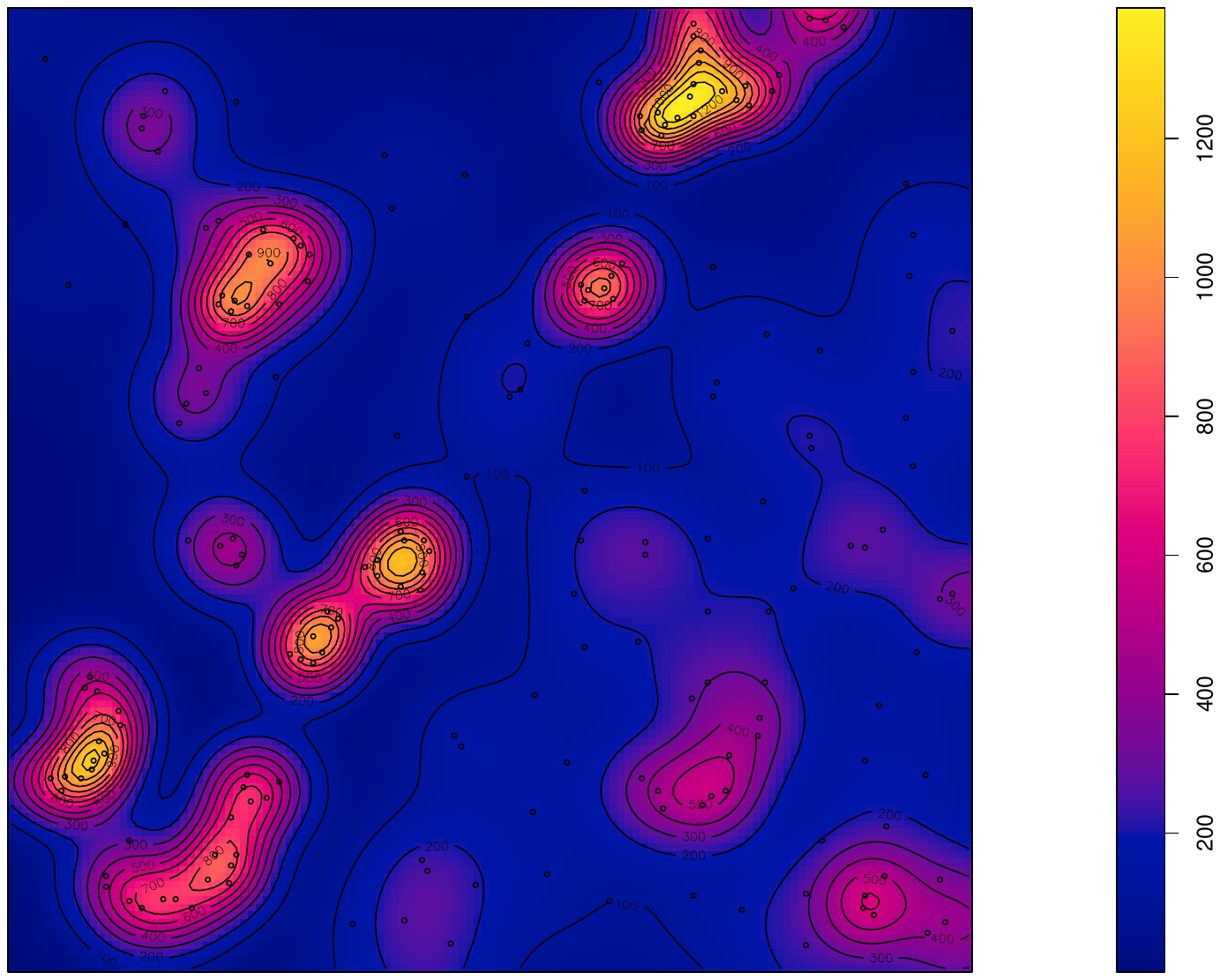}
    \caption{Adaptive-bandwidth estimator}
  \end{subfigure}
   \caption{Fixed-bandwidth chosen using likelihood cross-validation and adaptive-bandwidth kernel estimators.}
\label{fig8}
\end{figure}

 Finally, we use a data set available in the \textbf{R} package \textbf{spatstat} describing the locations of 195 trees in a square sampling region shown with dots in the figures below. \citet{adams2009tractable} analyzed the \textbf{redwoodfull} data using their recommended algorithm. We present the posterior mean and the posterior median obtained with our algorithm for different number of trees and the result of kernel estimators. Intensity inference via posterior mean (Figure \ref{Red_mean}) or posterior median (Figure \ref{Red_median}) for 10 Trees is similar to the fixed-bandwidth kernel estimator with edge correction and bandwidth selected using likelihood cross-validation (Figure \ref{Red_kernel}), and the inference from \citet{adams2009tractable}.

\section{Simulation Study on Synthetic Data}

\subsection{One dimensional Poisson Process with continuously varying intensity}
\label{sec:cosine_intenstity}
\begin{figure}[H] 
  \begin{subfigure}{8cm}
    \centering\includegraphics[width=8cm]{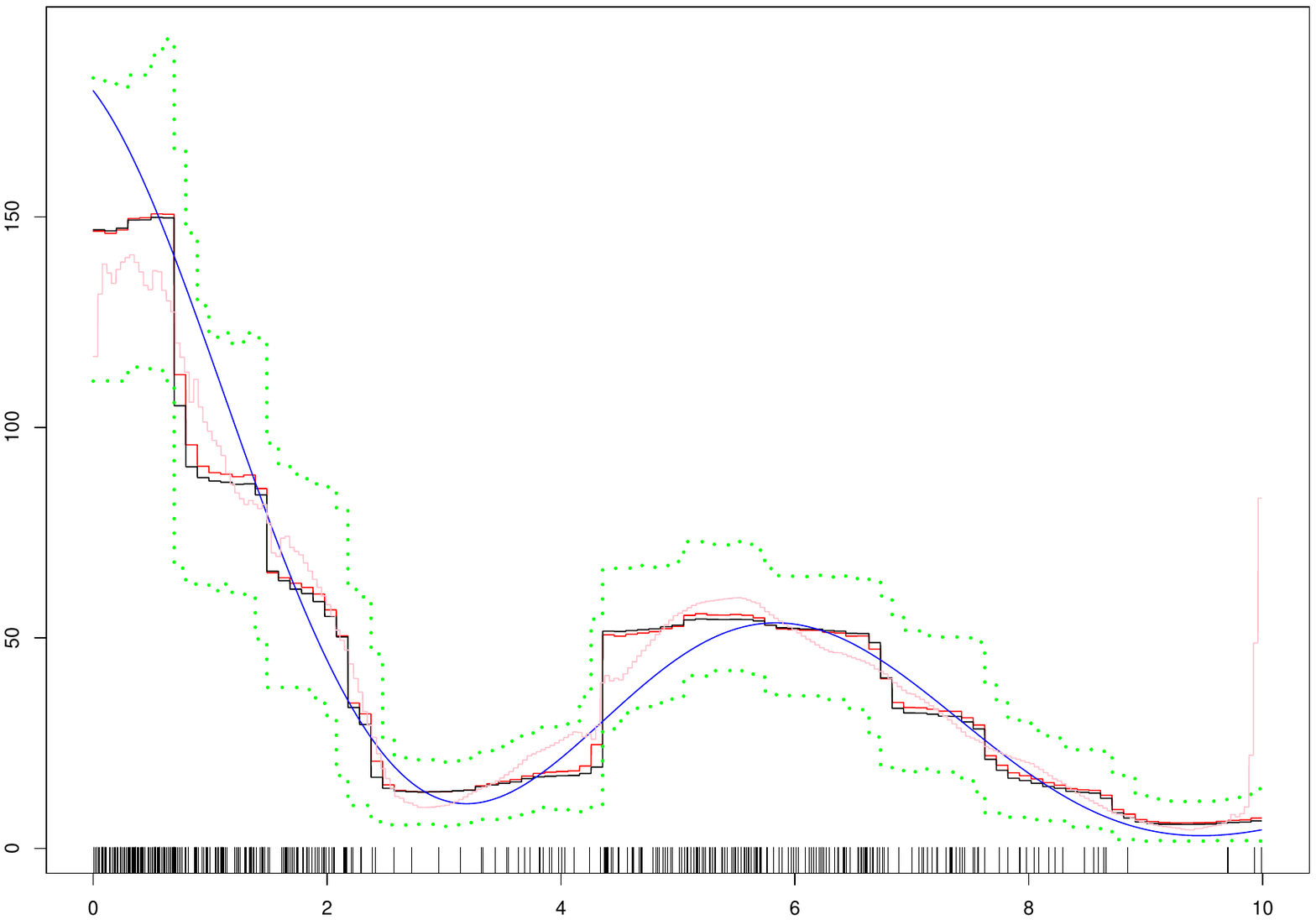}
    \caption{5 Trees}
  \end{subfigure}
 \begin{subfigure}{8cm}
    \centering\includegraphics[width=8cm]{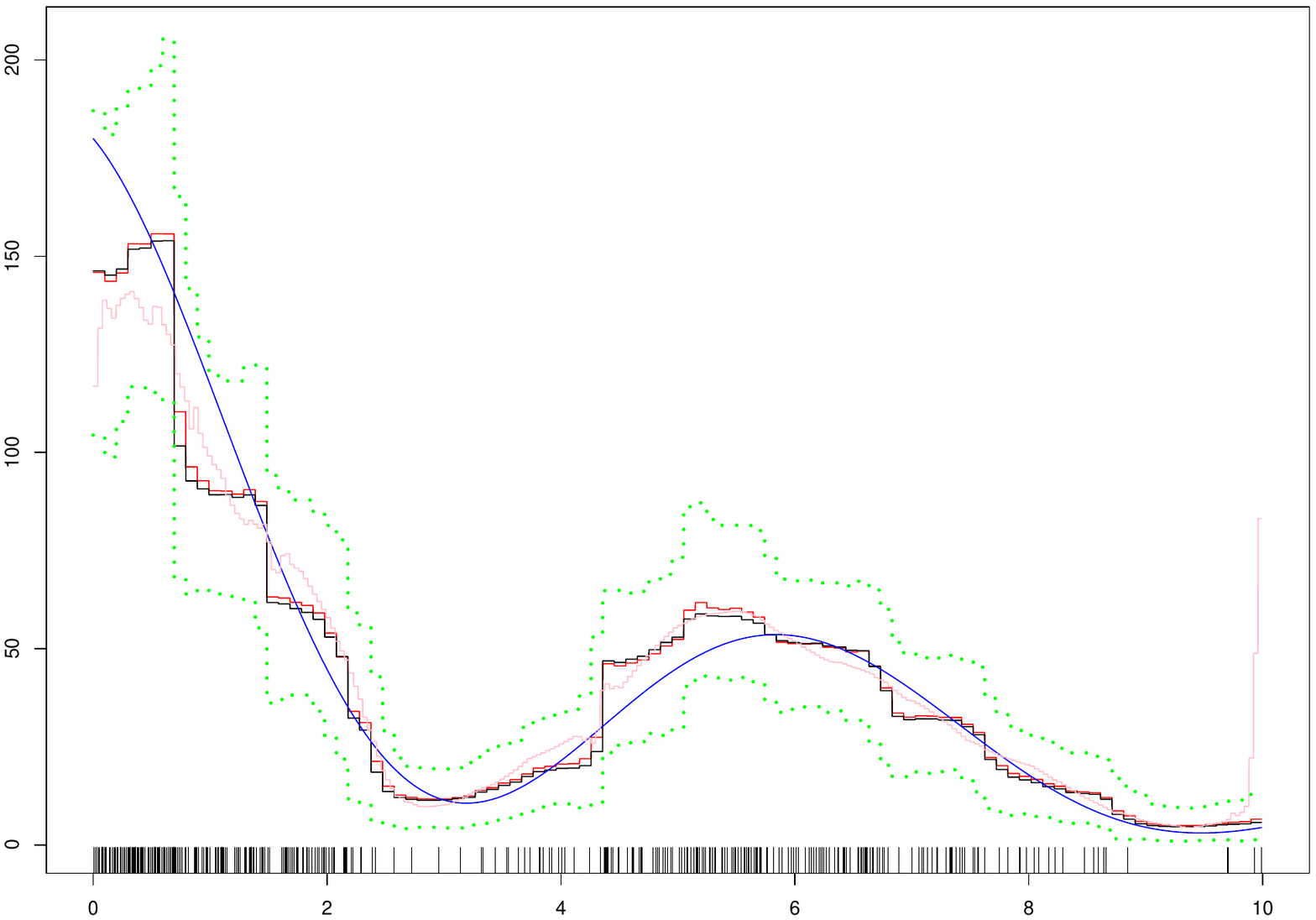}
    \caption{10 Trees}
  \end{subfigure}
\caption{Scenario 2:  The original intensity (blue curve), the posterior mean (red curve), the posterior median (black curve), the 95\% hdi interval of the estimated intensity illustrated by the dotted green lines and the Haar-Fisz estimator (pink curve). The rug plot on the bottom displays the 440 event times. }
\label{fig2}
\end{figure}

We have applied our algorithm to samples of a one dimensional Poisson process with intensity $\lambda(x)=20 e^{-x/5} (5 + 4 \cos(x))$ for $x \in [0,10]$.
Figure \ref{fig2} and Tables \ref{table2}-\ref{tableH2} show that  the algorithm works well on a smoothy varying intensity with fewer sample points and outperforms the Haar-Fisz Estimator for the majority the range. The convergence criteria indicate convergence of the simulated chains for 10 Trees and for the most testing points for 5 Trees (see supplementary material).

\begin{table}[H]
\begin{tabular}{ |p{2cm}||p{2.4cm}|p{2.4cm}|p{2.5cm}|p{2.5cm}| }
 \hline
 \multicolumn{5}{|c|}{Proposed Algorithm}  \\
 \hline
 Number of trees & AAE for Posterior Mean & AAE for Posterior Median & RISE for Posterior Mean & RISE for Posterior Median\\
 \hline
 5 &  6.14 &6.38&9.52 &10.17  \\
10 & 5.95& 6.01 & 9.39&9.8  \\
 \hline
\end{tabular}
\caption{Average Absolute Error and Root Integrated Square Error for  various number of trees for the data in Fig.~\ref{fig2}.}
\label{table2}

\end{table}
\begin{table}[H]
\begin{tabular}{ |p{3cm}||p{3cm}|  }
 \hline
 \multicolumn{2}{|c|}{Haar-Fisz Algorithm}  \\
 \hline
 AAE  &  RISE   \\
 \hline
7.16 & 11.67   \\
 \hline
\end{tabular}
\caption{Average Absolute Error and Root Integrated Square Error for Haar-Fisz estimator for the data in Fig.~\ref{fig2}}
\label{tableH2}
\end{table}
\subsection{Inhomogeneous two-dimensional Poisson Process with Gaussian intensity} 
\begin{figure}[H] 
  \begin{subfigure}{6cm}
    \centering\includegraphics[width=6cm]{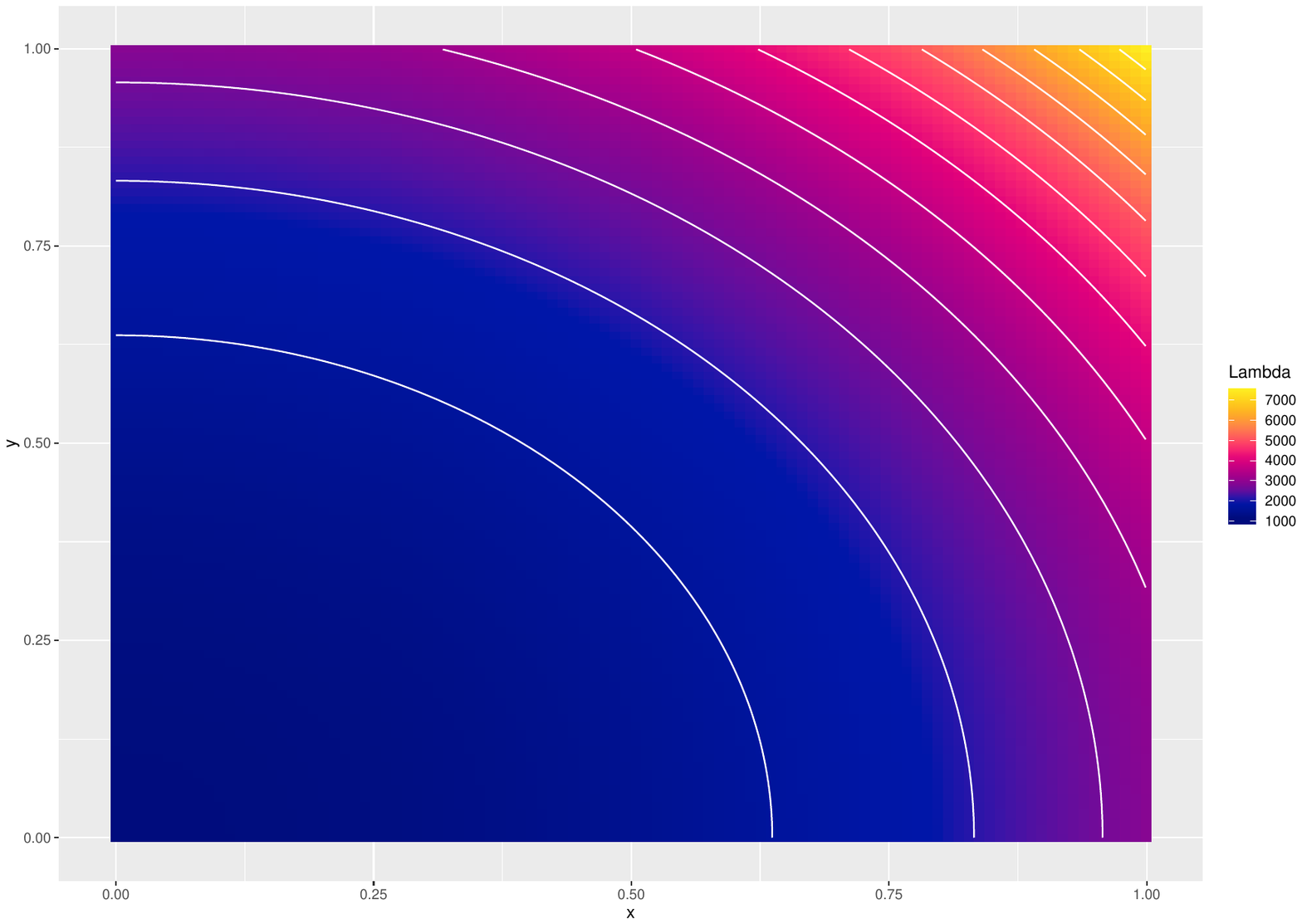}
    \caption{Original Intensity}
  \end{subfigure}
 \begin{subfigure}{6cm}
    \centering\includegraphics[width=6cm]{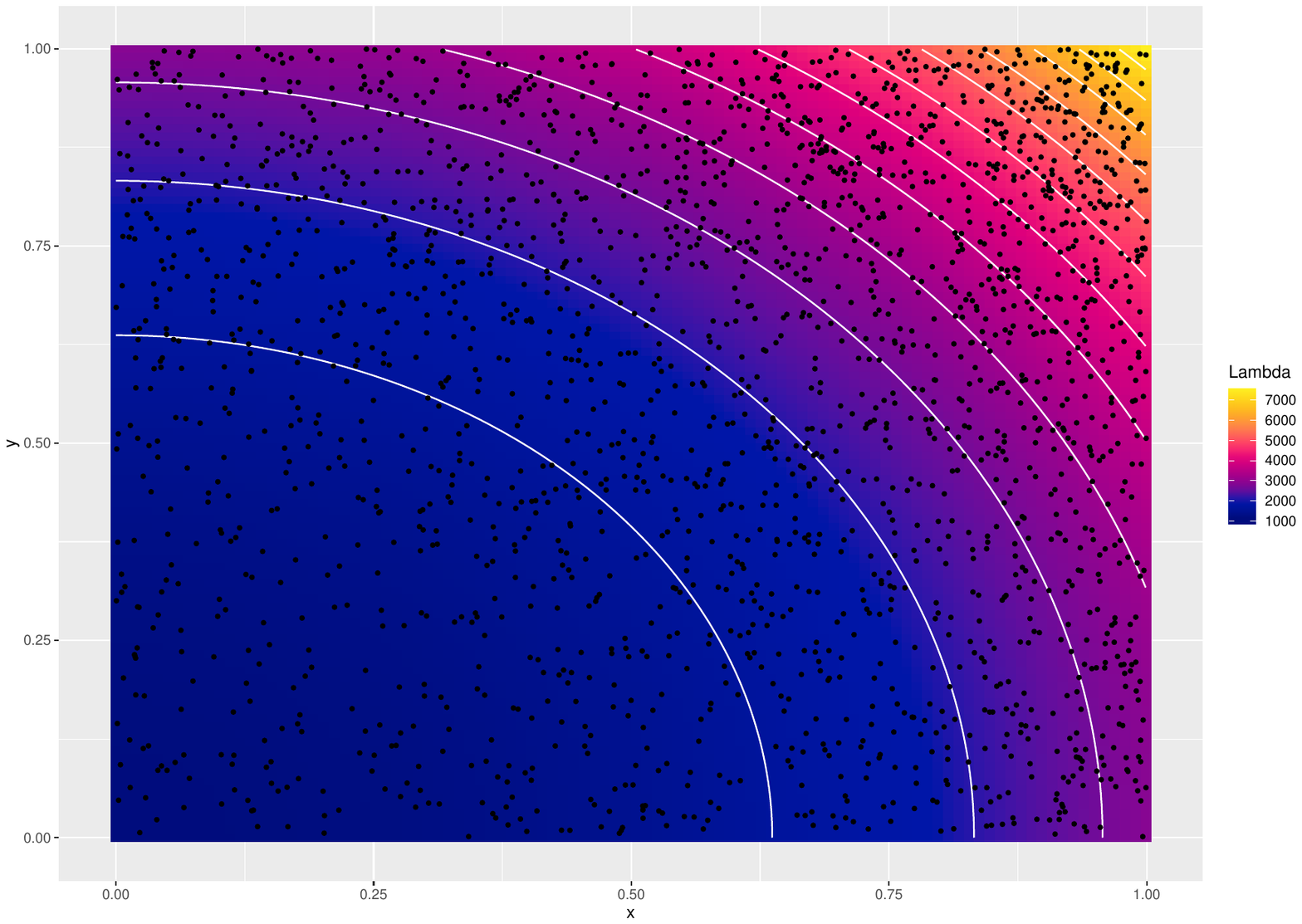}
    \caption{Realization of Process consisting of 2176 points}
  \end{subfigure}

  \begin{subfigure}{6cm}
    \centering\includegraphics[width=6cm]{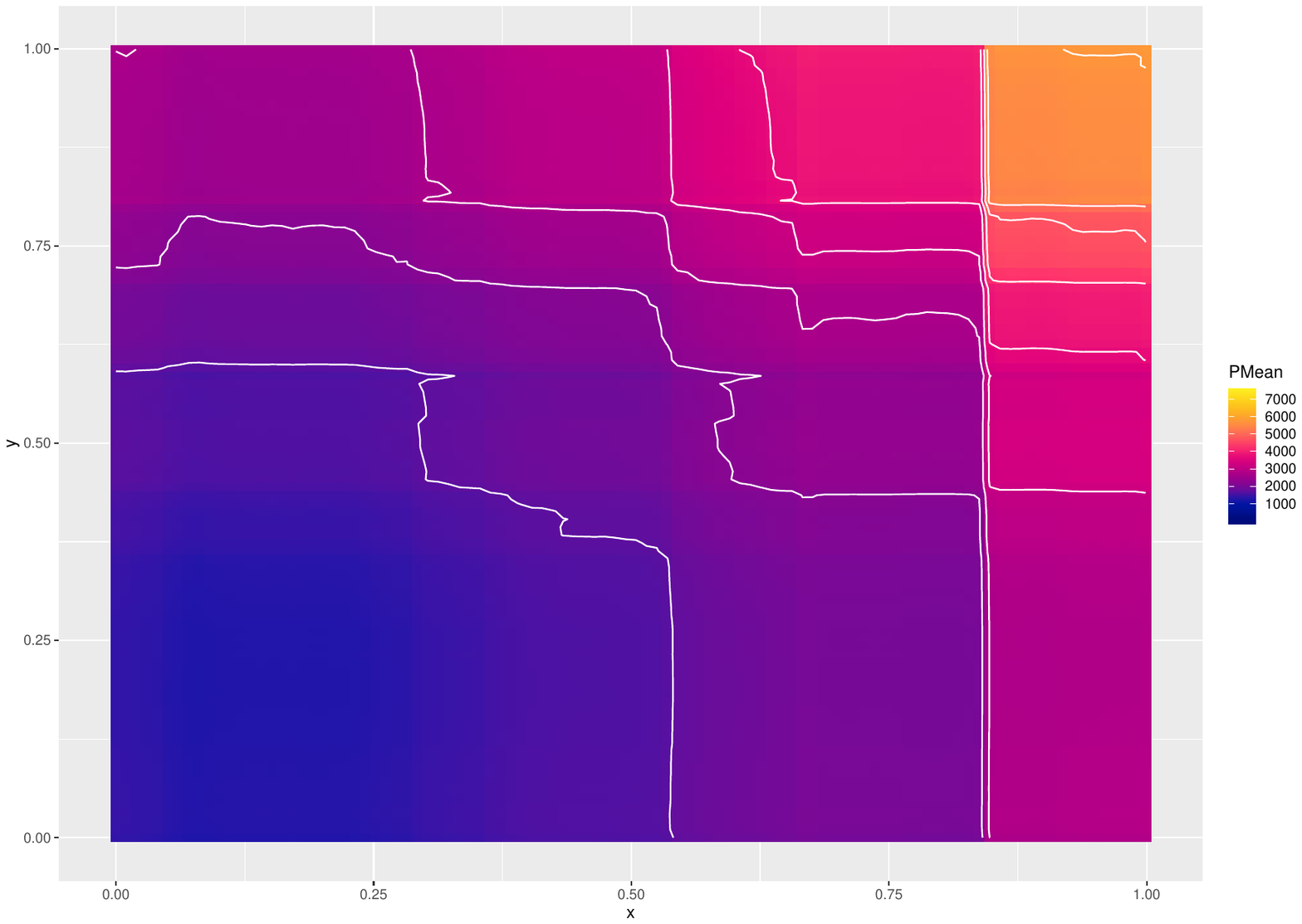}
    \caption{Posterior Mean for 8 Trees}
  \end{subfigure}
  \begin{subfigure}{6cm}
    \centering\includegraphics[width=6cm]{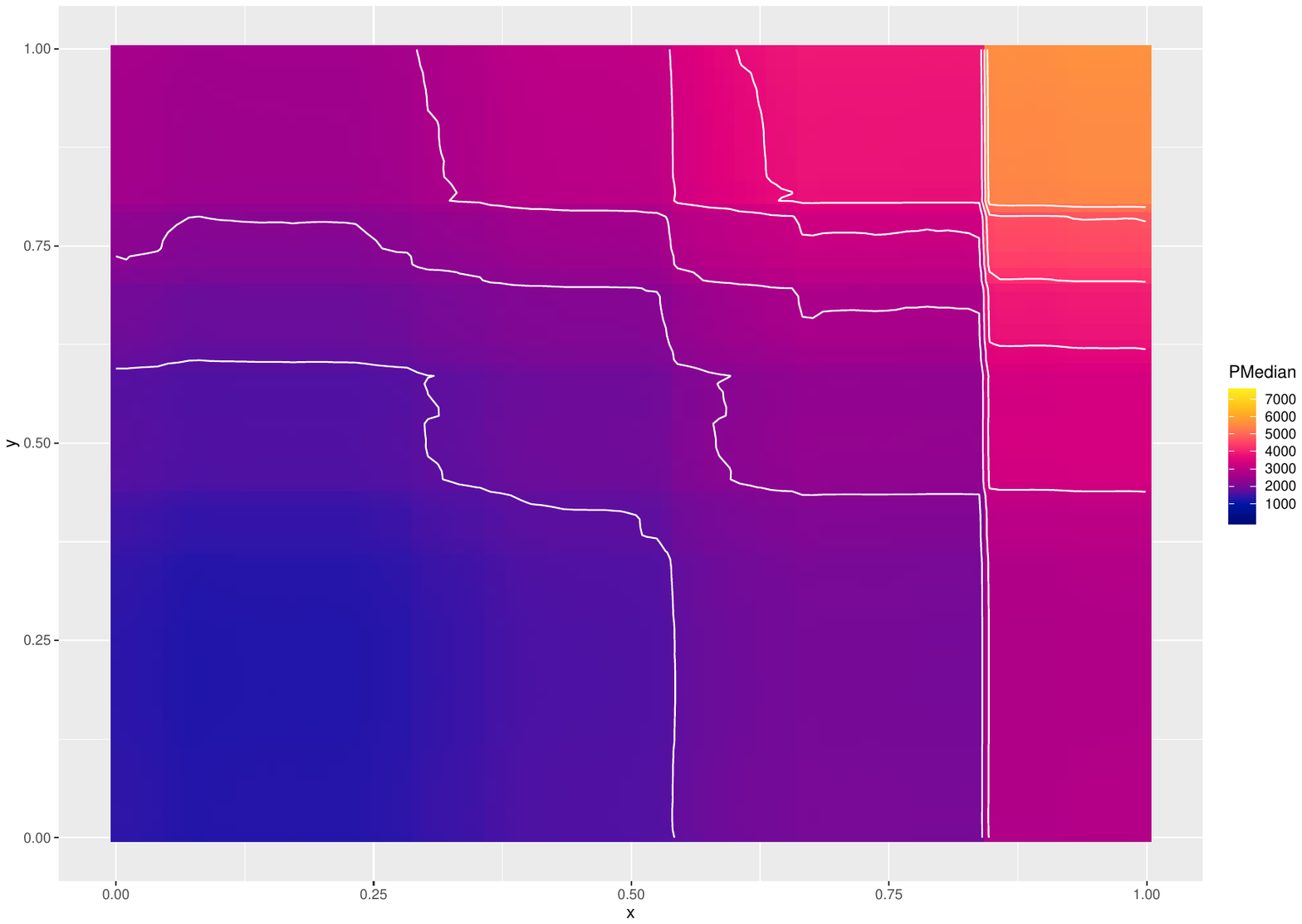}
    \caption{Posterior Median for 8 Trees}
  \end{subfigure}

  \begin{subfigure}{6cm}
    \centering\includegraphics[width=6cm]{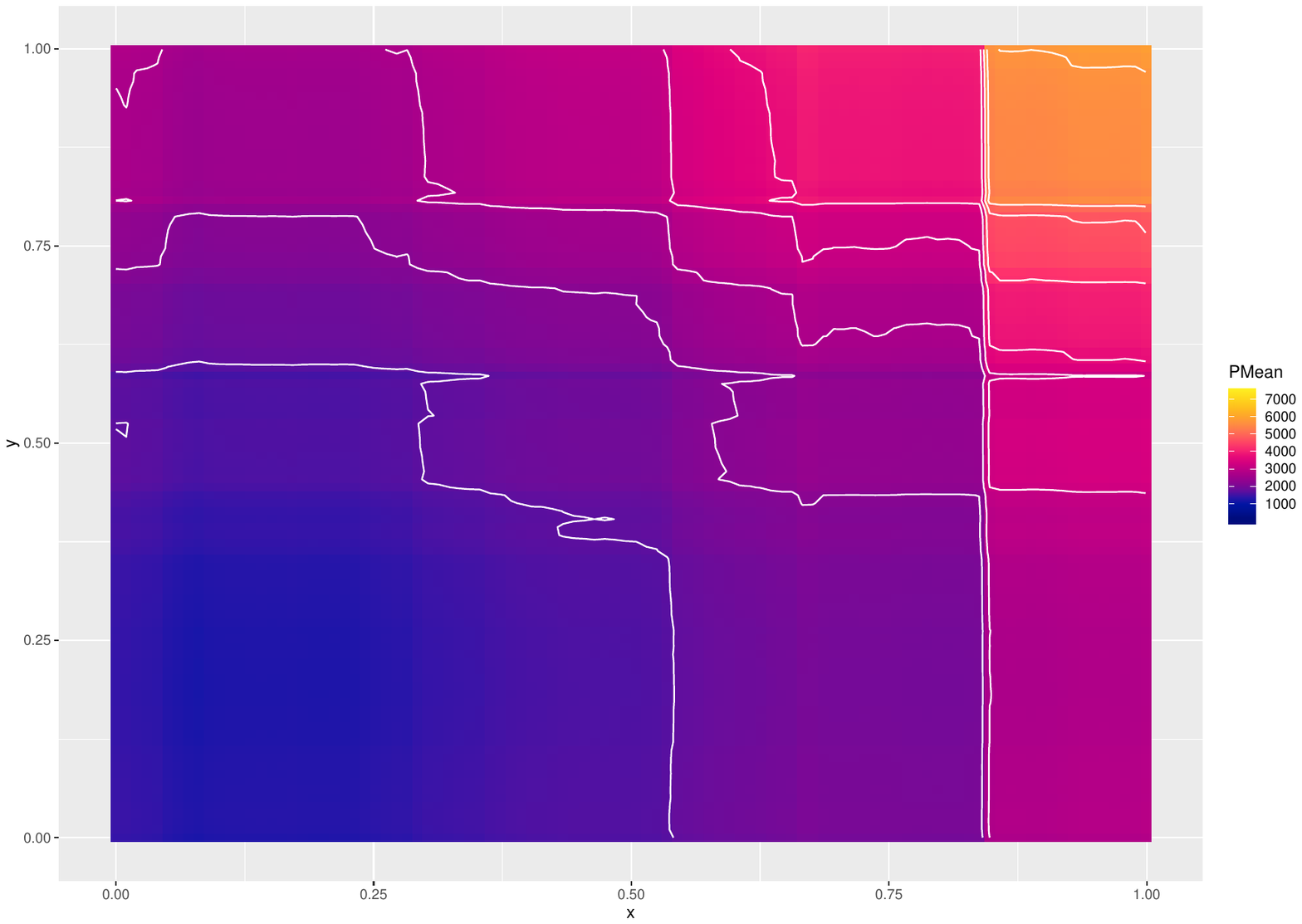}
    \caption{Posterior Mean for 10 Trees}
  \end{subfigure}
  \begin{subfigure}{6cm}
    \centering\includegraphics[width=6cm]{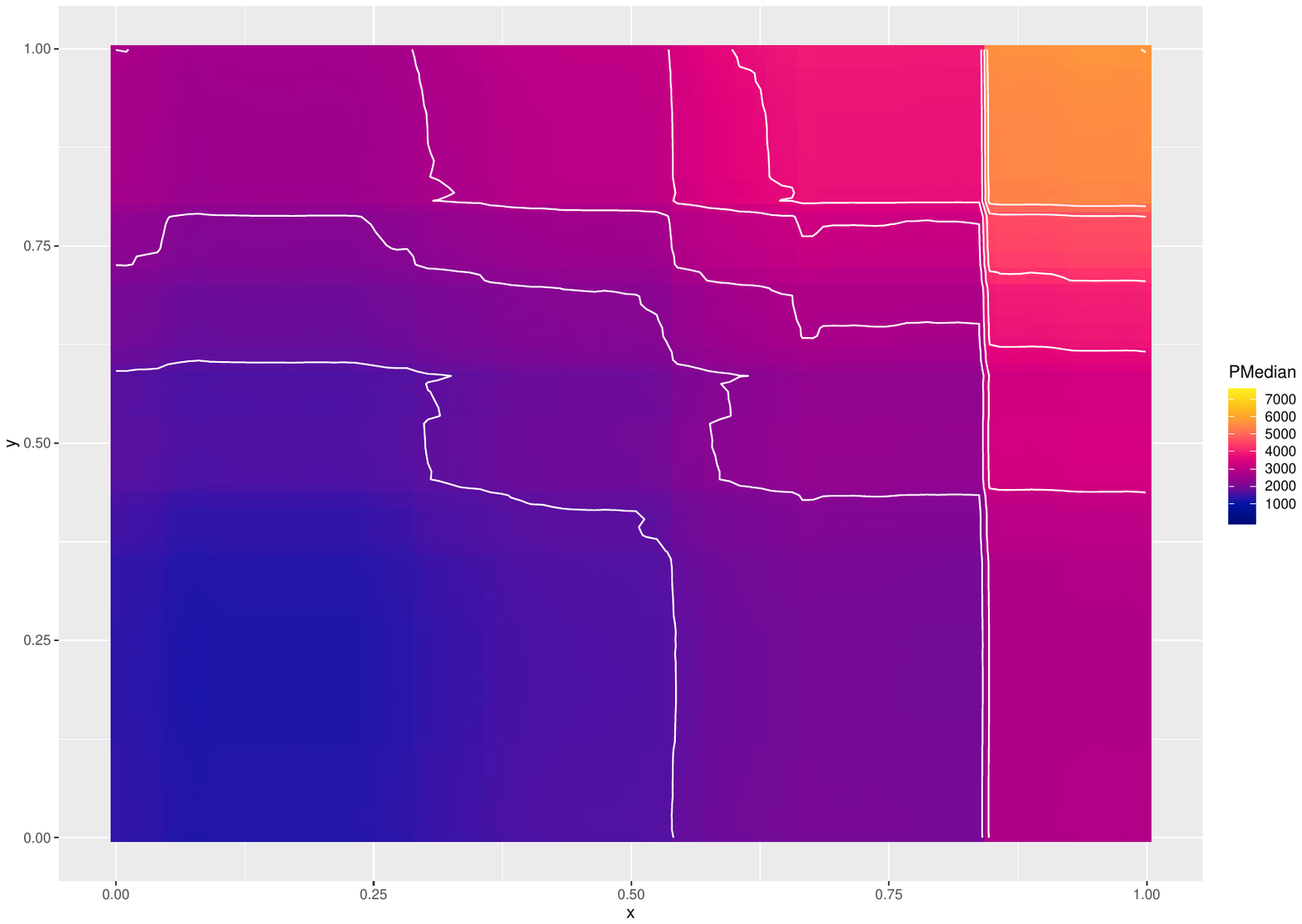}
    \caption{Posterior Median for 10 Trees}
  \end{subfigure}

 \begin{subfigure}{6cm}
    \centering\includegraphics[width=6cm]{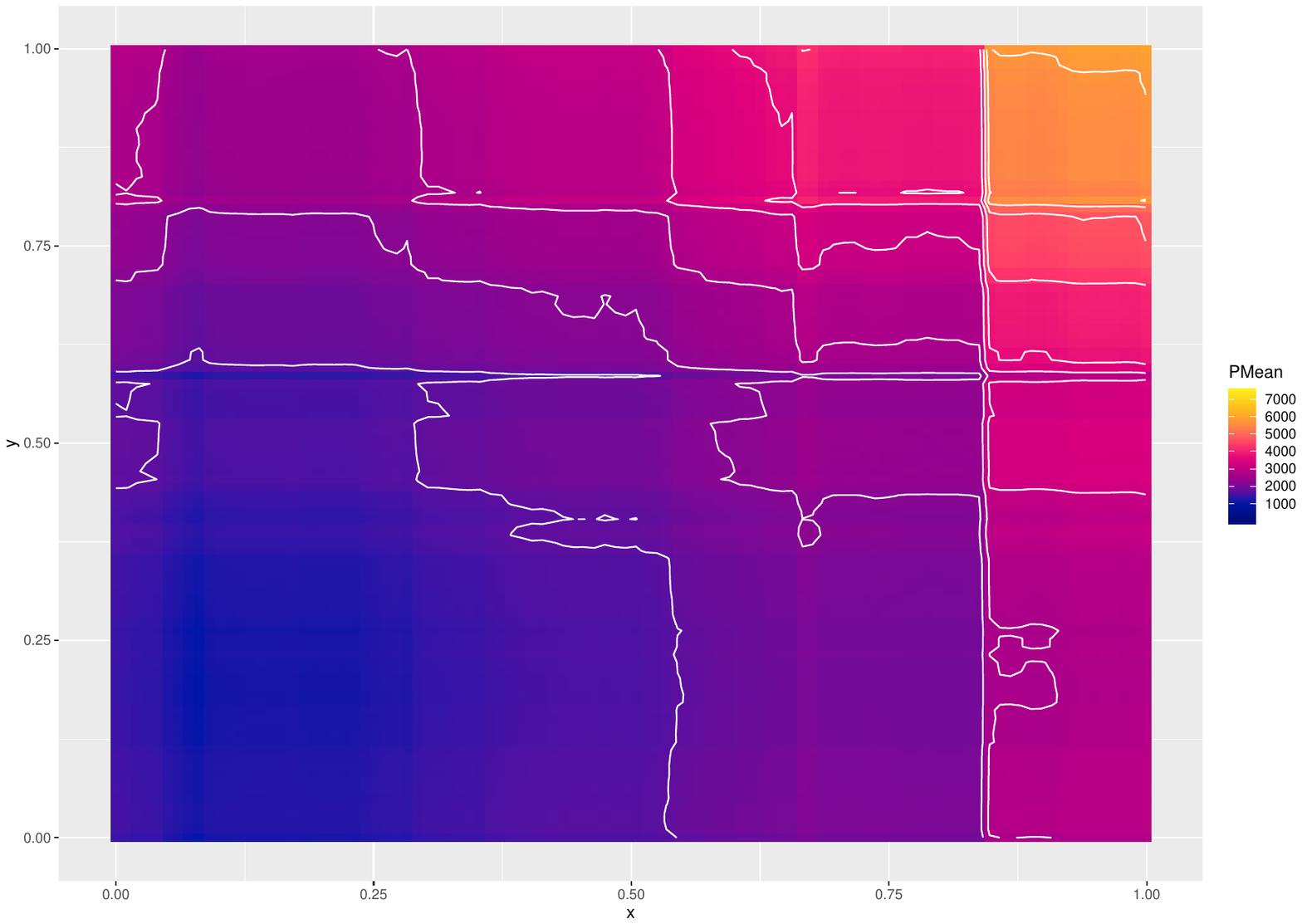}
    \caption{Posterior Mean for 15 Trees}
  \end{subfigure}
  \begin{subfigure}{6cm}
    \centering\includegraphics[width=6cm]{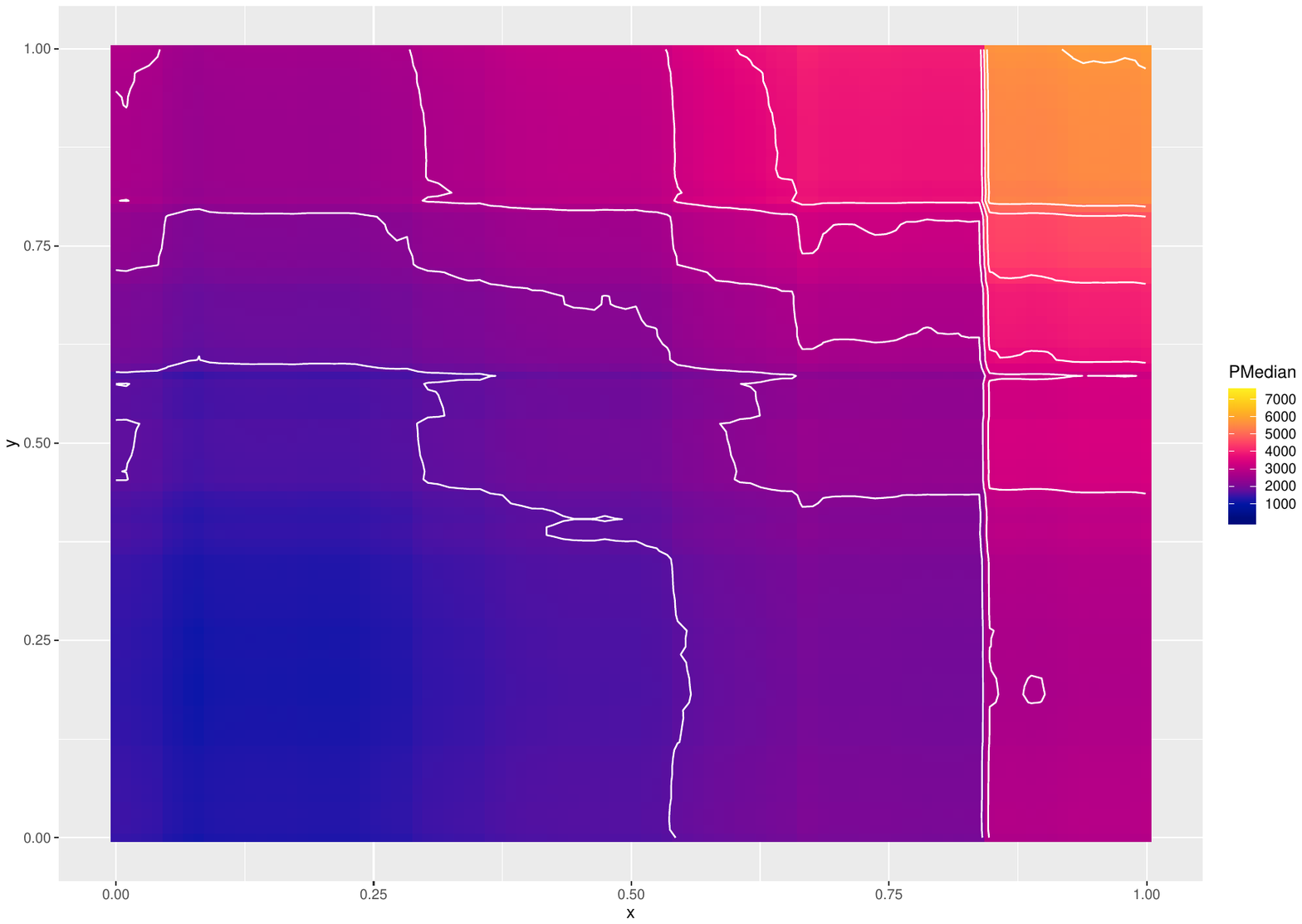}
    \caption{Posterior Median for 15 Trees}
  \end{subfigure}
\caption{Posterior Mean and Posterior Median for 8, 10 and 15 Trees}
\label{fig4}
\end{figure}

\begin{figure}[H] 
  \begin{subfigure}{8 cm}
    \centering\includegraphics[width=8cm]{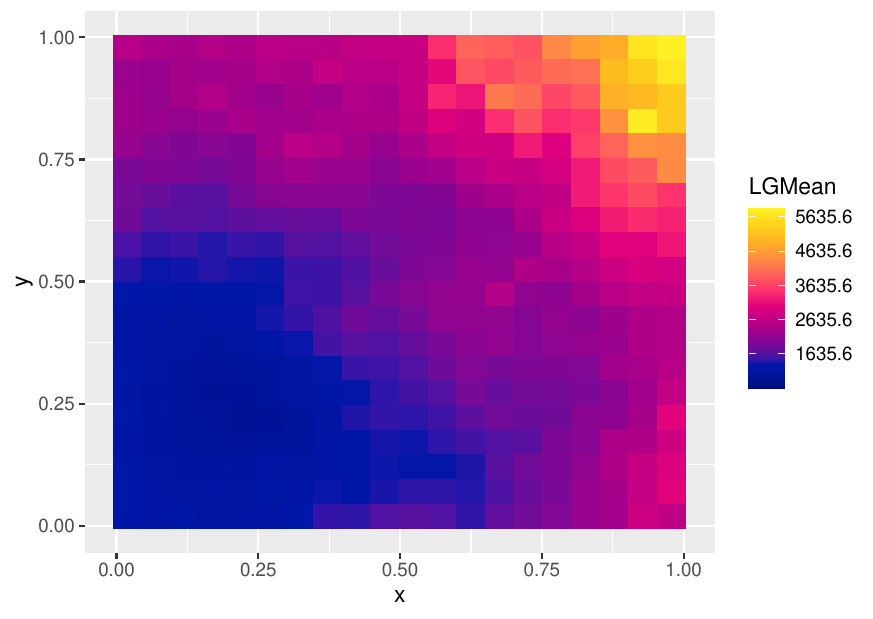}
    \caption{LGCP}
  \end{subfigure}
 \begin{subfigure}{8 cm}
    \centering\includegraphics[width=8cm]{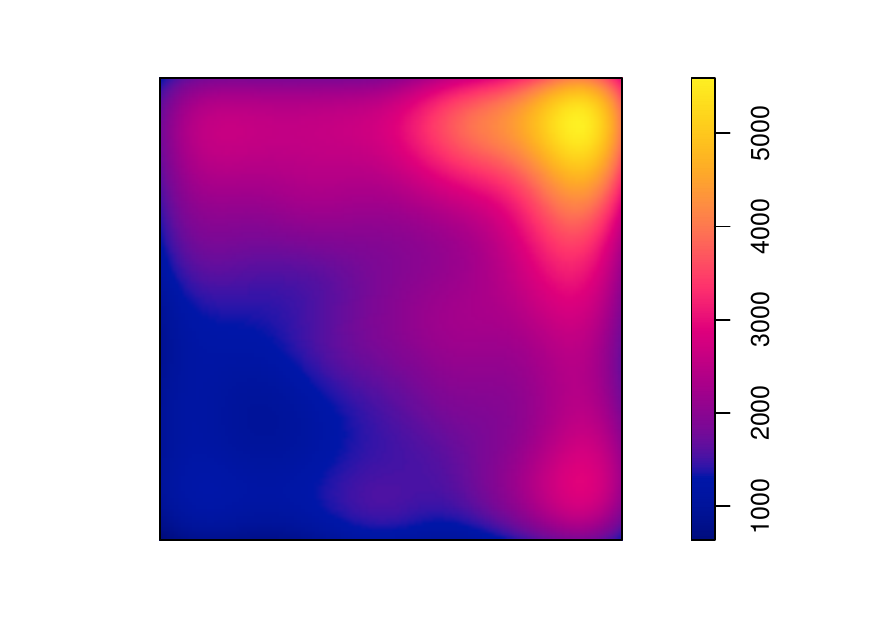}
    \caption{Kernel Smoothing with $h=0.087$}
  \end{subfigure}
\caption{Kernel estimator and inference with spatial log-Gaussian Cox processes.}
\label{fig4l}
\end{figure}

We also considered a two-dimensional Poisson process with intensity $\lambda(x,y)=1000 \, e^{x^2+y^2}$ for $x,y \in [0,1)$.
The outcomes of the algorithm, log-Gaussian Cox processes (LGCP) and kernel smoothing are illustrated in Figures~\ref{fig4}-\ref{fig4l} and Tables~\ref{table5l}-\ref{table6}.  The results demonstrate that the proposed algorithm performs well in this setting, is competitive with the kernel method, and spatial log-Gaussian Cox processes. In this scenario, the hyperparameter $\beta$ has been set equal to 1. In this scenario, the hyperparameter $\beta$ has been set equal to 1. The convergence criteria indicate convergence of the simulated chains (see also Section \ref{SS_AR}).

\begin{table}[H]
\begin{tabular}{ p{4cm} p{3cm} p{3cm}  }
 \hline
 \multicolumn{3}{c}{Inference with spatial log-Gaussian Cox processes}\\ 
 \hline
 grid & AAE  &  RISE  \\
 \hline
$10\times 10$ &195& 263 \\
$20\times 20$ & 182& 245 \\
\hline
\end{tabular}
\caption{Average Absolute Error and Root Integrated Square Error with LGCP for the data in Figure~\ref{fig4}.}
\label{table5l}
\end{table}

\begin{table}[H]
\begin{tabular}{ |p{2cm}||p{2.4cm}|p{2.4cm}|p{2.5cm}|p{2.5cm}| }
 \hline
 \multicolumn{5}{|c|}{Proposed Algorithm} \\
 \hline
 Number of  trees & AAE for Posterior Mean & AAE for Posterior Median &RISE for Posterior Mean & RISE for Posterior Median \\
 \hline
 8 &  177.44 &175.62  &255.23 &258.88 \\
 10&176.52 &174.02 &253.14 &255.92 \\
 15&177.48 &172.62 &254.22 &251.96 \\
 \hline
\end{tabular}
\caption{Average Absolute Error and Root Integrated Square Error for various number of trees for data in Fig.~\ref{fig4}. }
\label{table5}
\end{table}

\begin{table}[H]
\begin{tabular}{ |p{4cm}||p{3cm}|p{3cm}|  } \hline
 \multicolumn{3}{|c|}{Kernel Smoothing} \\
 \hline
 Bandwidth (sigma) & AAE  &  RISE  \\
 \hline
0.03 &  360.11 & 463.1  \\
0.04 &  277.89 &353.82   \\
0.087 (LCV)& 167.74 & 227.85  \\
0.095& 166.51 & 230.27  \\
 \hline
\end{tabular}
\caption{Average Absolute Error and Root Integrated Square Error for fixed bandwidth estimators for data in Fig.~\ref{fig4}.}
\label{table6}
\end{table}

 \subsection{Inhomogeneous five dimensional Poisson Process with Gaussian intensity} 

Our next example is a five dimensional Poisson process with intensity $\lambda(x)=50e^{x^Tx}$ for $x\in[0,1)^5$ and the generated process via thinning consists of 343 points. The statistics are presented in Tables~\ref{Table_MD3}~and~\ref{Table_MD4} for our algorithm and kernel smoothing, respectively.  We have checked that the Gelman-Rubin criterion indicates convergence of chains. 

\begin{table}[H]
\begin{tabular}{ |p{1.2cm}||p{1.5cm}|p{1.5cm}|p{1.8cm}|p{1.8cm}| }
 \hline
 \multicolumn{5} {|c|} {Proposed Algorithm}  \\
 \hline
 Number of  trees & AAE for Mean & AAE for Median & RISE for Mean & RISE for Median \\
 \hline
 8 & 65.6 & 66.7 & 104.2 &106.2  \\
 \hline
 10 & 66.94 & 67 & 106.09 &106.36   \\
 \hline
\end{tabular}
\caption{Average Absolute Error and Root Integrated Square Error with different number of trees.}
\label{Table_MD3}
\end{table}

\begin{table}[H]
\begin{tabular}{ |p{3.2cm}||p{2.2cm}|p{2.2cm}| }
 \hline
 \multicolumn{3} {|c|} {Kernel Smoothing}  \\
 \hline
 $h$ & AAE &  RISE  \\
 \hline
0.03 &631.1 & 5060.9  \\
\hline
0.06 &409.6 & 825.1   \\
\hline
0.08 & 287.8 & 419.8   \\
\hline
0.1 &213.5 & 295.6  \\
\hline
0.15 (LCV) &181.2 & 278.5  \\
\hline
0.3 &258.4 & 363 \\
\hline
0.5 &311.2 & 409.7  \\
\hline
\end{tabular}
\caption{Average Absolute Error and Root Integrated Square Error for  various isotropic variance matrices.}
\label{Table_MD4}
\end{table}

\
\section {One dimensional Poisson Process with stepwise intensity}

\subsection{5 Trees}
\noindent We run 3 parallel chains  each for 200000 iterations keeping every 100th sample. 
\begin{figure}[H] 
    \centering\includegraphics[width=8cm]{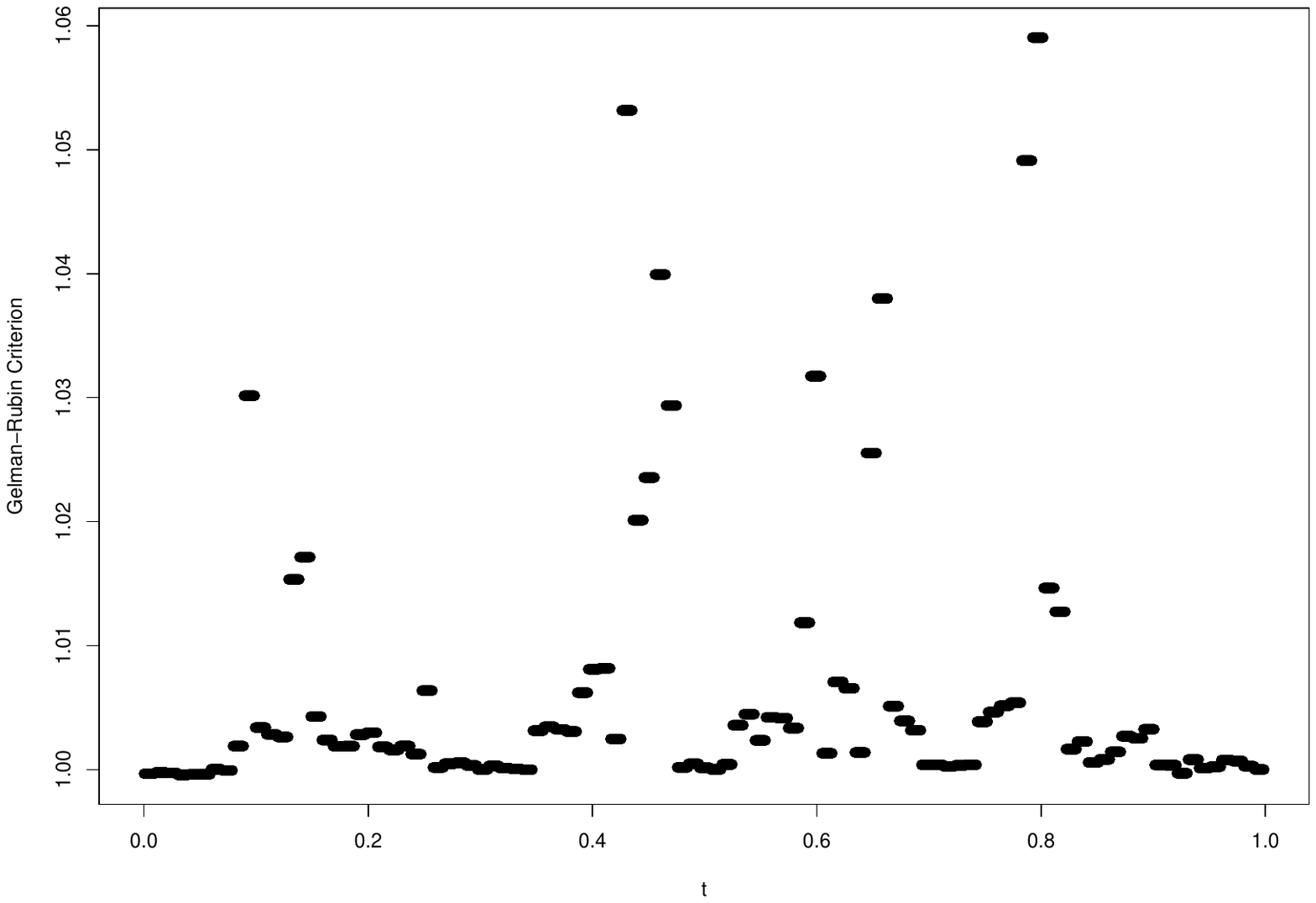}
\caption{The Gelman-Rubin Criterion for 5 Trees}
\label{SL_1DS1}
\end{figure}

\begin{figure}[H] 
    \centering\includegraphics[width=8cm]{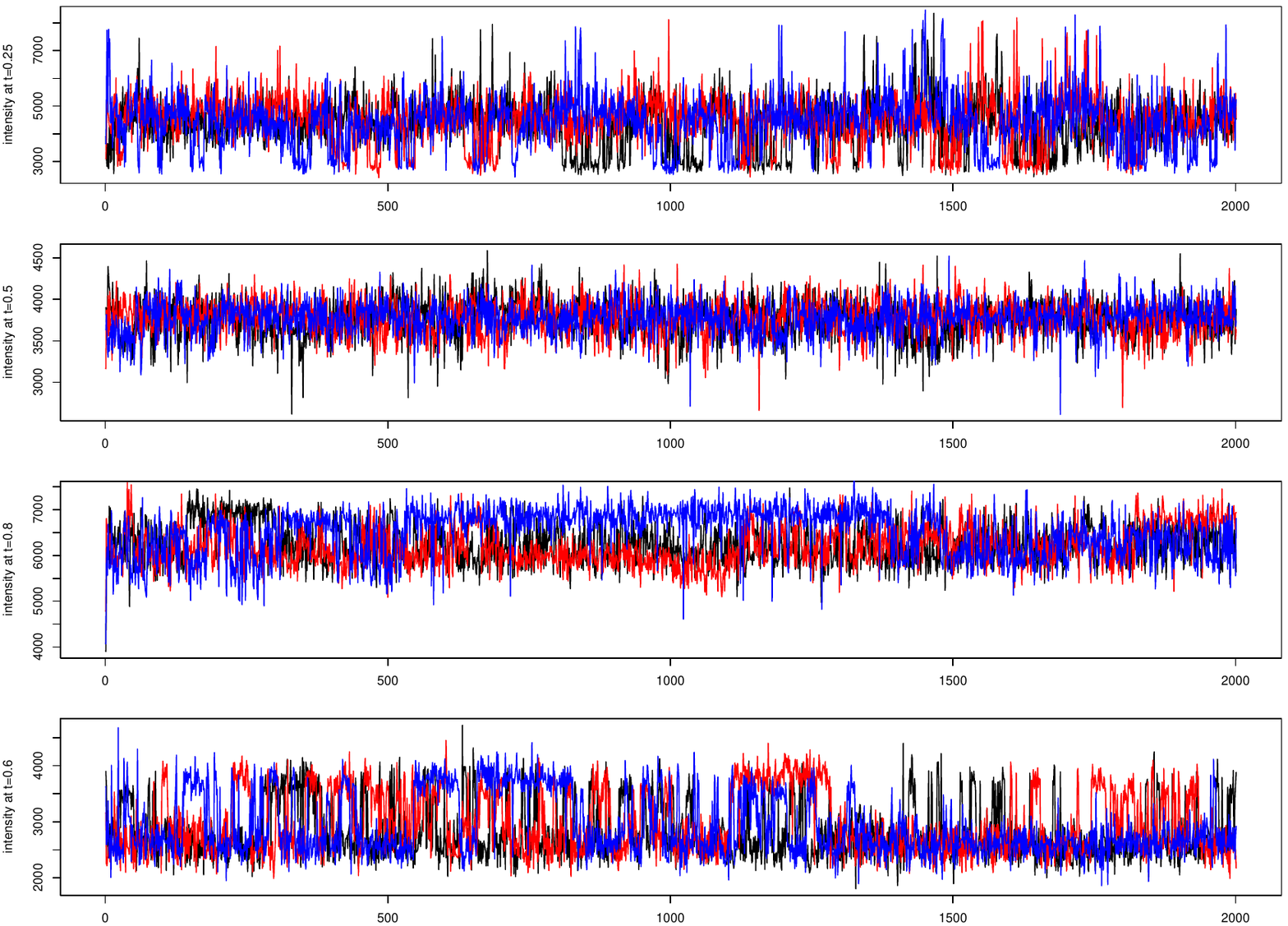}
    \caption{Trace plots for 5 Trees}
\label{SL_1DS2}
\end{figure}

\begin{figure}[H] 
    \centering\includegraphics[width=8cm]{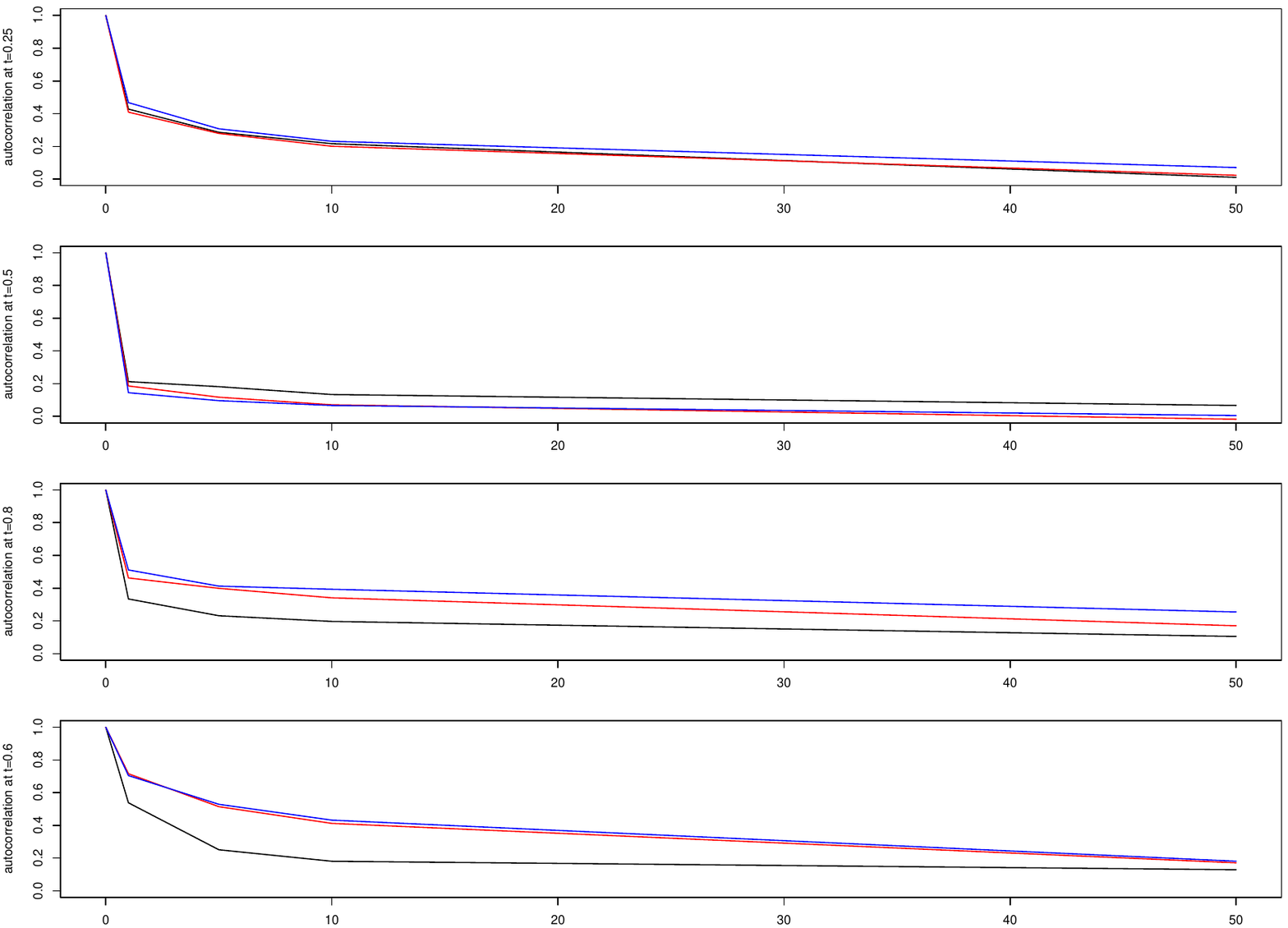}
    \caption{Autocorrelation plots for 5 Trees}
    \label{SL_1DS3}
\end{figure}

\begin{figure}[H] 
    \centering\includegraphics[width=8cm]{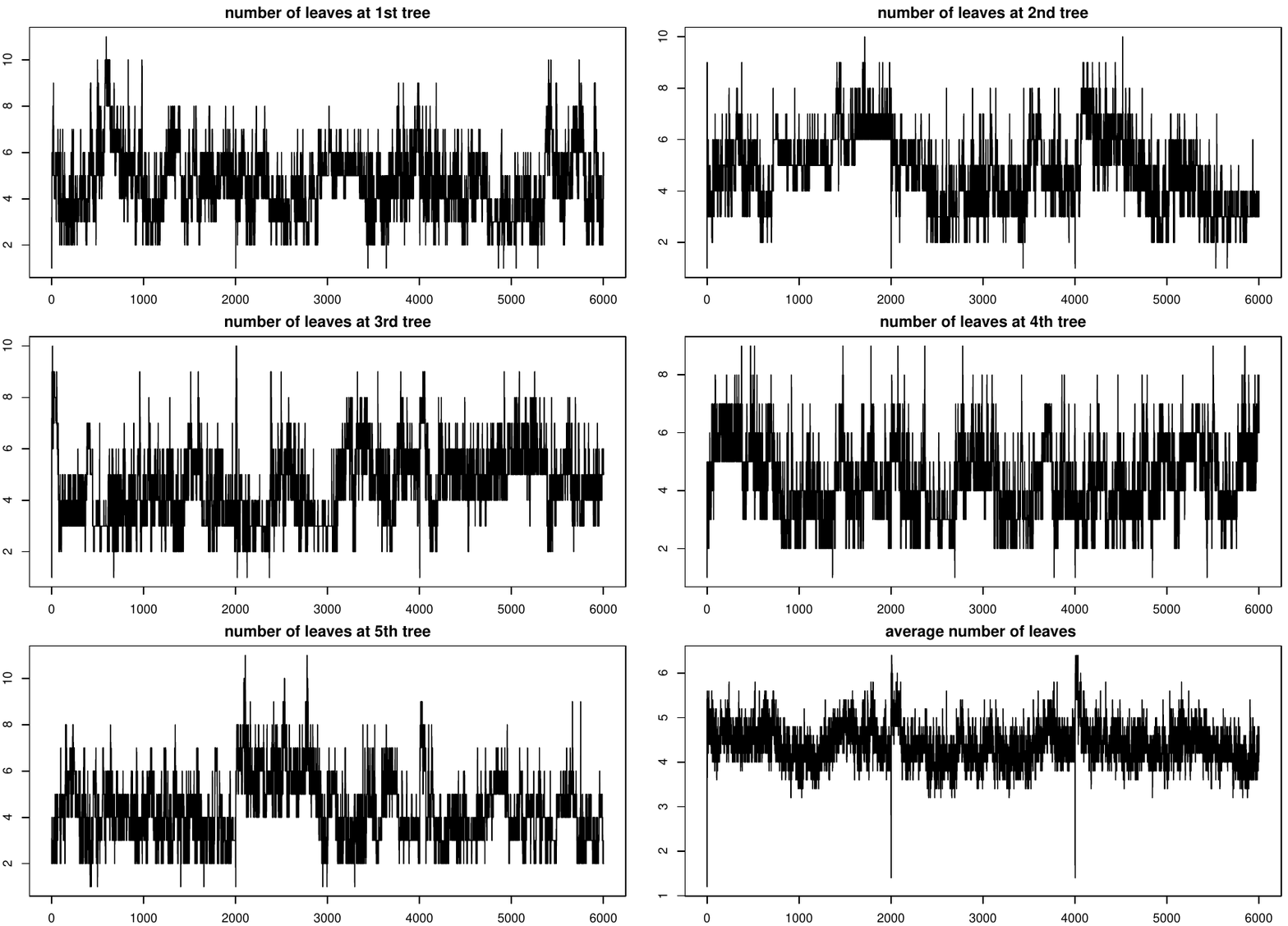}
    \caption{Average number of leaves at trees}
\end{figure}

\begin{figure}[H] 
    \centering\includegraphics[width=8cm]{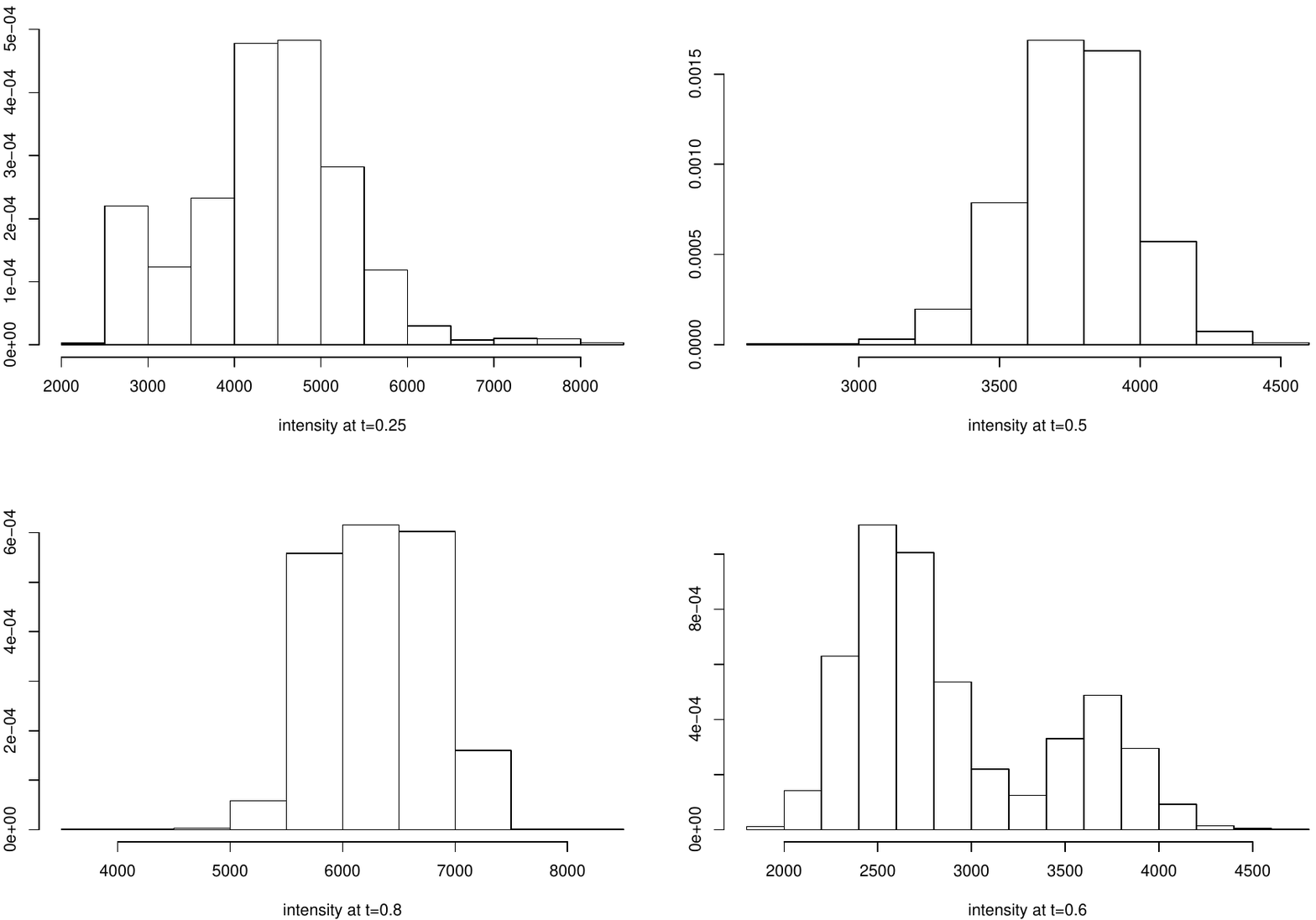}
    \caption{Density of the estimated intensity for 5 Trees}
\end{figure}

\begin{figure}[H] 
    \centering\includegraphics[width=8cm]{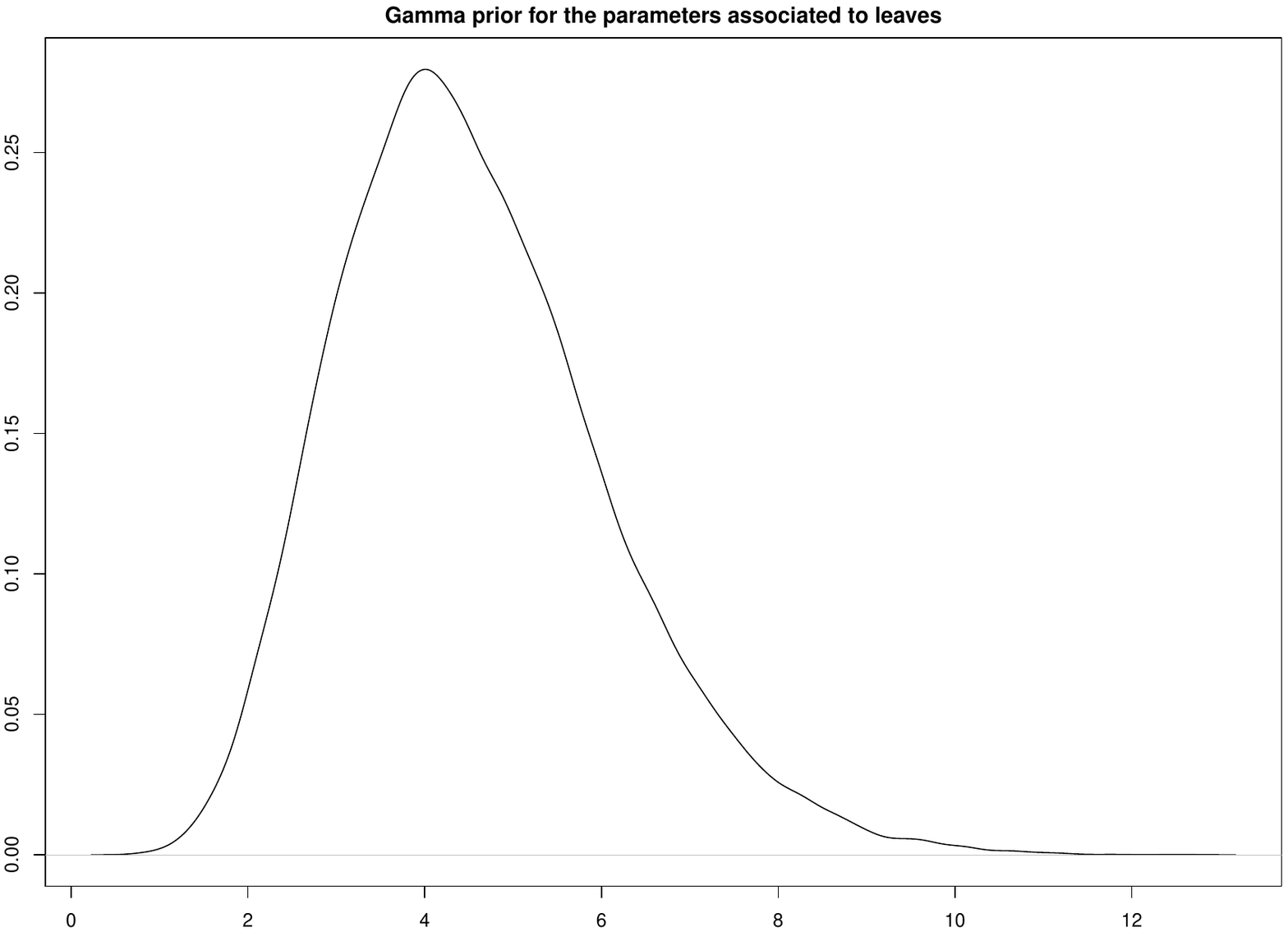}
    \caption{Prior for 5 Trees}
\end{figure}

\subsection{7 Trees}
\noindent We run 3 parallel chains  each for 200000 iterations keeping every 100th sample. 
\begin{figure}[H] 
    \centering\includegraphics[width=8cm]{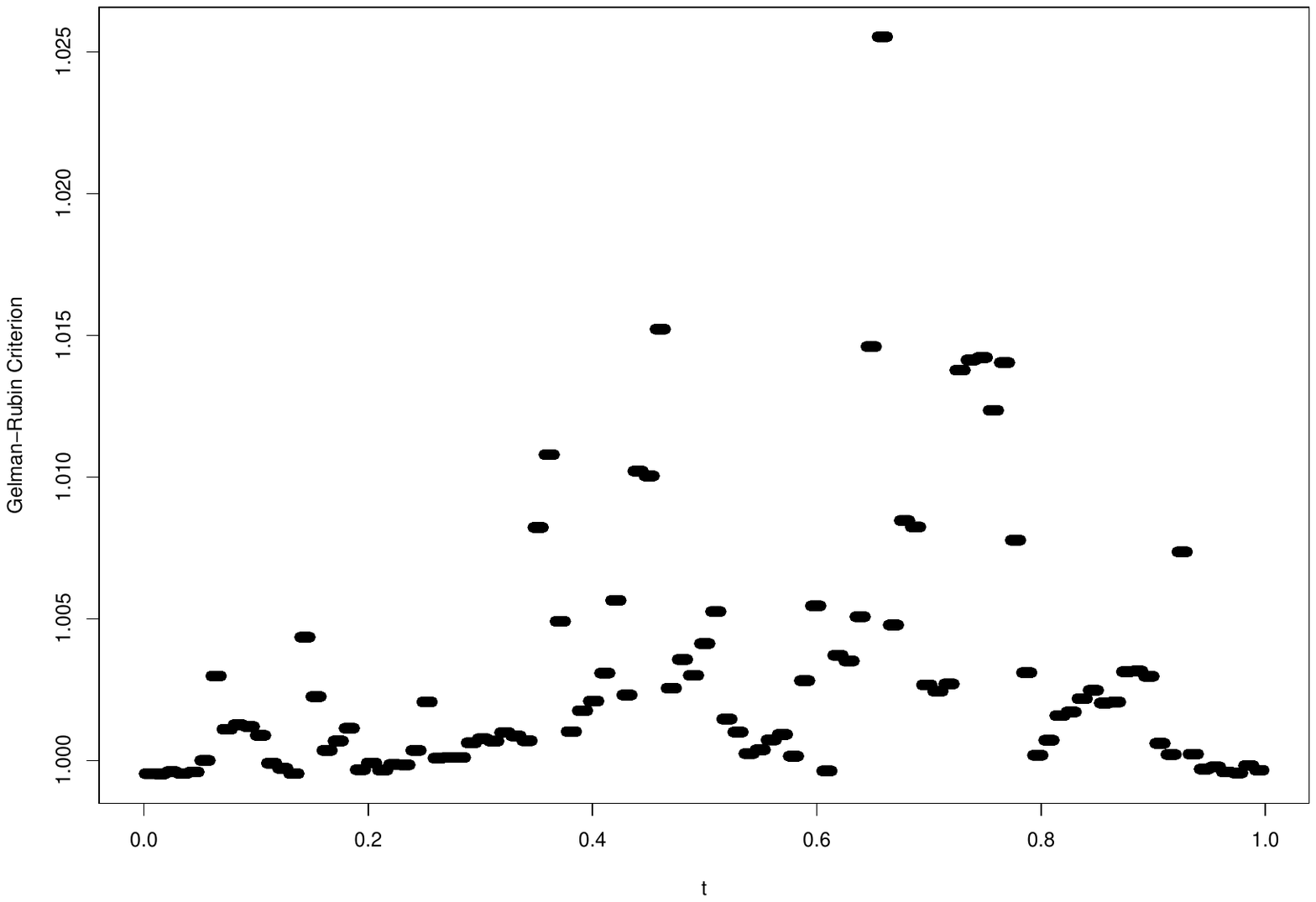}
\caption{The Gelman-Rubin Criterion for 7 Trees}
\label{SL_1DS4}
\end{figure}

\begin{figure}[H] 
    \centering\includegraphics[width=8cm]{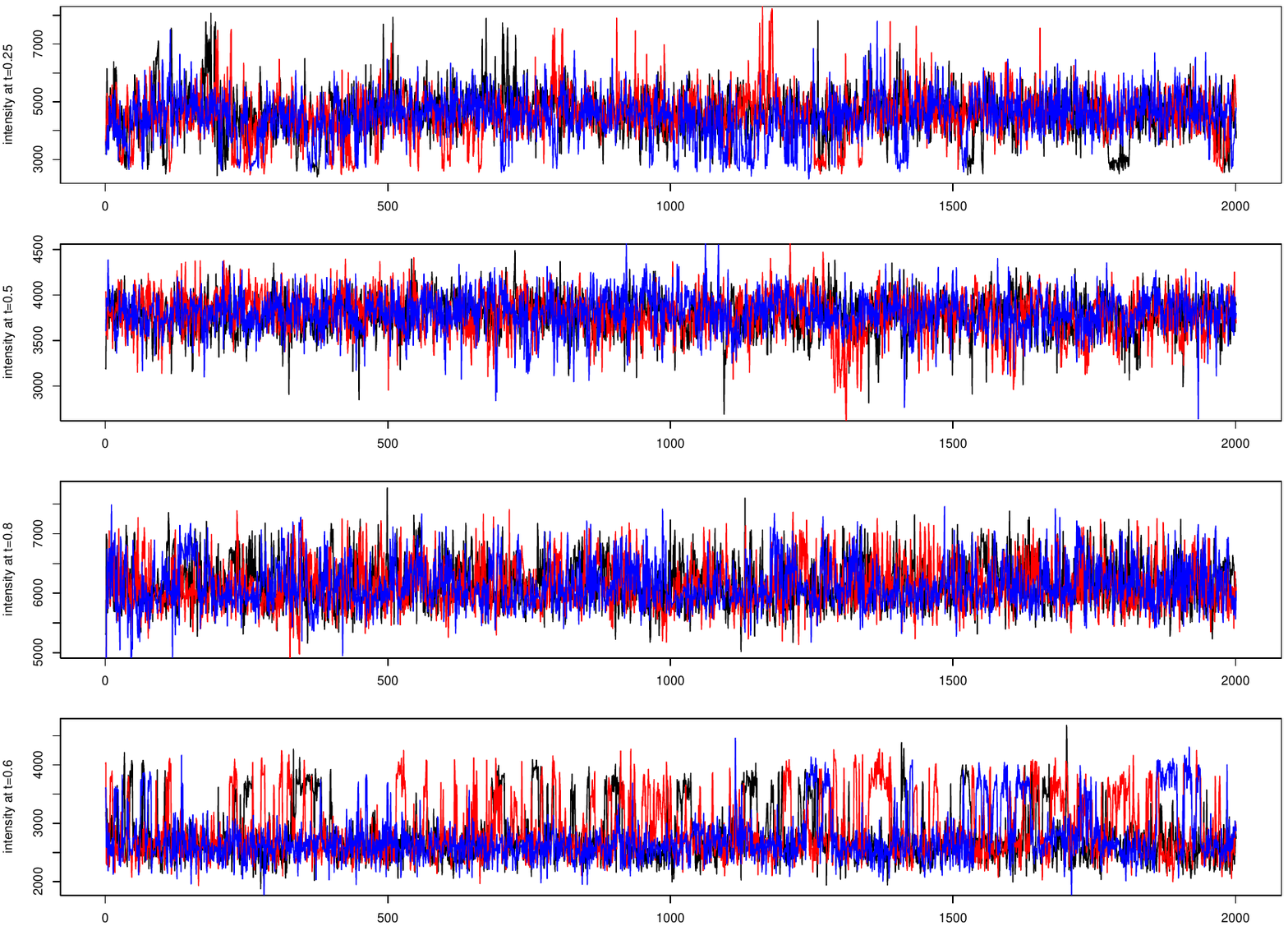}
    \caption{Trace plots for 7 Trees}
\label{SL_1DS5}
\end{figure}

\begin{figure}[H] 
    \centering\includegraphics[width=8cm]{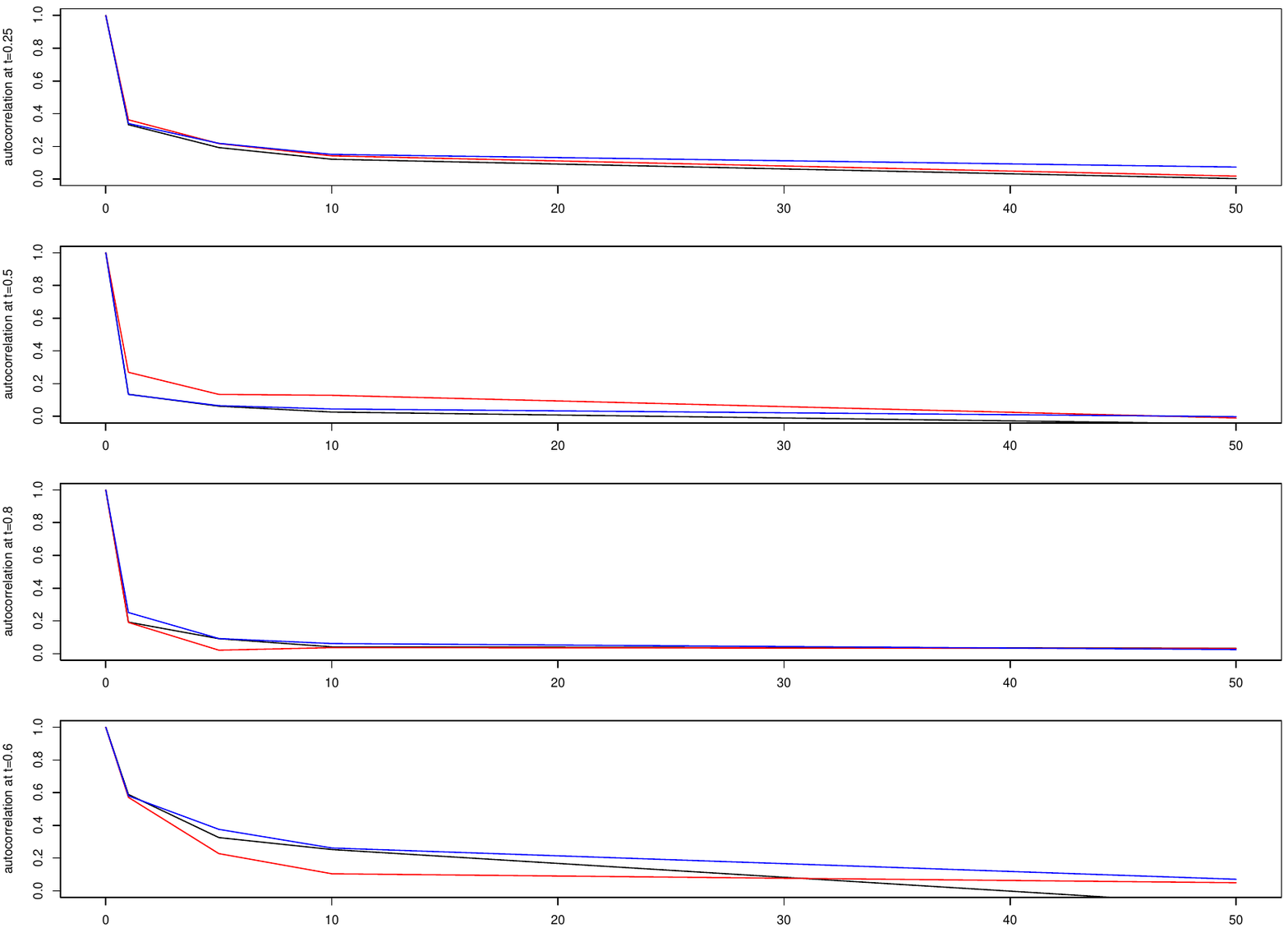}
    \caption{Autocorrelation plots for 7 Trees}
\label{SL_1DS6}
\end{figure}

\begin{figure}[H] 
    \centering\includegraphics[width=8cm]{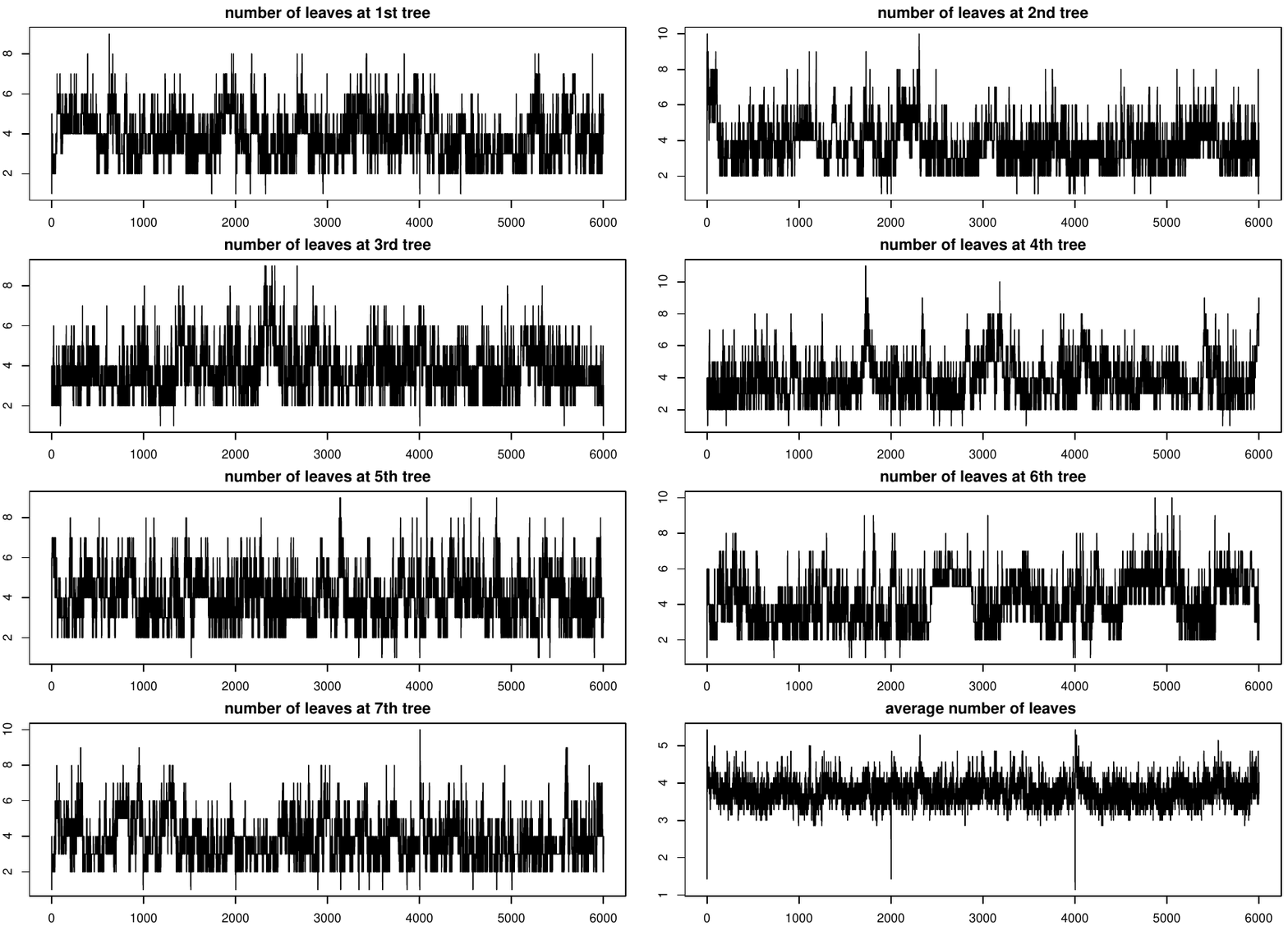}
    \caption{Average number of leaves at trees}
\end{figure}

\begin{figure}[H] 
    \centering\includegraphics[width=8cm]{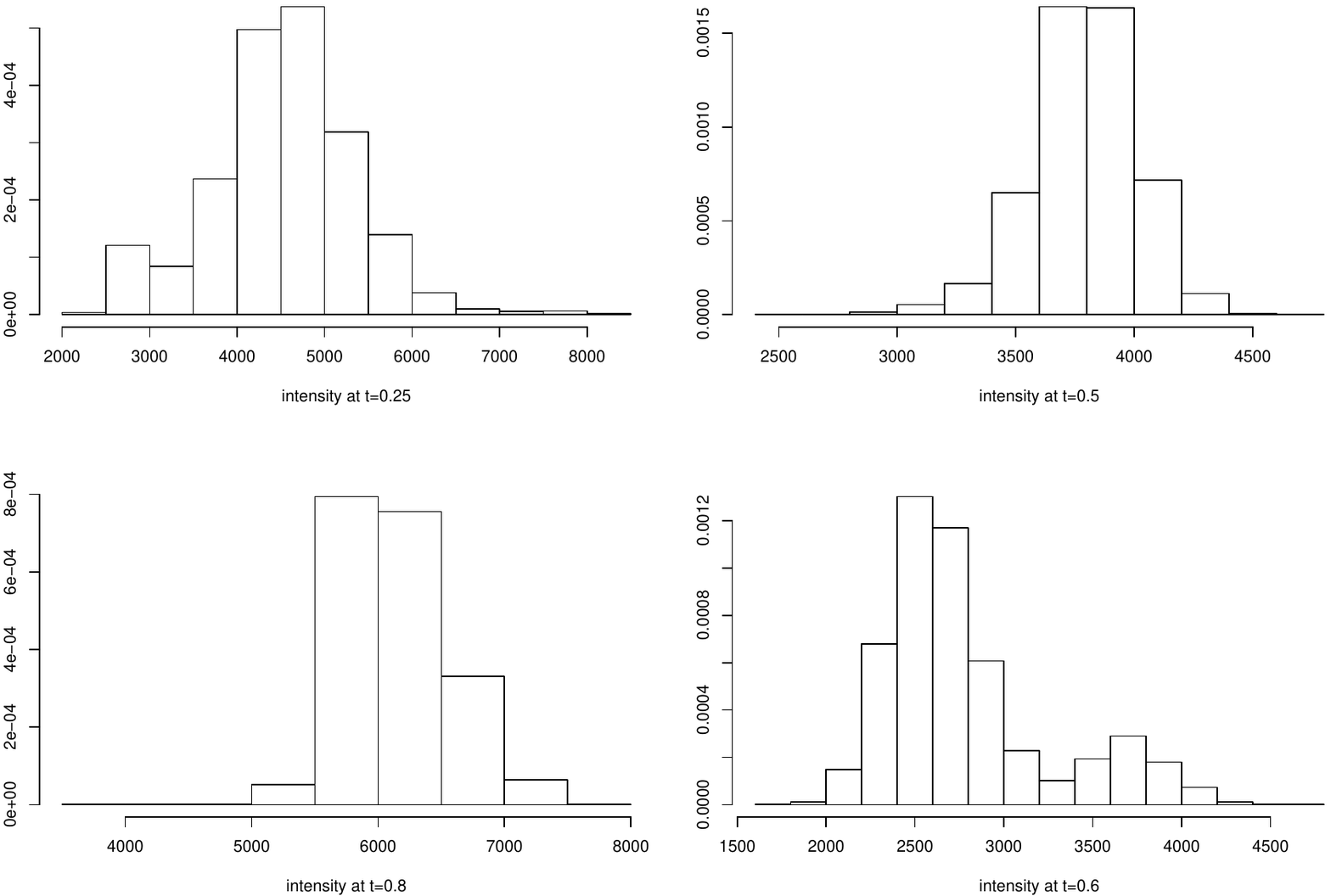}
    \caption{Density of the estimated intensity for 7 Trees}
\end{figure}

\begin{figure}[H] 
    \centering\includegraphics[width=8cm]{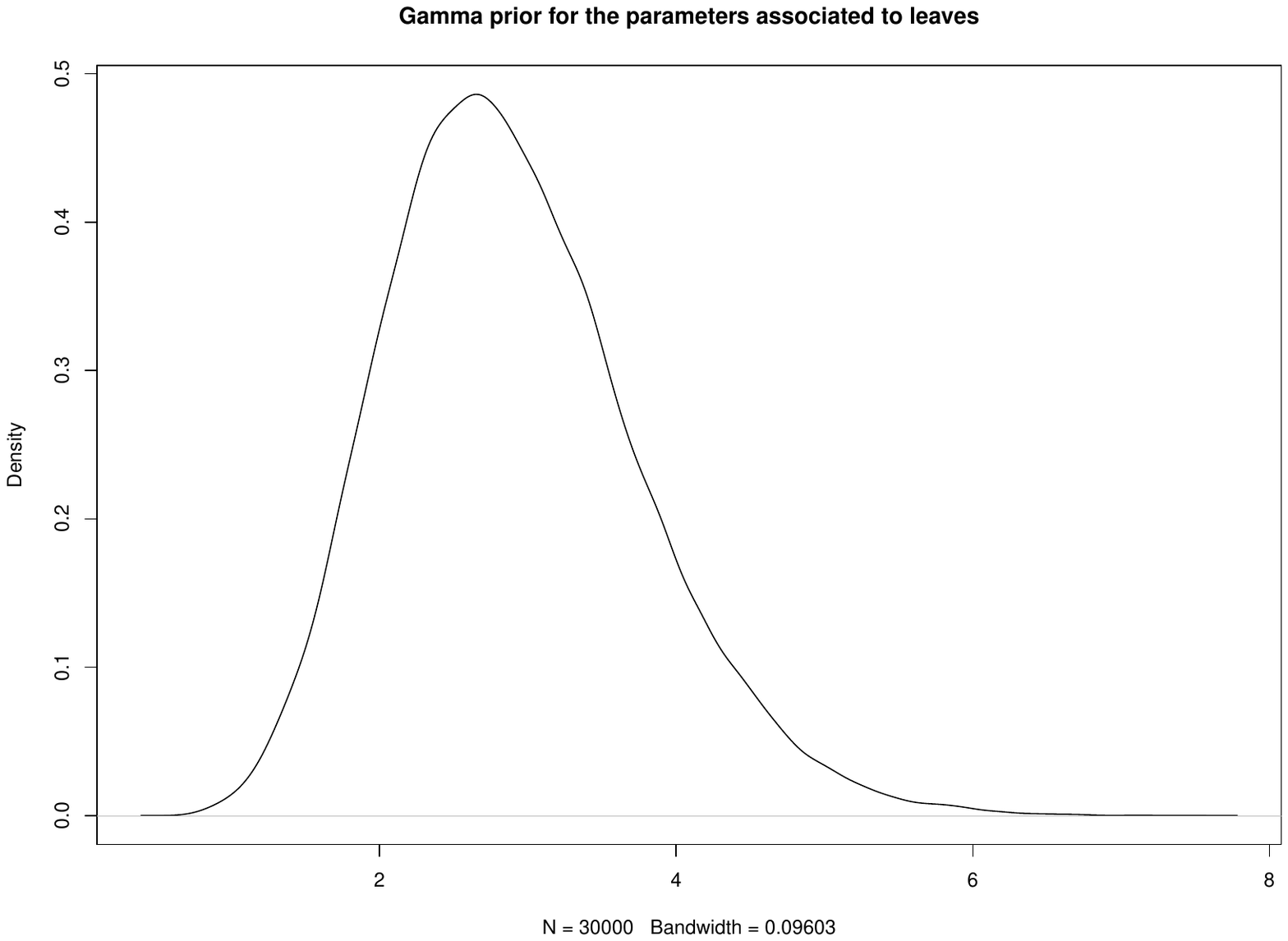}
    \caption{Prior for 7 Trees}
\end{figure}

\section{One dimensional Poisson Process with with continuously varying intensity}
\subsection{5 Trees}
\noindent We run 3 parallel chains  each for 100000 iterations keeping every 50th sample. 
\begin{figure}[H] 
    \centering\includegraphics[width=8cm]{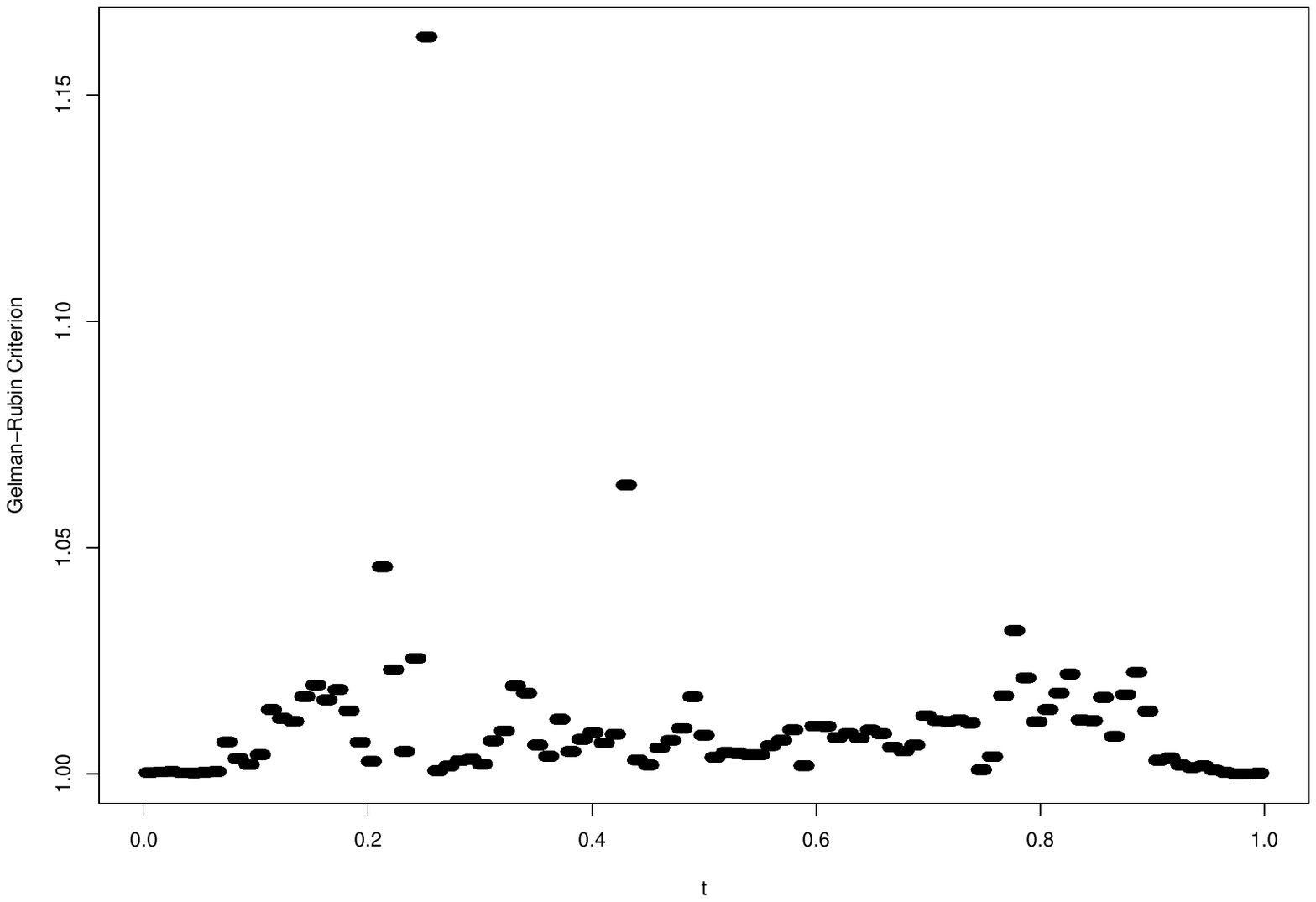}
\caption{The Gelman-Rubin Criterion for 5 Trees}
\label{SL_1DSM1}
\end{figure}

\begin{figure}[H] 
    \centering\includegraphics[width=8cm]{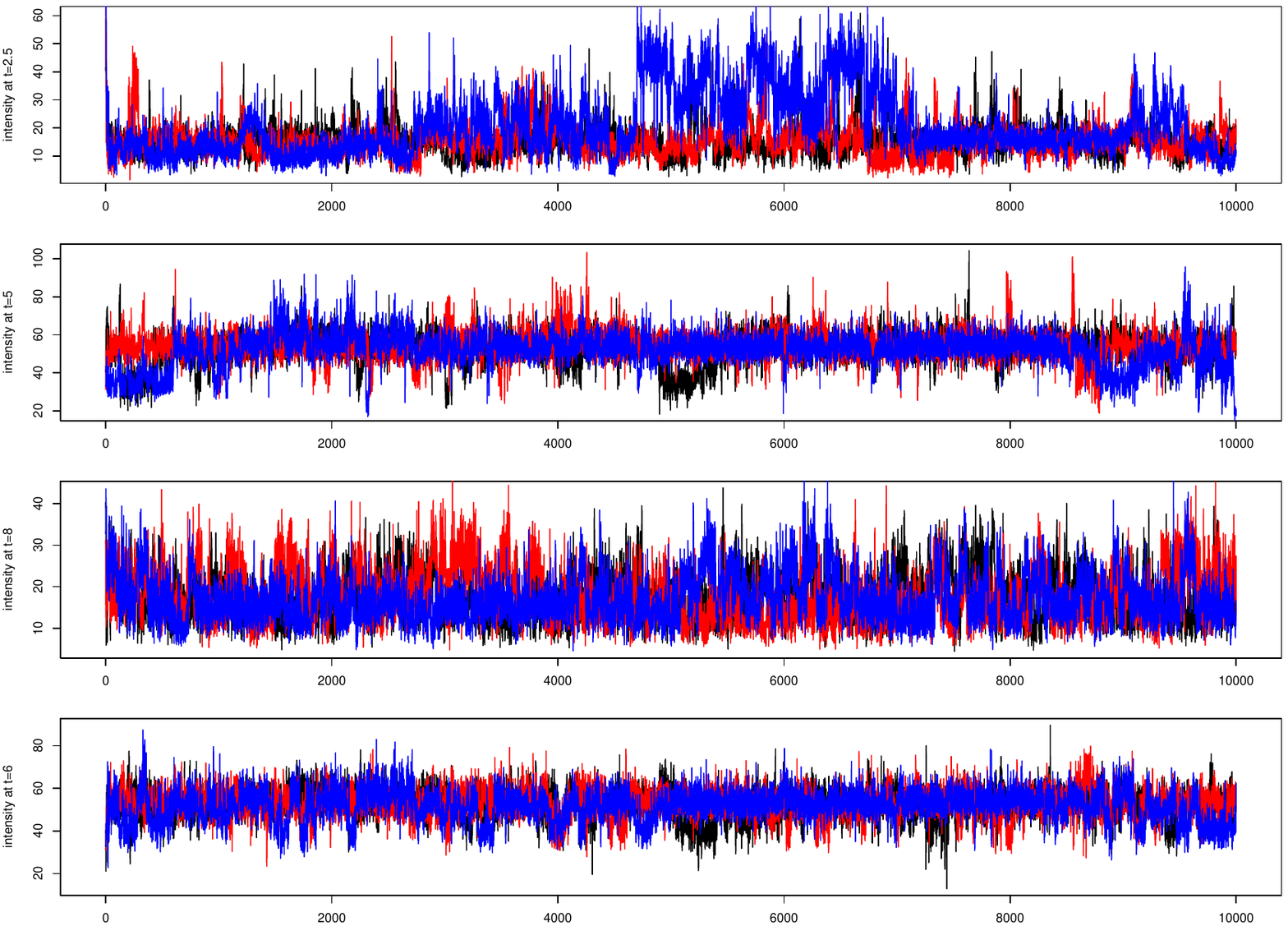}
    \caption{ Trace plots for 5 Trees}
\label{SL_1DSM2}
\end{figure}

\begin{figure}[H] 
    \centering\includegraphics[width=8cm]{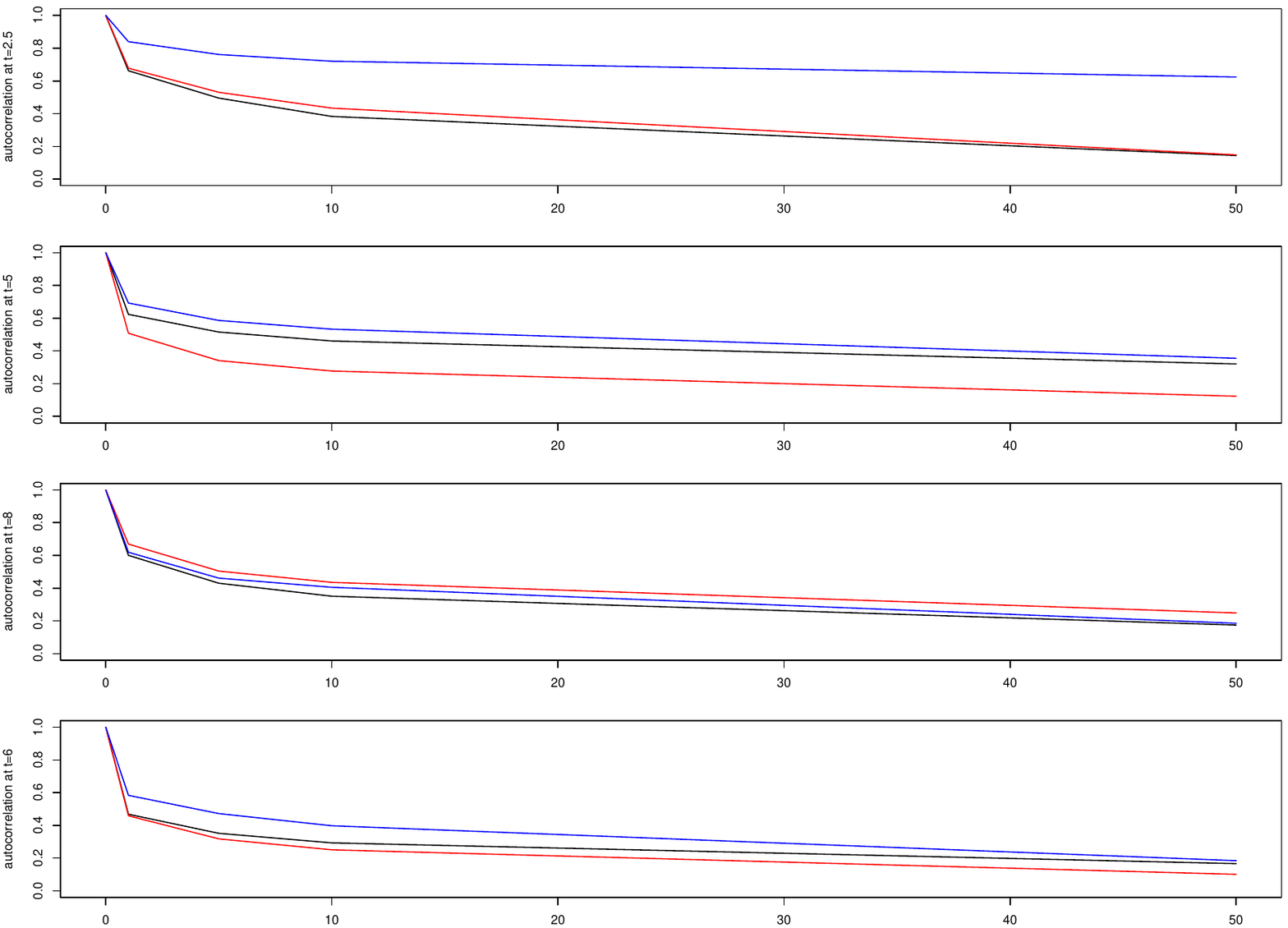}
    \caption{Autocorrelation plots for 5 Trees}
\label{SL_1DSM3}
\end{figure}


\begin{figure}[H] 
    \centering\includegraphics[width=8cm]{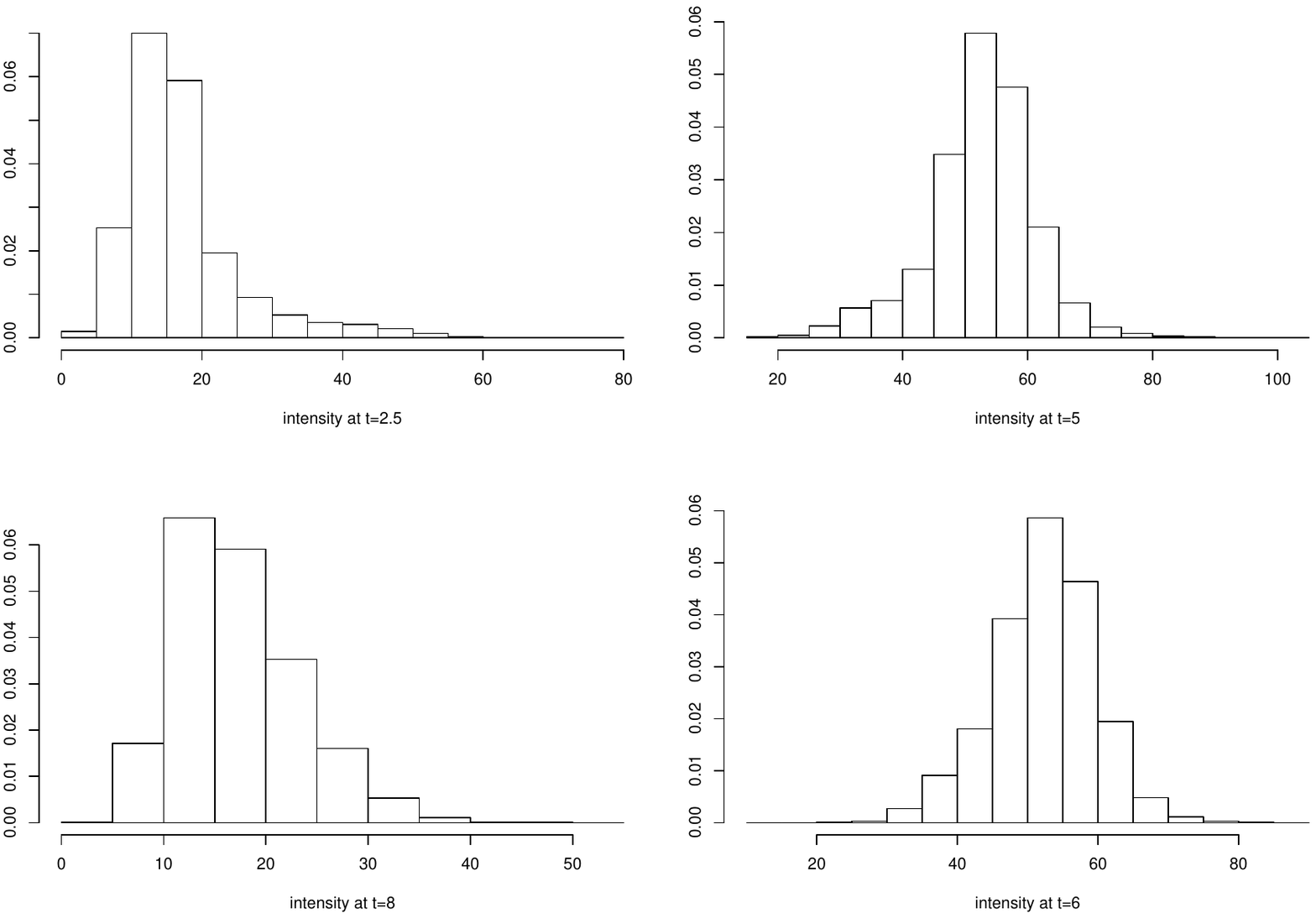}
    \caption{Density of the estimated intensity for 5 Trees}
\end{figure}

\begin{figure}[H] 
    \centering\includegraphics[width=8cm]{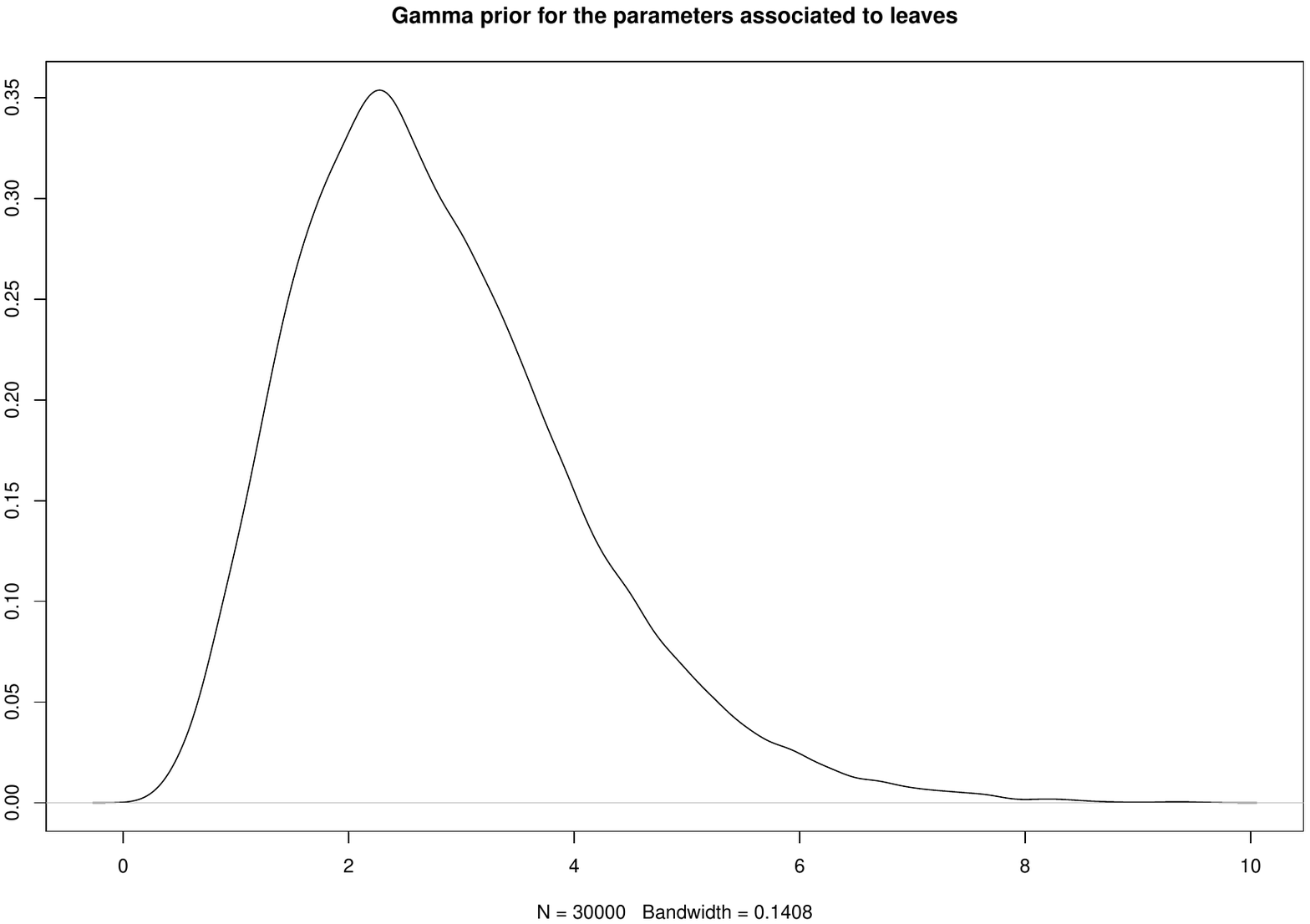}
    \caption{Prior for 5 Trees}
\end{figure}

\subsection{10 Trees}
\noindent We run 3 parallel chains  each for 100000 iterations keeping every  50th sample. 
\begin{figure}[H] 
    \centering\includegraphics[width=8cm]{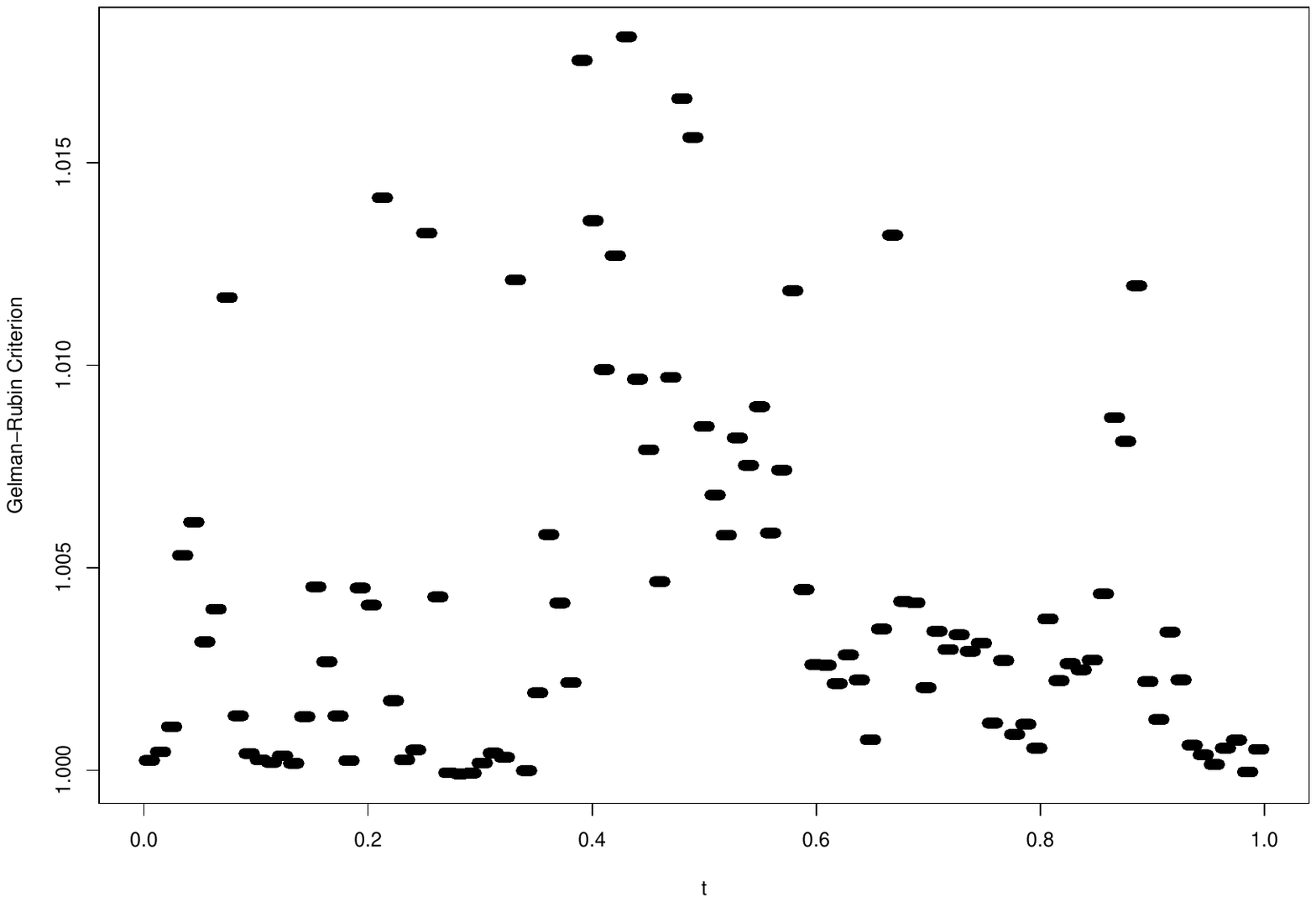}
\caption{The Gelman-Rubin Criterion for 10 Trees}
\label{SL_1DSM4}
\end{figure}

\begin{figure}[H] 
    \centering\includegraphics[width=8cm]{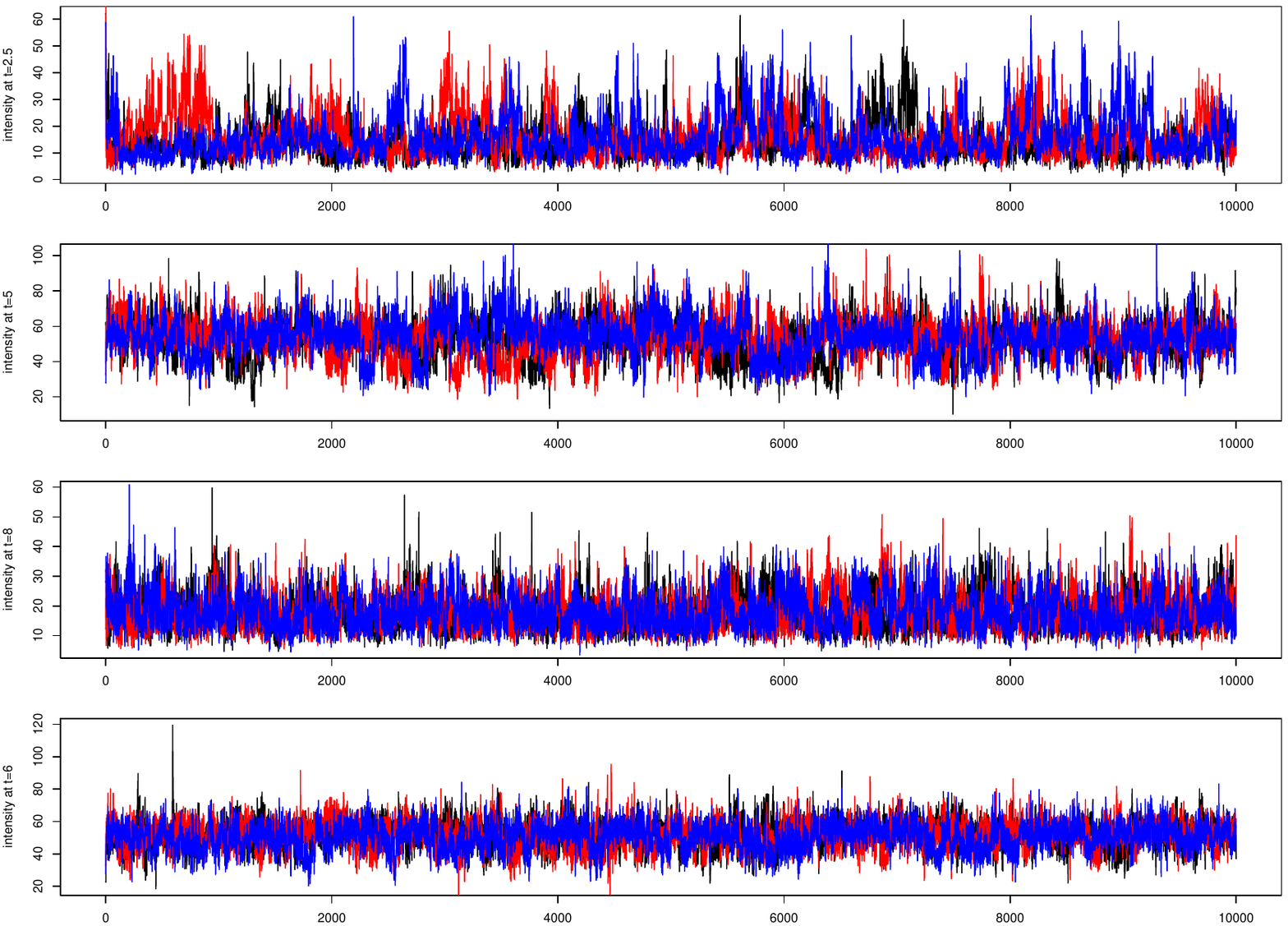}
    \caption{Trace plots for 10Trees}
\label{SL_1DSM5}
\end{figure}

\begin{figure}[H] 
    \centering\includegraphics[width=8cm]{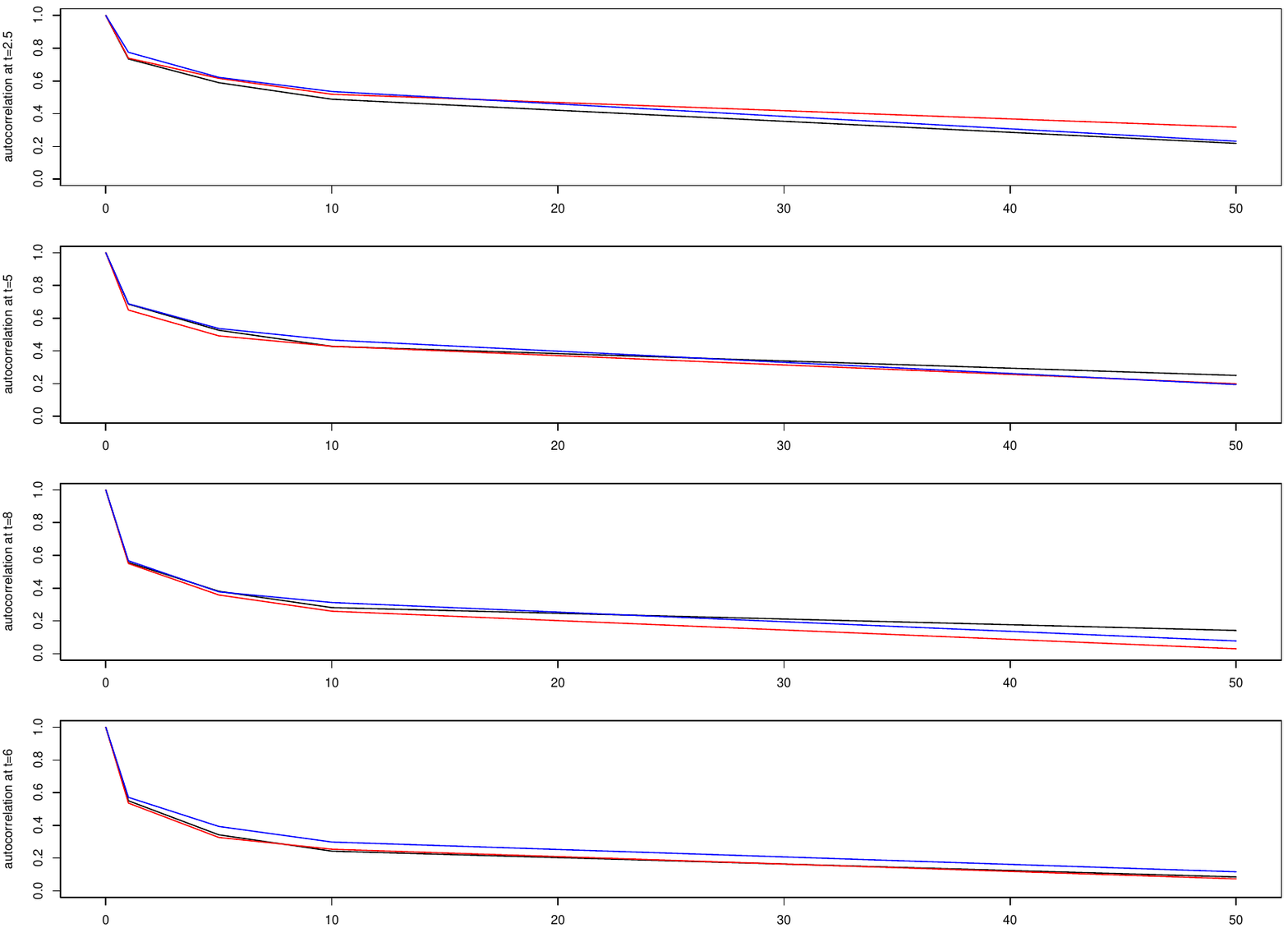}
    \caption{Autocorrelation plots for 10 Trees}
\label{SL_1DSM6}
\end{figure}

\begin{figure}[H] 
    \centering\includegraphics[width=8cm]{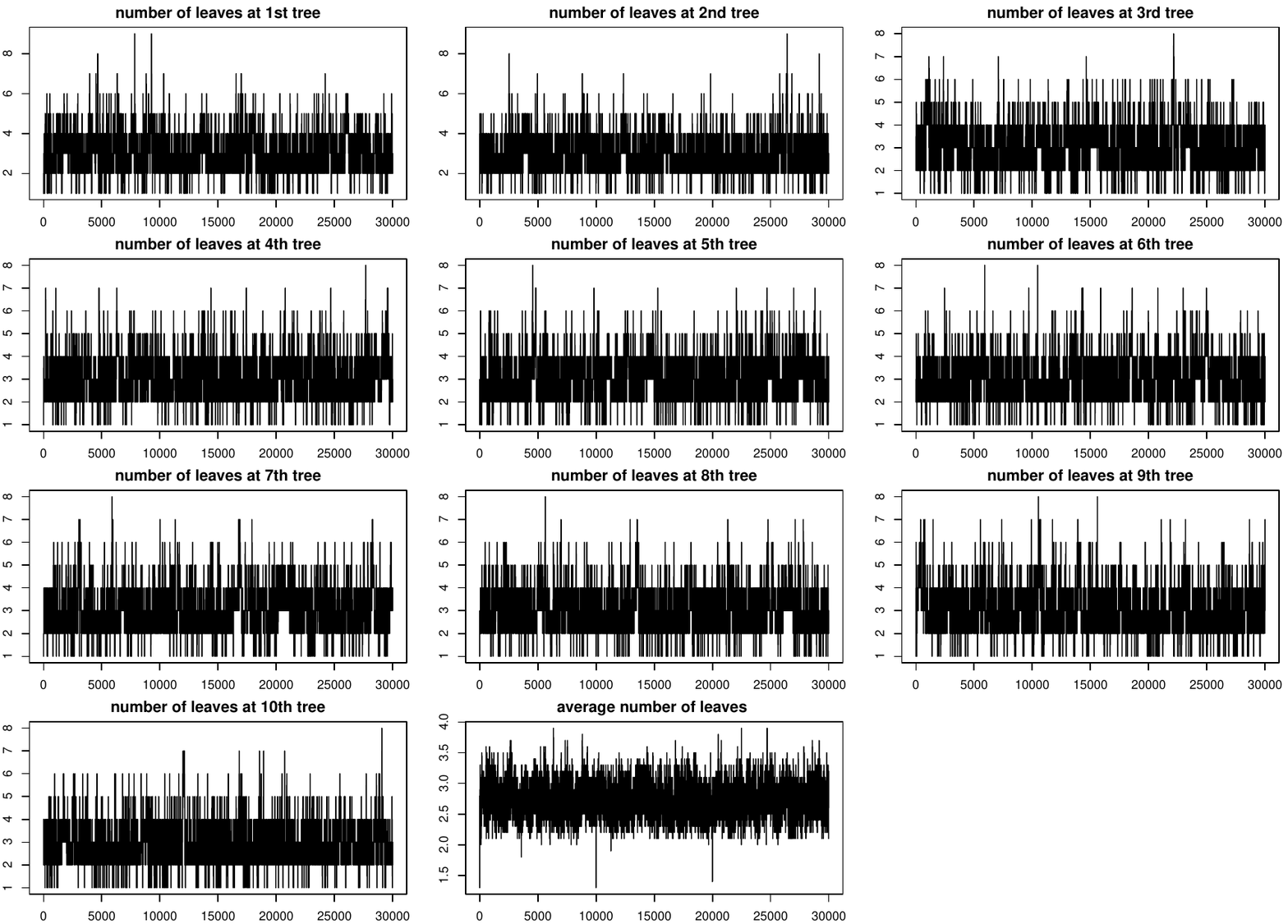}
    \caption{ Average number of leaves at trees}
\end{figure}

\begin{figure}[H] 
    \centering\includegraphics[width=8cm]{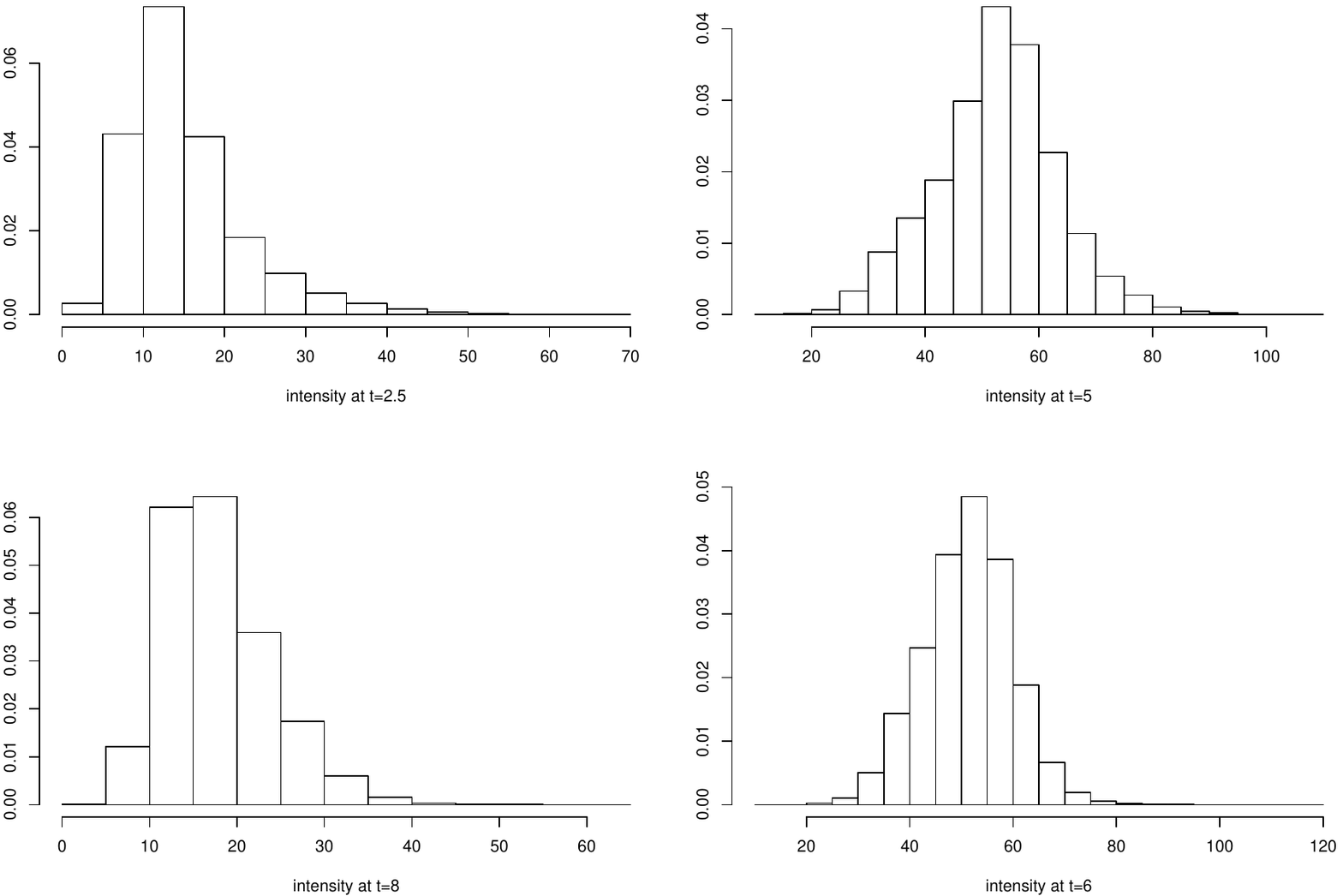}
    \caption{Density of the estimated intensity for 10 Trees}
\end{figure}

\begin{figure}[H] 
    \centering\includegraphics[width=8cm]{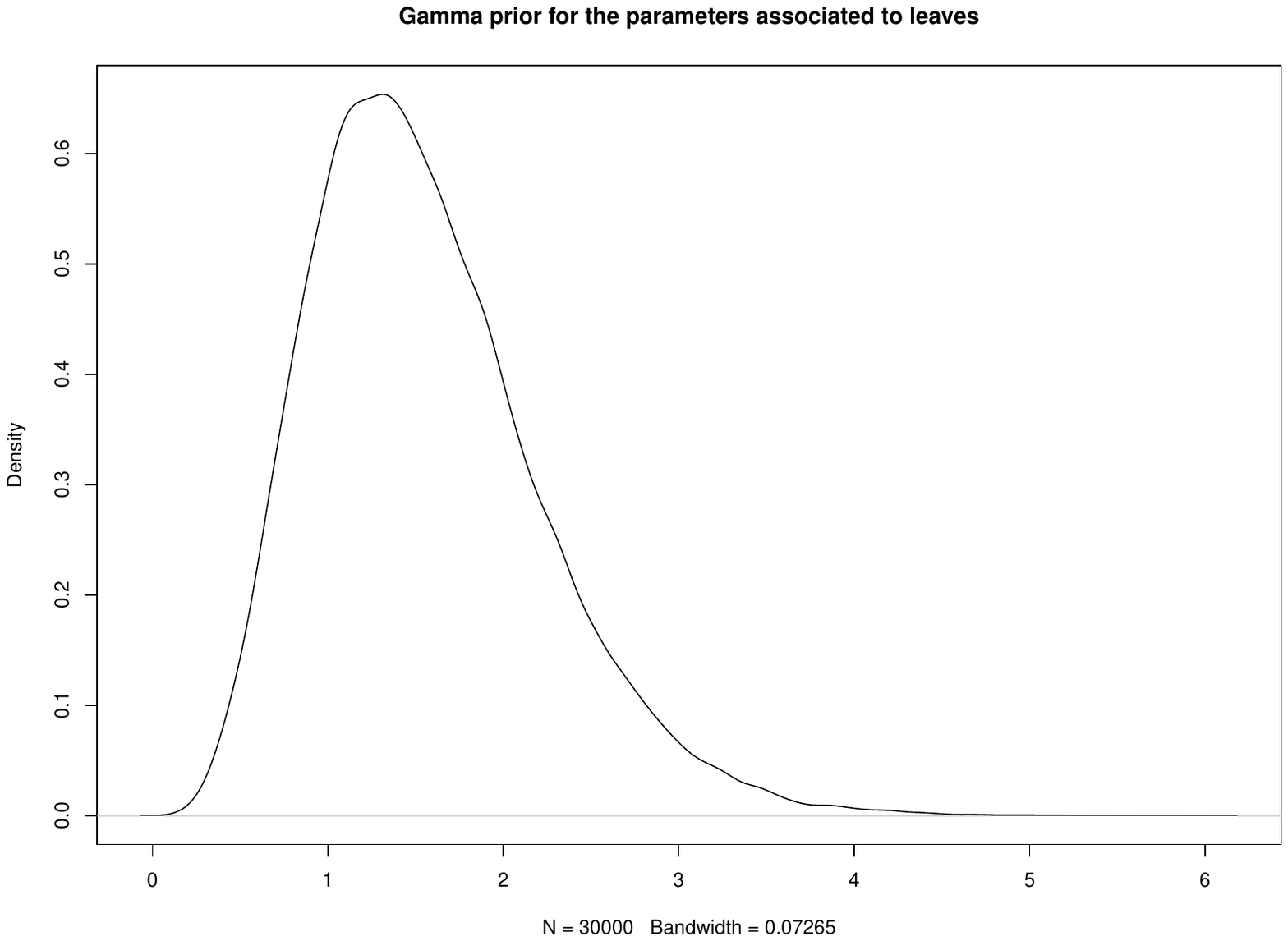}
    \caption{Prior for 10 Trees}
\end{figure}

\section{Inhomogeneous two-dimensional Poisson Process with Gaussian intensity} \label{SS_AR}

\subsection{8 Trees}
\noindent We run 3 parallel chains  each for 200000 iterations keeping every 100th sample. 
\begin{figure}[H] 
    \centering\includegraphics[width=8cm]{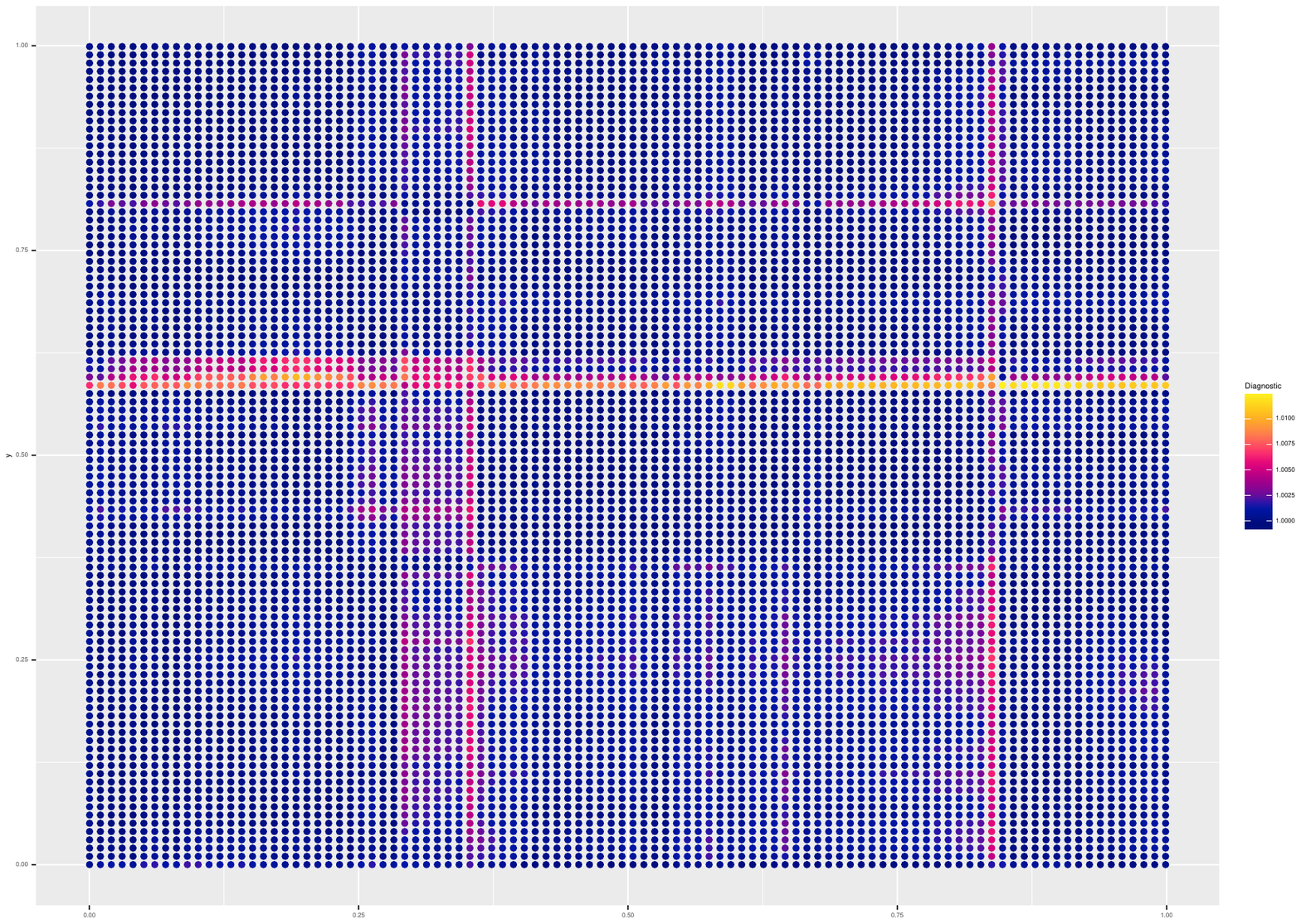}
\caption{The Gelman-Rubin Criterion for 8 Trees}
\label{SL_2DSM1}
\end{figure}

\begin{figure}[H] 
    \centering\includegraphics[width=8cm]{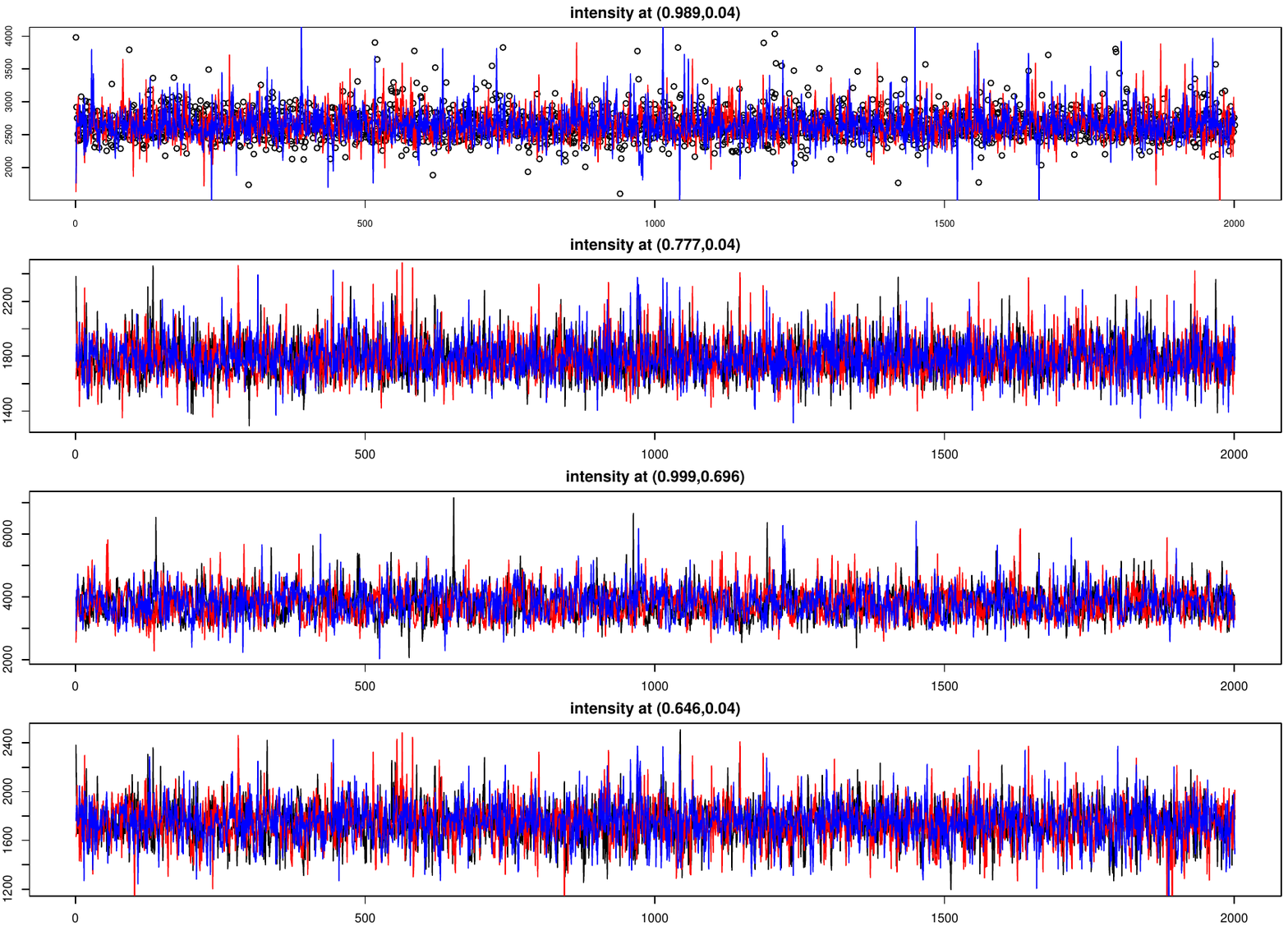}
    \caption{Trace plots for 8 Trees}
\label{SL_2DSM2}
\end{figure}

\begin{figure}[H] 
    \centering\includegraphics[width=8cm]{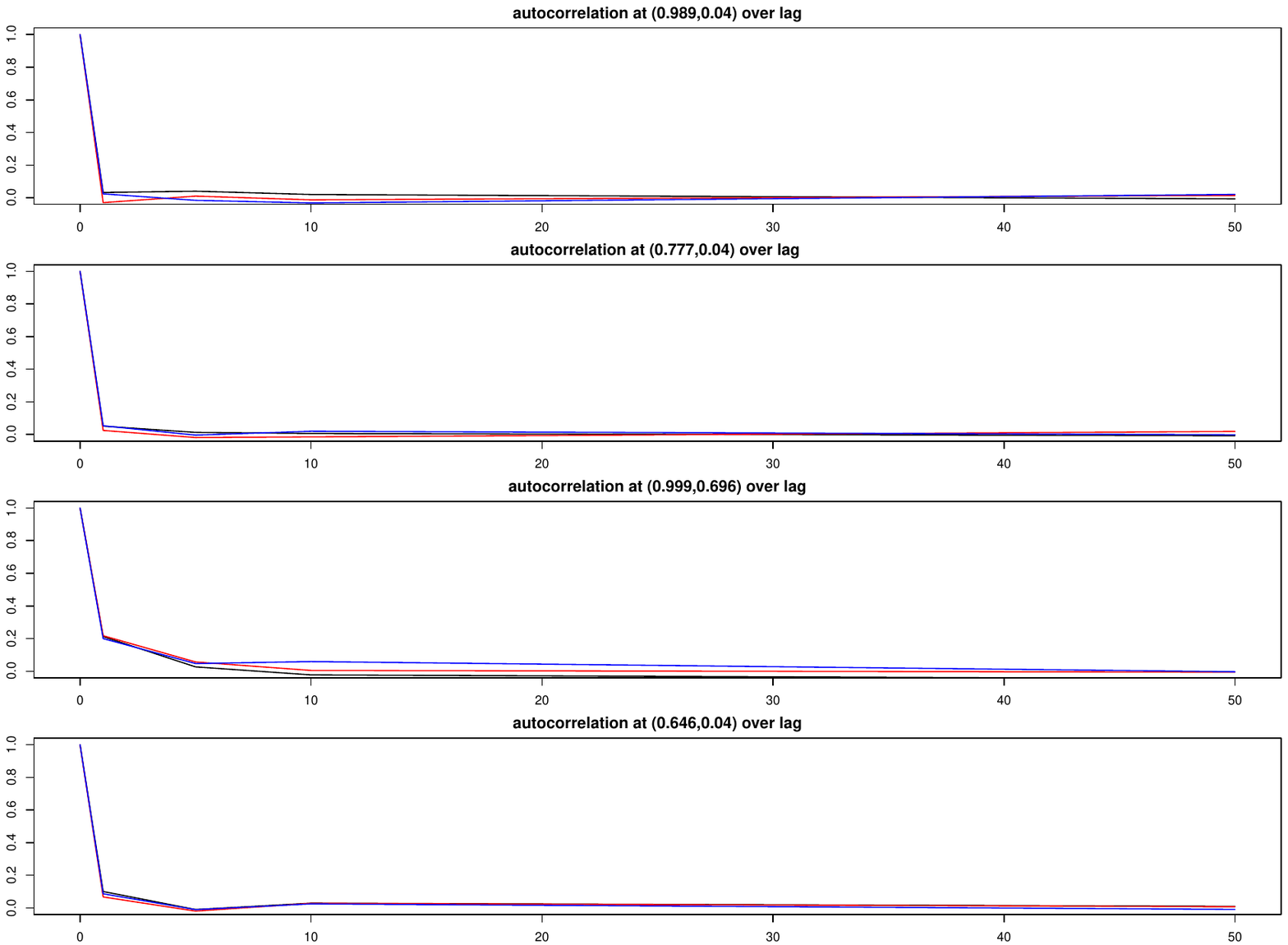}
    \caption{Autocorrelation plots for 8 Trees}
\label{SL_2DSM3}
\end{figure}

\begin{figure}[H] 
    \centering\includegraphics[width=8cm]{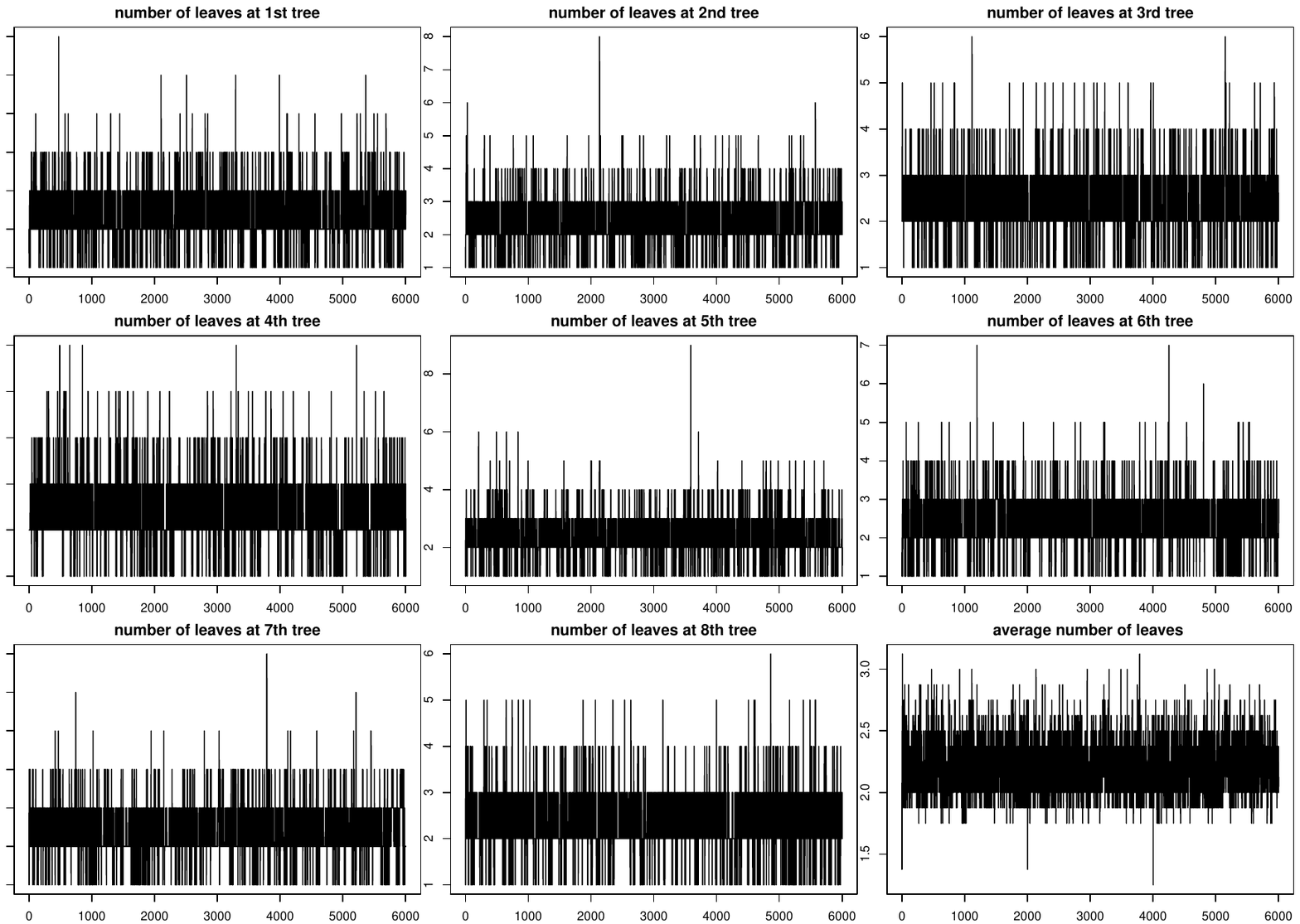}
    \caption{Average number of leaves at trees}
\end{figure}

\begin{figure}[H] 
    \centering\includegraphics[width=8cm]{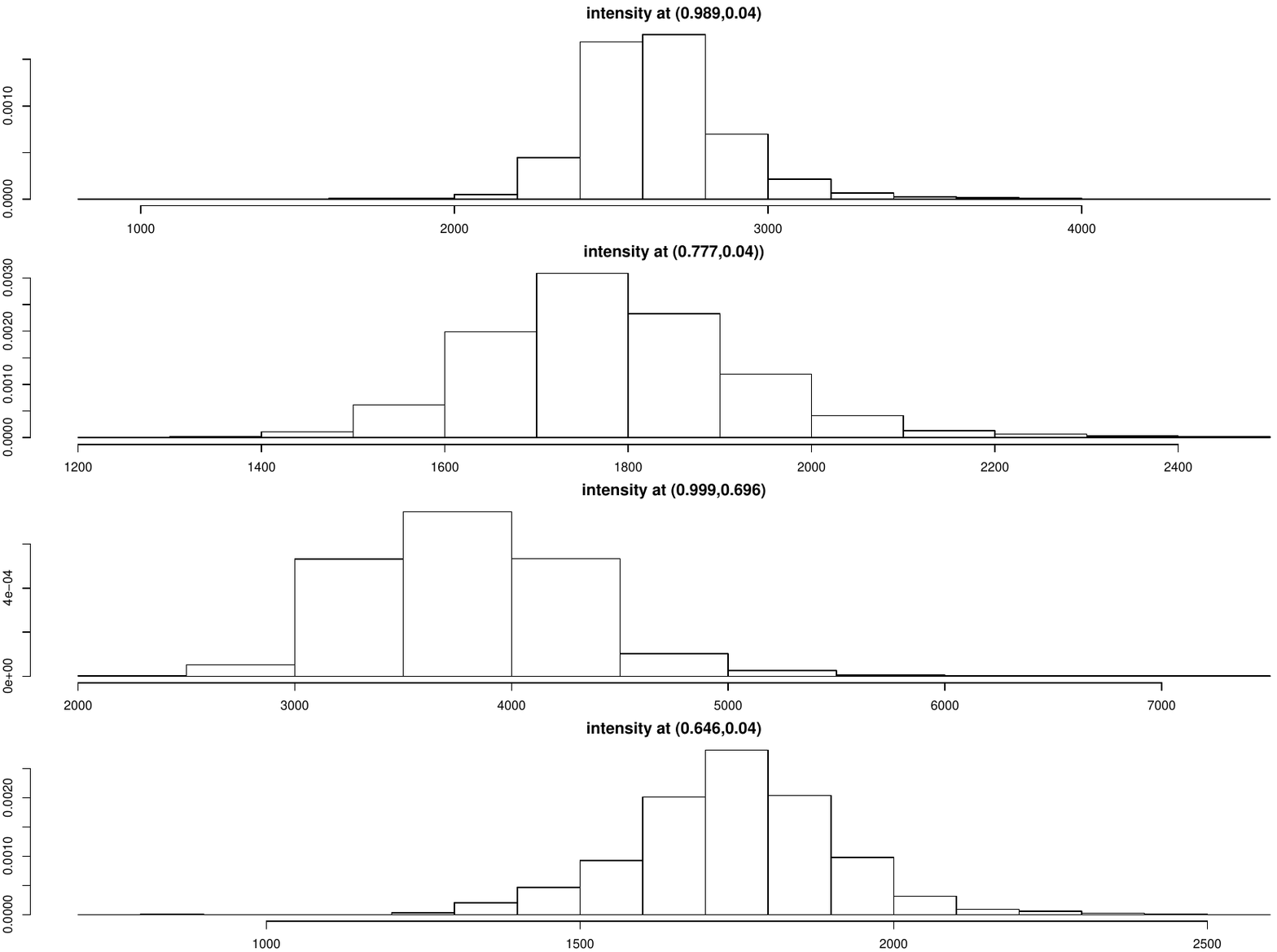}
    \caption{Density of the estimated intensity for 8 Trees}
\end{figure}

\begin{figure}[H] 
    \centering\includegraphics[width=8cm]{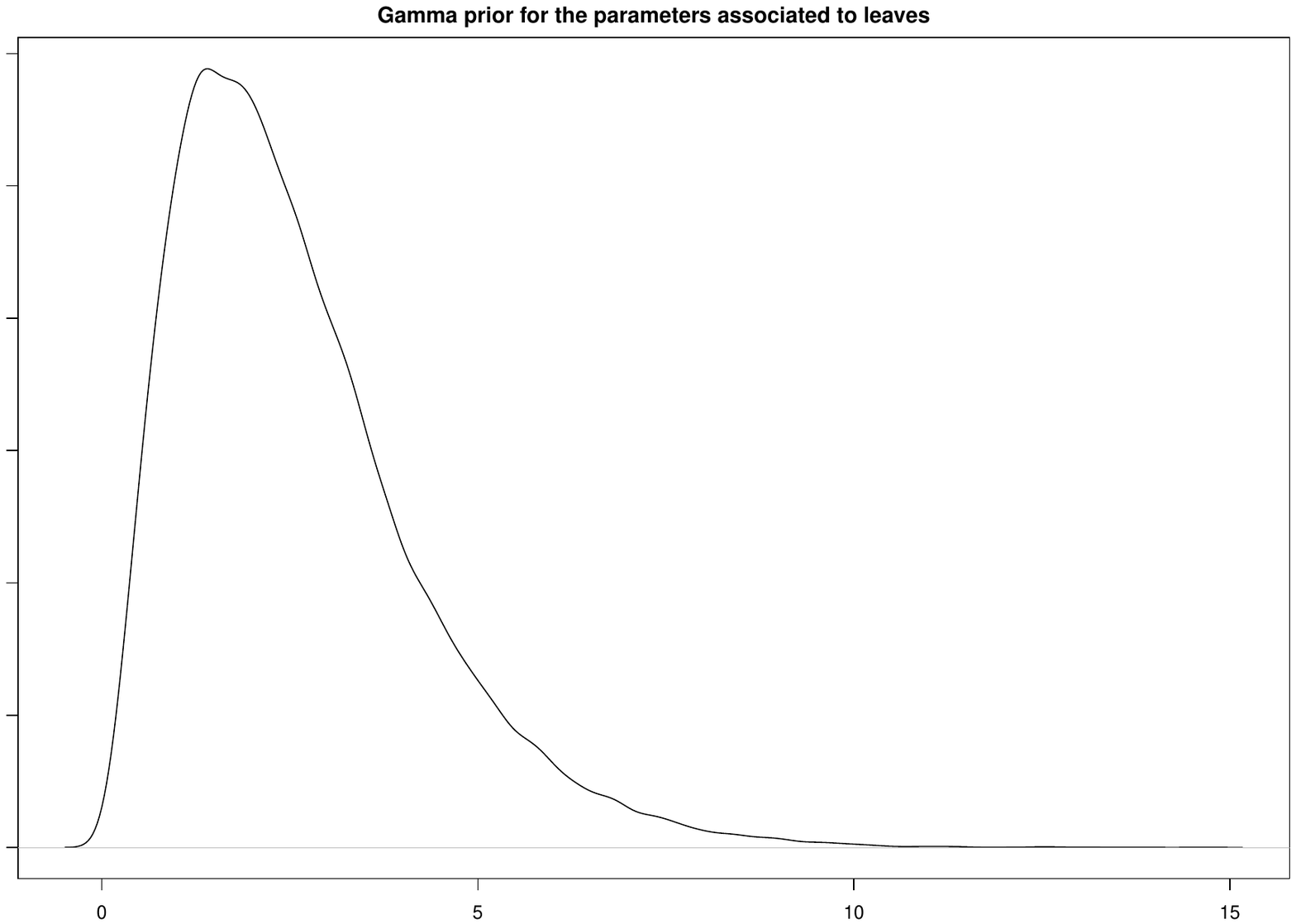}
    \caption{Prior for 8 Trees}
\end{figure}

\subsection{10 Trees}
\noindent We run 3 parallel chains  each for 200000 iterations keeping every 100th sample. 
\begin{figure}[H] 
    \centering\includegraphics[width=8cm]{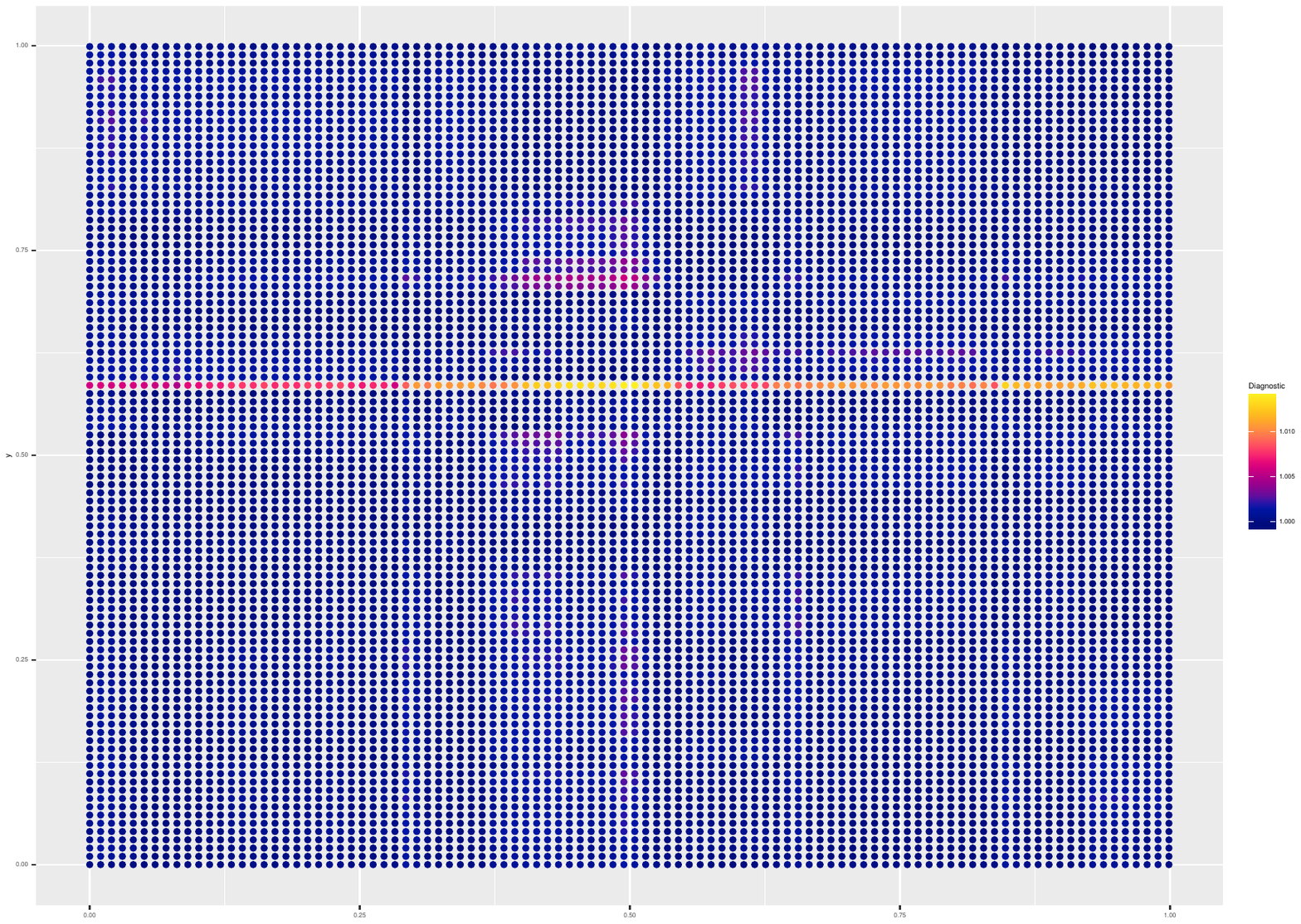}
\caption{The Gelman-Rubin Criterion for 10 Trees}
\label{SL_2DSM4}
\end{figure}

\begin{figure}[H] 
    \centering\includegraphics[width=8cm]{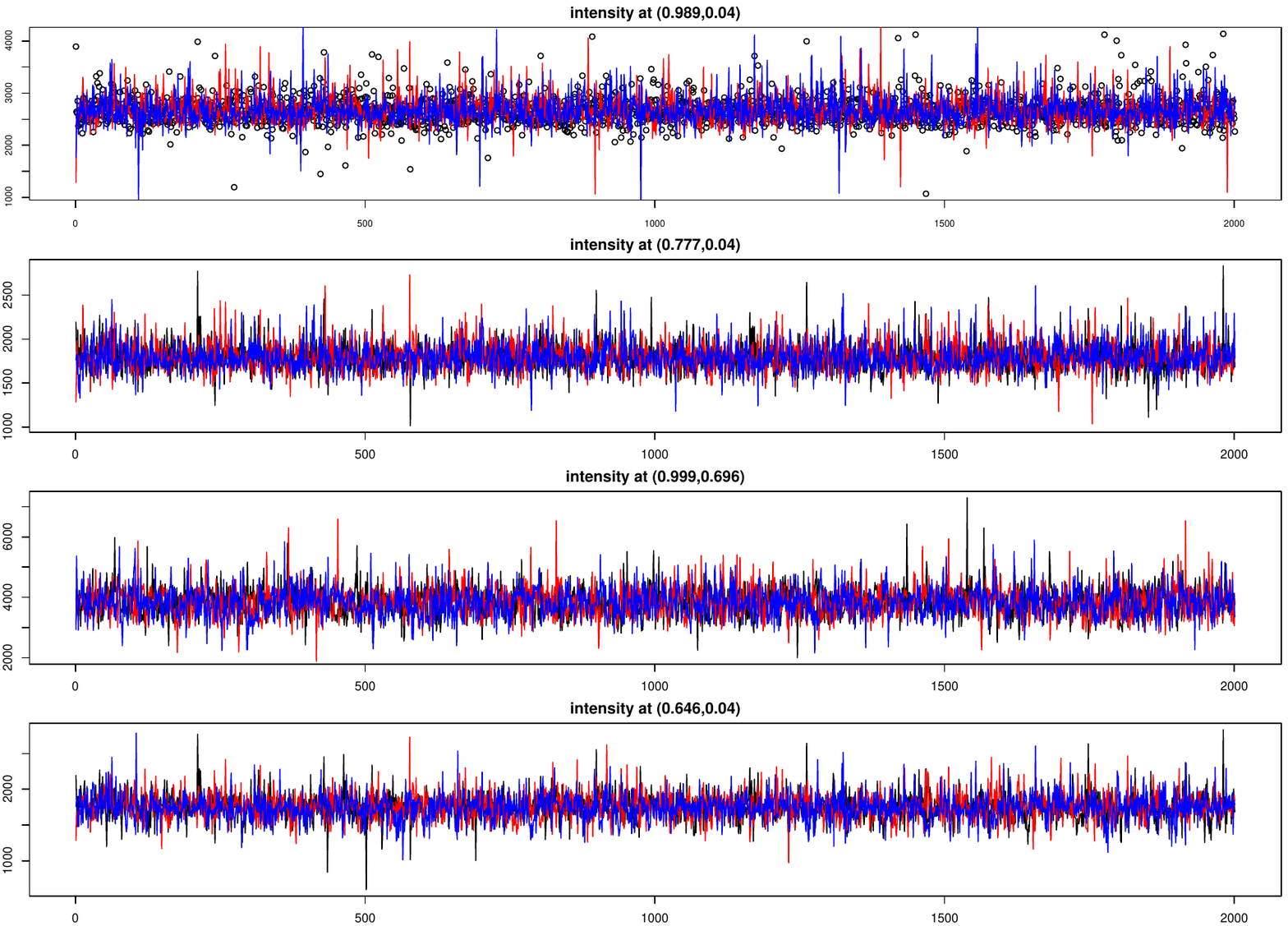}
    \caption{Trace plots for 10 Trees}
    \label{SL_2DSM5}
\end{figure}

\begin{figure}[H] 
    \centering\includegraphics[width=8cm]{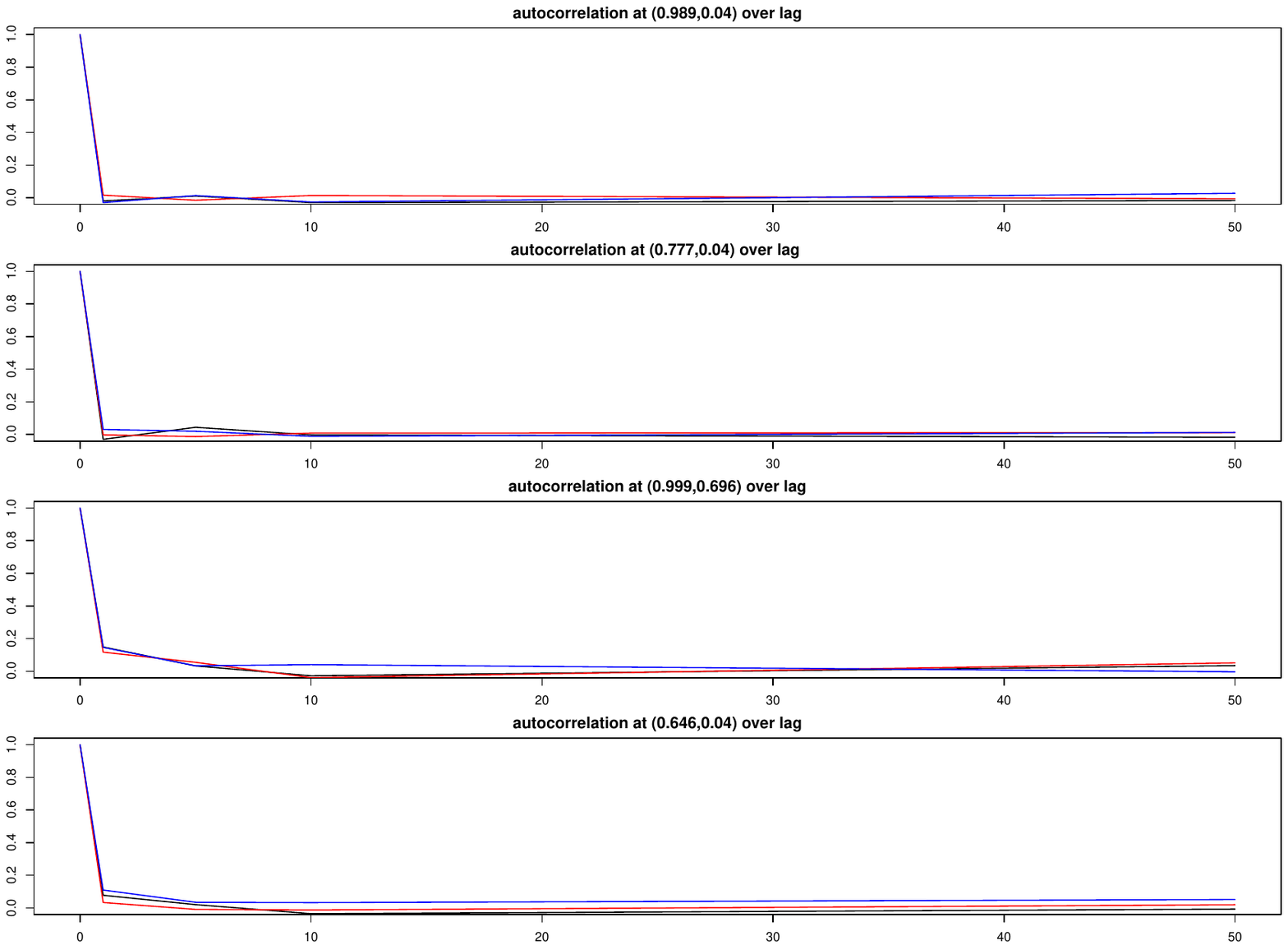}
    \caption{Autocorrelation plots for 10 Trees}
\label{SL_2DSM6}
\end{figure}

\begin{figure}[H] 
    \centering\includegraphics[width=8cm]{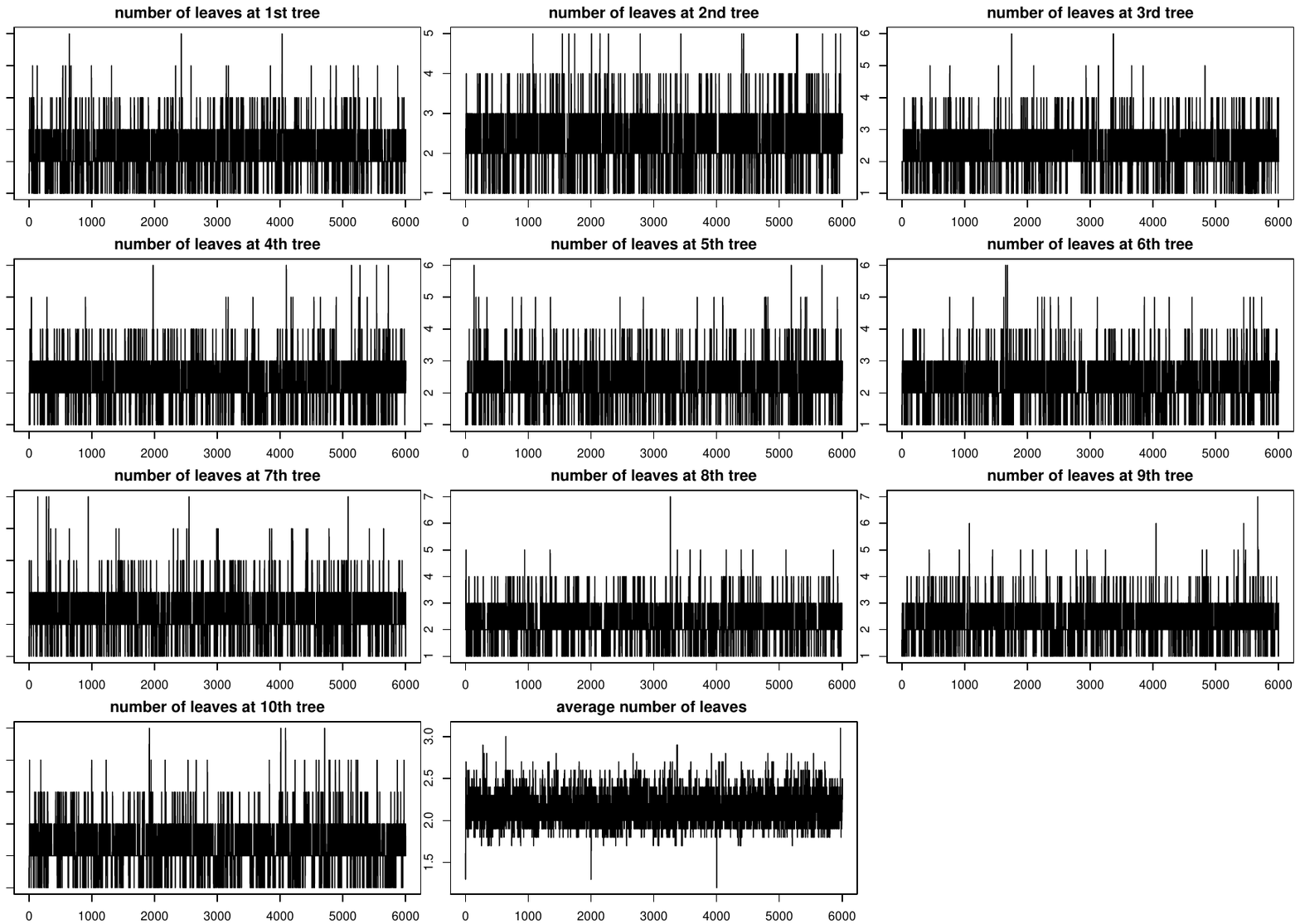}
    \caption{Average number of leaves at trees}
\end{figure}

\begin{figure}[H] 
    \centering\includegraphics[width=8cm]{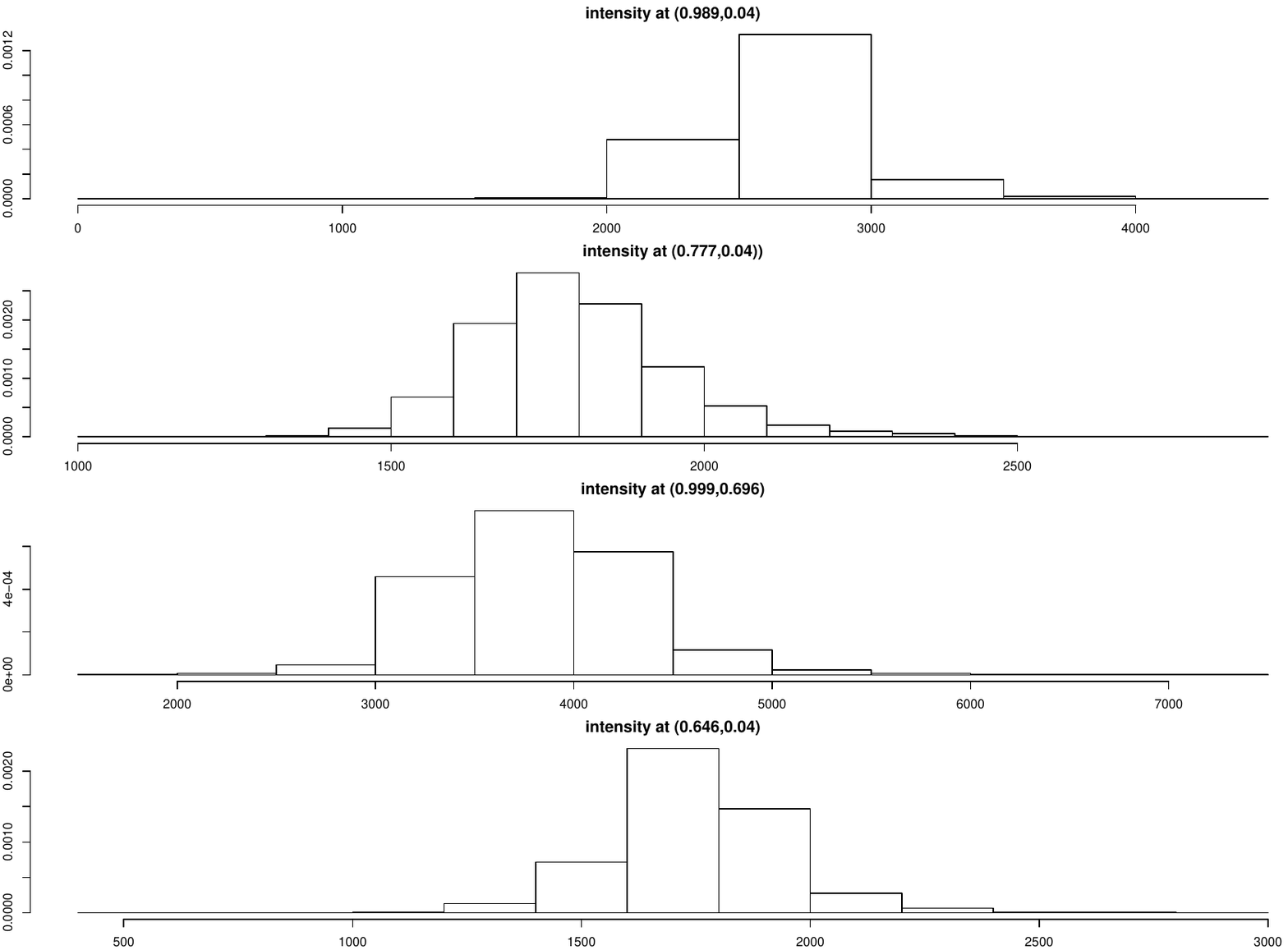}
    \caption{Density of the estimated intensity for 10 Trees}
\end{figure}

\begin{figure}[H] 
    \centering\includegraphics[width=8cm]{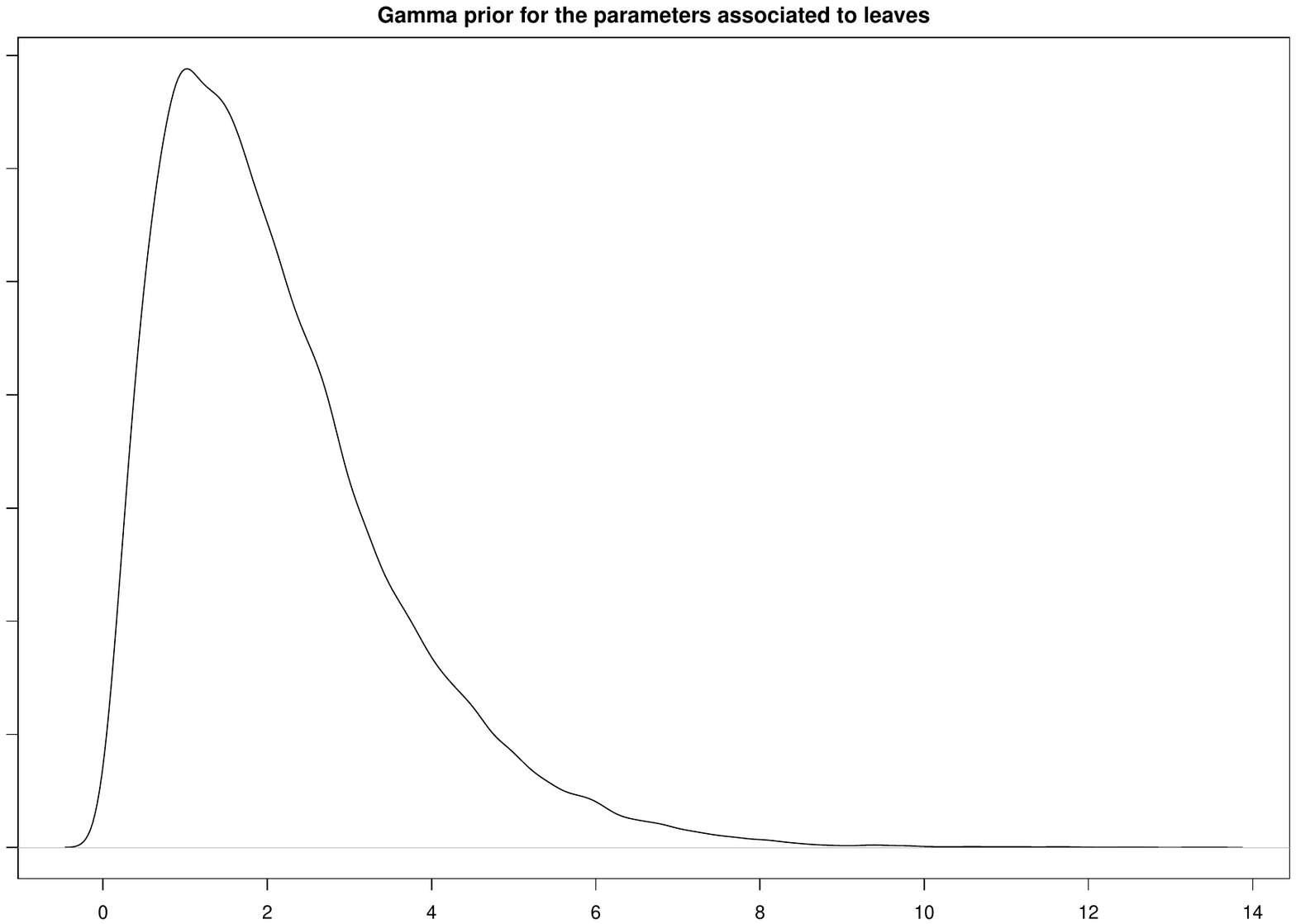}
    \caption{Prior for 10 Trees}
\end{figure}

\section{Two dimensional Poisson process with stepwise intensity function}
\subsection{4 Trees}
\noindent We run 3 parallel chains  each for 100000 iterations keeping every 50th sample. 
\begin{figure}[H] 
  \begin{subfigure}{6cm}
    \centering\includegraphics[width=6cm]{./Figures/step5_Rhat_4Trees}
    \caption{Gelman-Rubin Criterion for 4 Trees}
  \end{subfigure}
  \begin{subfigure}{6cm}
    \centering\includegraphics[width=6cm]{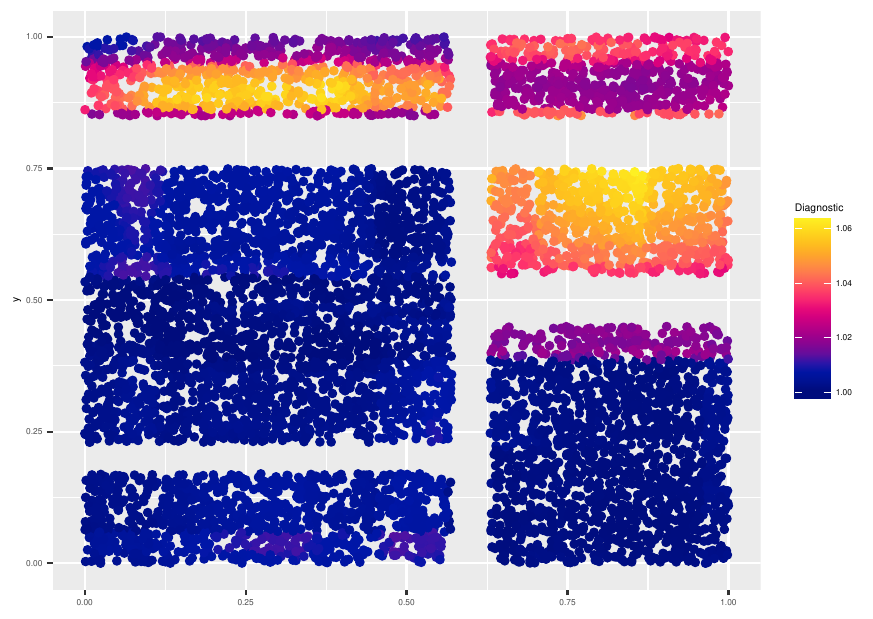}
    \caption{Gelman-Rubin Criterion for 4 Trees with removal of points close to jumps}
  \end{subfigure}
\caption{Gelman-Rubin Criterion}
\label{SL_2DM1}
\end{figure}

\begin{figure}[H] 
    \centering\includegraphics[width=8cm]{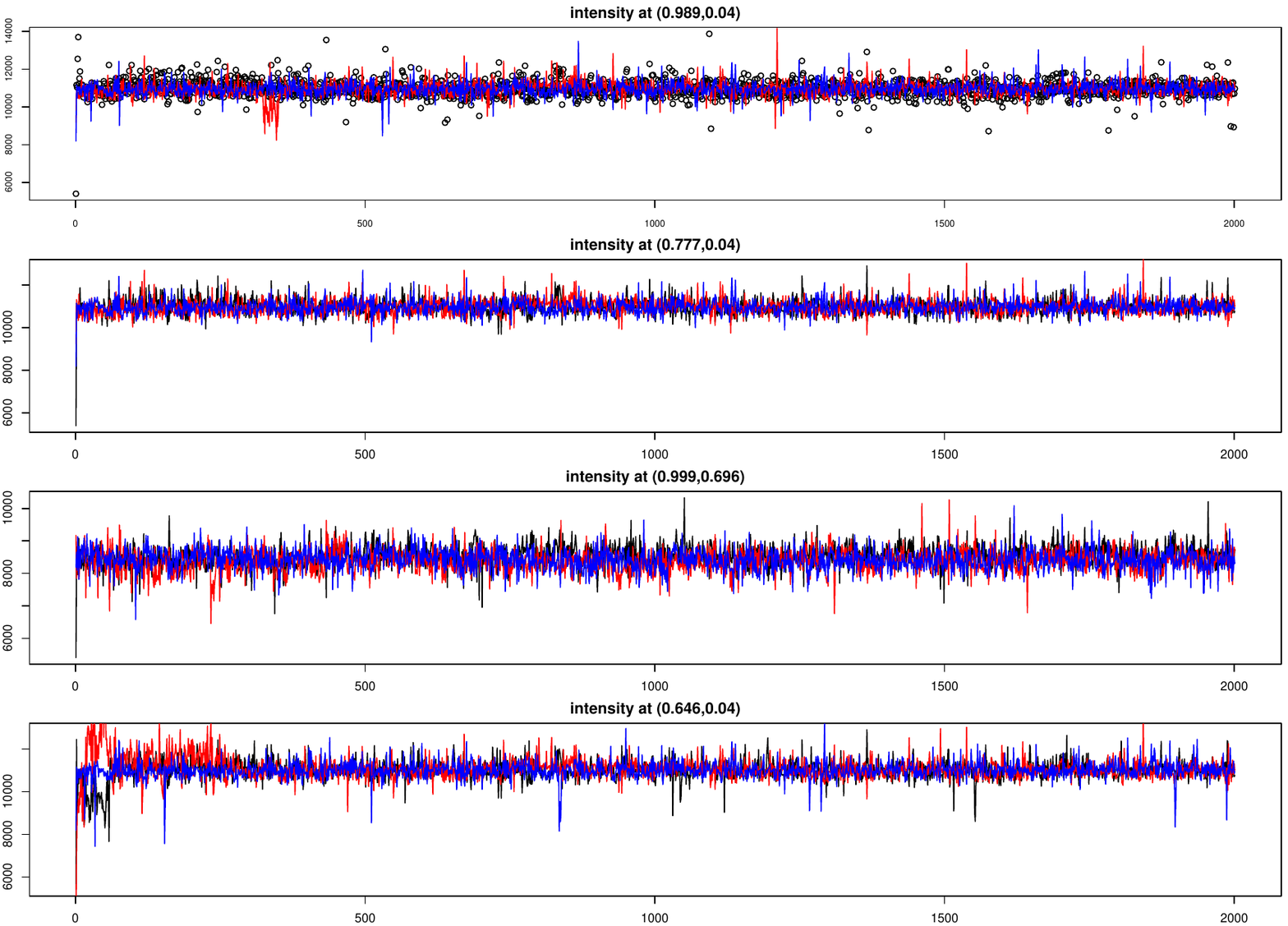}
    \caption{Trace plots for 4 Trees}
\label{SL_2DM2}
\end{figure}

\begin{figure}[H] 
    \centering\includegraphics[width=8cm]{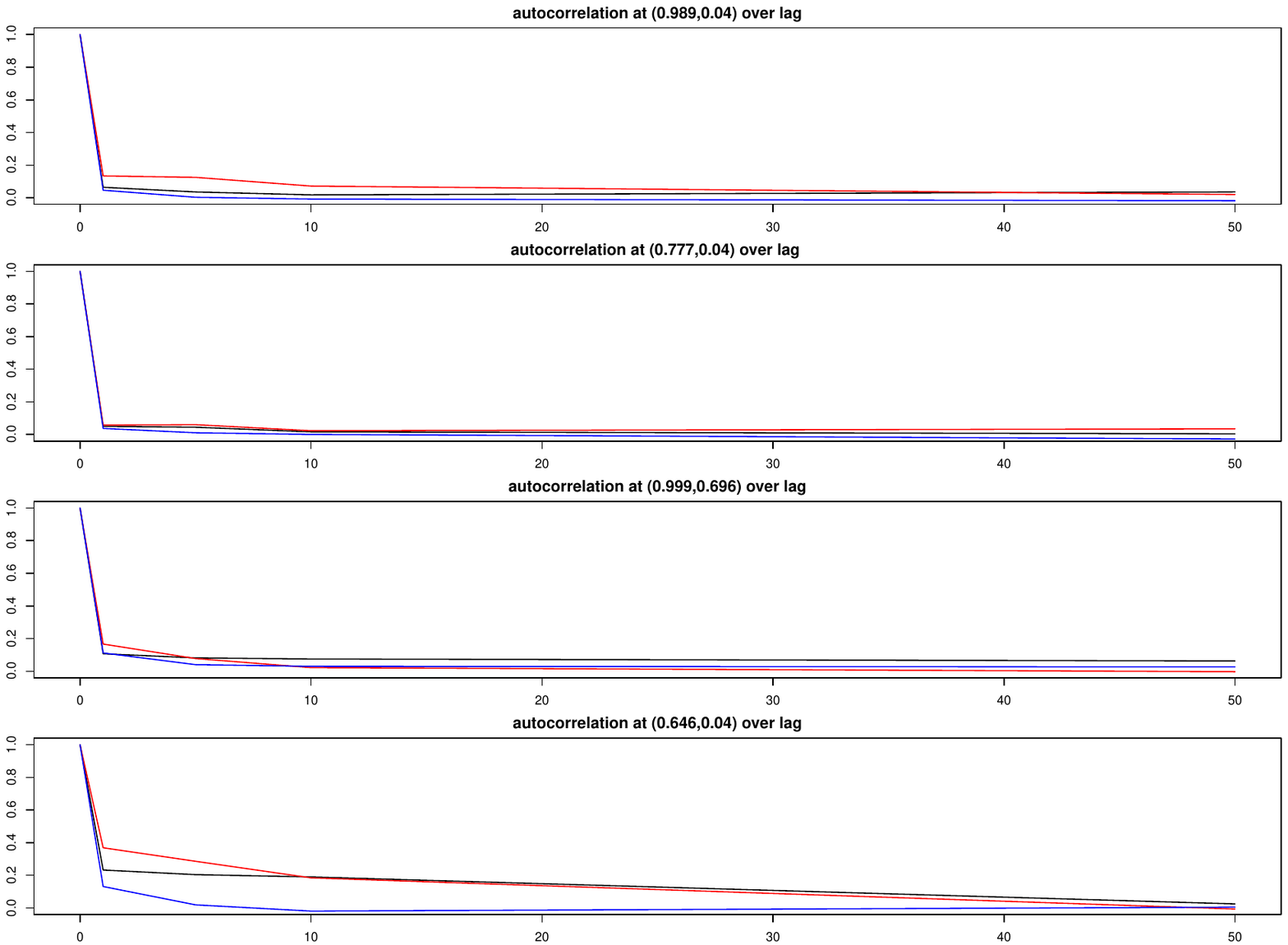}
    \caption{ Autocorrelation plots for 4 Trees}
\label{SL_2DM3}
\end{figure}

\begin{figure}[H] 
    \centering\includegraphics[width=8cm]{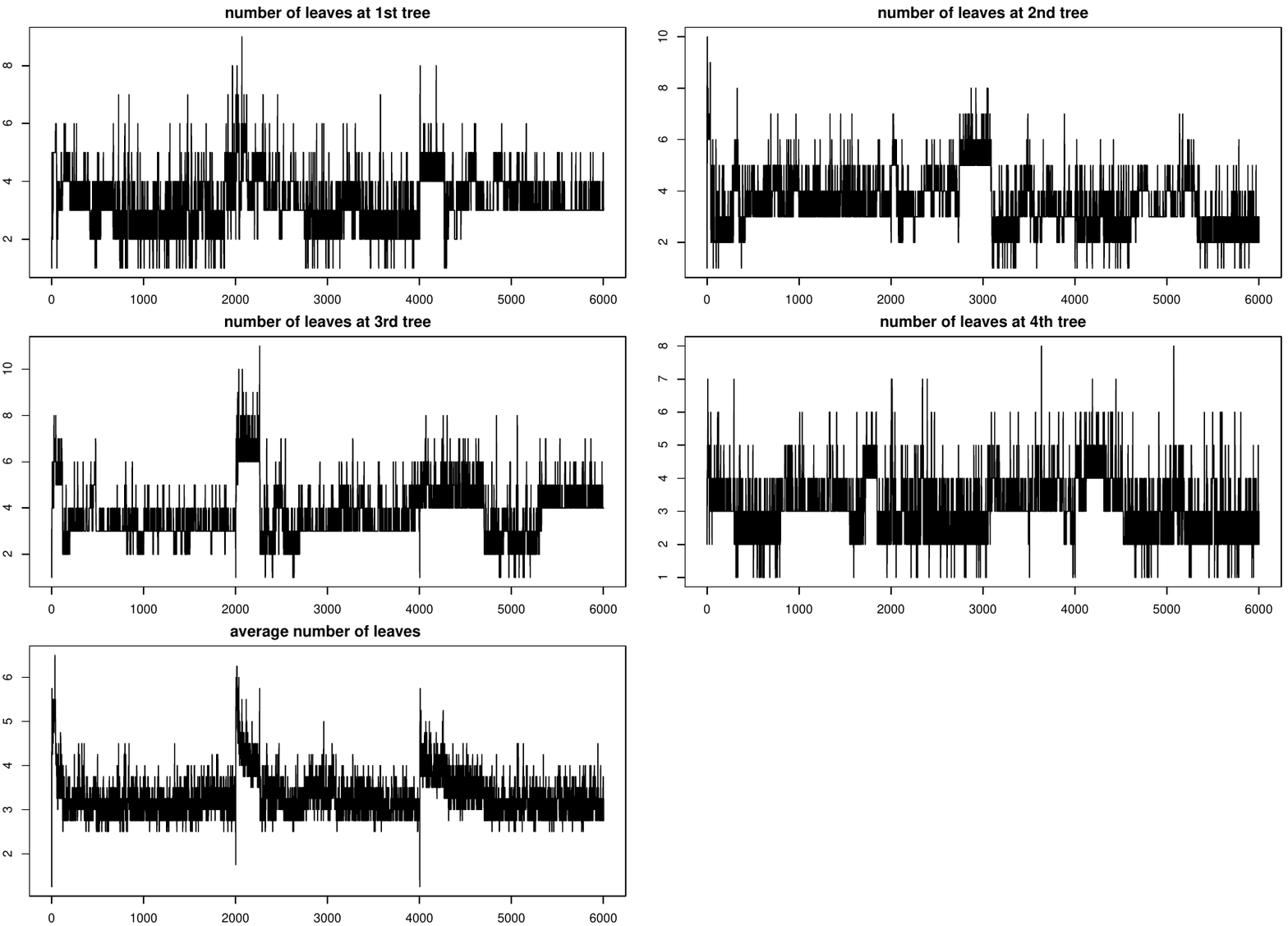}
    \caption{Average number of leaves at trees}
\end{figure}

\begin{figure}[H] 
    \centering\includegraphics[width=8cm]{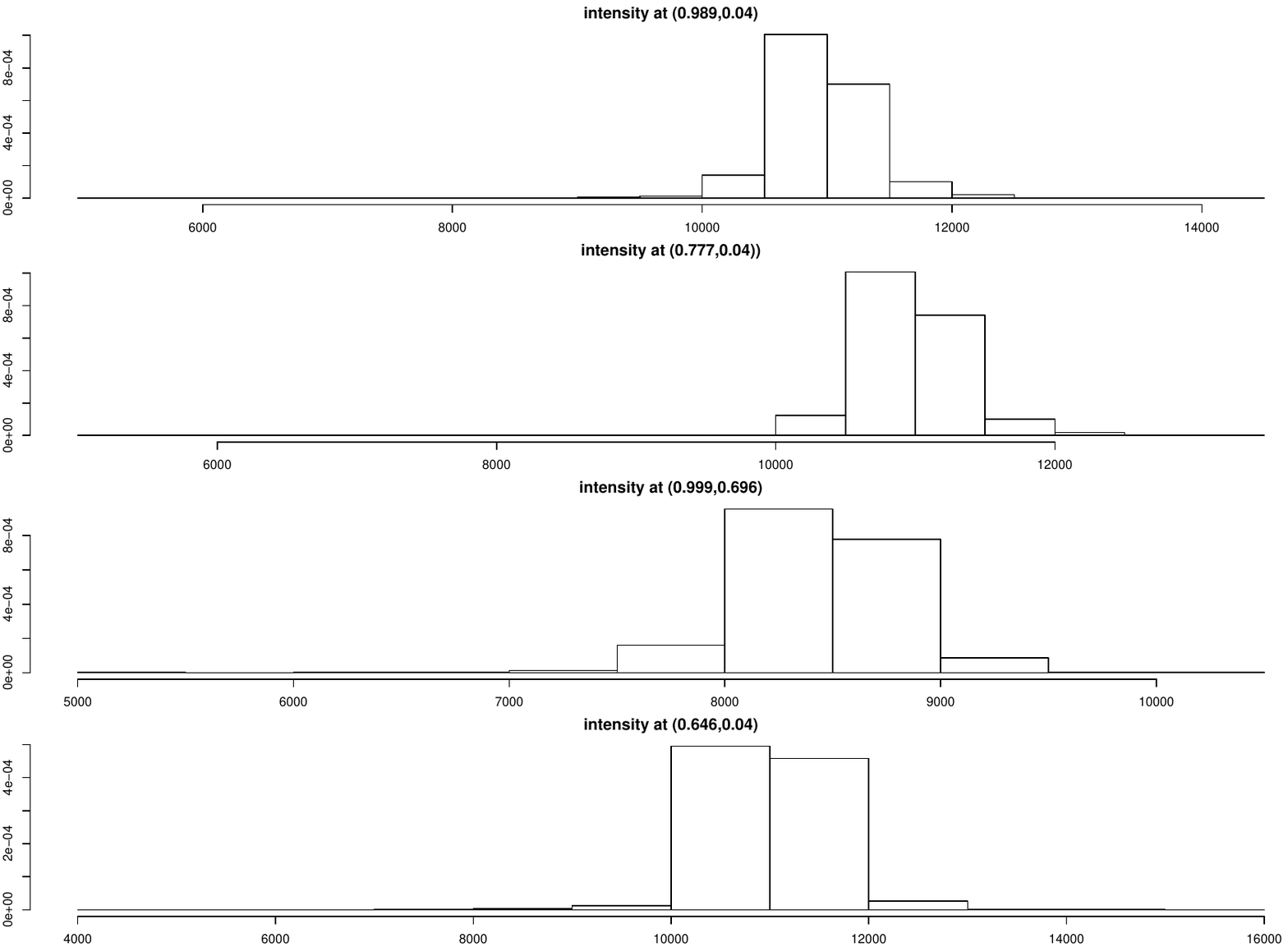}
    \caption{Density of the estimated intensity for 4 Trees}
\end{figure}


\section{Coal Data}

\subsection{8 Trees}
\noindent We run 3 parallel chains  each for 200000 iterations keeping every 100th sample. 
\begin{figure}[H] 
    \centering\includegraphics[width=8cm]{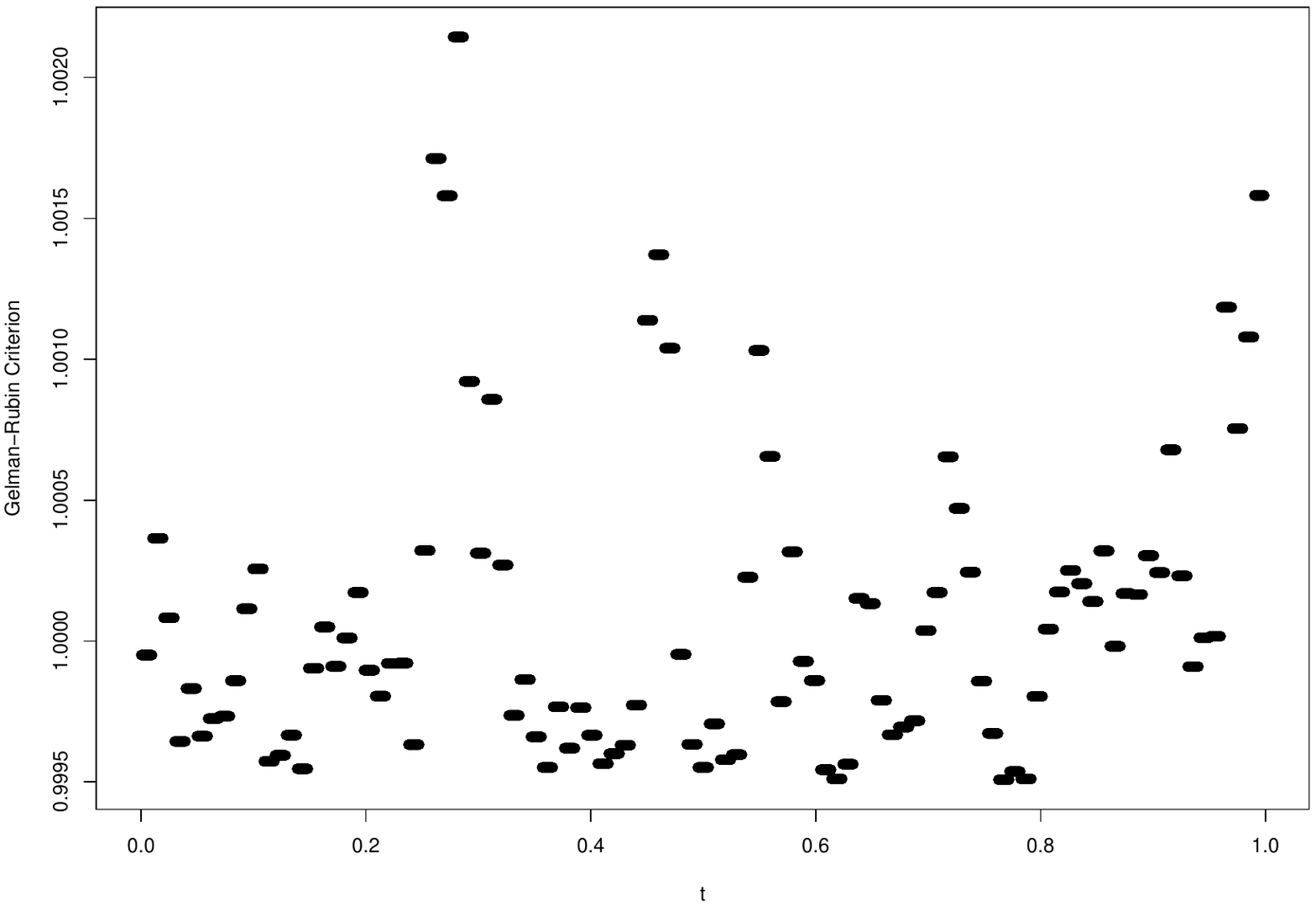}
\caption{The Gelman-Rubin Criterion for 8 Trees}
\end{figure}

\begin{figure}[H] 
    \centering\includegraphics[width=8cm]{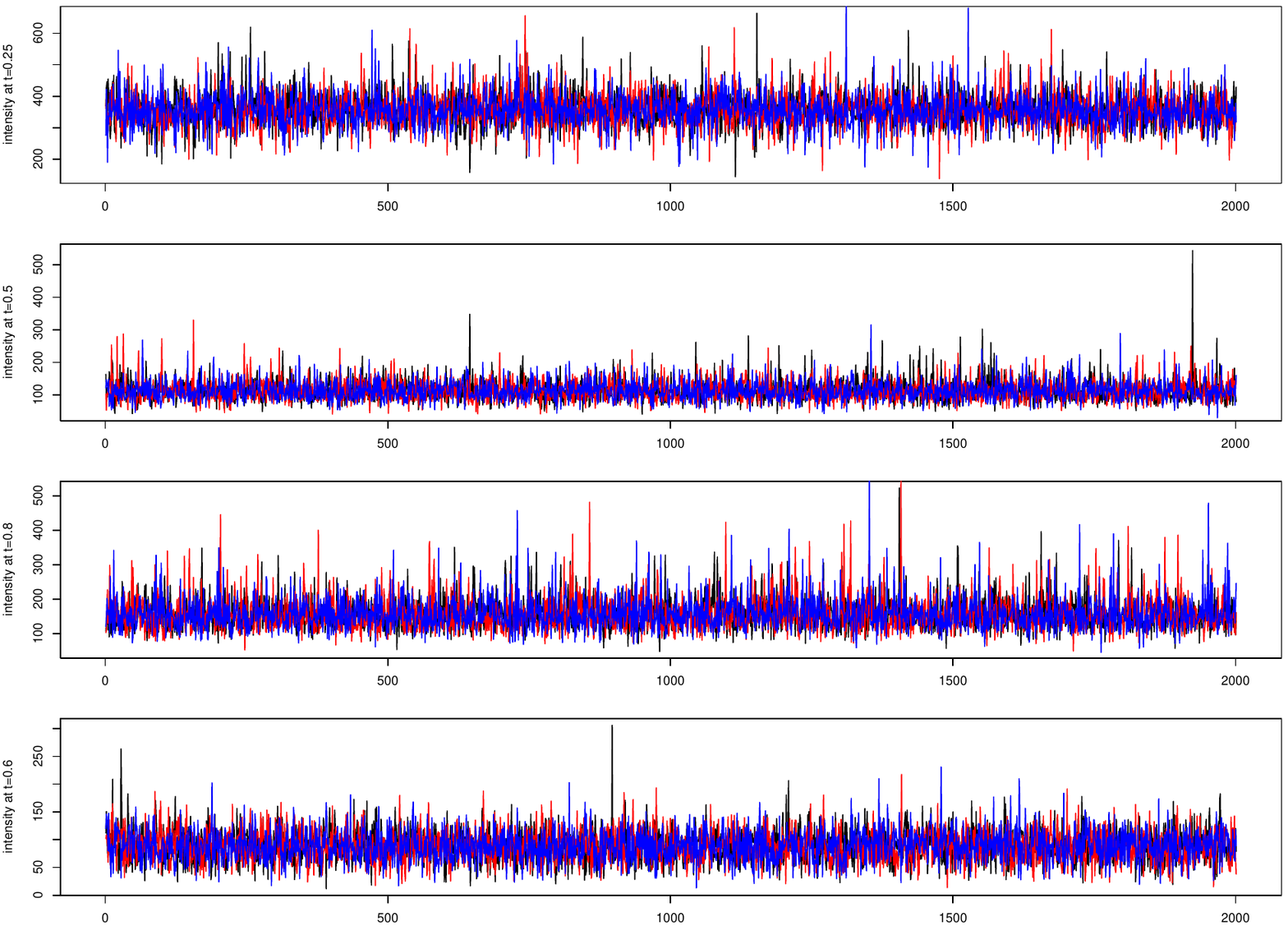}
    \caption{Trace plots for 8 Trees}
\end{figure}

\begin{figure}[H] 
    \centering\includegraphics[width=8cm]{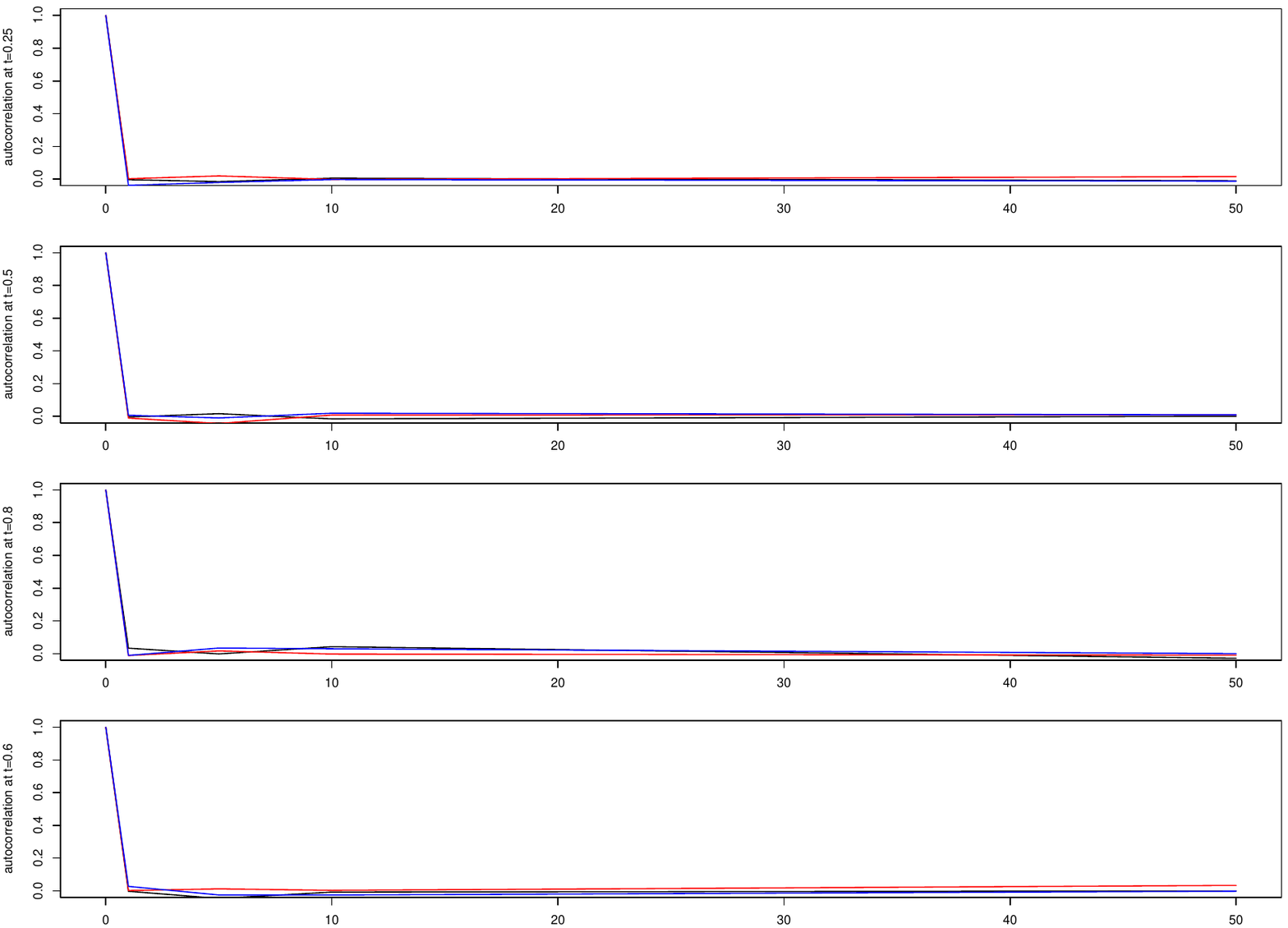}
    \caption{Autocorrelation plots for 8 Trees}
\end{figure}

\begin{figure}[H] 
    \centering\includegraphics[width=8cm]{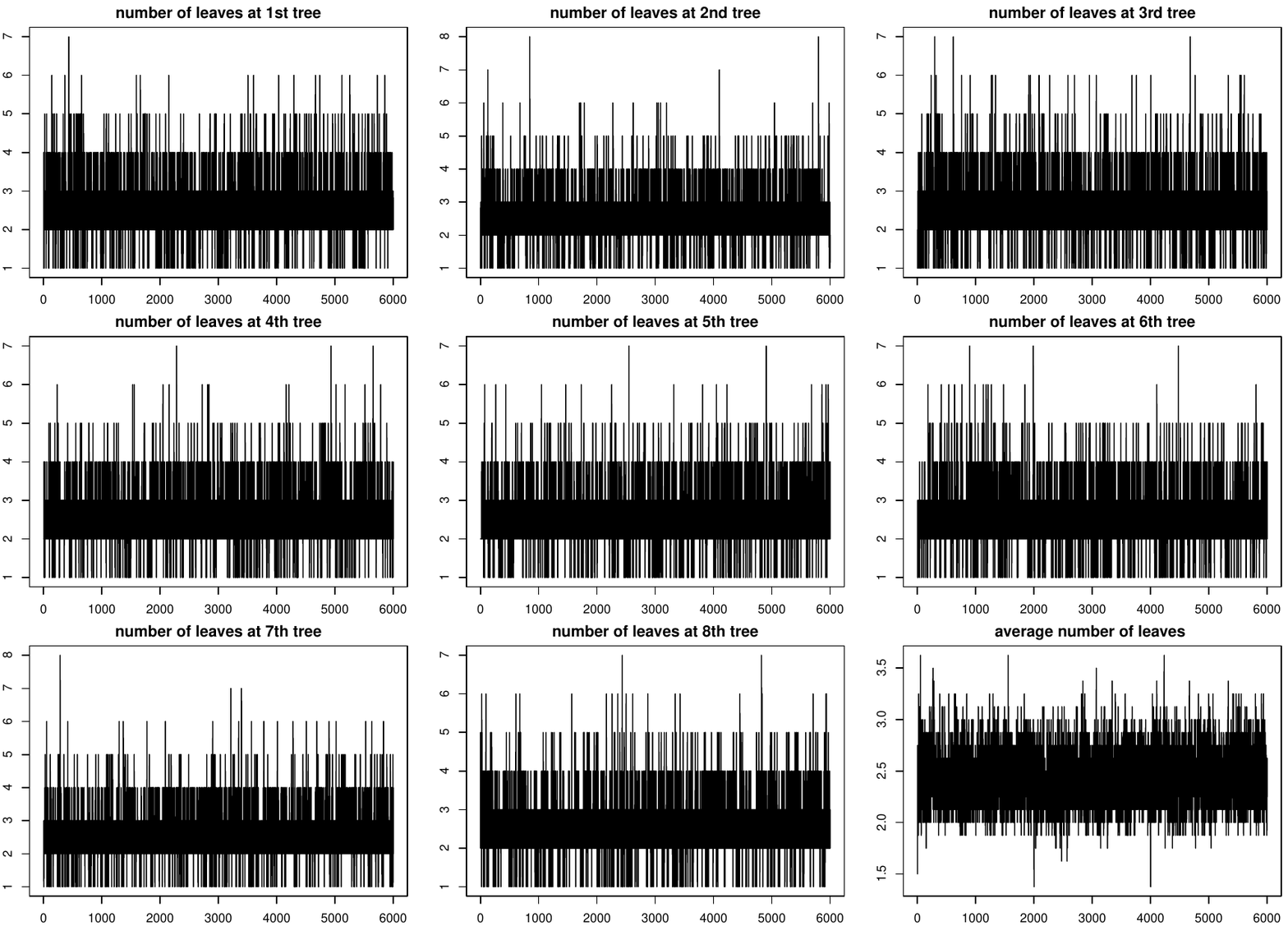}
    \caption{Average number of leaves at trees}
\end{figure}

\begin{figure}[H] 
    \centering\includegraphics[width=8cm]{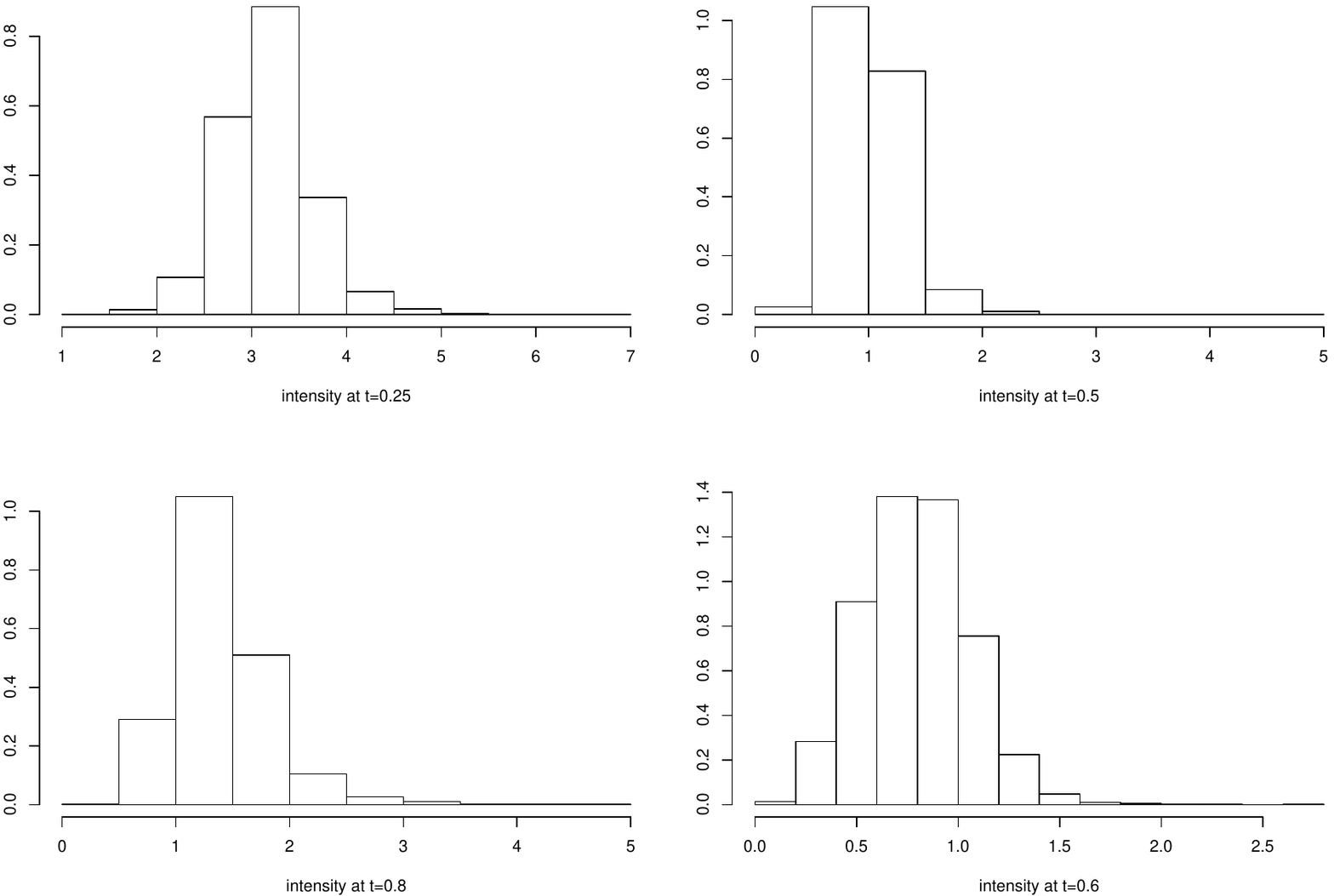}
    \caption{Density of the estimated intensity for 8 Trees}
\end{figure}

\begin{figure}[H] 
    \centering\includegraphics[width=8cm]{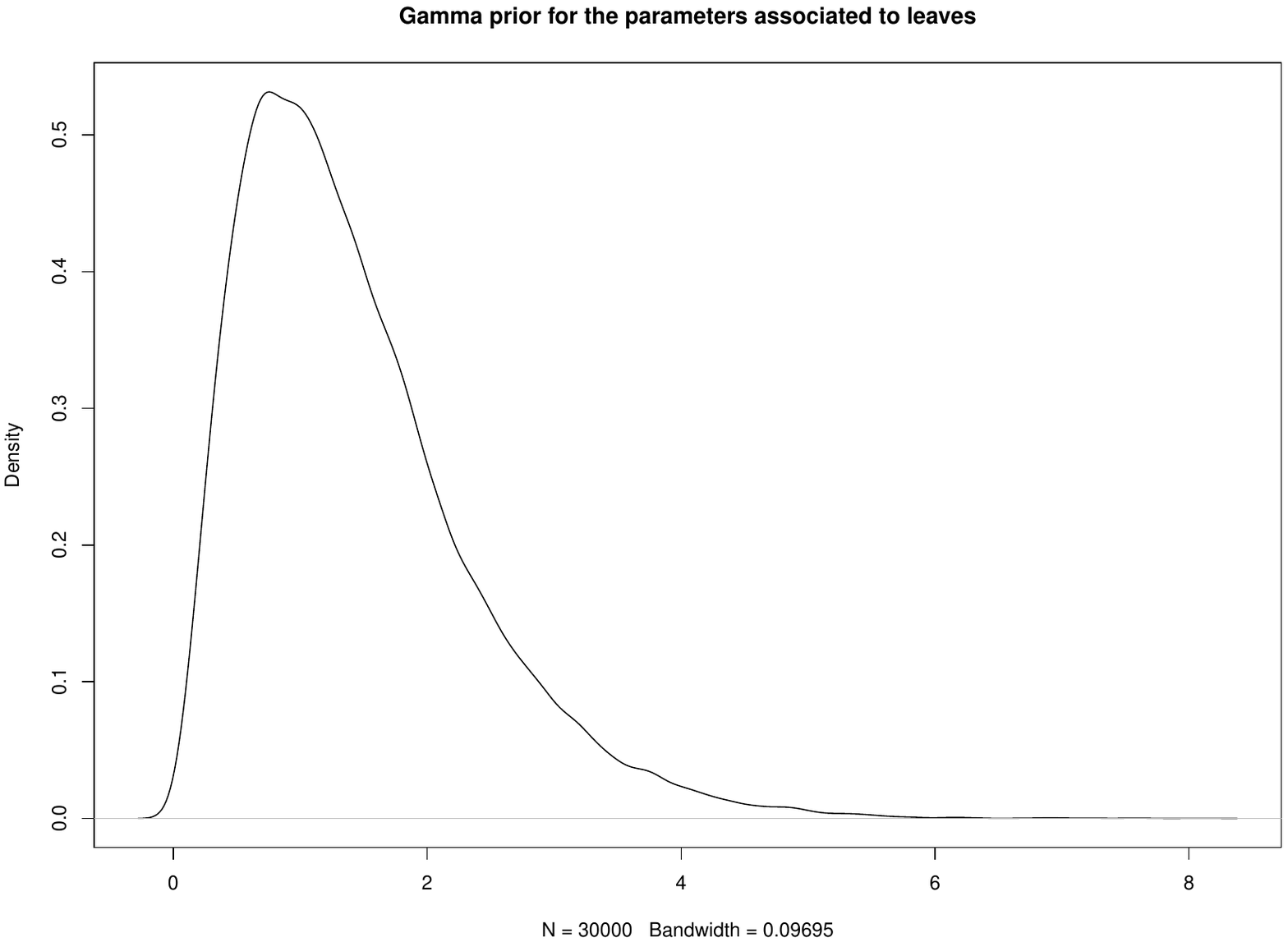}
    \caption{Prior for 8 Trees}
\end{figure}

\subsection{10 Trees}
\noindent We run 3 parallel chains  each for 200000 iterations keeping every 100th sample. 
\begin{figure}[H] 
    \centering\includegraphics[width=8cm]{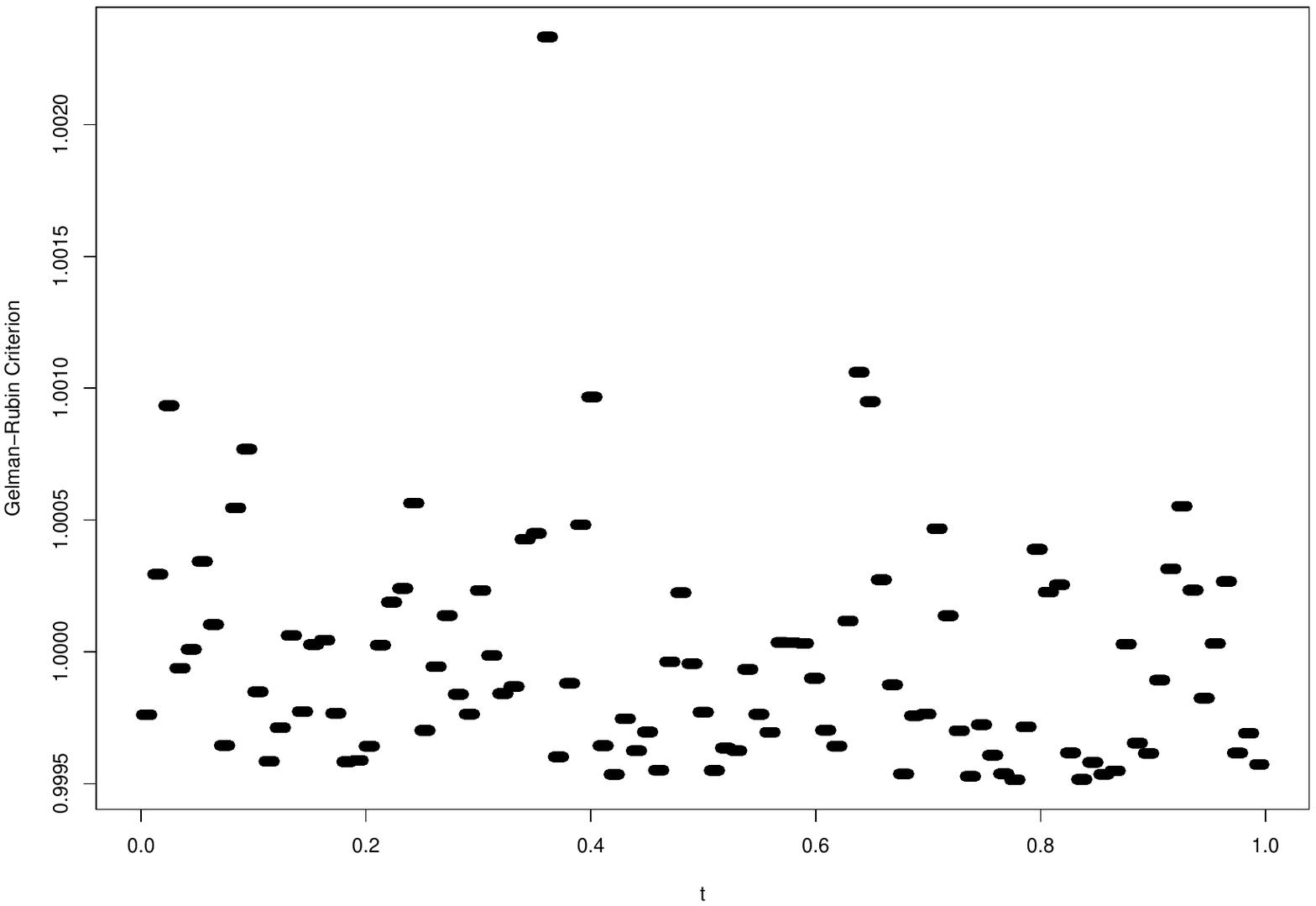}
\caption{The Gelman-Rubin Criterion for 10 Trees}
\end{figure}

\begin{figure}[H] 
    \centering\includegraphics[width=8cm]{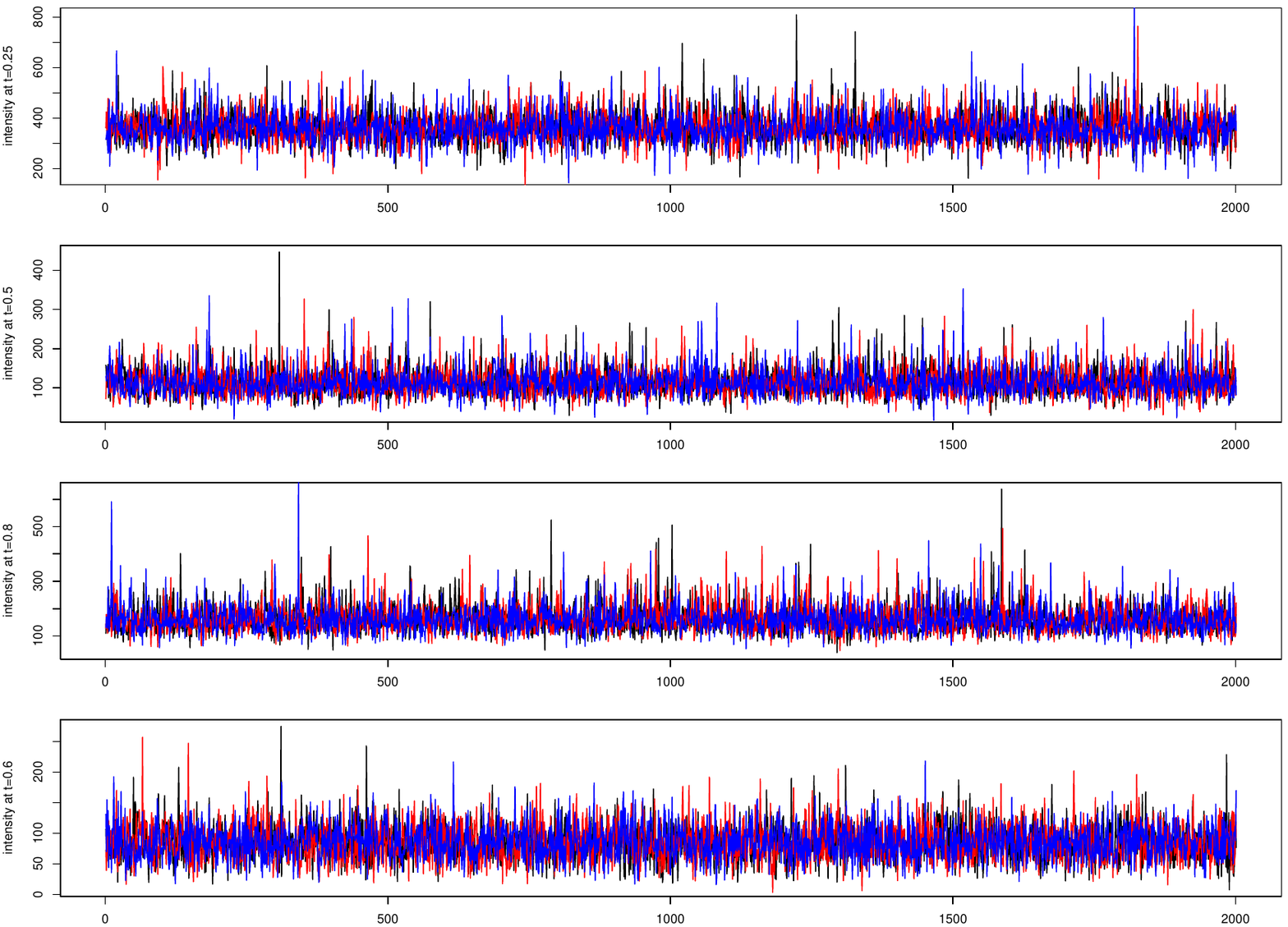}
    \caption{Trace plots for 10 Trees}
\end{figure}

\begin{figure}[H] 
    \centering\includegraphics[width=8cm]{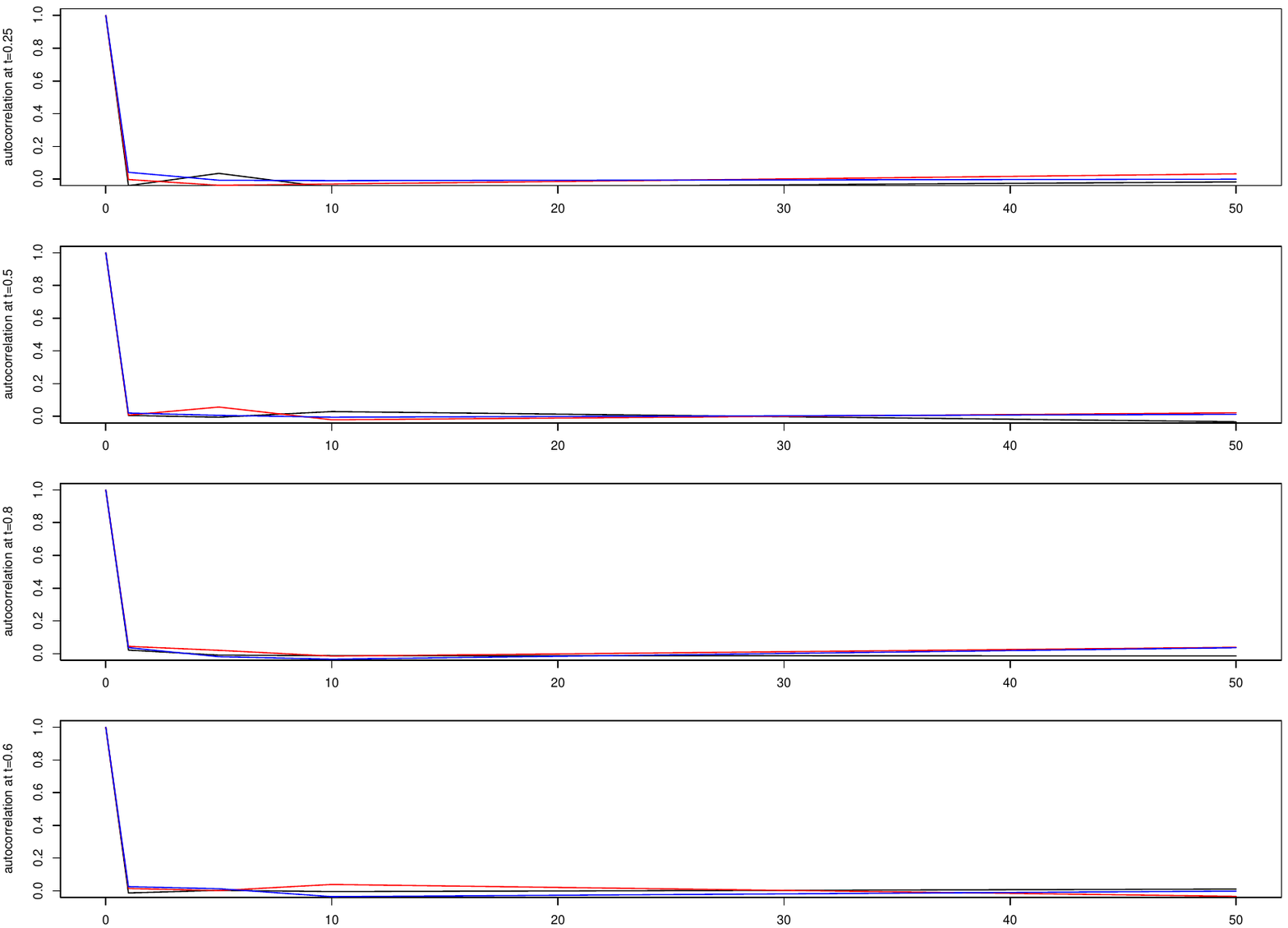}
    \caption{Autocorrelation plots for 10 Trees}
\end{figure}

\begin{figure}[H] 
    \centering\includegraphics[width=8cm]{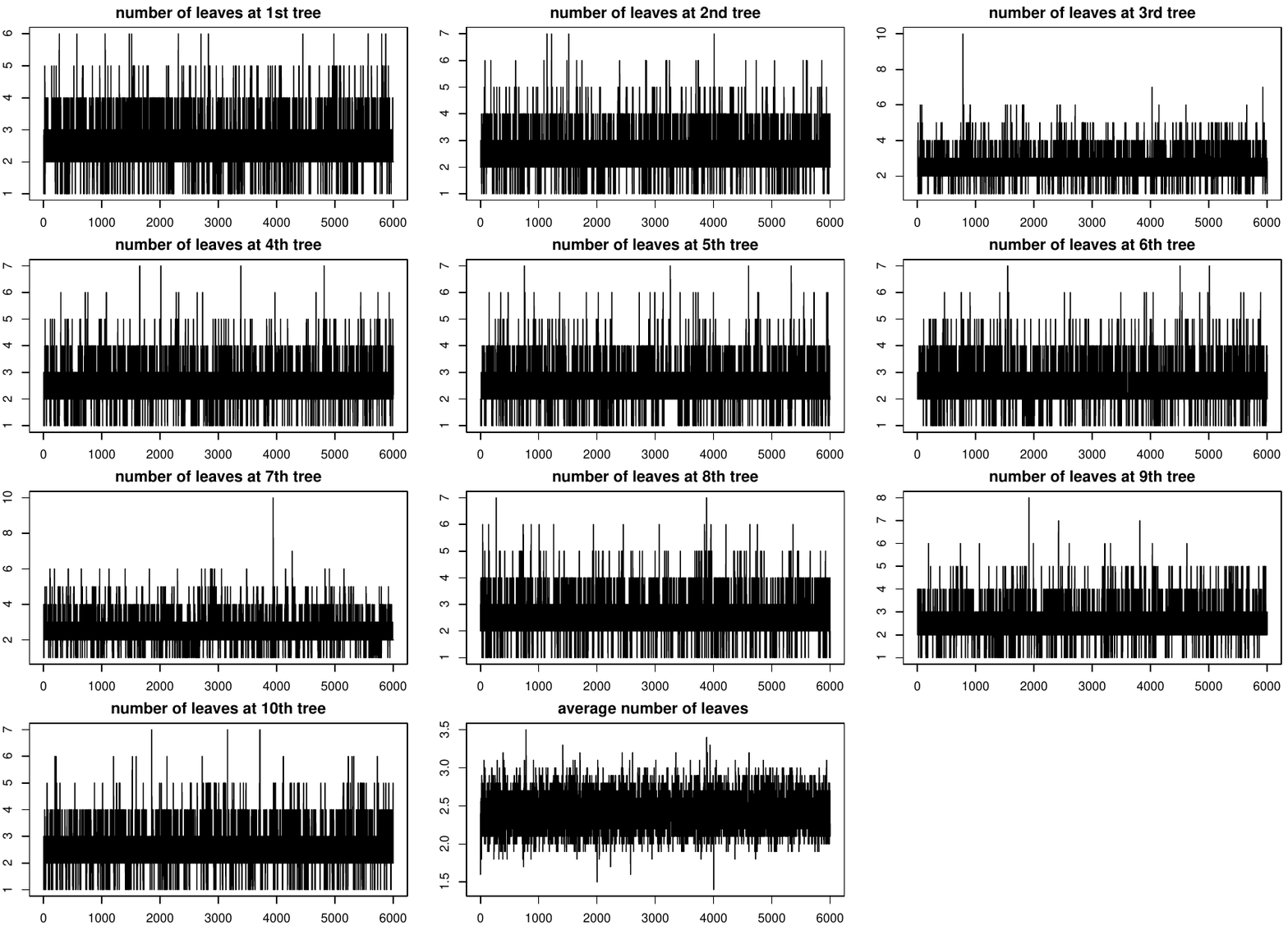}
    \caption{Average number of leaves at trees}
\end{figure}

\begin{figure}[H] 
    \centering\includegraphics[width=8cm]{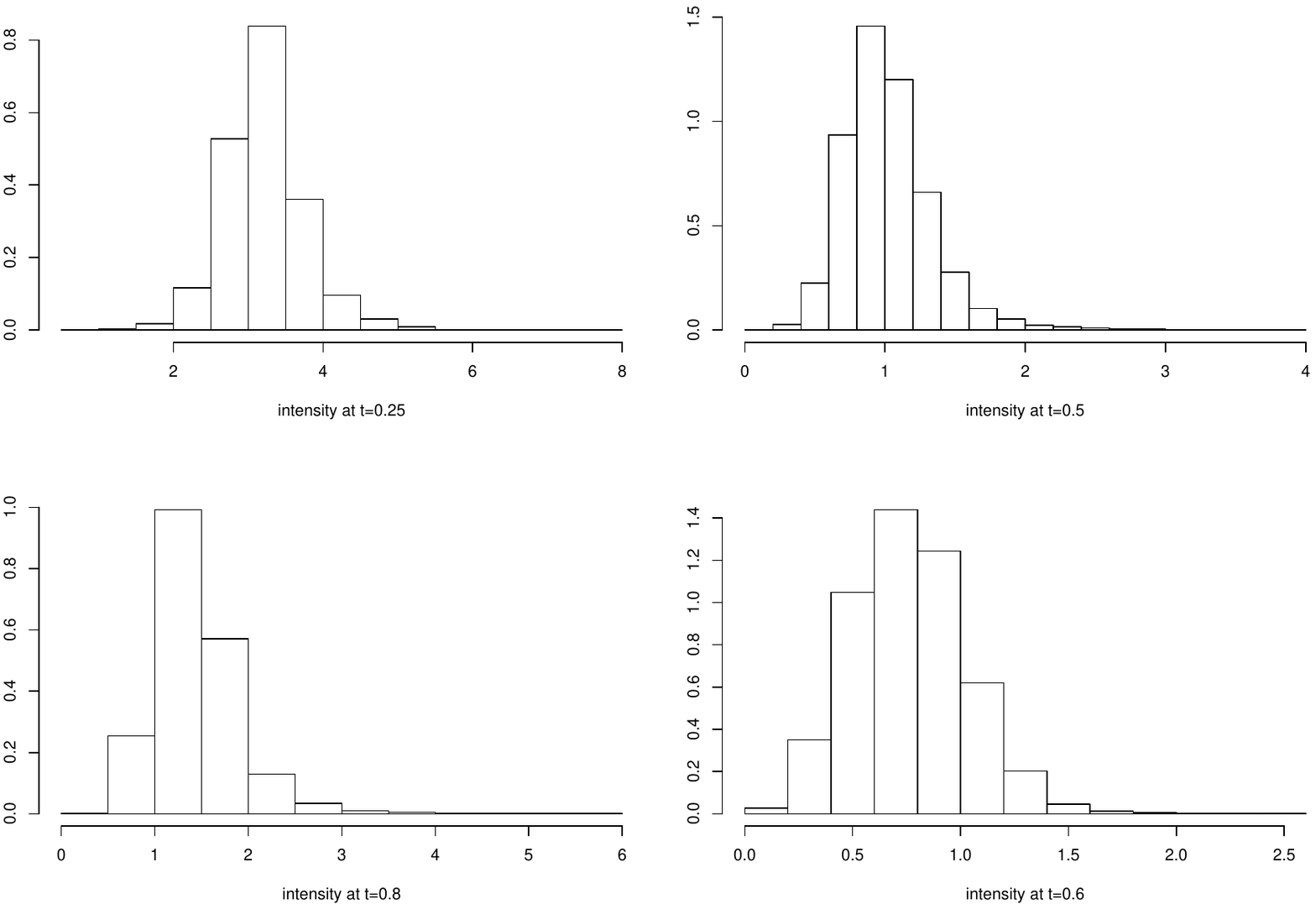}
    \caption{Density of the estimated intensity for 10 Trees}
\end{figure}

\begin{figure}[H] 
    \centering\includegraphics[width=8cm]{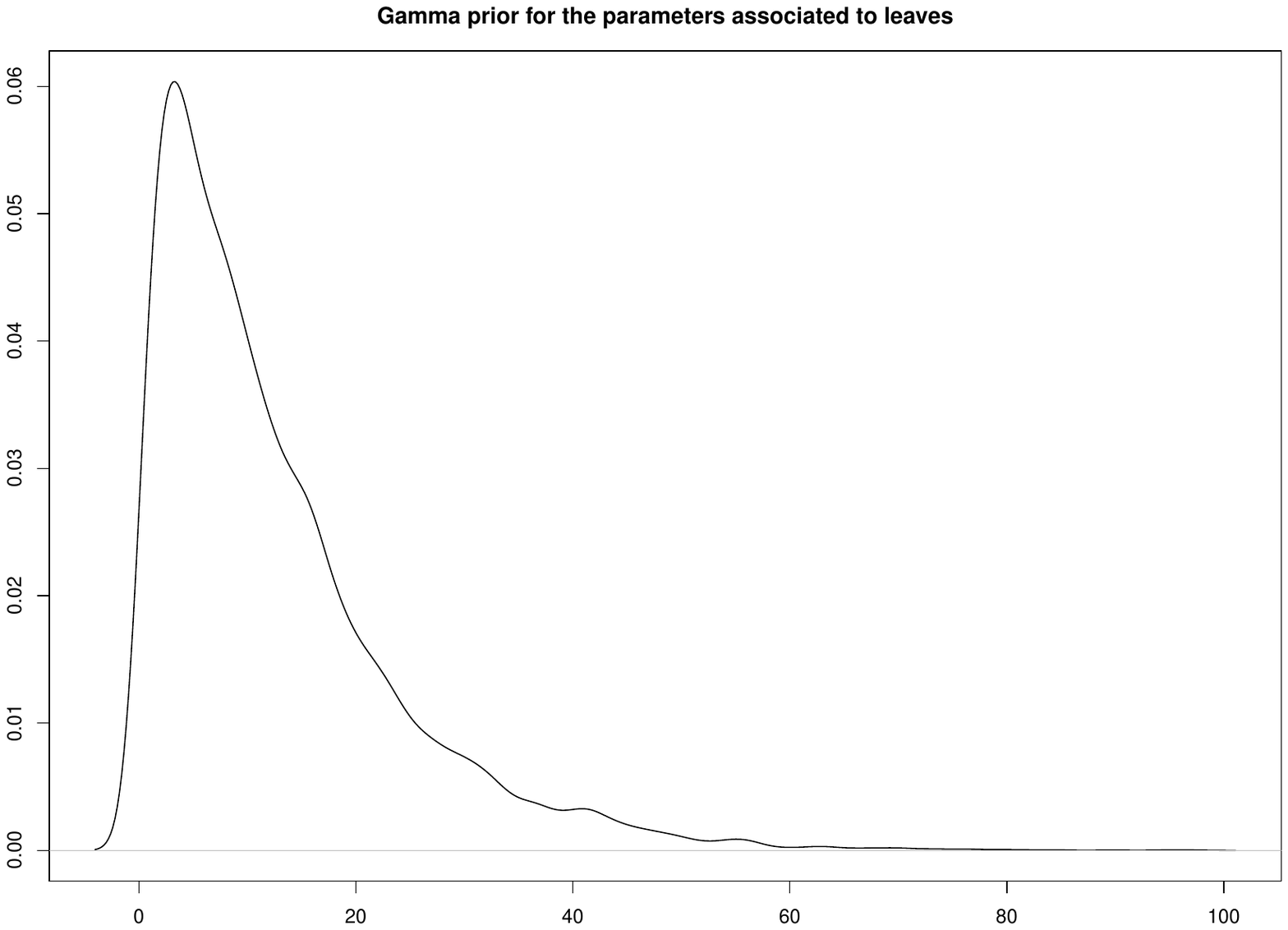}
    \caption{Prior for 10 Trees}
\end{figure}

\section{Earthquakes Data}

\subsection{10 Trees}
\noindent We run 3 parallel chains  each for 100000 iterations keeping every 50th sample. 
\begin{figure}[H] 
    \centering\includegraphics[width=8cm]{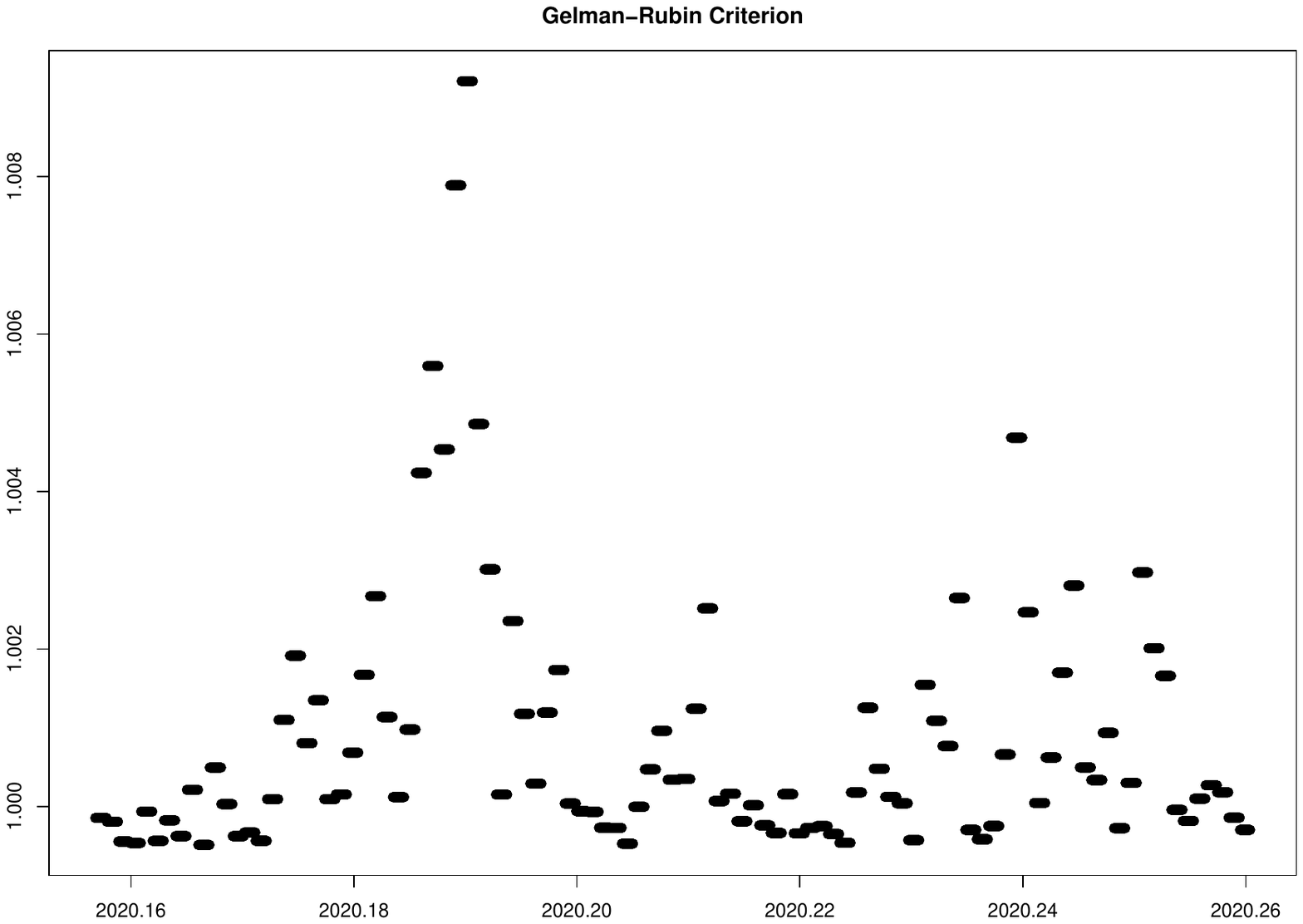}
\caption{The Gelman-Rubin Criterion for 10 Trees}
\end{figure}

\begin{figure}[H] 
    \centering\includegraphics[width=8cm]{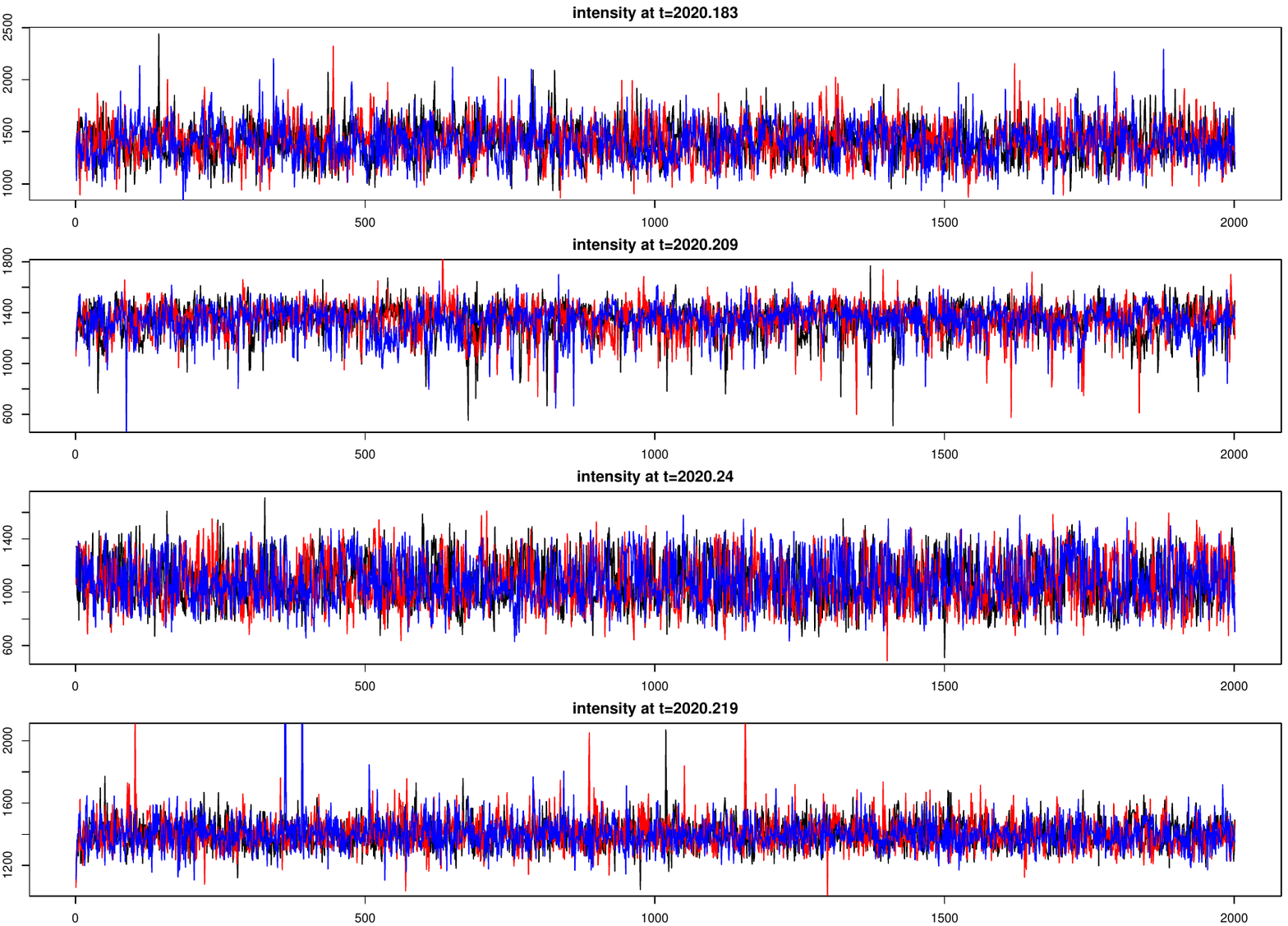}
    \caption{Trace plots for 10 Trees}
\end{figure}

\begin{figure}[H] 
    \centering\includegraphics[width=8cm]{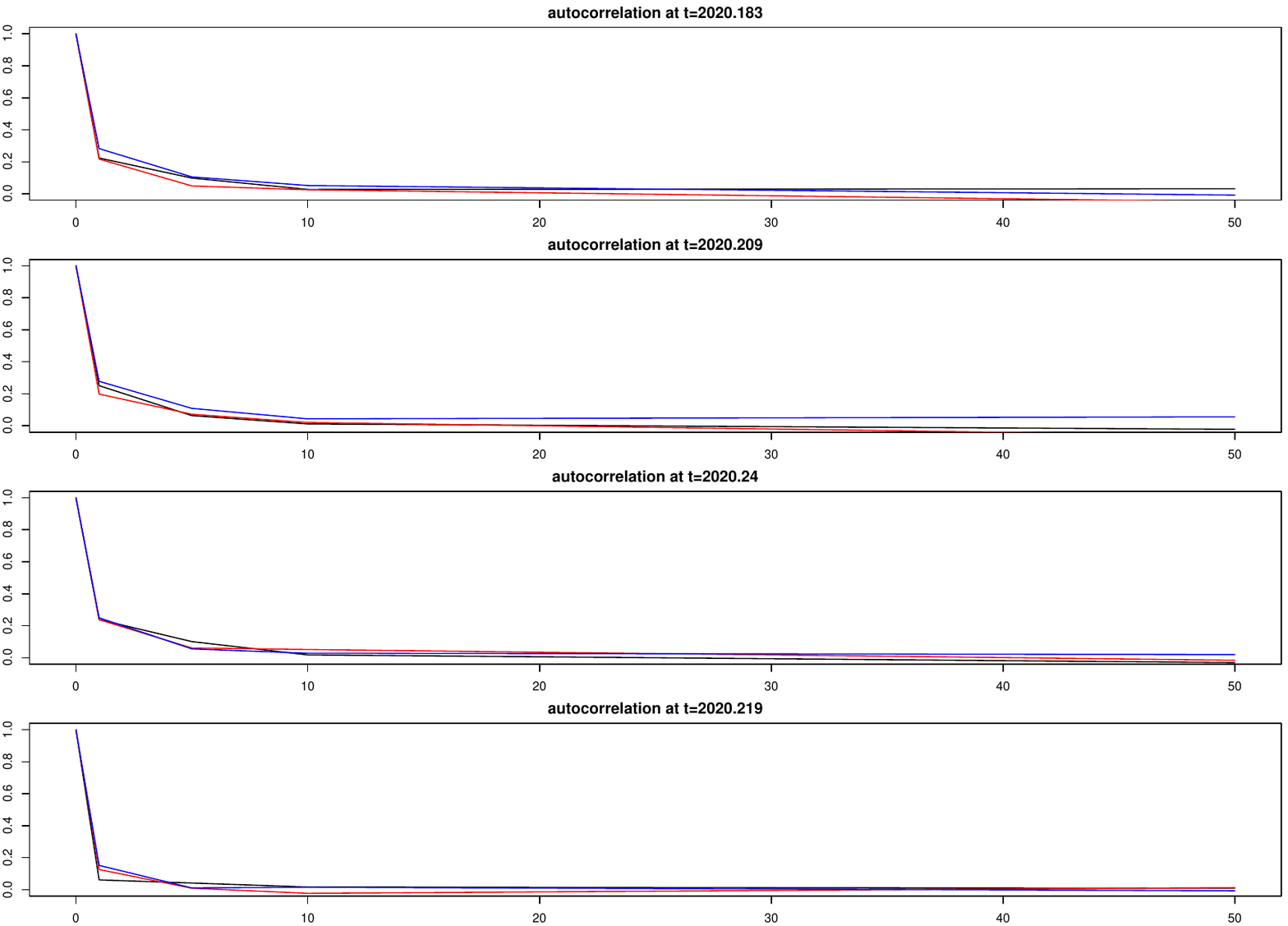}
    \caption{Autocorrelation plots for 10 Trees}
\end{figure}

\begin{figure}[H] 
    \centering\includegraphics[width=8cm]{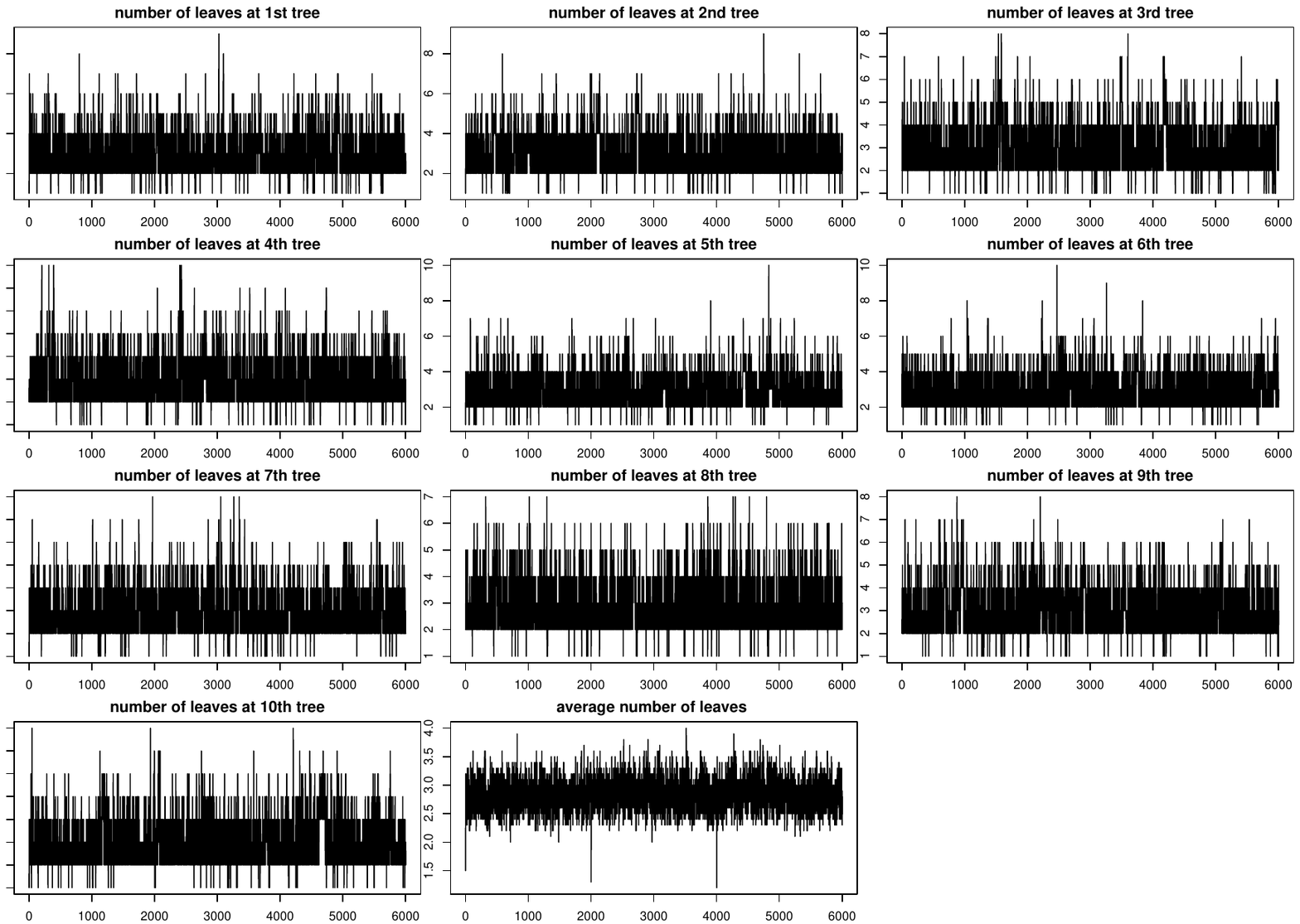}
    \caption{Average number of leaves at trees}
\end{figure}

\begin{figure}[H] 
    \centering\includegraphics[width=8cm]{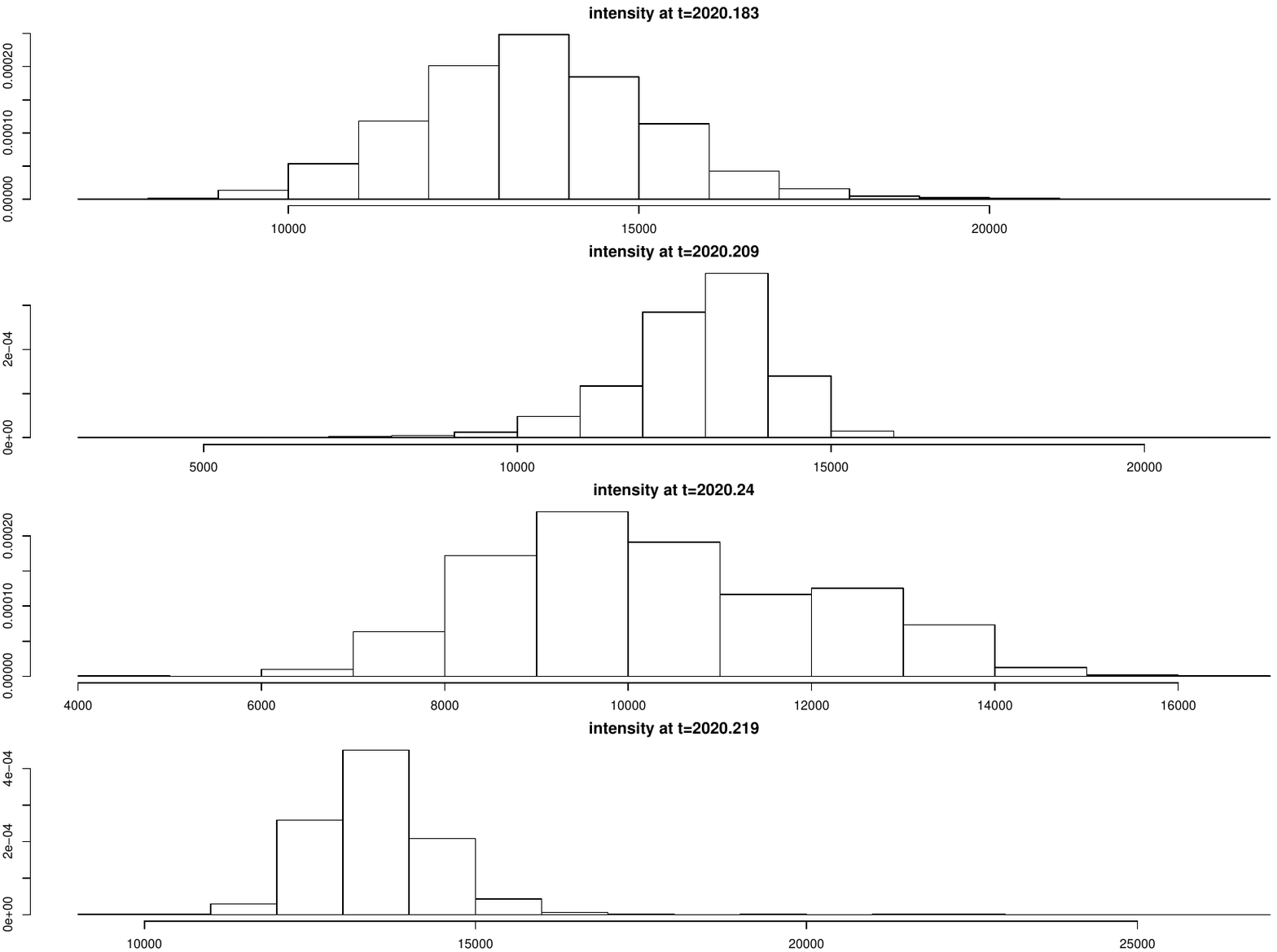}
    \caption{Density of the estimated intensity for 10 Trees}
\end{figure}

\begin{figure}[H] 
    \centering\includegraphics[width=8cm]{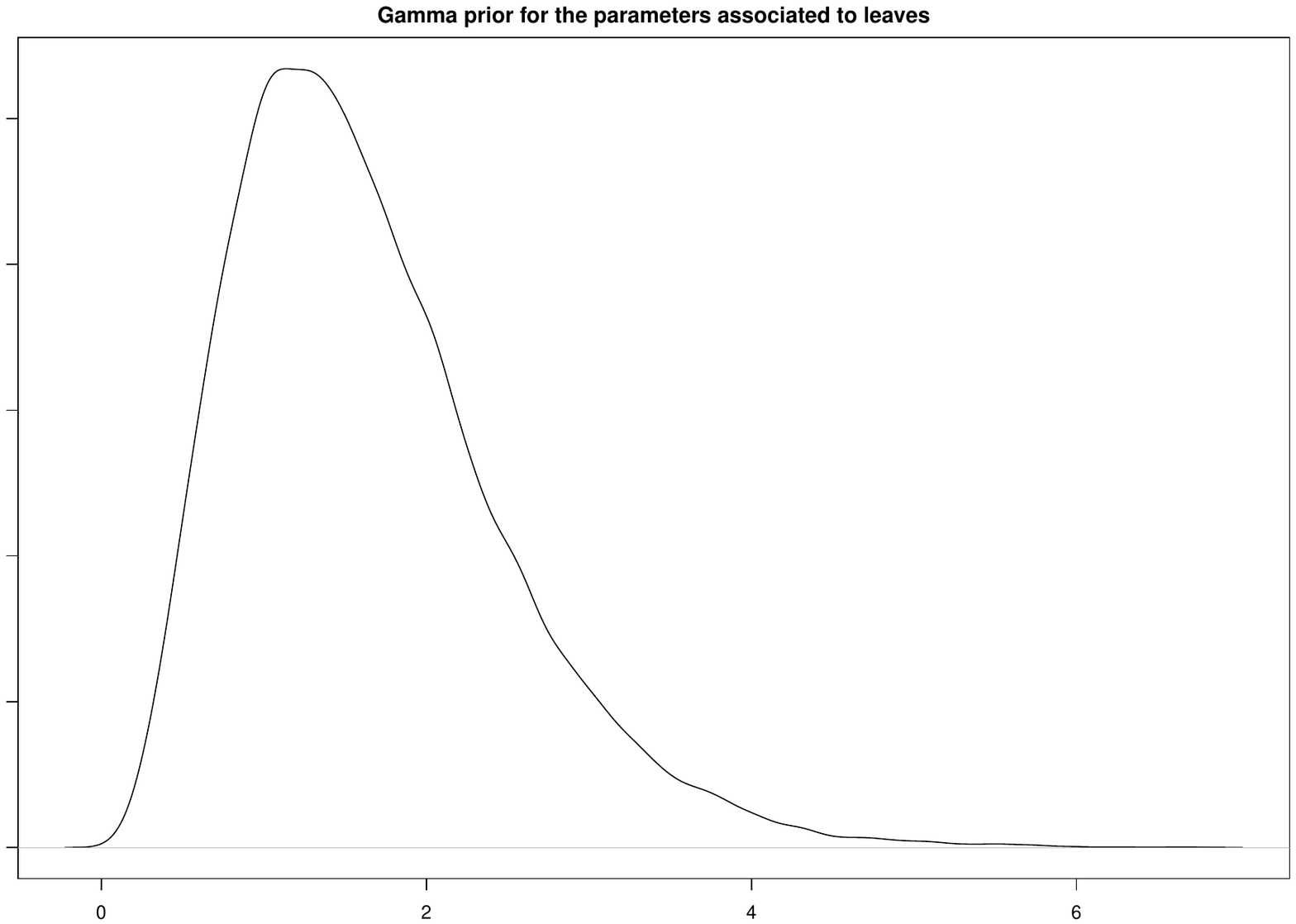}
    \caption{Prior for 10 Trees}
\end{figure}

\section{Mapples}

\subsection{5Trees}
\noindent We run 3 parallel chains  each for 300000 iterations keeping every 150th sample. 
\begin{figure}[H] 
    \centering\includegraphics[width=8cm]{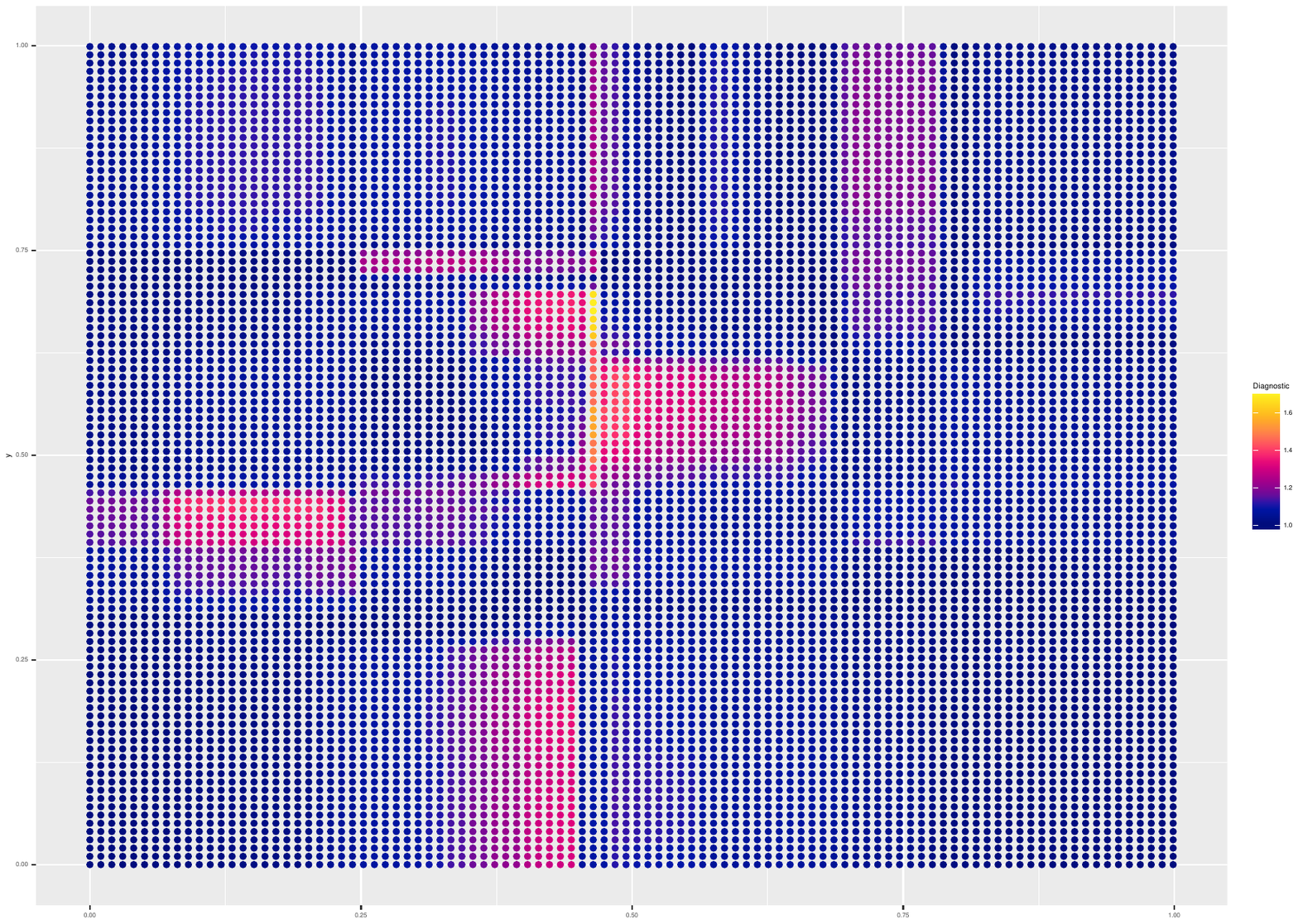}
\caption{The Gelman-Rubin Criterion for 5 Trees}
\end{figure}

\begin{figure}[H] 
    \centering\includegraphics[width=8cm]{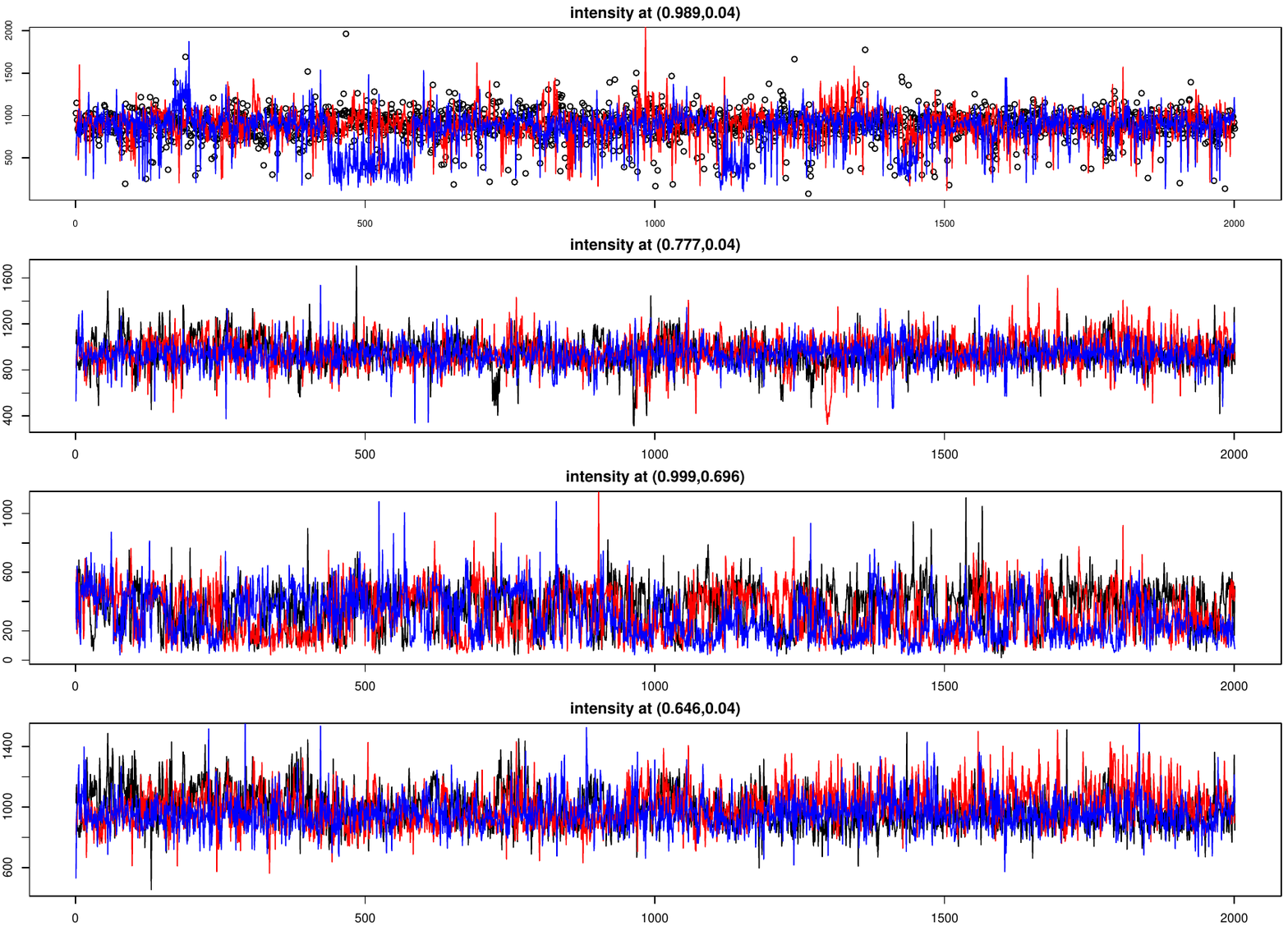}
    \caption{Trace plots for 5 Trees}
\end{figure}

\begin{figure}[H] 
    \centering\includegraphics[width=8cm]{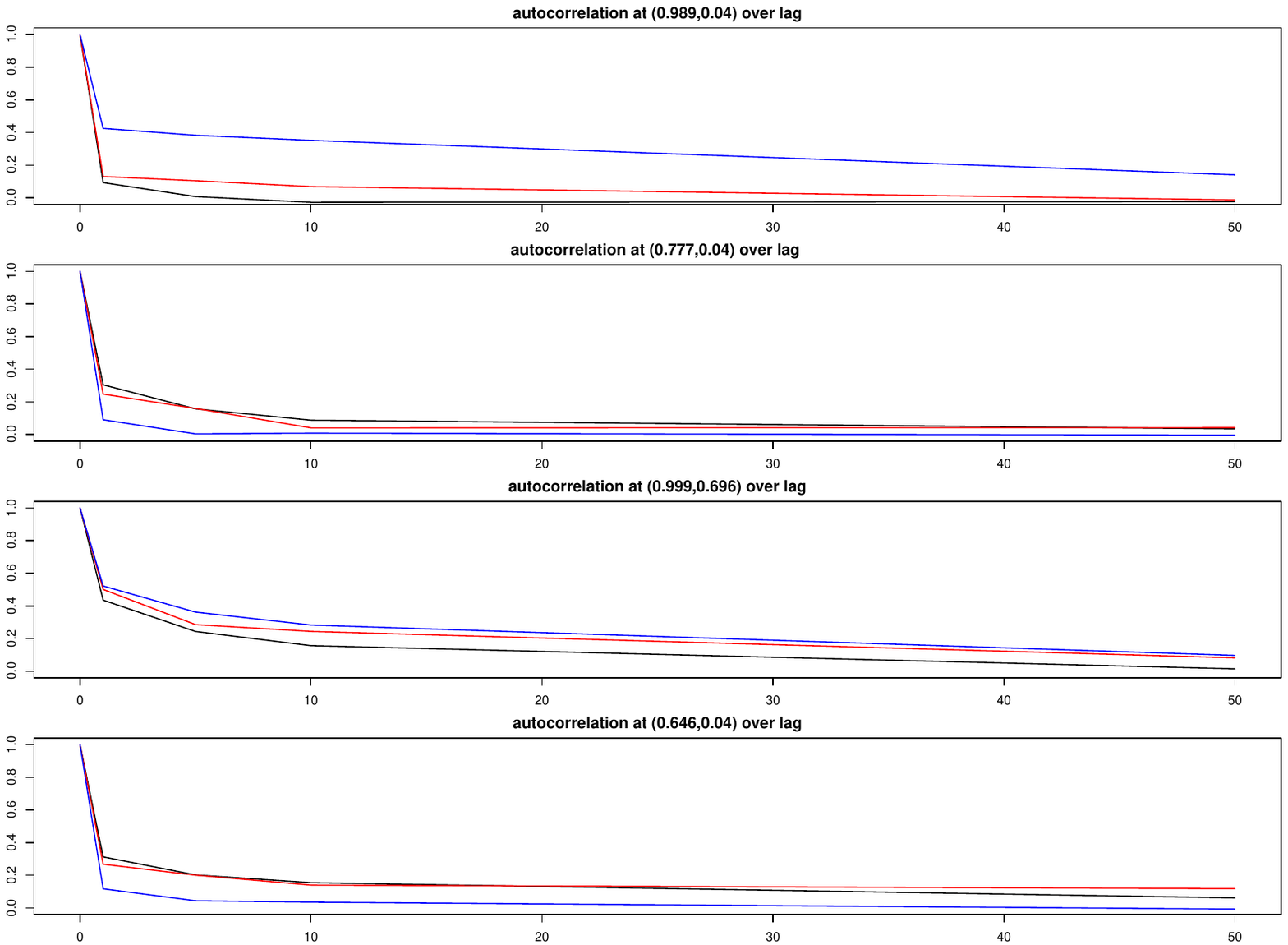}
    \caption{Autocorrelation plots for 5 Trees}
\end{figure}

\begin{figure}[H] 
    \centering\includegraphics[width=8cm]{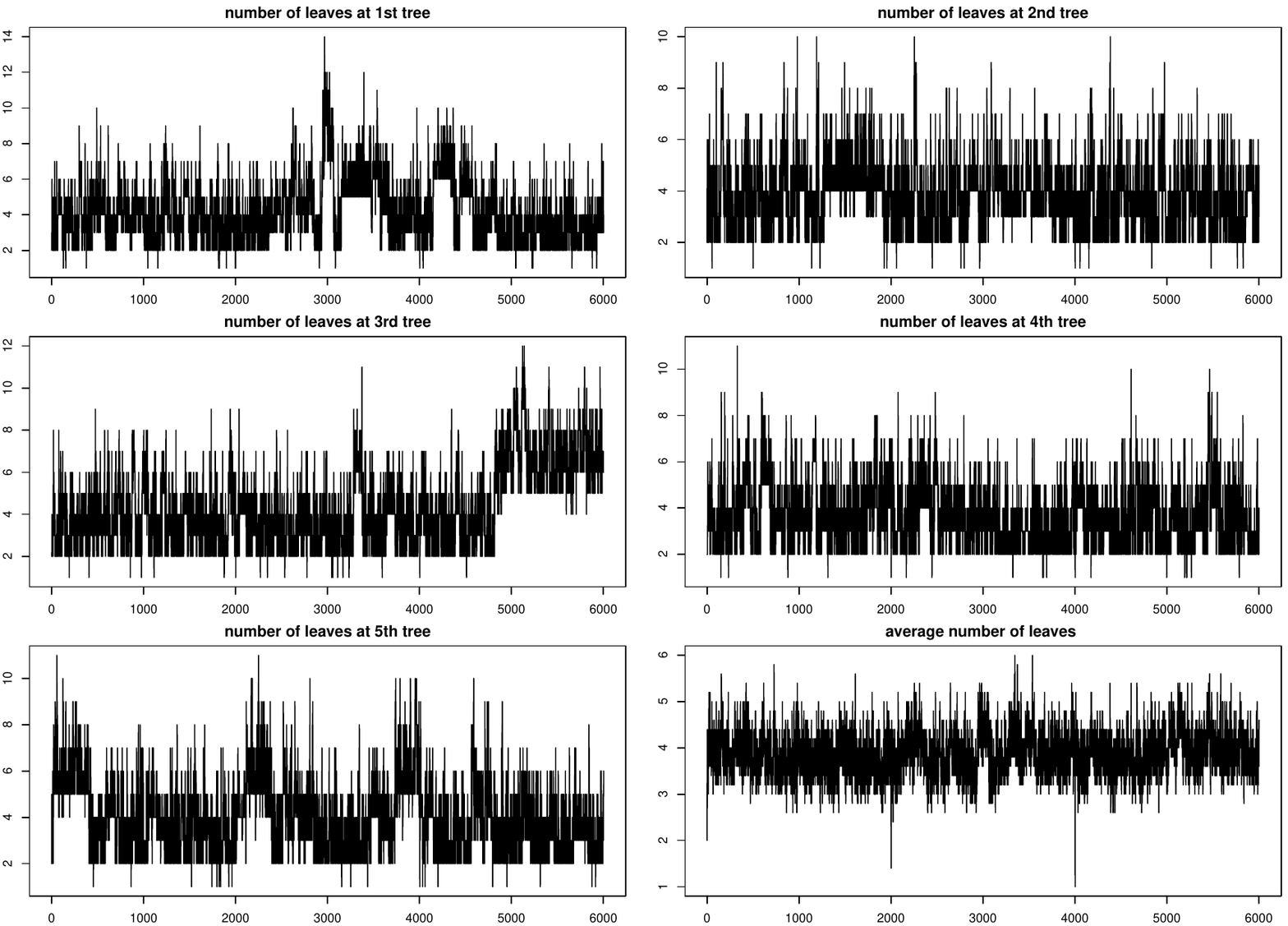}
    \caption{Average number of leaves at trees}
\end{figure}

\begin{figure}[H] 
    \centering\includegraphics[width=8cm]{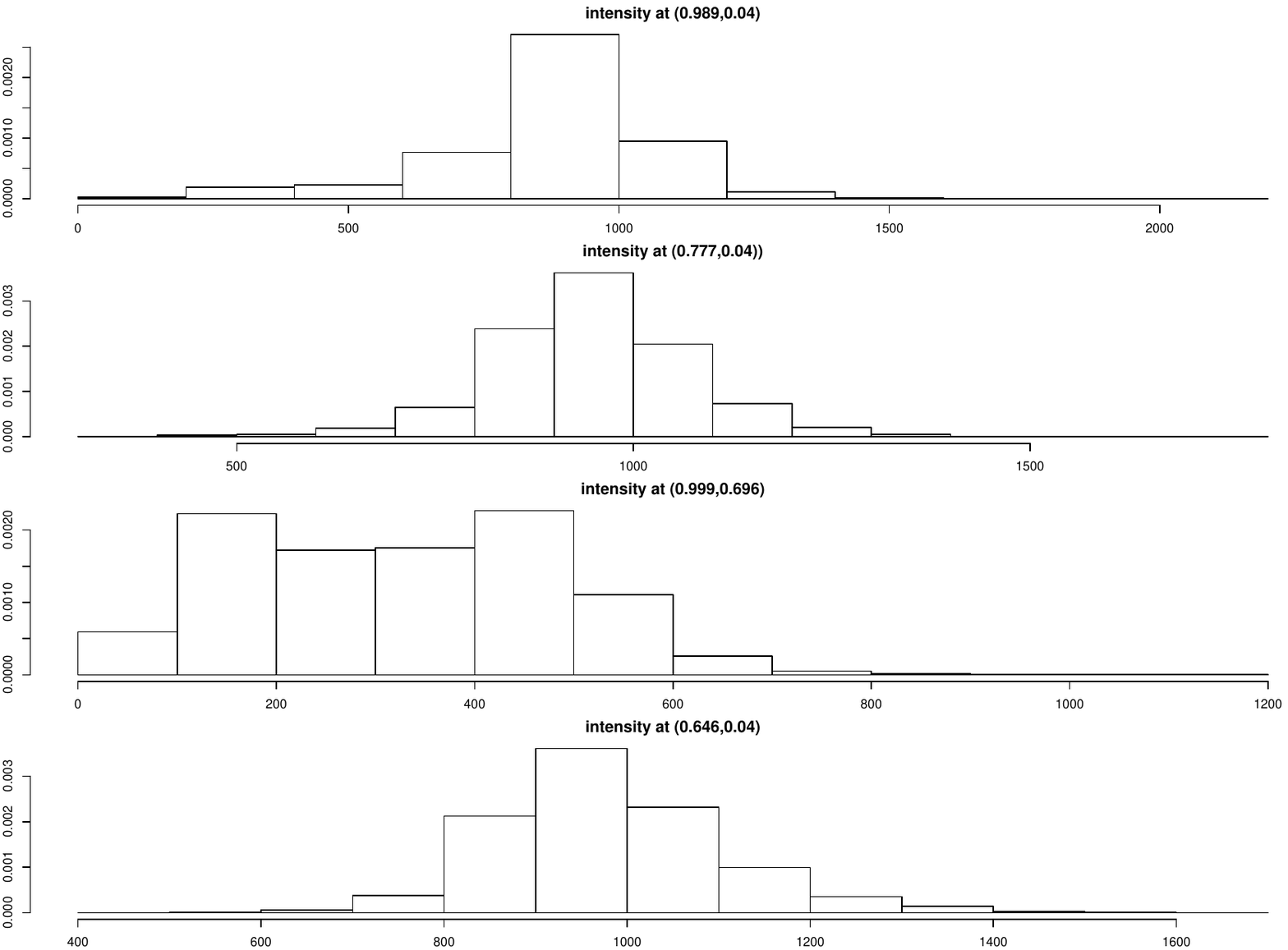}
    \caption{Density of the estimated intensity for 5 Trees}
\end{figure}

\begin{figure}[H] 
    \centering\includegraphics[width=8cm]{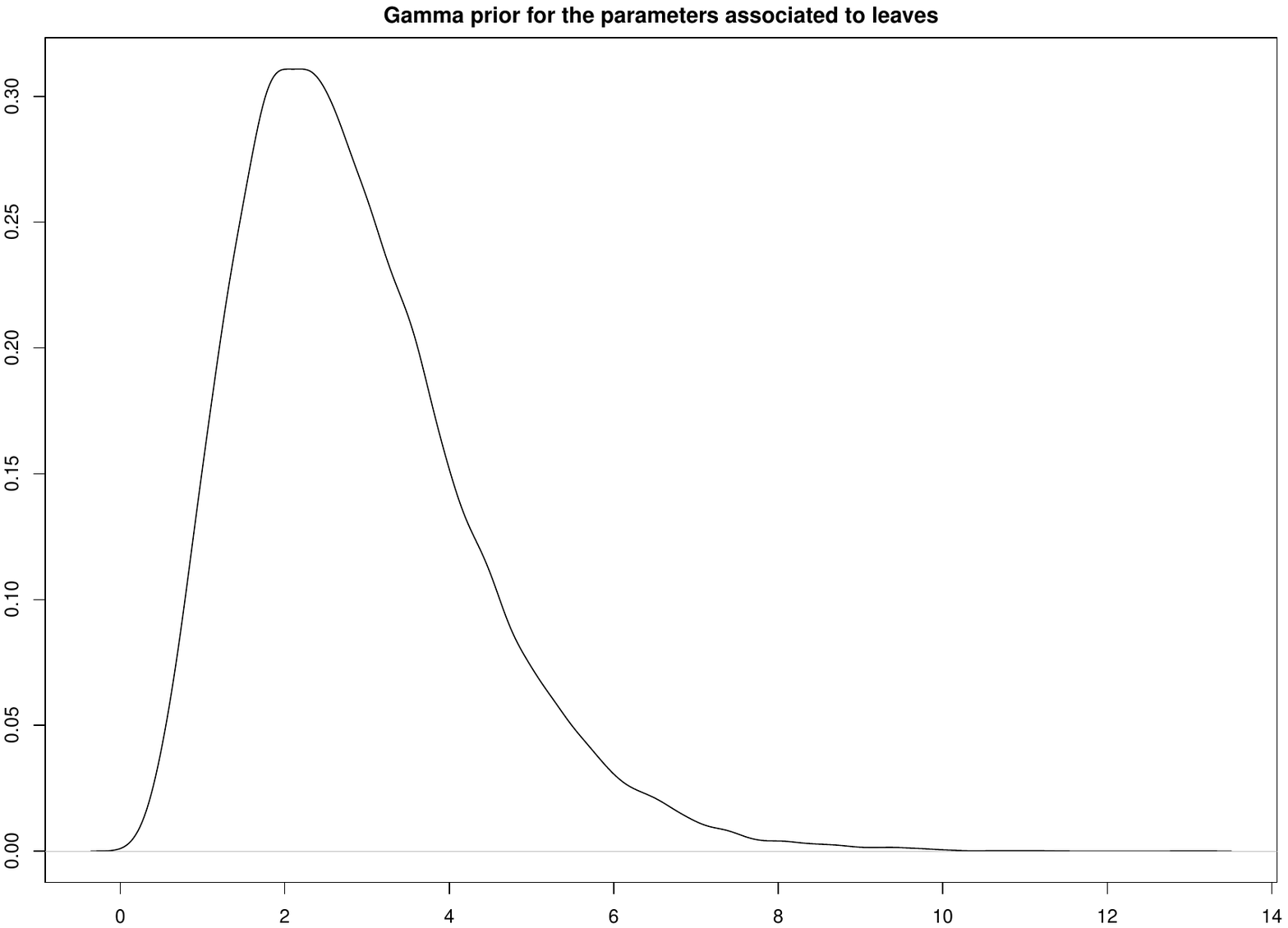}
    \caption{Prior for 5 Trees}
\end{figure}

\subsection{10 Trees}
\noindent We run 3 parallel chains  each for 300000 iterations keeping every 150th sample. 
\begin{figure}[H] 
    \centering\includegraphics[width=8cm]{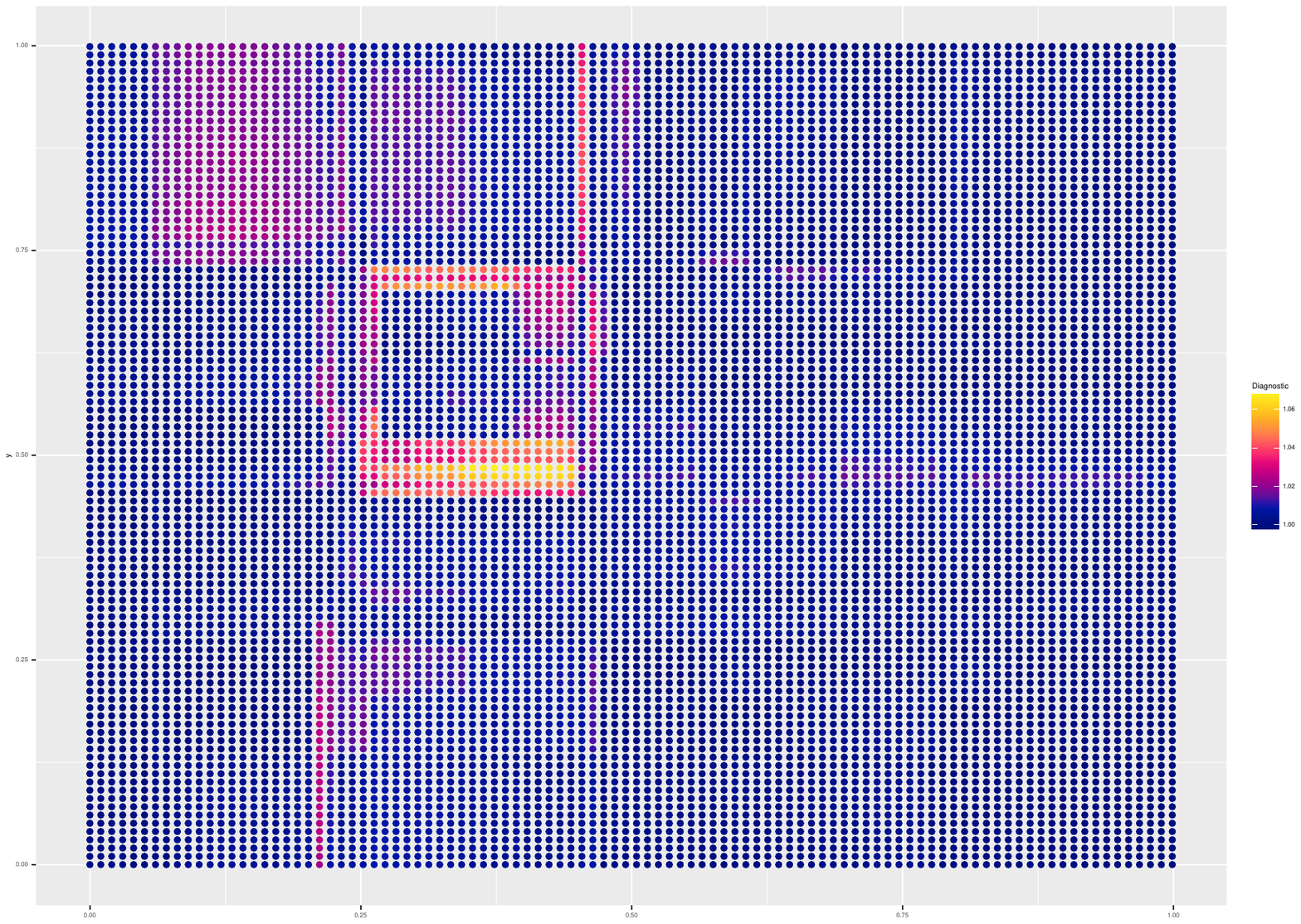}
\caption{The Gelman-Rubin Criterion for 10 Trees}
\end{figure}

\begin{figure}[H] 
    \centering\includegraphics[width=8cm]{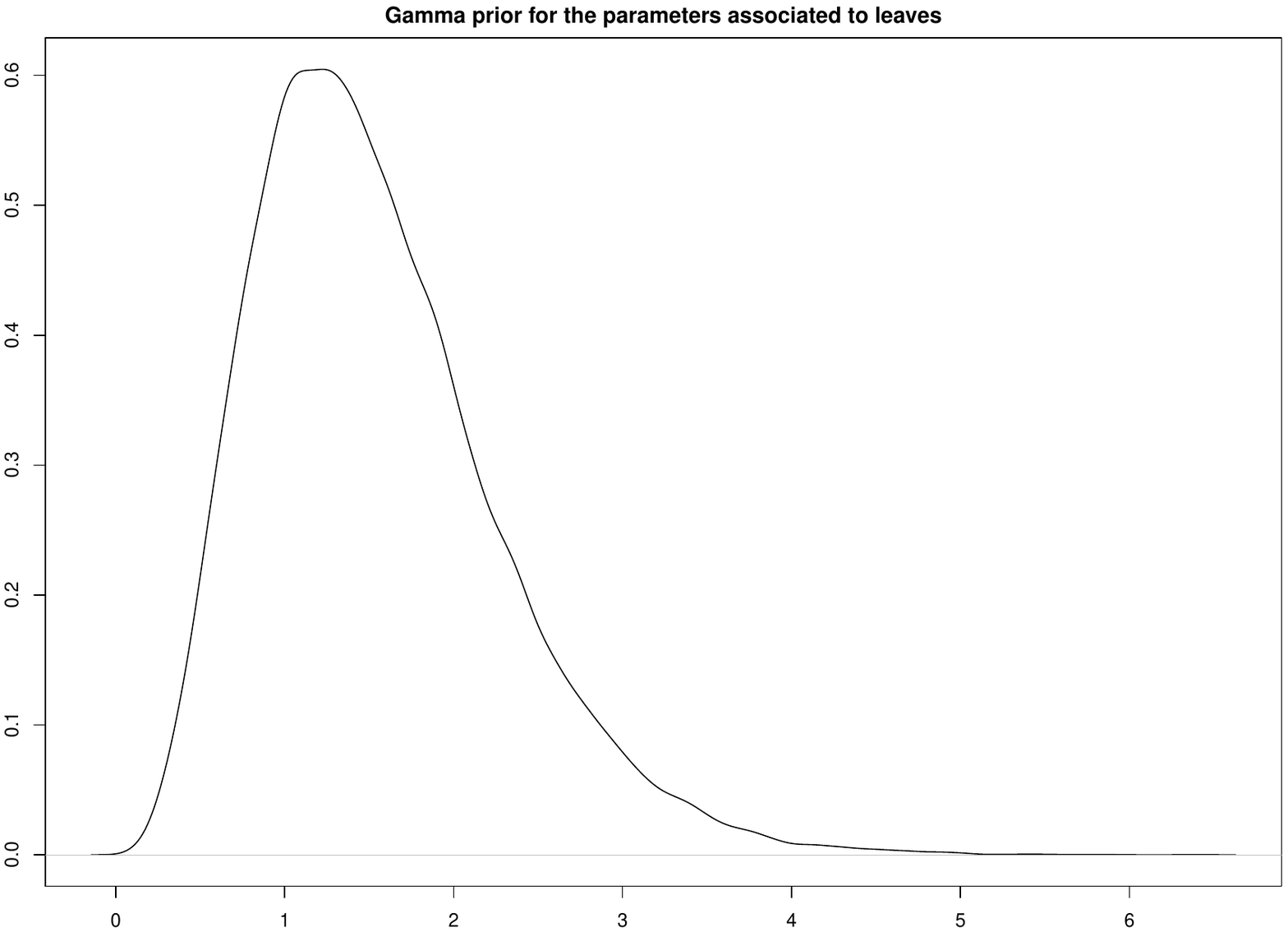}
    \caption{Prior for 10 Trees}
\end{figure}

\section{Redwood}

\subsection{5 Trees}
\noindent We run 3 parallel chains  each for 300000 iterations keeping every 150th sample. 
\begin{figure}[H] 
    \centering\includegraphics[width=8cm]{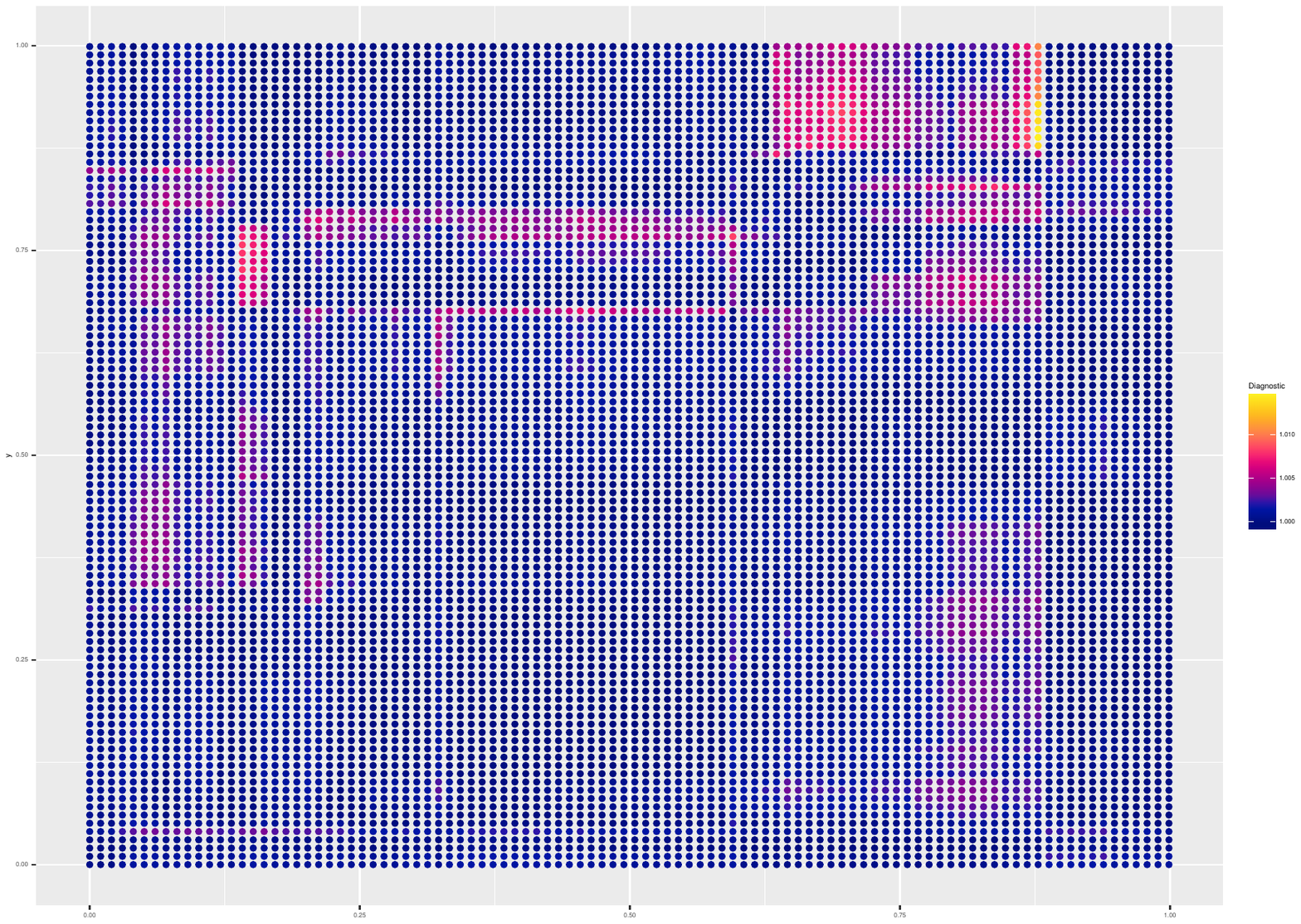}
\caption{The Gelman-Rubin Criterion for 5 Trees}
\end{figure}

\begin{figure}[H] 
    \centering\includegraphics[width=8cm]{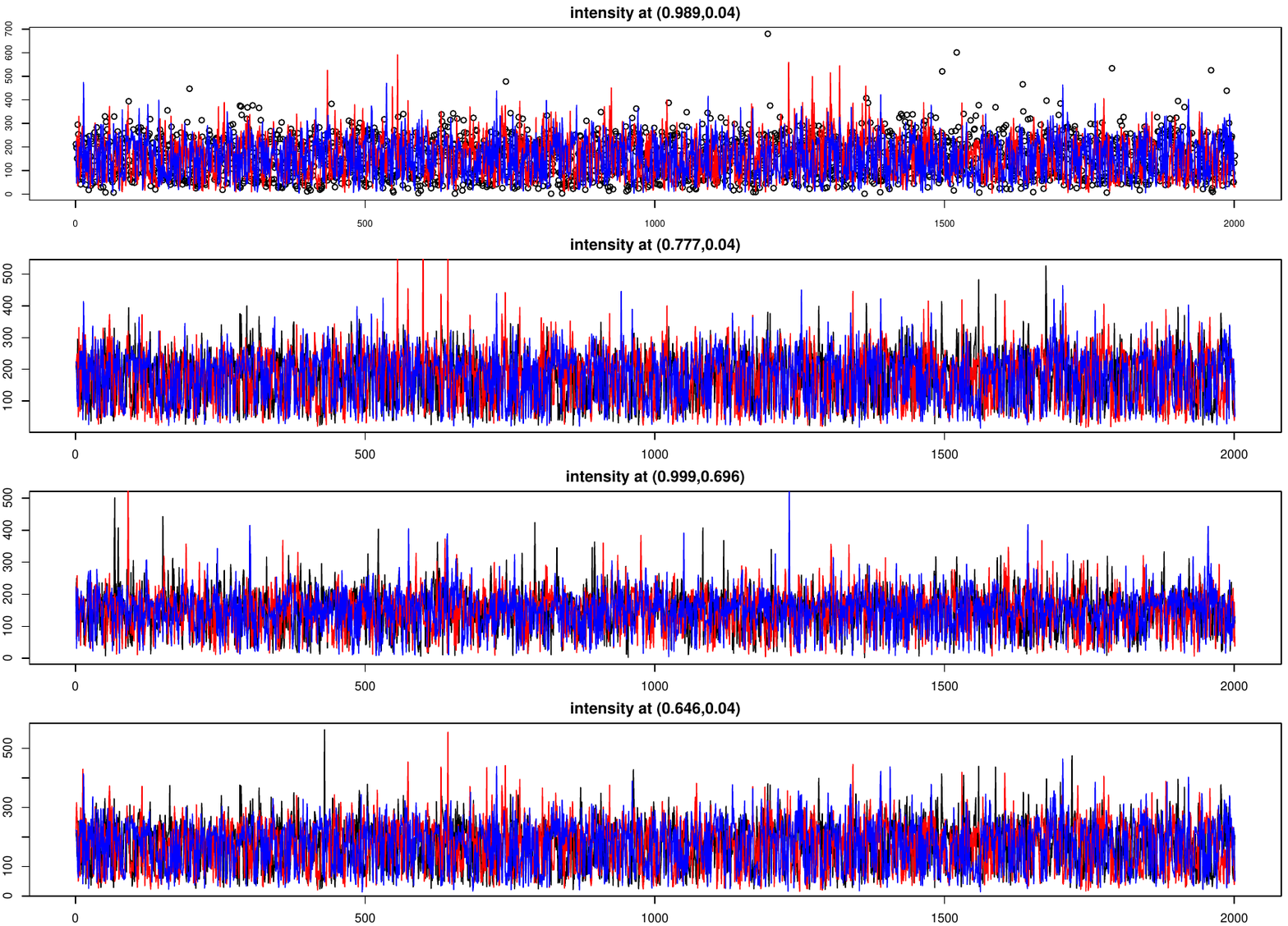}
    \caption{Trace plots for 5 Trees}
\end{figure}

\begin{figure}[H] 
    \centering\includegraphics[width=8cm]{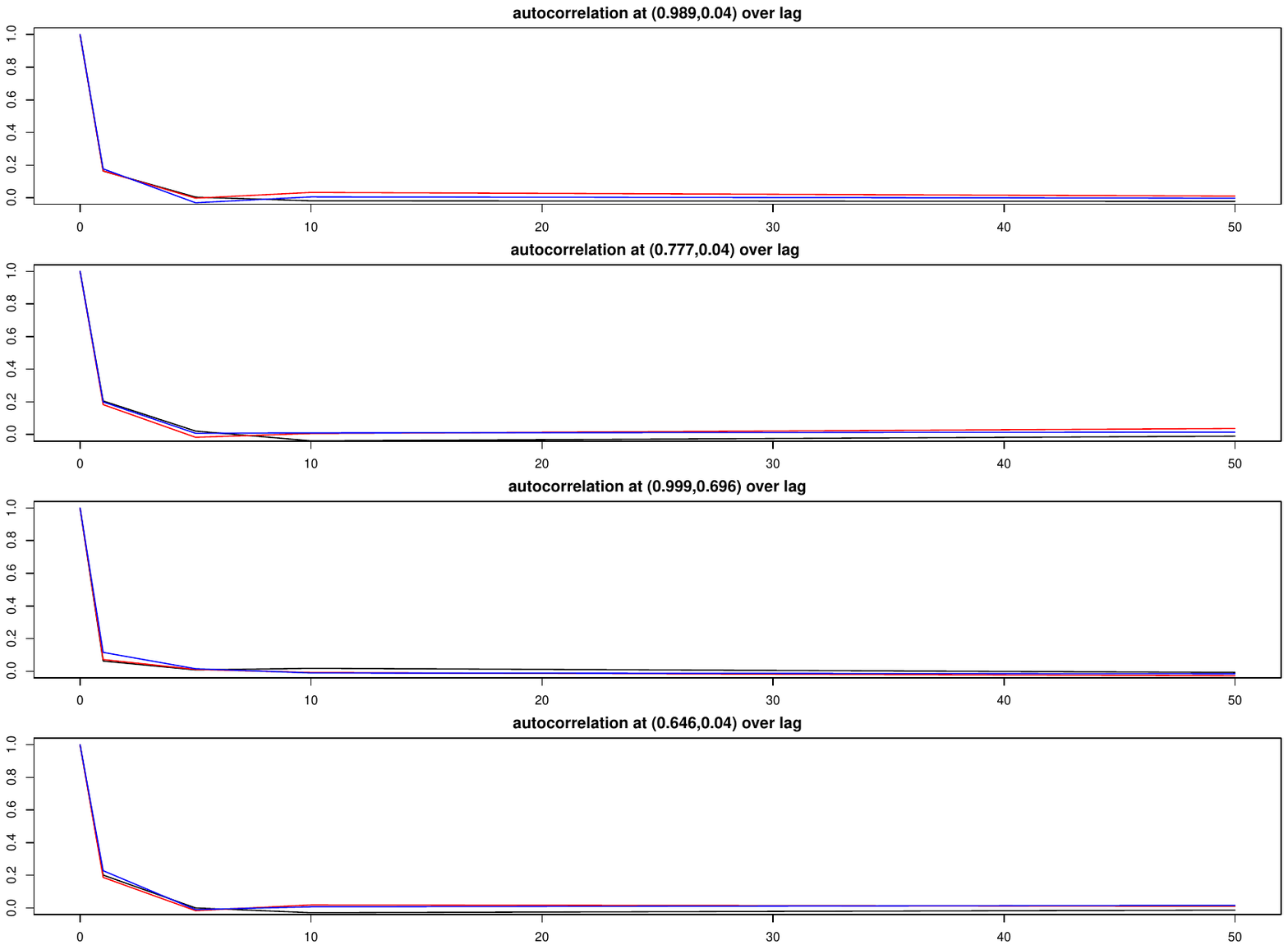}
    \caption{Autocorrelation plots for 5 Trees}
\end{figure}

\begin{figure}[H] 
    \centering\includegraphics[width=8cm]{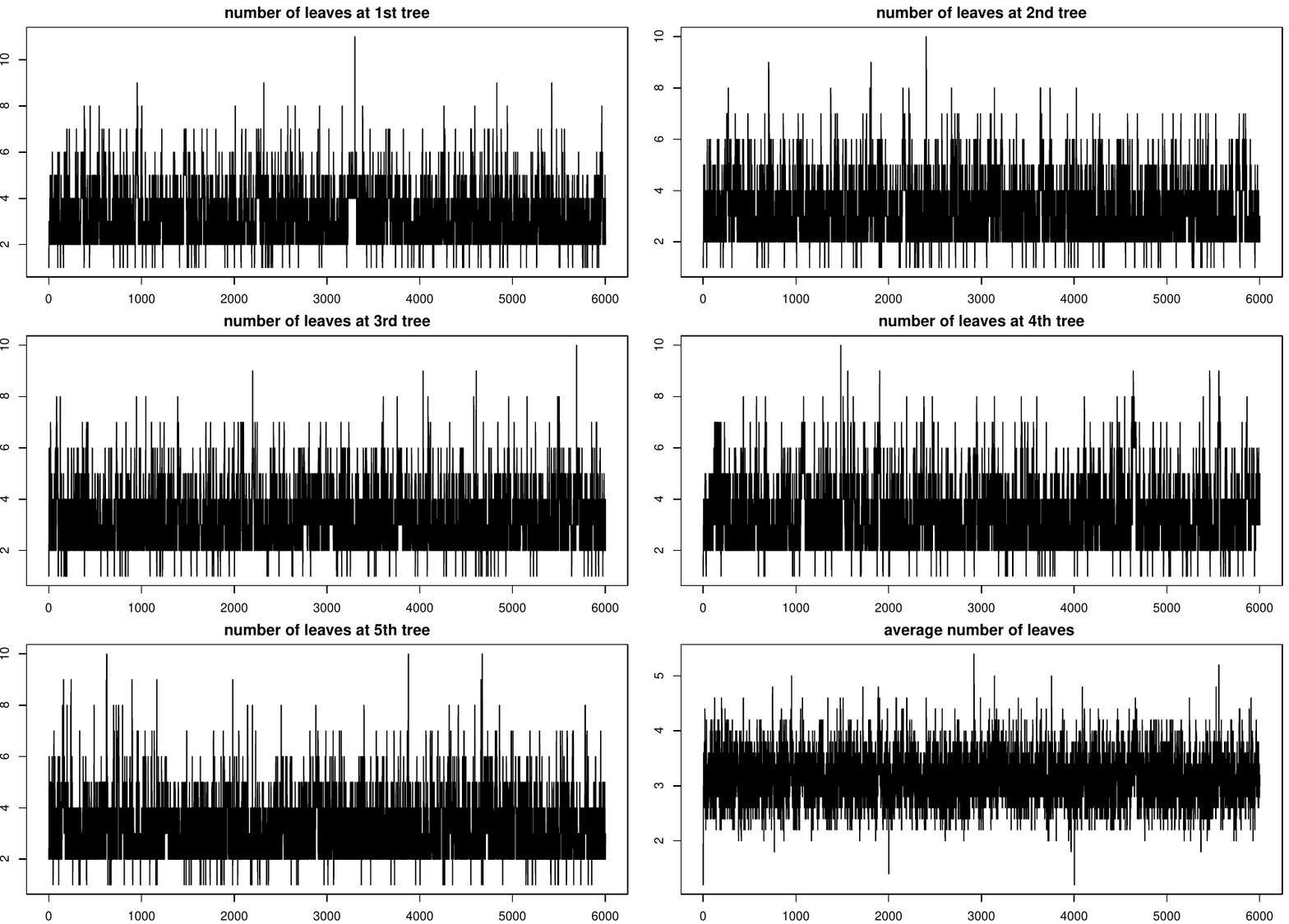}
    \caption{Average number of leaves at trees}
\end{figure}

\begin{figure}[H] 
    \centering\includegraphics[width=8cm]{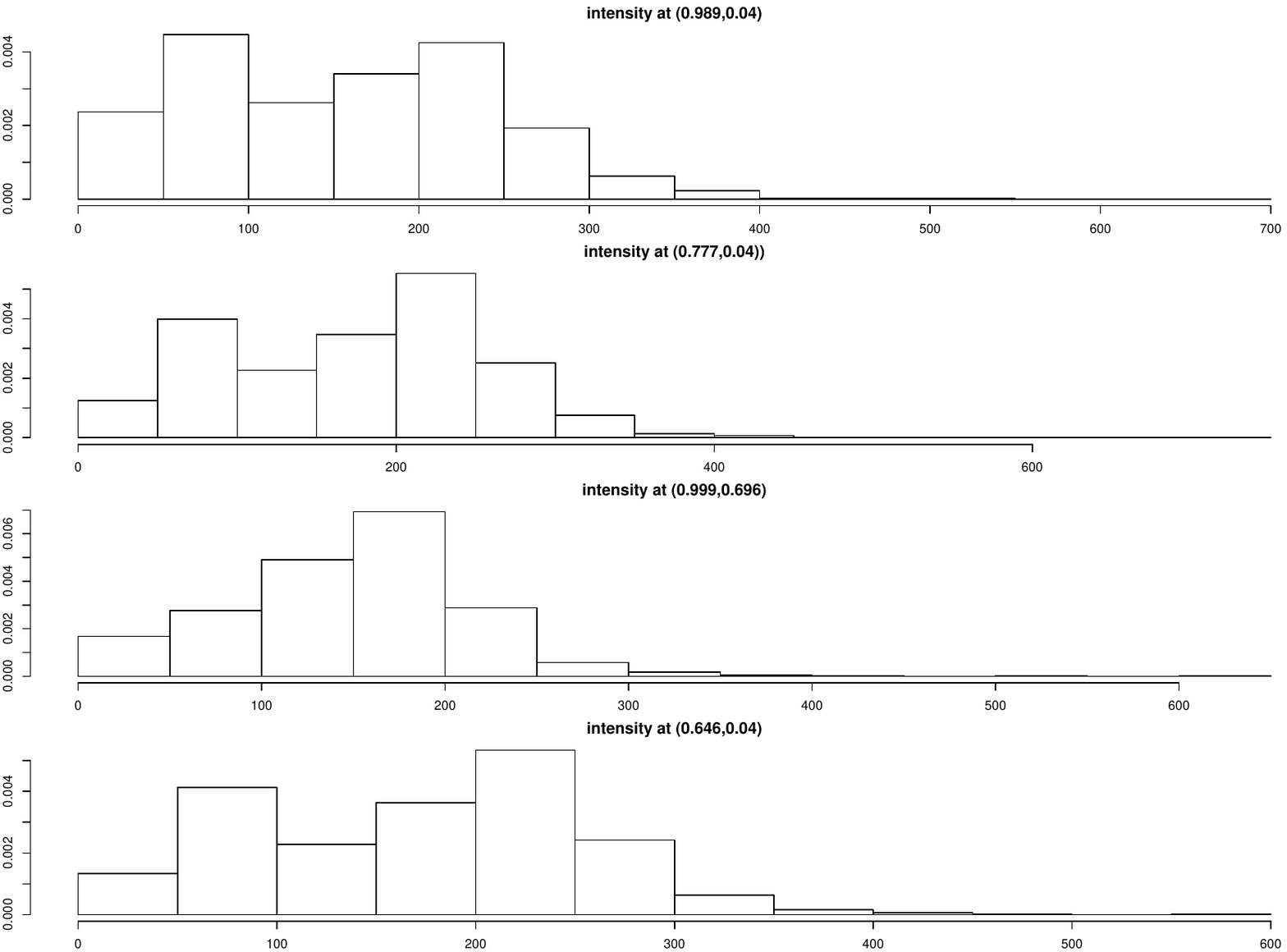}
    \caption{Density of the estimated intensity for 5 Trees}
\end{figure}

\begin{figure}[H] 
    \centering\includegraphics[width=8cm]{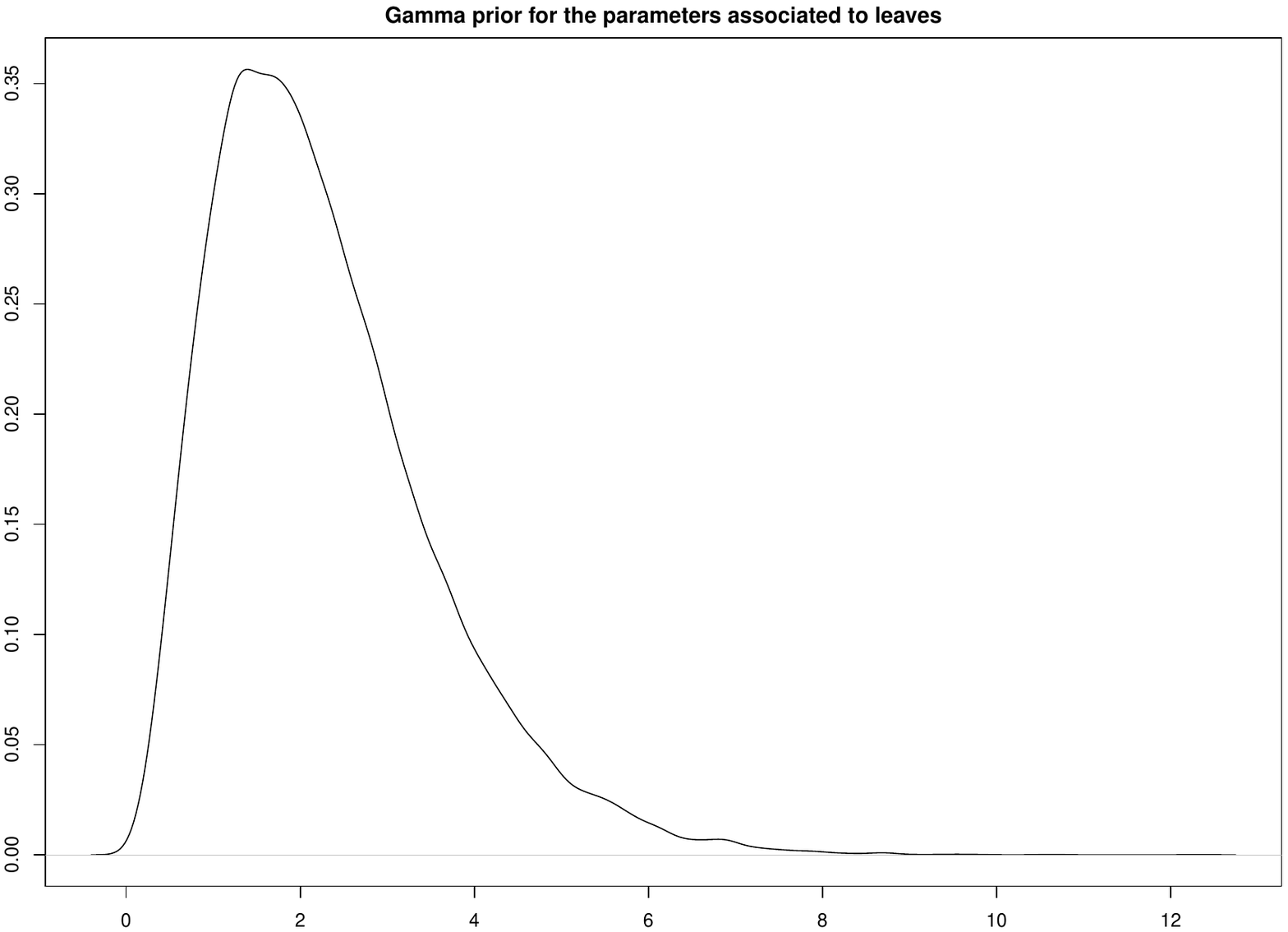}
    \caption{Prior for 5 Trees}
\end{figure}

\subsection{10 Trees}
\noindent We run 3 parallel chains  each for 300000 iterations keeping every 150th sample. 
\begin{figure}[H] 
    \centering\includegraphics[width=8cm]{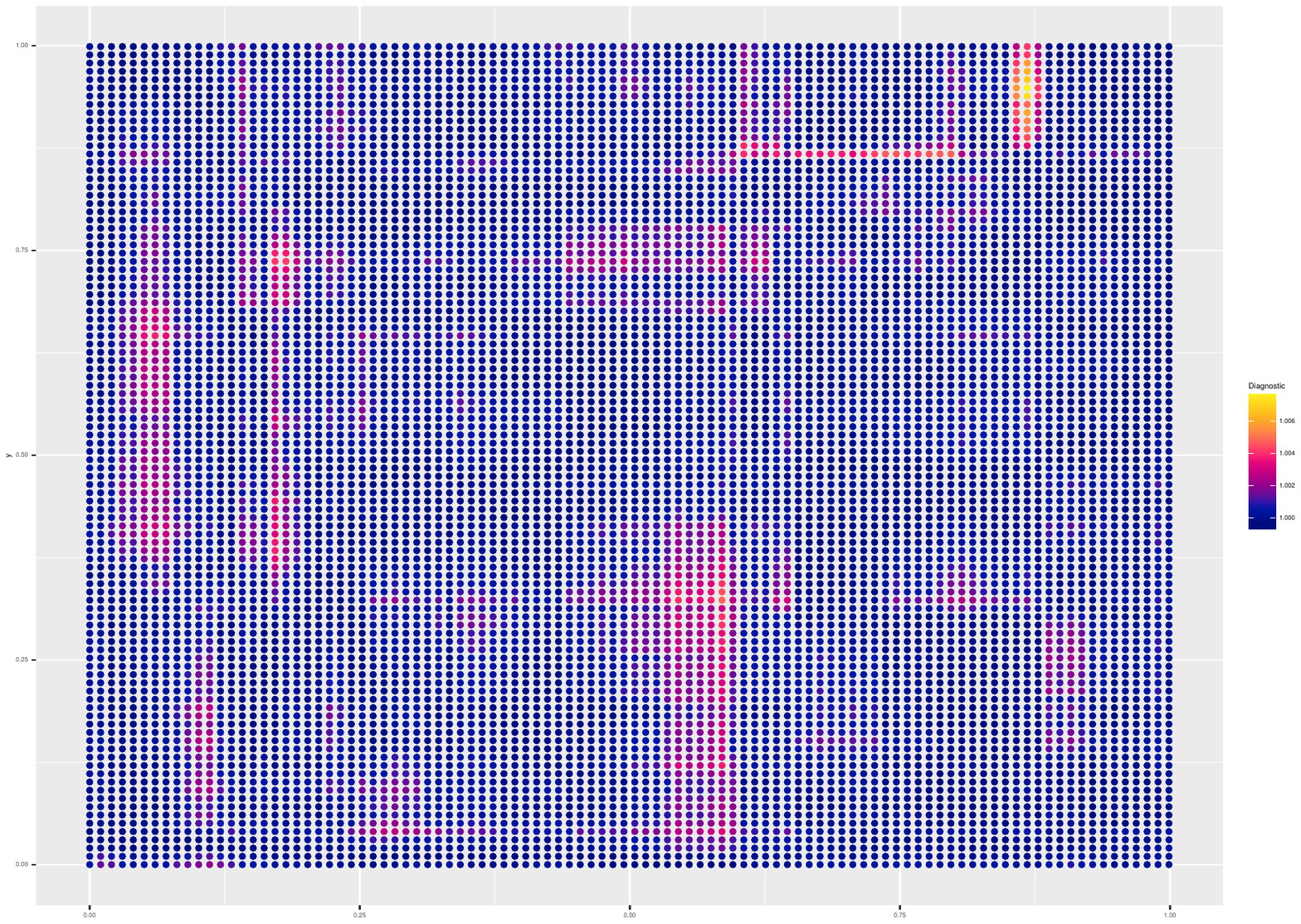}
\caption{The Gelman-Rubin Criterion for 10 Trees}
\end{figure}

\begin{figure}[H] 
    \centering\includegraphics[width=8cm]{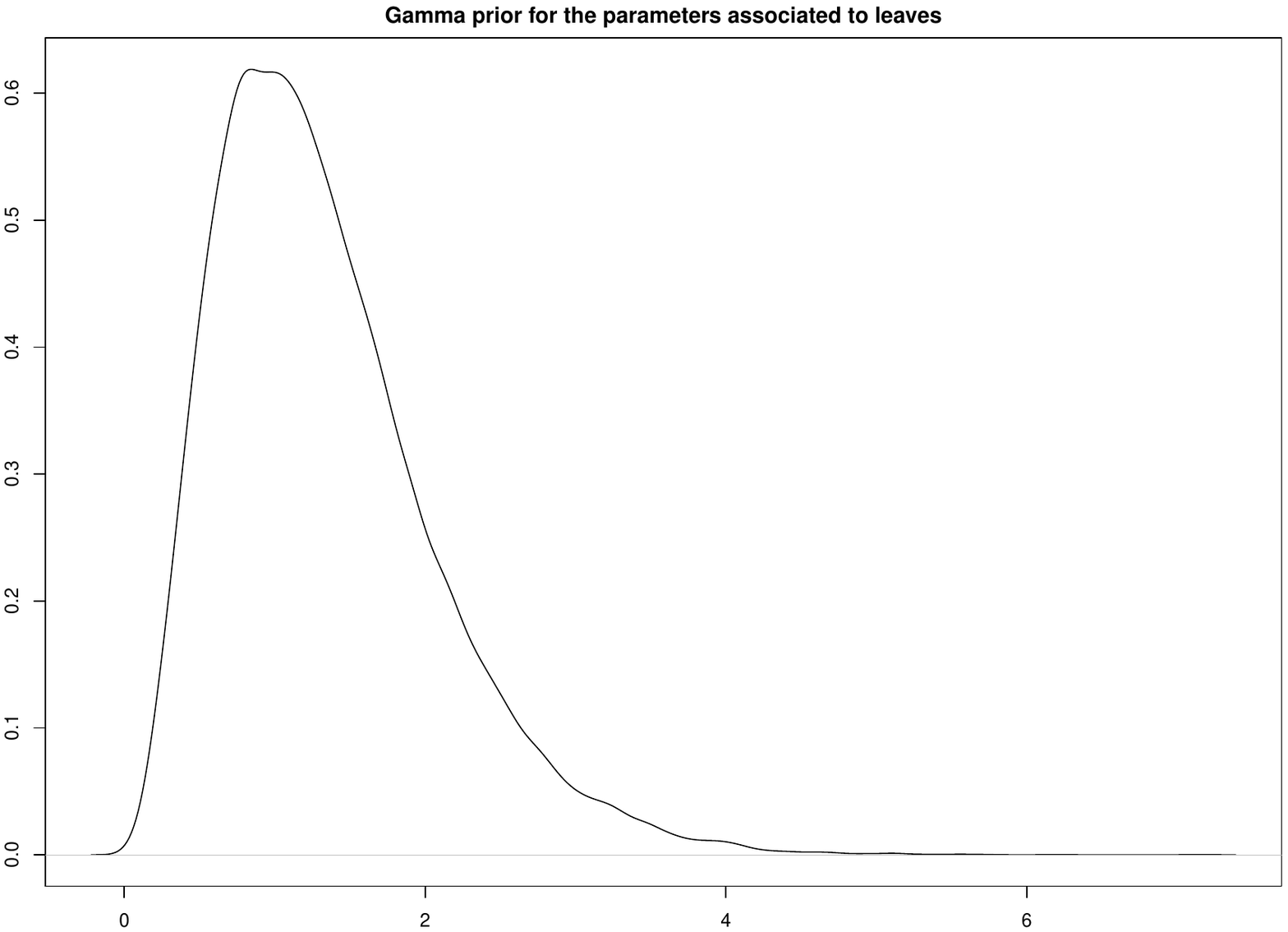}
    \caption{Prior for 10 Trees}
\end{figure}

\end{document}